\documentclass[12pt,reqno]{amsart}

\usepackage{amssymb,amsmath}

\usepackage[arrow,matrix,curve]{xy}

\newcommand{\al}{\alpha}
\newcommand{\bt}{\beta}
\newcommand{\gm}{\gamma}
\newcommand{\dl}{\delta}
\newcommand{\sg}{\sigma}

\newcommand{\om}{\omega}
\newcommand{\vf}{\varphi}

\newcommand{\bR}{\mathbb{R}}

\newcommand{\cA}{\mathcal{A}}
\newcommand{\cB}{\mathcal{B}}
\newcommand{\cC}{\mathcal{C}}
\newcommand{\cD}{\mathcal{D}}
\newcommand{\cF}{\mathcal{F}}
\newcommand{\FF}{\mathcal{F}}
\newcommand{\cG}{\mathcal{G}}
\newcommand{\cN}{\mathcal{N}}
\newcommand{\cP}{\mathcal{P}}
\newcommand{\cS}{\mathcal{S}}
\newcommand{\cT}{\mathcal{T}}
\newcommand{\cU}{\mathcal{U}}
\newcommand{\cV}{\mathcal{V}}

\newcommand{\cW}{\mathcal{W}}
\newcommand{\Cl}{\mathcal{C}l}
\newcommand{\Bd}{\mathcal{B}d}
\newcommand{\sD}{s\mathcal{D}}
\newcommand{\sBd}{s\mathcal{B}d}
\newcommand{\Catg}{\mathcal{C}atg}
\newcommand{\wND}{w\mathcal{N}\mathcal{D}}

\newcommand{\dd}{\text{-}\!}
\newcommand{\sbs}{\subset}
\newcommand{\sbsq}{\subseteq}
\newcommand{\ol}{\overline}
\newcommand{\ds}{\displaystyle}
\newcommand{\vnth}{\varnothing}
\newcommand{\lra}{\Longrightarrow}
\newcommand{\llra}{\Longleftrightarrow}

\newcommand{\RR}{\operatorname{R}}
\newcommand{\sS}{\operatorname{S}}
\newcommand{\TT}{\operatorname{T}}
\newcommand{\WW}{\operatorname{W}}

\newcommand{\Sp}{\operatorname{S}p}
\newcommand{\Wp}{\operatorname{W}p}

\newcommand{\I}{\operatorname{I}}
\newcommand{\II}{\operatorname{II}}

\newcommand{\BS}{\operatorname{BS}}
\newcommand{\TS}{\operatorname{TS}}
\newcommand{\SC}{\operatorname{SC}}
\newcommand{\SW}{\operatorname{SW}}
\newcommand{\WS}{\operatorname{WS}}
\newcommand{\WSC}{\operatorname{WSC}}
\newcommand{\BsS}{\operatorname{BsS}}
\newcommand{\BrS}{\operatorname{BrS}}
\newcommand{\TsS}{\operatorname{TsS}}

\newcommand{\cl}{\operatorname{cl}}
\newcommand{\co}{\operatorname{co}}
\newcommand{\nt}{\operatorname{int}}
\newcommand{\ord}{\operatorname{ord}}
\newcommand{\ind}{\operatorname{ind}}
\newcommand{\Ind}{\operatorname{Ind}}
\newcommand{\sind}{\operatorname{sind}}
\newcommand{\sInd}{\operatorname{sInd}}

\newcommand{\Fr}{\operatorname{Fr}}

\newenvironment{theorem}[1]{\vskip+0.2cm \textbf{Theorem #1.}\em }{\vskip+0.2cm}
\newenvironment{lemma}[1]{\vskip+0.2cm \textbf{Lemma #1.}\em }{\vskip+0.2cm}
\newenvironment{corollary}[1]{\vskip+0.2cm \textbf{Corollary #1.}\em }{\vskip+0.2cm}
\newenvironment{proposition}[1]{\vskip+0.2cm \textbf{Proposition #1.}\em }{\vskip+0.2cm}

\newenvironment{definition}[1]{\vskip+0.2cm \textbf{Definition #1.}}{\vskip+0.2cm}
\newenvironment{example}[1]{\vskip+0.2cm \textbf{Example #1.}}{\vskip+0.2cm}
\newenvironment{remark}[1]{\vskip+0.2cm \textbf{Remark #1.}}{\vskip+0.2cm}

\newenvironment{pf}{\noindent {\it Proof.}}{\ \hfill $\square$}

\begin{document}

\begin{center}
    \textbf{\large On the Theory of Relative Bitopological and Topological
    Properties}\\[0.4cm]
    \textbf{B. P. Dvalishvili} \\[0.2cm]
    {\small \textit{Department of Algebra-Geometry of the Institute
        of Mechanics and Mathematics of the Faculty
        of Exact and Natural Sciences
        of I.~Javakhishvili Tbilisi State University, \\
        2, University St., 0143 Tbilisi, Georgia}}
\end{center}

\section*{\textbf{1. Introduction}}
\vskip+0.2cm

\markboth{B. P. DVALISHVILI}{ON THE THEORY OF RELATIVE
BITOPOLOGICAL AND $\ldots$}

By analogy with [4], the leading idea, permeating the work, can be
briefly described as follows: with each (bi)topological property
$\cP$ one can associate a relative version of it formulated in
terms of location of a (bi)topological subspace $Y$ in a
(bi)topological space $X$ in such a natural way that, when $Y$
coincides with $X$, then the relative property coincides with
$\cP$. Moreover, please note that all bitopological versions are
constructed in the commonly accepted manner so that if the
topologies coincide, we obtain the original or new topological
notions and results.

Observe also here that if after bitopological results their
topological counterparts are also given, then this means that the
counterparts are new too.

All useful notions have been collected and the following
abbreviations are used throughout the work: $\TS$ for a
topological space, $\TsS$ for a topological subspace, $\BS$ for a
bitopological space and $\BsS$ for a bitopological subspace.
Always. $i,j\in\{1,2\}$, $i\neq j$, unless stated other\-wise.

Let $(X,\tau_1,\tau_2)$ be a $\BS$ and $\cP$ be some topological
property. Then $(i,j)\dd\cP$ denotes  the analogue of this
property for $\tau_i$ with respect to $\tau_j$, and $p\,\dd\cP$
denotes the conjunction $(1,2)\dd\cP\wedge(2,1)\dd\cP$, that is,
$p\,\dd\cP$ denotes an ``absolute'' bitopological analogue of
$\cP$, where ``$p$'' is the abbreviation for ``pairwise''. Also
note that $(X,\tau_i)$ has a property $\cP$ if and only if
$(X,\tau_1,\tau_2)$  has the property $i\dd\cP$, and $d\dd\cP$ is
equivalent to $1\dd\cP\wedge 2\dd\cP$, where ``$d$'' is the
abbreviation for ``double''. A $\BsS$ $(Y,\tau_1',\tau_2')$ of a
$\BS$ $(X,\tau_1,\tau_2)$ has a property $\cP$ if and only if
$(Y,\tau_1',\tau_2')$ has the property $\cP$ in itself.

As usual, the symbol $2^X$ is used for the power set of the set
$X$, and for a family $\cA=\{A_s\}_{s\in S}\sbs 2^X$, $\co\cA$
denotes the conjugate family $\{X\setminus A_s,\;A_s\in\cA\}_{s\in
S}$. If $A\sbs X$, then $\tau_i\nt A$ and $\tau_i\cl A$ denote
respectively interiors and closures of $A$ in the topologies
$\tau_i$ (for a $\TS$ $(X,\tau)$ the interior and the closure of a
subset $A\sbs X$ is denoted by $\tau\nt A$ and $\tau\cl A$,
respectively). A set $A\sbs X$ is $p$-closed in
$(X,\tau_1,\tau_2)$ if $A=\tau_1\cl A\cap\tau_2\cl A$ and the
family of all $p$-closed subsets of $(X,\tau_1,\tau_2)$ is denoted
by $p\,\dd\,\Cl(X)$. It is clear, that $
\co\tau_1\cup\co\tau_2\sbs p\,\dd\,\Cl(X)$ and for
$\tau_1\sbs\tau_2$, that is, for a $\BS$ $(X,\tau_1<\tau_2)$,
$p\,\dd\,\Cl(X)=\co\tau_2$.

A family $\boldsymbol{\cU}=\{U_s\}_{s\in S}$ of subsets of $X$ is
a $p$-open covering of $X$ if
$\boldsymbol{\cU}\sbs\tau_1\cup\tau_2$, $X=\bigcup\limits_{s\in S}
U_s$ and $\boldsymbol{\cU}\cap \tau_i$ contains a nonempty set
[11].

In our further discussion $(\bR,\om)$ is the natural $\TS$, and
the natural $\BS$ $(\bR,\om_1,\om_2)$ is the real line $\bR$ with
the lower
$\om_1=\{\vnth,\bR\}\cup\,\big\{(a,+\infty):\;a\in\bR\big\}$ and
upper $\om_2=\{\vnth,\bR\}\cup\,\big\{(-\infty,a):\;a\in\bR\big\}$
topologies.

A function $f:(X,\tau_1,\tau_2)\to (Y,\gm_1,\gm_2)$ is said to be
$d$-continuous if the induced functions $f:(X,\tau_i)\to
(Y,\gm_i)$ are continuous. Following [8], ``lower (upper)
semicontinuous'' is abbreviated to l.(u.).s.c. and by
Proposition~0.1.4 in [8], a function $f: (X,\tau_1,\tau_2)\to
(I,\om_1',\om_2')$ is       \linebreak    $d$-continuous if and
only if $f:(X,\tau_1,\tau_2)\to (I,\om')$ is $(i,j)$-l.u.s.c.,
where $f:(X,\tau_1,\tau_2)\to (I,\om')$ is $(i,j)$-l.u.s.c. if $f$
is $i$-l.s.c. and $j$-u.s.c.

We shall also use the following double indexation
$$  A_j^i=\big\{x\in A:\;x\;\text{is a $j$-isolated point of}\;A\big\}, $$
where ``$j$'' denotes the belonging to the topology, while ``$i$''
is fixed as the isolation symbol. By [6], for any subset $A$ of a
$\BS$ $(X,\tau_1,\tau_2)$ the $(i,j)$-boundaries (i.e.,
bitopological boundaries) of $A$ are the $p$-closed sets
$(i,j)\dd\Fr A=\tau_i\cl A\cap\tau_j\cl(X\setminus A)$.

Furthermore, in a $\BS$ $(X,\tau_1,\tau_2)$ the topology $\tau_1$
is coupled to the topology $\tau_2$ (briefly, $\tau_1C\tau_2)$ if
$\tau_1\cl U\sbsq\tau_2\cl U$ for any set $U\in\tau_1$ [27],
$\tau_1$ is near $\tau_2$ (briefly, $\tau_1N\tau_2)$ if $\tau_1\cl
U\sbsq\tau_2\cl U$ for any set $U\in\tau_2$ [8], and $\tau_1$ and
$\tau_2$ are $S$-related (briefly, $\tau_1S\tau_2)$ if
$$  \tau_1\nt A\sbs\tau_1\cl\tau_2\nt A\wedge
            \tau_2\nt A\sbs\tau_2\cl\tau_1\nt A     $$
for any subset $A\sbs X$ [25].

In the sequel it will be assumed that
\begin{eqnarray*}
    &\ds \big(\tau_1C\tau_2\wedge\tau_1\sbs\tau_2\big)\llra
                \tau_1<_C\tau_2, \\
    &\ds \big(\tau_1N\tau_2\wedge\tau_1\!\sbs\!\tau_2\big)\!\llra\!
            \tau_1<_N\tau_2 \;\;\text{and}\;\;
        \big(\tau_1S\tau_2\wedge\tau_1\!\sbs\!\tau_2\big)\!\llra\!
                \tau_1<_S\tau_2.
\end{eqnarray*}

By [8], $i\dd\,\cD(X)=\big\{A\sbs X:\;\tau_i\cl A=X\big\}$,
$i\dd\Bd(X)=\co i\dd\,\cD(X)=\{A\sbs X:\;\tau_i\nt A=\vnth\}$ and
$$  (i,j)\dd\,\cC\cD(X)=\big\{A\sbs X:\;A=\tau_i\cl\tau_j\nt A\big\}. $$

\begin{definition}{1.1}
let $(X,\tau_1,\tau_2)$ be $\BS$. Then
\begin{enumerate}
\item[(1)] $(X,\tau_1,\tau_2)$ is $\RR\dd\,p\,\dd \TT_1$ (i.e.,
$p\,\dd\TT_1$ in the sense of Reilly) if it is $d\dd\TT_1$ [20].

\item[(2)] $(X,\tau_1,\tau_2)$ is $(i,j)$-regular if for each
point $x\in X$ and each       \linebreak        $i$-closed set
$F$, $x\,\overline{\in}\,F$, there exist an $i$-open set $U$ and a
$j$-open set $V$ such that $x\in U$, $F\subset V$ and $U\cap
V=\varnothing$ [13] and by (2) of Proposition~0.1.7 in [8],
$(X,\tau_1,\tau_2)$ is $(i,j)$-regular if and only if for each
point $x\in X$ and any neighborhood $U(x)\in\tau_i$ there exists a
neighborhood $V(x)\in\tau_i$ such that $\tau_j\cl V(x)\subset
U(x)$.

\item[(3)] $(X,\tau_1,\tau_2)$ is $(i,j)$-Tychonoff if it is
$\RR\dd\,p\,\dd\TT_1$ and $(i,j)$-completely regular, that is, for
every $i$-closed set $F\sbs X$ and any point $x\in X$,
$x\,\ol{\in}\,F$, there is an $(i.j)$-l.u.s.c. function
$f:(X,\tau_1,\tau_2)\to (I,\om')$ such that $f(F)=0$ and $f(x)=1$
[14].

\item[(4)] $(X,\tau_1,\tau_2)$ is $p$-normal if for every pair of
disjoint sets $A$, $B$ in $X$, where $A$ is $1$-closed and $B$ is
$2$-closed, there exist a 2-open set $U$ and a $1$-open set $V$
such that $A\subset U$, $B\subset V$ and $U\cap V=\varnothing$
[13] and by (4) of Proposition~0.1.7 in [8], $(X,\tau_1,\tau_2)$
is $p$-normal if and only if for each $2$-closed $(1$-closed$)$
set $F$ and each $1$-open $(2$-open$)$ set $U$ with $F\subset U$,
there exists a $1$-open $(2$-open$)$ set $V$ such that $F\subset
V\subset\tau_2\cl V\subset U$ $(F\subset V\subset\tau_1\cl
V\subset U)$.

\item[(5)] $(X,\tau_1,\tau_2)$ is hereditarily $p$-normal if every
one of its $\BsS$ is     \linebreak        $p$-normal [6] and by
Theorem~0.2.2 in [8], $(X,\tau_1,\tau_2)$ is hereditarily
$p$-normal if and only if whenever $A,B\sbs X$, $(\tau_1\cl A\cap
B)\cup (A\cap\tau_2\cl B)=\vnth$ there are disjoint sets
$U\in\tau_2$, $V\in\tau_1$ such that $A\sbs U$ and $B\sbs V$.

\item[(6)] $(X,\tau_1,\tau_2)$ is $(i,j)\dd\RR\RR$-paracompact if
for each $i$-open covering $\boldsymbol{\cV}$ of $X$ there is an
$i$-open covering $\boldsymbol{\cU}$ of $X$ which refines
$\boldsymbol{\cV}$ and is $j$-locally finite at each point of $X$
[18].

\item[(7)] A $\BS$ $(X,\tau_1,\tau_2)$ is $p$-connected if $X$
cannot be expressed as the union of two nonempty sets $A$ and $B$
such that $(A\cap \tau_2\cl B)\cup(\tau_1\cl A\cap B)=\vnth$ [17]
and by (c) of Theorem~A and Theorem~C in [17], respectively, the
$p$-connectedness of a $\BS$ $(X,\tau_1,\tau_2)$ is equivalent to
each of the following two conditions: $X$ contains no nonempty
subset which is both $1$-open and $2$-closed (hence, none which is
$1$-closed and $2$-open), and every $d$-continuous function $f:
(X,\tau_1,\tau_2)\to (R,\om_1,\om_2)$ has the Darboux property,
that is, if $f(x_1)<c<f(x_2)$, then there is $x\in X$ such that
$f(x)=c$.
\end{enumerate}
\end{definition}

In the first part of the work (Sections 2--6) a special attention
is given to relative separation axioms and relative connectedness,
in particular, many relative versions of $p\,\dd\TT_0$,
$p\,\dd\TT_1$, $p\,\dd\TT_2$, $(i,j)$- and $p$-regularities,
$(i,j)$- and $p$-complete regularities, $p$-real normality and
$p$-normality are discussed. Moreover, relative properties of
$(i,j)$- and $p$-compactness types, including relative versions of
$(i,j)$- and $p$-paracompactness, $(i,j)$- and
$p$-Lin\-de\-l\"{o}fness, $(i,j)$- and  $p$-pseudocompactness are
also introduced and investigated. The second part (Sections 7--12)
is devoted, on the one hand, to relative bitopological inductive
and covering dimension functions and, on the other hand, to
relative versions of Baire spaces for both the topological and the
bitopological case.

At the end, note that relative (bi)topological properties play a
special role not only in the development of respective theories,
but also in the strengthening of the previously known results.

\vskip+0.7cm
\section*{\textbf{2. Relative Separation Axioms}}
\vskip+0.3cm

First of all, note that topological versions of relative
bitopological properties, introduced and studied in Sections 2--6,
are considered, in particular, in [2]--[4].

Below, the letters ``$W$'' and ``$S$'' abbreviate, respectively,
the words ``weakly'' and ``strongly''.

\begin{definition}{2.1}
Let $(Y,\tau_1',\tau_2')$ be a $\BsS$ of a $\BS$
$(X,\tau_1,\tau_2)$. Then
\begin{enumerate}
\item[(1)] $Y$ is $p\,\dd\TT_0$ in $X$ if for every pair of
distinct point of $Y$ there exists a $1$-neighborhood or a
$2$-neighborhood in $X$ of one point, not containing the other.

\item[(2)] $Y$ is $\WS\dd\TT_0$ in $X$ if for each pair of
distinct points $y\in Y$, $x\in X$ there exists a $1$-neighborhood
or a $2$-neighborhood of $y$ in $Y$, not containing $x$, or there
exists a $1$- or $2$-neighborhood of $x$ in $X$, not
con\-ta\-i\-ning~$y$.

\item[(3)] $Y$ is $\Wp\dd\RR_0$ in $X$ of for every set
$U\in\tau_1\setminus\{\vnth\}$ it follows from $x\in U\cap Y$ that
$\tau_2'\cl\{x\}\sbs U$ and for every set
$V\in\tau_2\setminus\{\vnth\}$ it follows from $y\in V\cap Y$ that
$\tau_1'\cl\{y\}\sbs V$.

\item[(4)] $Y$ is $p\,\dd\RR_0$ in $X$ of for every set
$U\in\tau_1\setminus\{\vnth\}$ it follows from $x\in U\cap Y$ that
$\tau_2\cl\{x\}\sbs U$ and for every set
$V\in\tau_2\setminus\{\vnth\}$ it follows from $y\in V\cap Y$ that
$\tau_1\cl\{y\}\sbs V$.

\item[(5)] $Y$ is $\Wp\dd\TT_1$ in $X$ if for any pair of distinct
points of $Y$ one point has a $1$-neighborhood in $X$ not
containing the other, while the other point has a $2$-neighborhood
in $X$, not containing the first.

\item[(6)] $Y$ is $p\,\dd\TT_1$ in $X$ if for every pair of
distinct points $x,\,y\in Y$ there exists a $1$- or a
$2$-neighborhood of $x$ in $X$, not con\-ta\-i\-ning~$y$.

\item[(7)] $Y$ is $\Sp\dd\TT_1$ in $X$ if for every pair of
distinct points $x,\,y\in Y$ there exists a $1$- and a
$2$-neighborhood of $x$ in $X$, not con\-ta\-i\-ning~$y$.

\item[(8)] $Y$ is $(i,j)\dd\WS\dd\TT_1$ in $X$ of for every pair
of distinct points $y\in Y$, $x\in X$ there exist an
$i$-neighborhood of $y$ in $Y$, not containing $x$, and a
$j$-neighborhood of $x$ in $X$, not con\-ta\-i\-ning~$y$.

\item[(9)] $Y$ is $(i,j)\dd\WW\dd\RR_1$ in $X$ if for each pair of
distinct points $x,\,y\in Y$ such that
$x\,\ol{\in}\,\tau_i\cl\{y\}$ there are disjoint neighborhoods
$U(x)\in\tau_i$ and $U(y)\in\tau_j$.

\item[(10)] $Y$ is $(i,j)\dd\RR_1$ in $X$ if for each pair of
distinct points $x,\,y\in Y$ such that
$x\,\ol{\in}\,\tau_i'\cl\{y\}$ there are disjoint neighborhoods
$U(x)\in\tau_i$ and $U(y)\in\tau_j$.

\item[(11)] $Y$ is $(i,j)\dd\WS\dd\RR_1$ in $X$ if for each pair
of distinct points $y\in Y$, $x\in X$ such that
$y\,\ol{\in}\,\tau_i\cl\{x\}$ there are disjoint neighborhoods
$U(y)\in\tau_i'$ and $U(x)\in\tau_j$.

\item[(12)] $Y$ is $(i,j)\dd\SW\dd\RR_1$ in $X$ if for each pair
of distinct points $y\in Y$, $x\in X$ such that
$x\,\ol{\in}\,\tau_i'\cl\{y\}$ there are disjoint neighborhoods
$U(y)\in\tau_j'$ and $U(x)\in\tau_i$.

\item[(13)] $Y$ is $\Wp\dd\RR_1$ in $X$ if for each pair of
distinct points $x,\,y\in Y$ such that
$\tau_1\cl\{x\}\neq\tau_2\cl\{y\}$ there are disjoint
neighborhoods $U(x)\in\tau_2$ and $U(y)\in\tau_1$.

\item[(14)] $Y$ is $\Sp\dd\RR_1$ in $X$ if for each pair of
distinct points $x,\,y\in Y$ such that
$\tau_1'\cl\{x\}\neq\tau_2'\cl\{y\}$ there are disjoint
neighborhoods $U(x)\in\tau_2$ and $U(y)\in\tau_1$.

\item[(15)] $Y$ is $\Wp\dd\TT_2$ in $X$ if for each pair of
distinct points $x,\,y\in Y$ there are disjoint sets $U\in\tau_1$,
$V\in\tau_2$ such that either $x\in U$, $y\in V$ or $x\in V$,
$y\in U$.

\item[(16)] $Y$ is $p\,\dd\TT_2$ in $X$ if for each pair of
distinct points $x,\,y\in Y$ there are disjoint sets $U\in\tau_1$,
$V\in\tau_2$ such that $x\in U$ and $y\in V$.

\item[(17)] $Y$ is $(i,j)\dd\WS\dd\TT_2$ in $X$ if for each pair
of distinct points $y\in Y$, $\!x\in\!X$ there are disjoint
neighborhoods $U(y)\!\in\!\tau_i'$ and $U(x)\!\in~\!\!\tau_j$.

\item[(18)] $Y$ is $(i,j)\dd\sS\dd\TT_2$ in $X$ if for each pair
of distinct points $y\in Y$, $x\!\in\!X$ there are disjoint
neighborhoods $U(y)\!\in\!\tau_i$ and $U(x)\!\in~\!\!\tau_j$.
\end{enumerate}
\end{definition}

Furthermore, we pass to the bitopological modifications of various
types of relative regularity, relative complete regularity and
relative normality.

\begin{definition}{2.2}
Let $(Y,\tau_1',\tau_2')$ be a $\BsS$ of a $\BS$
$(X,\tau_1,\tau_2)$. Then
\begin{enumerate}
\item[(1)] $Y$ is $(i,j)\dd\WS$-quasi regular in $X$ if for $x\in
Y$, $F\in\co\tau_i$ and $x\,\ol{\in}\,F$, there are disjoint sets
$U\in\tau_i'$, $V\in\tau_j'$ such that $x\in U$ and $F\cap Y\sbs
V$.

\item[(2)] $Y$ is $(i,j)$-regular in $X$ if for $x\in Y$,
$F\in\co\tau_i$ and $x\,\ol{\in}\,F$, there are disjoint sets
$U\in\tau_i$, $V\in\tau_j$ such that $x\in U$ and $F\cap Y\sbs V$.

\item[(3)] $Y$ is $(i,j)$-strongly regular in $X$ if for $x\in Y$,
$F\in\co\tau_i'$ and $x\,\ol{\in}\,F$, there are disjoint sets
$U\in\tau_i$, $V\in\tau_j$ such that $x\in U$ and $F\sbs V$.

\item[(4)] $Y$ is $(i,j)$-superregular in $X$ if for $x\in Y$,
$F\in\co\tau_i$ and $x\,\ol{\in}\,F$, there are disjoint sets
$U\in\tau_i$, $V\in\tau_j$ such that $x\in U$ and $F\sbs V$.

\item[(5)] $Y$ is $(i,j)\dd\WS$-regular in $X$ if for $x\in Y$,
$F\in\co\tau_i$ and $x\,\ol{\in}\,F$, there are disjoint sets
$U\in\tau_i'$, $V\in\tau_j$ such that $x\in U$ and $F\cap Y\sbs
V$.

\item[(6)] $Y$ is $(i,j)\dd\WS$-superregular in $X$ if for $x\in
Y$, $F\in\co\tau_i$ and $x\,\ol{\in}\,F$, there are disjoint sets
$U\in\tau_i'$, $V\in\tau_j$ such that $x\in U$ and $F\sbs V$.

\item[(7)] $Y$ is $(i,j)$-free regular in $X$ if for $x\in X$,
$F\in\co\tau_i$ and $x\,\ol{\in}\,F$, there are disjoint sets
$U\in\tau_i$, $V\in\tau_j$ such that $x\in U$ and $F\cap Y\sbs V$.
\end{enumerate}
\end{definition}

It is obvious that taking into account (2) of Definition~1.1, for
a $\BsS$ $(Y,\tau_1',\tau_2')$ of a $\BS$ $(X,\tau_1,\tau_2)$ the
following implications hold: {\small
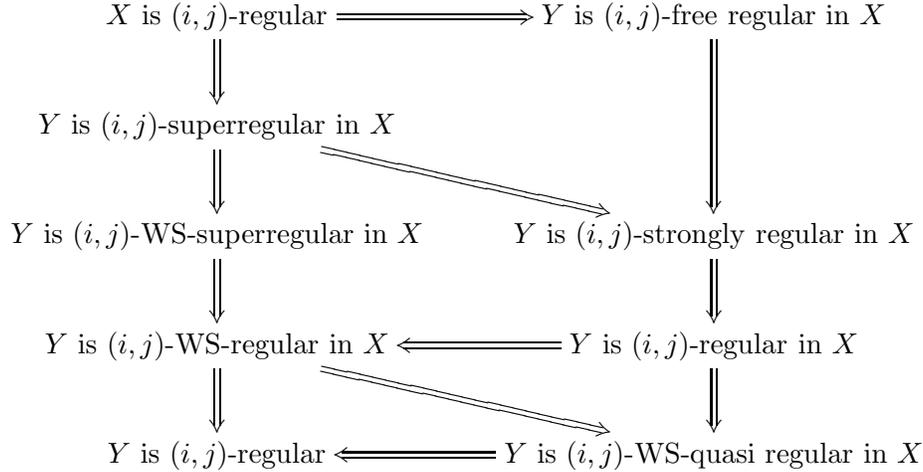
$$  \xymatrix{
    \text{$X$ is $(i,j)$-regular} \ar@{=>}[d] \ar@{=>}[r]
        & \text{$Y$ is $(i,j)$-free regular in $X$} \ar@{=>}[dd] \\
    \text{$Y$ is $(i,j)$-superregular in $X$}\ar@{=>}[d] \ar@{=>}[dr] & \\
    \text{$Y$ is $(i,j)\dd\WS$-superregular in $X$} \ar@{=>}[d] &
        \text{$Y$ is $(i,j)$-strongly regular in $X$} \ar@{=>}[d] \\
    \text{$Y$ is $(i,j)\dd\WS$-regular in $X$}  \ar@{=>}[d] \ar@{=>}[dr]
        & \text{$Y$ is $(i,j)$-regular in $X$} \ar@{=>}[l] \ar@{=>}[d] \\
    \text{$Y$ is $(i,j)$-regular}
        & \text{$Y$ is $(i,j)\dd\WS$-quasi regular in $X$} \ar@{=>}[l] }      $$
}

Using (2) of Definition~1.1 one can reformulate three axioms of
bitopological relative regularity.

\begin{proposition}{2.3}
Let $(Y,\tau_1',\tau_2')$ be a $\BsS$ of a $\BS$
$(X,\tau_1,\tau_2)$. Then
\begin{enumerate}
\item[(1)] $Y$ is $(i,j)$-strongly regular in $X$ if and only if
for each point $x\in Y$ and any neighborhood $U'(x)\in\tau_i'$
there is a neighborhood $V(x)\in\tau_i$ such that $\tau_j\cl
V(x)\cap Y\sbs U'(x)$.

\item[(2)] $Y$ is $(i,j)\dd\WS$-superregular in $X$ if and only if
for each point $x\in Y$ and each neighborhood $U(x)\in\tau_i$
there is a neighborhood $V'(x)\in\tau_i'$ such that $\tau_j\cl
V'(x)\sbs U(x)$.

\item[(3)] $Y$ is $(i,j)$-superregular in $X$ if and only if for
each point $x\!\in~\!\!Y$ and each neighborhood $U(x)\in\tau_i$
there is a neighborhood    \linebreak      $V(x)\in\tau_i$ such
that $\tau_j\cl V(x)\sbs U(x)$.
\end{enumerate}
\end{proposition}

\begin{pf}
(1) Let, first, $Y$ be $(i,j)$-strongly regular in $X$, $x\in Y$
and $U'(x)\in\tau_i'$. Then $x\,\ol{\in}\,A=Y\setminus U'(x)$,
$A\in \co\tau_i'$. Hence, there are $V(x)\in\tau_i$,
$V(A)\in\tau_j$ such that $\tau_j\cl V(x)\cap V(A)=\vnth$,
Clearly,
$$  \tau_j\cl V(x)\sbs X\setminus V(A)\sbs X\setminus A=
        X\setminus (Y\setminus U'(x))=U'(x)\cup (X\setminus Y)    $$
and thus $\tau_j\cl V(x)\cap Y\sbs U'(x)$.

Conversely, let the condition be satisfied, $x\in Y$,
$A\in\co\tau_i'$ and $x\,\ol{\in}\,A$. Then $x\in Y\setminus
A=U'(x)\in\tau_i'$ and hence, there is $V(x)\in\tau_i$ such that
$\tau_j\cl V(x)\cap Y\sbs U'(x)$. Therefore
$$  X\setminus U'(x)\sbs X\setminus(\tau_j\cl V(x)\cap Y)=
        (X\setminus\tau_j\cl V(x))\cup(X\setminus Y)        $$
so that $(X\setminus U'(x))\cap Y\sbs X\setminus\tau_j\cl V(x)$.
Thus $A\sbs X\setminus\tau_j\cl V(x)=V(A)\in\tau_j$ and $V(x)\cap
V(A)=\vnth$.

(2) Let, first, $Y$ be $(i,j)\dd\WS$-superregular in $X$, $x\in Y$
and $U(x)\in\tau_i$. Then $x\,\ol{\in}\,F=X\setminus
U(x)\in\co\tau_i$ and so, there are $V(x)\in\tau_i'$,
$V(F)\in\tau_j$ such that $\tau_j\cl V(x)\cap V(F)=\vnth$.
Therefore,
$$  \tau_j\cl V(x)\sbs X\setminus V(F)\sbs X\setminus F=U(x).  $$

Conversely, let the condition be satisfied, $x\in Y$,
$F\in\co\tau_i$ and $x\,\ol{\in}\,F$. Then $x\in U(x)=X\setminus
F\in\tau_i$ and by condition, there is $V(x)\in\tau_i'$ such that
$\tau_j\cl V(x)\sbs U(x)$. Hence
$$  F=X\setminus U(x)\sbs X\setminus\tau_j\cl V(x)=
        V(F)\in\tau_j       $$
and $V(x)\cap V(F)=\vnth$.

(3) Follows directly from (2) of Definition~1.1, taking into
account that $x\in Y$.

Hence, $Y$ is $(i,j)$-superregular in $X$ if and only if $X$ is
$(i,j)$-regular at each point of $Y$.
\end{pf}
\vskip+0.2cm

It is evident that a $\TsS$ $(Y,\tau')$ of a $\TS$ $(X,\tau)$ is
$\WS$-regular $(\WS$-superregular) in $X$ if for $x\in Y$,
$F\in\co\tau$ and $x\,\ol{\in}\,F$ there are disjoint sets
$U\in\tau'$, $V\in \tau$ such that $x\in U$ and $F\cap Y\sbs
V(F\sbs V)$.

Therefore, the topological version of Proposition~2.3 says that
$(Y,\tau')$ is strongly regular (respectively, $\WS$-superregular,
superregular) in $X$ if and only if for each point $x\in Y$ and
any neighborhood $U'(x)\in\tau'$ (respectively, $U(x)\in\tau)$
there is a neighborhood $V(x)\in\tau$ (respectively,
$V'(x)\in\tau',$ $V(x)\in\tau)$ such that $\tau\cl V(x)\cap Y\sbs
U'(x)$ (respectively, $\tau\cl V'(x)\sbs U(x)$, $\tau\cl V(x)\sbs
U(x))$.

\begin{definition}{2.4}
Let $(Y,\tau_1',\tau_2')$ be a $\BsS$ of a $\BS$
$(X,\tau_1,\tau_2)$. Then
\begin{enumerate}
\item[(1)] $Y$ is $(i,j)$-regular in $X$ from inside if every
$i$-closed in $X$ $\BsS$ of $Y$ is $(i,j)$-regular.

\item[(2)] $Y$ is $(i,j)$-internally regular in $X$ if for $x\in
Y$, $F\in\co\tau_i$, $F\sbs Y$ and $x\,\ol{\in}\,F$ there are
disjoint sets $U\in\tau_i$, $V\in\tau_j$ such that $x\in U$ and
$F\sbs V$.

\item[(3)] $X$ is $(i,j)$-regular on $Y$ if for each set
$F\in\co\tau_i$, $i$-concentrated on $Y$ (that is
$F\sbs\tau_i\cl(F\cap Y)$ [3]), and each point $x\in Y\setminus F$
there are disjoint sets $U\in\tau_i$, $V\in\tau_j$ such that $x\in
U$ and $F\sbs V$.

\item[(4)] $X$ is $(i,j)$-strongly regular on $Y$ if for each
$F\in\co\tau_i$, $i$-con\-cen\-tra\-ted on $Y$, and each point
$x\in X\setminus F$ there are disjoint sets $U\in\tau_i$,
$V\in\tau_j$ such that $x\in U$ and $F\sbs V$.
\end{enumerate}
\end{definition}

\begin{proposition}{2.5}
Let $(Y,\tau_1',\tau_2')$ be a $\BsS$ of a $\BS$
$(X,\tau_1,\tau_2)$. Then
\begin{enumerate}
\item[(1)] If $Y$ is $(i,j)$-regular, then $Y$ is $(i,j)$-regular
from inside in $X$ and so, in every larger $\BS$.

\item[(2)] If $Y$ is $(i,j)$-internally regular in $X$, then $Y$
is $(i,j)$-regular in $X$ from inside.
\end{enumerate}
\end{proposition}

\begin{pf}
(1) The condition is obvious since if $Y$ is $(i,j)$-regular, then
any $\BsS$ of $Y$ is also $(i,j)$-regular.

(2) Let $F\in\co\tau_i$, $F\sbs Y$ and let us prove that
$(F,\tau_1'',\tau_2'')$ is       \linebreak $(i,j)$-regular. If
$x\in F$, $\Phi\in\co\tau_i''$, $x\,\ol{\in}\,\Phi$, then $x\in Y$
and $\Phi\in\co\tau_i$ as $F\in\co\tau_i$. Since $x\in Y$,
$\Phi\sbs Y$, $\Phi\in\co\tau_i$, $x\,\ol{\in}\,\Phi$ and $Y$ is
$(i,j)$-internally regular in $X$, there are $U'\in\tau_i$,
$V'\in\tau_j$ such that $x\in U'$, $\Phi\sbs V'$ and $U'\cap
V'=\vnth$. Let $U=U'\cap F$, $V=V'\cap F$. Then $U\in \tau_i''$,
$V\in\tau_j''$, $x\in U$, $\Phi\sbs U$ and $U\cap V=\vnth$. Thus
$(F,\tau_1'',\tau_2'')$ is $(i,j)$-regular.~\end{pf} \vskip+0.2cm

Note here that a real-valued function $f: (X,\tau_1,\tau_2)\to
(I,\om')$ is $Y\dd\,(i,j)$-l.u.s.c. if it is $(i,j)$-l.u.s.c. at
each point $y\in Y$, where $(Y,\tau_1',\tau_2')$ is a $\BsS$ of a
$\BS$ $(X,\tau_1,\tau_2)$.

\begin{definition}{2.6}
Let $(Y,\tau_1',\tau_2')$ be a $\BsS$ of a $\BS$
$(X,\tau_1,\tau_2)$. Then
\begin{enumerate}
\item[(1)] $Y$ is $(i,j)$-almost completely regular in $X$ if for
each point $x\in Y$ and each set $F\in\co\tau_i'$,
$x\,\ol{\in}\,F$, there is a $Y\dd\,(j,i)$-l.u.s.c. function
$f:(X,\tau_1,\tau_2)\!\to\!(I,\om)$ such that $f(x)\!=\!0$ and
$f(F)\!\sbs~\!\!\!\{1\}$.

\item[(2)] $Y$ is $(i,j)$-completely regular in $X$ if for each
point $x\in Y$ and each set $F\in\co\tau_i$, $x\,\ol{\in}\,F$,
there is a $(j,i)$-l.u.s.c. function $f:(X,\tau_1,\tau_2)\to
(I,\om)$ such that $f(x)=0$ and $f(F\cap Y)\sbs\{1\}$.

\item[(3)] $Y$ is $(i,j)$-strongly completely regular in $X$ if
for each point      \linebreak       $x\in Y$ and each set
$F\in\co\tau_i$, $x\,\ol{\in}\,F$, there is a  $(j,i)$-l.u.s.c.
function $f:(X,\tau_1,\tau_2)\to (I,\om)$ such that $f(x)\!=\!0$
and $f(F)\!\sbs~\!\!\{1\}$.
\end{enumerate}
\end{definition}

Therefore, $Y$ is $(i,j)$-strongly completely regular in $X\lra Y$
is $(i,j)$-completely regular in $X\lra Y$ is $(i,j)$-almost
completely regular in $X$.

\begin{proposition}{2.7}
A $\BsS$ $(Y,\tau_1',\tau_2')$ of a $\BS$ $(X,\tau_1,\tau_2)$ is
$(i,j)$-almost completely regular in $X$ if and only if for each
point $x\in Y$ and any set $U\in\tau_i$, $x\in U$ there is a
$Y\dd\,(i,j)$-l.u.s.c. function $g:(X,\tau_1,\tau_2)\to (I,\om)$
such that $g(x)=1$ and $g(X\setminus U)\sbs\{0\}$.
\end{proposition}

\begin{pf}
If $P=X\setminus U$, then $P\in \co\tau_i$ and $F=P\cap
Y\in\co\tau_i'$. Since $x\,\ol{\in}\,F$ and $Y$ is $(i,j)$-almost
completely regular, there is a $Y\dd\,(j,i)$-l.u.s.c. function
$f:(X,\tau_1,\tau_2)\to (I,\om)$ such that $f(x)=0$ and
$f(F)\sbs\{1\}$. If $g'=1-f$, then $g':(X,\tau_1,\tau_2)\to
(I,\om)$ is $Y\dd\,(i,j)$-l.u.s.c., $g'(x)=1$ and
$g'(F)\sbs\{0\}$. Let $g(x)=g'(x)$ at each point $x\in U$ and
$g(x)=0$ at each point $x\in P=X\setminus U$. Then
$g:(X,\tau_1,\tau_2)\to (I,\om)$ is $Y\dd\,(i,j)$-l.u.s.c.,
$g(x)=1$ and $g(X\setminus U)\sbs\{0\}$.

Conversely, let the condition be satisfied and $x\in Y$,
$F\in\co\tau_i'$, $x\,\ol{\in}\,F$. Then $x\,\ol{\in}\,\tau_i\cl
F=P\in\co\tau_i$ and so $x\in U=X\setminus P$. Hence, by
condition, there is a $Y\dd\,(i,j)$-l.u.s.c. function
$f:(X,\tau_1,\tau_2)\to (I,\om)$ such that $f(x)=1$ and
$f(X\setminus U)\sbs\{0\}$. Clearly $g=1-f$ is
$Y\dd\,(j,i)$-l.u.s.c., $g(x)=0$ and
$$  g(F)\sbs g(P)=1-f(P)=1-f(X\setminus U)=\{1\},       $$
that is, $Y$ is $(i,j)$-almost completely regular in $X$.
\end{pf}

\begin{definition}{2.8}
Let $(Y,\tau_1',\tau_2')$ be a $\BsS$ of a $\BS$
$(X,\tau_1,\tau_2)$. Then
\begin{enumerate}
\item[(1)] $Y$ is $p$-quasi normal in $X$ if for each pair of
disjoint sets       \linebreak        $A\in\co\tau_1$,
$B\in\co\tau_2$ there are disjoint sets $U\in\tau_2'$,
$V\in\tau_1'$ such that $A\cap Y\sbs U$ and $B\cap Y\sbs V$.

\item[(2)] $Y$ is $p$-normal in $X$ if for each pair of disjoint
sets $A\in\co\tau_1$, $B\in\co\tau_2$ there are disjoint sets
$U\in\tau_2$, $V\in\tau_1$ such that $A\cap Y\sbs U$ and $B\cap
Y\sbs V$.

\item[(3)] $Y$ is $p$-strongly normal in $X$ if for each pair of
disjoint sets $A\in\co\tau_1'$, $B\in\co\tau_2'$ there are
disjoint sets $U\in\tau_2$, $V\in\tau_1$ such that $A\sbs U$ and
$B\sbs V$.

\item[(4)] $Y$ is $(i,j)$-supernormal in $X$ if for each pair of
disjoint sets $A\in\co\tau_i'$, $B\in\co\tau_j$ there are disjoint
sets $U\in\tau_j$, $V\in\tau_i$ such that $A\sbs U$ and $B\sbs V$.

\item[(5)] $Y$ is $(i,j)\dd\WS$-normal in $X$ if for each pair of
disjoint sets       \linebreak       $A\in\co\tau_i'$,
$B\in\co\tau_j$ there are disjoint sets $U\in\tau_j'$,
$V\in\tau_i$ such that $A\sbs U$ and $B\cap Y\sbs V$.

\item[(6)] $Y$ is $(i,j)\dd\WS$-supernormal in $X$ if for each
pair of disjoint sets $A\in\co\tau_i'$, $B\in\co\tau_j$ there are
disjoint sets $U\in\tau_j'$, $V\in\tau_i$ such that $A\sbs U$ and
$B\sbs V$ [9].
\end{enumerate}
\end{definition}

In general, the $p$-normality of $(Y,\tau_1',\tau_2')$ in
$(X,\tau_1,\tau_2)$ does not imply the $p$-normality of
$(Y,\tau_1',\tau_2')$. Nevertheless, if $(Y,\tau_1',\tau_2')$ is
$p$-normal in a $\RR\dd\,p\,\dd\TT_1$ $\BS$ $(X,\tau_1,\tau_2)$,
then $(Y,\tau_1',\tau_2')$ is $p$-regular. Indeed, let $x\in Y$,
$F\in\co\tau_i'$ and $x\,\ol{\in}\,F$. Then
$x\,\ol{\in}\,\tau_i\cl F$ and by condition there are disjoint
sets $U(x)\in \tau_i$ and $U(\tau_i\cl F)\in\tau_j$. Clearly,
$U'(x)=U(x)\cap Y\in\tau_i'$ and $U'(F)=U(\tau_i\cl F)\cap
Y\in\tau_j'$ are also disjoint, and hence, $(Y,\tau_1',\tau_2')$
is $p$-regular.

\begin{proposition}{2.9}
Let $(Y,\tau_1',\tau_2')$ be a $\BsS$ of a $\BS$
$(X,\tau_1,\tau_2)$. Then
\begin{enumerate}
\item[(1)] $Y$ is $p$-strongly normal in $X$ if and only if for
each set $F\in\co\tau_1'$ $(F\in\co\tau_2')$ and any neighborhood
$U'(F)\in\tau_2'$ $(U'(F)\in\tau_1')$ there is a neighborhood
$V(F)\in\tau_2$ $(V(F)\in\tau_1)$ such that $\tau_1\cl V(F)\cap
Y\sbs U'(F)$ $(\tau_2\cl V(F)\cap Y\sbs U'(F))$.

\item[(2)] $Y$ is $(i,j)\dd\WS$-supernormal in $X$ if and only if
for each set     \linebreak           $F\in\co\tau_i'$ and any
neighborhood $U(F)\in\tau_j$ there is a neighborhood
$V(F)\in\tau_j'$ such that $\tau_i\cl V(F)\sbs U(F)$.

\item[(3)] $Y$ is $(i,j)$-supernormal in $X$ if and only if for
each set $F\in\co\tau_i'$ and any neighborhood $U(F)\in\tau_j$
there is a neighborhood $V(F)\in\tau_j$ such that $\tau_i\cl
V(F)\sbs U(F)$.
\end{enumerate}
\end{proposition}

\begin{pf}
(1) Let, first, $Y$ be $p$-strongly normal in $X$,
$F\in\co\tau_1'$ and $U'(F)\in\tau_2'$. Then $\Phi=(X\setminus
U(F))\cap Y$, where $U(F)\in\tau_i$, $U(F)\cap Y=U'(F)$,
$\Phi\in\co\tau_2'$ and $F\cap \Phi=\vnth$. Hence there are
$V(F)\in\tau_2$, $V(\Phi)\in\tau_1$ such that $\tau_1\cl V(F)\cap
V(\Phi)=\vnth$ and so
$$  V(F)\sbs \tau_1\cl V(F)\sbs X\setminus
            V(\Phi)\sbs X\setminus\Phi=U(F)\cup(X\setminus Y).    $$
Therefore, $\tau_1\cl V(F)\cap Y\sbs U'(F)$.

Conversely, let the condition be satisfied, $A\in\co\tau_1'$,
$B\in\co\tau_2'$ and $A\cap B=\vnth$. Then $A\sbs Y\setminus
B=U'(A)\in\tau_2'$ and by condition there is $V(A)\in\tau_2$ such
that $\tau_1\cl V(A)\cap Y\sbs U'(A)$. Hence
$$  X\setminus U'(A)\sbs X\setminus (\tau_1\cl V(A)\cap Y)=
        (X\setminus\tau_1\cl V(A))\cup(X\setminus Y)       $$
and so $(X\setminus U'(A))\cap Y\sbs X\setminus\tau_1\cl V(A)$,
i.e., $B=Y\setminus U'(A)\sbs (X\setminus\tau_1\cl
V(A))=V(B)\in\tau_1$ and $V(A)\cap V(B)=\vnth$.

The case in the brackets can be proved by the similar manner.

(2) Let, first, $Y$ be $(i,j)\dd\WS$-supernormal in $X$,
$A\in\co\tau_i'$ and $U(A)\in\tau_j$. Then $A\cap B=\vnth$, where
$B=X\setminus U(A)\in\co\tau_j$. Hence there are $V(A)\in\tau_j'$,
$V(B)\in\tau_i$ such that $\tau_i\cl V(A)\cap V(B)=\vnth$ and thus
$$  V(A)\sbs\tau_i\cl V(A)\sbs X\setminus V(B)\sbs
            X\setminus B=U(A).      $$

Conversely, let the condition be satisfied, $A\in\co\tau_i'$,
$B\in\co\tau_j$ and $A\cap B=\vnth$. Then $A\sbs U(A)=X\setminus
B\in\tau_j$. Hence, by condition, there is $V(A)\in\tau_j'$ such
that $\tau_i\cl V(A)\sbs U(A)$. Therefore,
$$  B=X\setminus U(A)\sbs X\setminus\tau_i\cl V(A)=
            V(B)\in\tau_i       $$
and $V(A)\cap V(B)=\vnth$.

(3) This proof is similar to the proof of (2) and can be omitted.
\end{pf}
\vskip+0.2cm

Now, taking into account the topological version of $(i,j)$-strong
normality, given in [2], and the topological versions of
$(i,j)\dd\WS$-su\-per\-nor\-ma\-li\-ty and $(i,j)$-supernormality,
we come to

\begin{corollary}{2.10}
Let $(Y,\tau')$ be a $\TsS$ of a $\TS$ $(X,\tau)$. Then
\begin{enumerate}
\item[(1)] $Y$ is strongly normal in $X$ if and only if for each
set $F\in\co\tau'$ and any neighborhood $U(F)\in\tau'$ there is a
neighborhood $V(F)\in\tau$ such that $\tau\cl V(F)\cap Y\sbs
U(F)$.

\item[(2)] $Y$ is $\WS$-supernormal in $X$ if and only if for each
set $F\in\co\tau'$ and any neighborhood $U(F)\in\tau$ there is a
neighborhood $V(F)\in\tau'$ such that $\tau\cl V(F)\sbs U(F)$.

\item[(3)] $Y$ is supernormal in $X$ if and only if for each set
$F\in\co\tau'$ and any neighborhood $U(F)\in\tau$ there is a
neighborhood $V(F)\in\tau$ such that $\tau\cl V(F)\sbs U(F)$.
\end{enumerate}
\end{corollary}

\begin{definition}{2.11}
Let $(Y,\tau_1',\tau_2')$ be a $\BsS$ of a $\BS$
$(X,\tau_1,\tau_2)$. Then
\begin{enumerate}
\item[(1)] $Y$ is $p$-normal in $X$ from inside if every
$p$-closed in $X$ $\BsS$ of $Y$ is $p$-normal.

\item[(2)] $Y$ is $p$-internally normal in $X$ if for each pair of
disjoint subsets $A,B\sbs Y$, where $A\in\co\tau_1$ and
$B\in\co\tau_2$, there are disjoint sets $U\in\tau_2$,
$V\in\tau_1$ such that $A\sbs U$ and $B\sbs V$.

\item[(3)] $X$ is $p$-normal on $Y$ if for each pair of disjoint
sets $A\in\co\tau_1$, $B\in\co\tau_2$, where $A\sbs\tau_1\cl(A\cap
Y)$ and $B\sbs \tau_2\cl(B\cap Y)$, there are disjoint sets
$U\in\tau_2$, $V\in\tau_1$ such that $A\sbs U$ and $B\sbs V$.
\end{enumerate}
\end{definition}

It is obvious that taking into account (4) of Definition~1.1, for
a $\BsS$ $(Y,\tau_1',\tau_2')$ of a $\BS$ $(X,\tau_1,\tau_2)$ we
have the following implications {\fontsize{8}{12pt}\selectfont
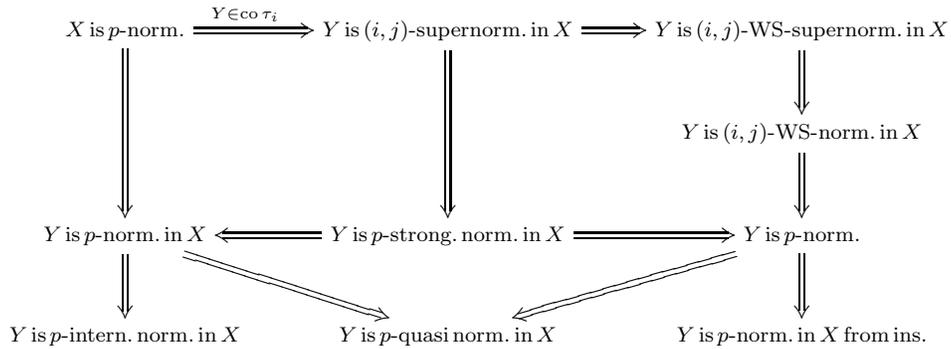
$$  \xymatrix{
  \text{$X$\,is\,$p$-norm.}\ar@{=>}[r]^{Y\in\co\tau_i\ \hskip+1cm }\ar@{=>}[dd]
        & \text{$Y$\,is\,$(i,j)$-supernorm.\,in\,$X$} \ar@{=>}[r]\ar@{=>}[dd]
        & \text{$Y$\,is\,$(i,j)\dd\WS$-supernorm.\,in\,$X$} \ar@{=>}[d] \\
    & & \text{$Y$\,is\,$(i,j)\dd\WS$-norm.\,in\,$X$} \ar@{=>}[d] \\
    \text{$Y$\,is\,$p$-norm.\,in\,$X$} \ar@{=>}[d] \ar@{=>}[dr]
        & \text{$Y$\,is\,$p$-strong.\,norm.\,in\,$X$} \ar@{=>}[r] \ar@{=>}[l]
        & \text{$Y$\,is\,$p$-norm.} \ar@{=>}[d] \ar@{=>}[dl] \\
    \text{$Y$\,is\,$p$-intern.\,norm.\,in\,$X$}
        & \text{$Y$\,is\,$p$-quasi\,norm.\,in\,$X$}
        & \text{$Y$\,is\,$p$-norm.\,in\,$X$\,from\,ins.}}   $$ }

Therefore, if $Y$ is $p$-normal, then $Y$ is $p$-quasi normal in
any larger $\BS$, but $Y$ is not $p$-normal in a larger $\BS$.
Besides, any $\BsS$ $Y$ of a hereditarily $p$-normal $\BS$ $X$ is
$p$-strongly normal in $X$ and hence, $Y$ is $p$-normal in $X$,
$Y$ is $p$-quasi normal in $X$, $Y$ is $p$-internally normal in
$X$ and $Y$ is $p$-normal in $X$ from inside.

Note also here that if $X$ is $\RR\dd\,p\,\dd\TT_1$ and $p$-normal
on $Y$, then $X$ is $p$-re\-gu\-lar on $Y$, but $X$ is not
$(1,2)$-stron\-gly regular on $Y$ nor       \linebreak
      $(2,1)$-stron\-gly regular on $Y$. Indeed, if $x\in (X\setminus
F)\cup(X\setminus Y)$, then $\{x\}$ is not $1$-concentrated on $Y$
nor $2$-concentrated on $Y$.

\begin{remark}{2.12}
First, note that if $F\in\co\tau_i$ in a $\BS$ $(X,\tau_1,\tau_2)$
and $Y\sbs X$, then $F$ is $i$-concentrated on $Y$ if and only if
there is a subset $A\sbs Y$ such that $F=\tau_i\cl A$ [3].
Similarly, it is easily to see that if $A\in p\,\dd\,\Cl(X)$, then
$A$ is $p$-concentrated on $Y$, that is,
$$  A\sbs p\,\dd\cl(A\cap Y)=\tau_1\cl(A\cap Y)\cap\tau_2\cl(A\cap Y)   $$
if and only if there is a subset $B\sbs Y$ such that $A=p\,\dd\cl
B$. Therefore, if $A$ is $i$-concentrated $(p$--concentrated) on
an $i$-closed (a $p$-closed) subset $Y\sbs X$, then $A\sbs Y$.
Indeed, $A\sbs \tau_i\cl(A\cap Y)$ $(A\sbs p\,\dd\cl(A\cap Y))$
implies that $A\sbs\tau_i\cl A\cap Y$ $(A\sbs p\,\dd\cl A\cap Y)$
and so $A\sbs Y$.

Furthermore, if $(Y,\tau_1',\tau_2')$ is a $\BsS$ of a $\BS$
$(X,\tau_1,\tau_2)$ and $Y\in p\,\dd\,\Cl(X)$, then
$p\,\dd\,\Cl(Y)\sbs p\,\dd\,\Cl(X)$. Hence, for $T\in
p\,\dd\,\Cl(X)$ and for any set $A\sbs T\cap Y$ we have $A\in
p\,\dd\,\Cl(T\cap Y)\llra A\in p\,\dd\,\Cl(Y)$.

Indeed, let $A=\tau_1'\cl A\cap\tau_2'\cl A$ and $Y=\tau_1\cl
Y\cap \tau_2\cl Y$. Then
\begin{eqnarray*}
    &\ds A=(\tau_1\cl A\cap Y)\cap(\tau_2\cl A\cap Y)= \\
    &\ds =\big(\tau_1\cl A\cap\tau_1\cl Y\cap \tau_2\cl Y\big)\cap
        \big(\tau_2\cl A\cap\tau_1\cl Y\cap \tau_2\cl Y\big)= \\
    &\ds =\tau_1\cl A\cap\tau_2\cl A,
\end{eqnarray*}
i.e., $p\,\dd\,\Cl(Y)\sbs p\,\dd\,\Cl(X)$.

Now, let $T\in p\,\dd\,\Cl(X)$. First, suppose that $A\in
p\,\dd\,\Cl(T\cap Y)$. Since $T\in p\,\dd\,\Cl(X)$, the set $T\cap
Y\in p\,\dd\,\Cl(Y)$. Hence, $A\in p\,\dd\,\Cl(T\cap Y)$ and
$T\cap Y\in p\,\dd\,\Cl(Y)$ imply that $A\in p\,\dd\,\Cl(Y)$.
Conversely, if $A\in p\,\dd\,\Cl(Y)$, then $A\sbs T\cap Y$ implies
that $A\in p\,\dd\,\Cl(T\cap Y)$.
\end{remark}

\begin{proposition}{2.13}
Let $(Y,\tau_1',\tau_2')$ be a $\BsS$ of a $\BS$
$(X,\tau_1,\tau_2)$. Then
\begin{enumerate}
\item[(1)] If $Y$ is $p$-normal, then $Y$ is $p$-normal in $X$
from inside and so it is $p$-normal from inside in every larger
$\BS$.

\item[(2)] If $Y\in\co\tau_1\cap\co\tau_2$ and $Y$ is
$p$-internally normal in $X$, then $Y$ is $p$-normal in $X$ from
inside.

\item[(3)] $X$ is $p$-normal on $Y$ if and only if for any pair of
subsets $A,B\!\sbs~\!\!Y$ such that $\tau_1\cl A\cap \tau_2\cl
B=\vnth$ there are disjoint sets $U\in\tau_2$, $V\in\tau_1$ such
that $\tau_1\cl A\sbs U$ and $\tau_2\cl B\sbs V$.
\end{enumerate}
\end{proposition}

\begin{pf}
(1) Let $F\in p\,\dd\,\Cl(X)$ and $F\sbs Y$. Then $F\in
p\,\dd\,\Cl(Y)$ and since $Y$ is $p$-normal, by Corollary~3.2.6 in
[8], $F$ is also $p$-normal.

(2) Let $Y\sbs X$, $F\in p\,\dd\,\Cl(X)$ and $F\sbs Y$. Then $F\in
p\,\dd\,\Cl(Y)$ and for $A\in\co\tau_1''$, $B\in\co\tau_2''$ in
$(F,\tau_1'',\tau_2'')$, $A\cap B=\vnth$, by Lemma~3.2.5 in [8],
there are $A'\in\co\tau_1'$, $B'\in\co\tau_2'$ such that $A'\cap
B'=\vnth$, $A'\cap F=A$ and $B'\cap F=B$. Since
$Y\in\co\tau_1\cap\co\tau_2$, it is obvious that $A'\in\co\tau_1$,
$B'\in\co\tau_2$ and hence, by condition, there are disjoint sets
$U'\in\tau_2$, $V'\in\tau_1$ such that $A'\sbs U'$ and $B'\sbs
V'$. Now, if $U=U'\cap F$, $V=V'\cap F$, then $U\in\tau_2''$,
$V\in\tau_1''$, $A\sbs U$ and $B\sbs V$. Thus
$(F,\tau_1'',\tau_2'')$ is $p$-normal.

(3) Let, first, $X$ be $p$-normal on $Y$, $A,B\sbs Y$ and
$\tau_1\cl A\cap\tau_2\cl B=\vnth$. Then, by Remark~2.12,
$F_1=\tau_1\cl A\in\co\tau_1$, $F_1$ is $1$-concentrated on $Y$,
$F_2=\tau_2\cl B\in\co\tau_2$ and $F_2$ is $2$-concentrated on
$Y$. Since $X$ is $p$-normal on $Y$, there are disjoint sets
$U\in\tau_2$, $V\in\tau_1$ such that $F_1=\tau_1\cl A\sbs U$ and
$F_2=\tau_2\cl B\sbs V$.

Conversely, let the condition be satisfied, $F_1\in\co\tau_1$ be
$1$-concentrated on $Y$, $F_2\in\co\tau_2$ be $2$-concentrated on
$Y$ and $F_1\cap F_2=\vnth$. Then by Remark~2.12 there are
$A,B\sbs Y$ such that $F_1=\tau_1\cl A$ and $F_2=\tau_2\cl B$.
Hence, it remains to use the condition.
\end{pf}

\begin{definition}{2.14}
Let $(Y,\tau_1',\tau_2')$ be a $\BsS$ of a $\BS$
$(X,\tau_1,\tau_2)$. Then
\begin{enumerate}
\item[(1)] $Y$ is $p$-weakly realnormal in $X$ if for every pair
of disjoint sets $A\in\co\tau_1$ and $B\in\co\tau_2$ there is a
$(1,2)$-l.u.s.c. function        \linebreak
$f:(Y,\tau_1',\tau_2')\!\to\!(I,\om)$ such that $f(A\cap
Y)\!\sbs\!\{0\}$ and $f(B\cap Y)\!\sbs~\!\!\{1\}$.

\item[(2)] $Y$ is $p$-realnormal in $X$ if for every pair of
nonempty disjoint sets $A\in\co\tau_1$ and $B\in\co\tau_2$ there
is a $Y\dd\,(1,2)$-l.u.s.c. function
$f:(X,\tau_1,\tau_2)\to(I,\om)$ such that $f(A)=\{0\}$ and
$f(B)=\{1\}$.

\item[(3)] $Y$ is $p$-strongly realnormal in $X$ if for every pair
of nonempty disjoint sets $A\in\co\tau_1'$ and $B\in\co\tau_2'$
there is a $Y\dd\,(1,2)$-l.u.s.c. function
$f:(X,\tau_1,\tau_2)\to(I,\om)$ such that $f(A)=\{0\}$ and
$f(B)=\{1\}$.

\item[(4)] $Y$ is $p$-super strongly realnormal in $X$ if for
every pair of non\-empty disjoint sets $A\in\co\tau_1'$ and
$B\in\co\tau_2'$ there is a $(1,2)$-l.u.s.c. function
$f:(X,\tau_1,\tau_2)\to(I,\om)$ such that $f(A)=\{0\}$ and
$f(B)=\{1\}$.

\item[(5)] $X$ is $p$-realnormal on $Y$ if for every pair of
disjoint sets $A\in\co\tau_1$ and $B\in\co\tau_2$ there is a
$(1,2)$-l.u.s.c. function $f:(X,\tau_1,\tau_2)\to(I,\om)$ such
that $f(A\cap Y)\sbs\{0\}$ and $f(B\cap Y)\sbs\{1\}$.
\end{enumerate}
\end{definition}

\begin{proposition}{2.15}
For a $\BsS$ $(Y,\tau_1',\tau_2')$ of a $\BS$ $(X,\tau_1,\tau_2)$
the following conditions are equivalent:
\begin{enumerate}
\item[(1)] $Y$ is $p$-realnormal in $X$.

\item[(2)] For every pair of nonempty subsets $A,B\sbs Y$, where
$\tau_1\cl A\cap\tau_2\cl B=\vnth$ there is a
$Y\dd\,(1,2)$-l.u.s.c. function $f:(X,\tau_1,\tau_2)\to(I,\om)$
such that $f(A)\sbs\{0\}$ and $f(B)\sbs\{1\}$.

\item[(3)] For every pair of nonempty disjoint sets
$A\in\co\tau_1$ and $B\in\co\tau_2$ there is a
$Y\dd\,(1,2)$-l.u.s.c. function $f:(X,\tau_1,\tau_2)\to(I,\om)$
such that $f(A\cap Y)\sbs\{0\}$ and $f(B\cap Y)\sbs\{1\}$.
\end{enumerate}
\end{proposition}

\begin{pf}
$(1)\lra (2)$ is obvious.

$(2)\lra (3)$. Let $A\in\co\tau_1\setminus\{\vnth\}$,
$B\in\co\tau_2\setminus\{\vnth\}$ and $A\cap B=\vnth$. Then
$\tau_1\cl(A\cap Y)\sbs \tau_1\cl A=A$, $\tau_2\cl(B\cap Y)\sbs
\tau_2\cl B=B        $ and by (2), there is a
$Y\dd\,(1,2)$-l.u.s.c. function $f:(X,\tau_1,\tau_2)\to(I,\om)$
such that $f(A\cap Y)\sbs\{0\}$ and $f(B\cap Y)\sbs\{1\}$.

$(3)\lra (1)$. Let $A\in\co\tau_1\setminus\{\vnth\}$,
$B\in\co\tau_2\setminus\{\vnth\}$ and $A\cap B=\vnth$. Then there
is a $Y\dd\,(1,2)$-l.u.s.c. function
$\vf:(X,\tau_1,\tau_2)\to(I,\om)$ such that $\vf(A\cap
Y)\sbs\{0\}$ and $\vf(B\cap Y)\sbs\{1\}$. Let us define a new
function $f:(X,\tau_1,\tau_2)\to(I,\om)$ as follows: $f(A)=\{0\}$,
$f(B)=\{1\}$ and $f\big|_{X\setminus(A\cup
B)}=\vf\big|_{X\setminus(A\cup B)}$. Then, it is clear, that $f$
is $Y\dd\,(1,2)$-l.u.s.c., and thus $(3)\lra(1)$.
\end{pf}
\vskip+0.2cm

Clearly, take place the following implications:
$$  
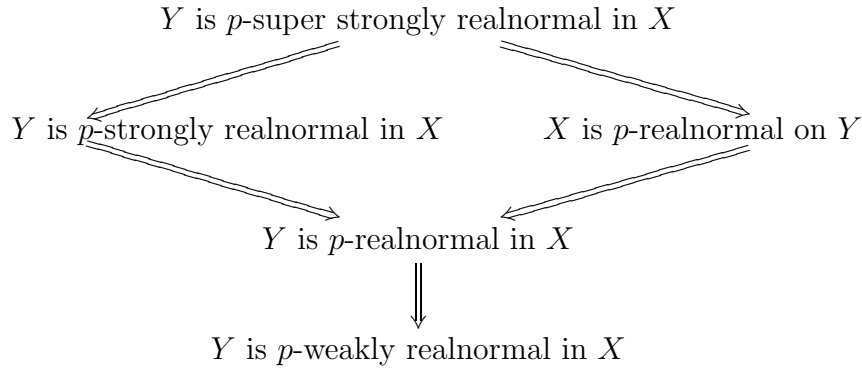
\xymatrix{
    & \text{$Y$ is $p$-super strongly realnormal in $X$}
            \ar@{=>}[dl] \ar@{=>}[dr] & \\
    \text{$Y$ is $p$-strongly realnormal in $X$}\ \hskip-5cm \ar@{=>}[dr] & &
        \ \hskip-3cm \text{$X$ is $p$-realnormal on $Y$} \ar@{=>}[dl] \\
    & \text{$Y$ is $p$-realnormal in $X$} \ar@{=>}[d] & \\
    & \text{$Y$ is $p$-weakly realnormal in $X$} }     $$

The rest of the section is devoted to some special notions and
results, which together with their applications have an
independent interest too.

First, let us prove an elementary, but very important

\begin{lemma}{2.16}
If $Y\in i\dd\,\cD(X)$ in a $\BS$ $(X,\tau_1<_S\tau_2)$, then for
any set $U\in\tau_i$ we have $\tau_j\cl U=\tau_j\cl(U\cap Y)$.

Moreover, if $Y\in d\dd\,\cD(X)$, then for each pair of disjoint
sets $U'\in\tau_2'$, $V'\in\tau_1'$, there is a pair of disjoint
sets $U\in\tau_2$, $V\in\tau_1$ such that $U\cap Y=U'$ and $V\cap
Y=V'$.
\end{lemma}

\begin{pf}
Since $\tau_1\sbs\tau_2$, it is obvious that for $Y\in
2\dd\,\cD(X)$ we have $\tau_1\cl U=\tau_1(U\cap Y)$ for any set
$U\in\tau_2$.

If $Y\in 1\dd\,\cD(X)$, $U\in\tau_1$ and $x\in\tau_2\cl U$, then
for each $U(x)\in\tau_2$ we have
$U(x)\cap\,U\in\tau_2\setminus\vnth$. Since $\tau_1S\tau_2$, by
(2) of Theorem~2.1.5 in [8], $\tau_1\nt(U(x)\cap\,U)\neq\vnth$ and
so $\tau_1\nt\big(U(x)\cap\,U\big)\cap Y\neq\vnth$ as $Y\in
1\dd\,\cD(X)$. Hence, $U(x)\cap (U\cap Y)\neq\vnth$ and thus,
$x\in\tau_2\cl(U\cap Y)$ so that $\tau_2\cl U=\tau_2\cl(U\cap Y)$
for any $U\in\tau_1$.

Now, let $U'\in\tau_2'$, $V'\in\tau_1'$ and $U'\cap V'=\vnth$. Let
$U''\in\tau_2$, $V''\in\tau_1$ and $U''\cap Y=U'$, $V''\cap Y=V'$.
By the first part,
\begin{eqnarray*}
    \vnth &&\hskip-0.6cm =U'\cap\tau_2'\cl V'=
            U'\cap(\tau_2\cl V'\cap Y)= \\
    &&\hskip-0.6cm =U'\cap(\tau_2\cl(V''\cap Y)\cap Y)=
            U'\cap(\tau_2\cl V''\cap Y)
\end{eqnarray*}
and similarly, $\vnth=V'\cap(\tau_1\cl U''\cap Y)$. Let
$U=U''\setminus\tau_2\cl V''$ and $V=V''\setminus\tau_1\cl U''$.
Then $U\in\tau_2$, $V\in\tau_1$, $U\cap V=\vnth$ and $U\cap Y=U'$,
$V\cap Y'=V'$.
\end{pf}
\vskip+0.2cm

Note also here, that if $\vnth\neq A\in(i,j)\dd\,\cC\cD(X)$ in
$(X,\tau_1<_S\tau_2)$ and $Y\in j\dd\,\cD(X)$, then $A$ is
$i$-concentrated on $Y$. Indeed, let $A=\tau_i\cl\tau_j\nt
A\neq\vnth$. Then $\tau_j\nt A\neq\vnth$ and since $Y\in
j\dd\,\cD(X)$, by Lemma~2.16, $\tau_i\cl\tau_j\nt
A\!=\!\tau_i\cl\big(\tau_j\nt A\cap Y\big)$, i.e.,
$A\!=\!\tau_i\cl(\tau_j\nt A\cap~\!\!Y)$. Hence
$$  A\sbs\tau_i\cl\big(\tau_i\cl\tau_j\nt A\cap Y\big)=
                    \tau_i\cl(A\cap Y).     $$

\begin{definition}{2.17}
A $\BS$ $(X,\tau_1,\tau_2)$ is $(d,p)$-densely normal if there is
a $\BsS$ $(Y,\tau_1',\tau_2')$ of $X$ such that $Y\in
d\dd\,\cD(X)$ and $X$ is $p$-normal on $Y$.

Following [22], a $\BS$ $(X,\tau_1,\tau_2)$ is $p$-middly normal
if for every pair of disjoint sets $A\in(1,2)\dd\,\cC\cD(X)$ and
$B\in(2,1)\dd\,\cC\cD(X)$ there are disjoint sets $U\in\tau_2$,
$V\in\tau_1$ such that $A\sbs U$ and $B\sbs V$.
\end{definition}

\begin{proposition}{2.18}
Every $(d,p)$-densely normal $\BS$ $(X,\tau_1<_S\tau_2)$ is
$p$-middly normal.
\end{proposition}

\begin{pf}
Let $A=\tau_1\cl\tau_2\nt A$, $B=\tau_2\cl\tau_1\nt B$ and $A\cap
B=\vnth$. Following condition, there is a subset $Y\in
d\dd\,\cD(X)$ such that $X$ is $p$-normal on $Y$. By the remark
before Definition~2.17, $A \sbs\tau_1\cl(A\cap Y)$ and
$B\sbs\tau_2\cl(B\cap Y)$, i.e., $A$ and $B$ are respectively $1$-
and $2$-concentrated on $Y$. Since $X$ is $p$-normal on $Y$, there
are $U\in\tau_2$, $V\in\tau_1$ such that $A\sbs U$ and $B\sbs V$.
\end{pf}

\begin{remark}{2.19}
If $Y\in i\dd\FF_\sg(X)$ in a $\BS$ $(X,\tau_1,\tau_2)$ and
$F\in\co\tau_i'$ in $(Y,\tau_1',\tau_2')$, then $F\in
i\dd\FF_\sg(X)$. Indeed, let
$Y=\bigcup\limits_{k=1}^\infty\Phi_k$, where $\Phi_k\in\co\tau_i$.
Since $F\in\co\tau_i'$, there is $T\in\co\tau_i$ such that $T\cap
Y=F$. Hence $F=\bigcup\limits_{k=1}^\infty (\Phi_k\cap T)\in
i\dd\FF_\sg(X)$.
\end{remark}

\begin{lemma}{2.20}
Let $A_1$, $A_2$ are subsets of a $\BS$ $(X,\tau_1,\tau_2)$ and
there are countable families
$\boldsymbol{\cU}_1=\{U_k^1\}_{k=1}^\infty\sbs \tau_1$,
$\boldsymbol{\cU}_2=\{U_k^2\}_{k=1}^\infty\sbs \tau_2$ such that
\begin{gather}
    A_1\sbsq \bigcup\limits_{k=1}^\infty U_k^1, \;\;\;
        A_2\sbsq \bigcup\limits_{k=1}^\infty U_k^2, \nonumber \\
    A_1\cap\tau_1\cl U_k^2=\vnth \;\;\text{and}\;\;
        A_2\cap\tau_2\cl U_k^1=\vnth, \;\; \text{for each} \;\;
                    k=\ol{1,\infty}.  \tag{$*$}
\end{gather}
Then there are $U_1\in\tau_1$, $U_2\in\tau_2$ such that $A_i\sbs
U_i$ and $U_1\cap\,U_2=\vnth$.
\end{lemma}

\begin{pf}
Let $V_1^1=U_1^1$, $V_k^1=U_k^1\setminus\bigcup\limits_{n=1}^{k-1}
            \tau_1\cl U_n^2$, $k=\ol{2,\infty} $
and
$$  V_k^2=U_k^2\setminus\bigcup\limits_{n=1}^k
            \tau_2\cl U_n^1, \;\; k=\ol{1,\infty}.     $$
The sets $V_k^1$ and $V_k^2$ are respectively $1$-open and
$2$-open for each $k=\ol{1,\infty}$. Hence
$U_1=\bigcup\limits_{k=1}^\infty V_k^1\in\tau_1$,
$U_2=\bigcup\limits_{k=1}^\infty V_k^2\in\tau_2$ and by $(*)$, we
have $A_1\sbs U_1$ and $A_2\sbs U_2$.

Finally, by construction, if $\ell\geq k$, then
\begin{eqnarray*}
    V_k^1\cap V_\ell^2=
        &&\hskip-0.6cm \big(U_k^1\setminus\big(\tau_1\cl U_1^2\cup\cdots\cup
                \tau_1\cl U_{k-1}^2\big)\big)\cap \\
    \cap &&\hskip-0.6cm \big(U_\ell^2\setminus\big(\tau_2\cl U_1^1\cup\cdots\cup
        \tau_2\cl U_k^1\cup\cdots\cup\tau_2\cl U_\ell^1\big)\big)\sbs \\
    \sbs &&\hskip-0.6cm U_k^1\cap(U_\ell^2\setminus U_k^1)=\vnth,
\end{eqnarray*}
and if $\ell<k$, then
\begin{eqnarray*}
    V_k^1\cap V_\ell^2= &&\hskip-0.6cm
        \big(U_k^1\setminus\big(\tau_1\cl U_1^2\cup\cdots\cup
            \tau_1\cl U_\ell^2\cup\cdots\cup\tau_1\cl U_{k-1}^2)\big)\cap \\
    \cap &&\hskip-0.6cm \big(U_\ell^2\setminus\big(\tau_2\cl U_1^1\cup\cdots\cup
            \tau_2\cl U_\ell^1\big)\big)\sbs \\
    \sbs &&\hskip-0.6cm (U_k^1\setminus\,U_\ell^2)\cap\,U_\ell^2=\vnth.
\end{eqnarray*}
Thus, $V_k^1\cap V_\ell^2=\vnth$ for each $k$ and $\ell$ and so
$U_1\cap\,U_2=\vnth$.
\end{pf}

\begin{lemma}{2.21}
If $(X,\tau_1,\tau_2)$ is a $p$-normal $\BS$, $P\in
1\dd\FF_\sg(X)$, $Q\in 2\dd\FF_\sg(X)$ and $(\tau_1\cl P\cap
Q)\cup(P\cap\tau_2\cl Q)=\vnth$, then there are disjoint sets
$U\in\tau_2$, $V\in\tau_1$ such that $P\sbs U$ and $Q\sbs V$.
\end{lemma}

\begin{pf}
Let $P=\bigcup\limits_{k=1}^\infty F_k$, where $F_k\in\co\tau_1$
for each $k=\ol{1,\infty}$. By condition, $P\cap\tau_2\cl Q=\vnth$
and so, $F_k\cap\tau_2\cl Q=\vnth$ for each $k=\ol{1,\infty}$.
Since $(X,\tau_1,\tau_2)$ is $p$-normal, by (4) of Definition~1.1,
for each $k=\ol{1,\infty}$ there is $U_k^2\in\tau_2$ such that
$F_k\sbs U_k^2$ and $ \tau_1\cl U_k^2\cap\tau_2\cl Q=\vnth$. It is
obvious that $P\sbs\bigcup\limits_{k=1}^\infty U_k^2$. Similarly,
one can construct a countable family
$\{U_k^1\}_{k=1}^\infty\sbs\tau_1$ such that
$Q\sbs\bigcup\limits_{k=1}^\infty U_k^1$ and $\tau_2\cl U_k^1\cap
\tau_1\cl P=\vnth$ for each $k=\ol{1,\infty}$. Hence, by
Lemma~2.20, there are disjoint sets $U\in\tau_2$, $V\in\tau_1$
such that $P\sbs U$ and $Q\sbs V$.
\end{pf}

\begin{proposition}{2.22}
If $Y\in 1\dd\FF_\sg(X)\cap 2\dd\FF_\sg(X)$ and
$(X,\tau_1,\tau_2)$ is $p$-normal, then $(Y,\tau_1',\tau_2')$ is
$p$-strongly normal in $X$.
\end{proposition}

\begin{pf}
Let $A\in\co\tau_1'$, $B\in\co\tau_2'$ and $A\cap B=\vnth$. Then
$$  (\tau_1\cl A\cap B)\cup (A\cap\tau_2\cl B)=\vnth.       $$
Moreover, by Remark~2.19, $A\in 1\dd\FF_\sg(X)$ and $B\in
2\dd\FF_\sg(X)$. Hence, by Lemma~2.21, there are disjoint sets
$U\in\tau_2$, $V\in\tau_1$ such that $A\sbs U$ and $B\sbs V$.
\end{pf}
\vskip+0.2cm

Let us show, that if $\tau_i\sbs j\dd\FF_\sg(X)$ in a $\BS$
$(X,\tau_1,\tau_2)$, then for any $\BsS$ $(Y,\tau_1',\tau_2')$ we
have $\tau_i'\sbs j\dd\FF_\sg(Y)$. Indeed, let $U'\in\tau_i'$.
Then there is $U\in \tau_i$ such that
$$  U'=U\cap Y=\bigcup\limits_{k=1}^\infty F_k\cap Y=
            \bigcup\limits_{k=1}^\infty (F_k\cap Y),    $$
where $F_k\cap Y\in\co\tau_j'$ for each $k=\ol{1,\infty}$, and so
$\tau_i'\sbs j\dd\FF_\sg(Y)$. Hence, take place

\begin{proposition}{2.23}
Every $\BsS$ of an $(i,j)$-perfectly normal $\BS$ is
$(i,j)$-perfectly normal, where by $[16]$, $(X,\tau_1,\tau_2)$ is
$(i,j)$-perfectly normal if it is $p$-normal and $\tau_i\sbs
j\dd\cF_\sg(X)$.
\end{proposition}

The reasonings below we come to the new type of bitopological
separation axioms.

A double family is a pair of families $\cA=\{\cA_1,\cA_2\}$, where
$\cA_i\sbs 2^X$.

\begin{definition}{2.24}
In a $\BS$ $(X,\tau_1,\tau_2)$ a double family
$\cA=\{\cA_1,\cA_2\}$ is said to be $d\wedge p$-discrete if each
point $x\in X$ has as a neighborhood $U(x)\in\tau_1$ so a
neighborhood $V(x)\in\tau_2$ each of which intersects at most one
set from $\cA$. A double family $\cA=\{\cA_1,\cA_2\}$ is
$\Wp$-discrete if each point $x\in X$ has a neighborhood
$U(x)\in\tau_1$ that intersects at most one set of $\cA_2$ or has
a neighborhood $V(x)\in\tau_2$ that intersects at most one set of
$\cA_1$.
\end{definition}

It is obvious that if $\cA=\{\cA_1,\cA_2\}$, is $d\wedge
p$-discrete, then $\cA_1$ and $\cA_2$ are discrete in the usual
sense and $\cA$ is $d\wedge p$-disjoint, that is, if
$A_1\in\cA_i$, $A_2\in\cA_1\cup\cA_2$ and $A_1\neq A_2$, then
$A_1\cap A_2=\vnth$.

\begin{definition}{2.25}
A family $\boldsymbol{\cF}=\{F_s\}_{s\in S}$
$(\boldsymbol{\cU}=\{U_s\}_{s\in S})$ of subsets of a $\BS$
$(X,\tau_1,\tau_2)$ is said to be $p$-closed $(p$-open) if
$\boldsymbol{\cF}\sbs\co\tau_1\cup\co\tau_2$
$(\boldsymbol{\cU}\sbs\tau_1\cup\tau_2)$ and
$\boldsymbol{\cF}\cap\big(\co\tau_i\setminus\{\vnth\}\big)\neq\vnth$
$(\boldsymbol{\cU}\cap\big(\tau_i\setminus\{\vnth\}\big)\neq\vnth)$.
\end{definition}

Hence, $\boldsymbol{\cF}=\{F_s\}_{s\in S}$
$(\boldsymbol{\cU}=\{U_s\}_{s\in S})$ is a double family
$\boldsymbol{\cF}=\{\boldsymbol{\cF}_1,\boldsymbol{\cF}_2\}$
$(\boldsymbol{\cU}=\{\boldsymbol{\cU}_1,\boldsymbol{\cU}_2\})$,
where $\boldsymbol{\cF}_1=\{F_s\}_{s\in S_1}\sbs\co\tau_1$,
$\boldsymbol{\cF}_2=\{F_s\}_{s\in S_2}\sbs\co\tau_2$
$\big(\,\boldsymbol{\cU}_1=\{U_s\}_{s\in S_1}\sbs\tau_1$,
$\boldsymbol{\cU}_2=\{U_s\}_{s\in S_2}\sbs\tau_2\,\big)$ and
$S_1\cup S_2=S$.

\begin{definition}{2.26}
A $\BS$ $(X,\tau_1,\tau_2)$ is said to be $p$-collectionwise
normal if it is $\RR\dd\,p\,\dd\TT_1$ and for every $p$-closed
$d\wedge p$-discrete family
$\boldsymbol{\cF}=\{\boldsymbol{\cF}_1,\boldsymbol{\cF}_2\}$,
where $\boldsymbol{\cF}_i=\{F_s\}_{s\in S_i}\sbs\co\tau_i$, there
exists a $p$-open $d\wedge p$-disjoint family
$\boldsymbol{\cU}=\{\boldsymbol{\cU}_1,\boldsymbol{\cU}_2\}$ such
that $\boldsymbol{\cU}_i=\{U_s\}_{s\in S_j}\sbs\tau_i$ and
$F_s\sbs U_s$ for each $s\in S$.
\end{definition}

Clearly, every $p$-collectionwise normal $\BS$ is $p$-normal.

\begin{theorem}{2.27}
In a $p$-collectionwise normal $\BS$ $(X,\tau_1,\tau_2)$ for any
$p$-closed $d\wedge p$-discrete family
$\boldsymbol{\cF}=\{\boldsymbol{\cF}_1,\boldsymbol{\cF}_2\}$,
where $\boldsymbol{\cF}_i=\{F_s\}_{s\in S_i}\sbs\co\tau_i$, there
exists a $p$-open $\Wp$-discrete family
$\boldsymbol{\cV}=\{\boldsymbol{\cV}_1,\boldsymbol{\cV}_2\}$ such
that $\boldsymbol{\cV}_i=\{V_s\}_{s\in S_j}\sbs\tau_i$ and
$F_s\sbs V_s$ for each $s\in S$.
\end{theorem}

\begin{pf}
We have
$\boldsymbol{\cF}=\{\boldsymbol{\cF}_1,\boldsymbol{\cF}_2\}$,
where $\boldsymbol{\cF}_1=\{F_s\}_{s\in S_1}\sbs\co\tau_1$,
$\boldsymbol{\cF}_2=\{F_s\}_{s\in S_2}\sbs\co\tau_2$, and by
Definition~2.24, $\boldsymbol{\cF}_1$ is $1$-discrete,
$\boldsymbol{\cF}_2$ is
         \linebreak  $2$-discrete. By Definition~2.26, there exists a $p$-open $d\wedge
p$-disjoint family
$\boldsymbol{\cU}=\{\boldsymbol{\cU}_1,\boldsymbol{\cU}_2\}$,
where $\boldsymbol{\cU}_i=\{U_s\}_{s\in S_j}\sbs \tau_i$ and
$F_s\sbs U_s$ for each $s\in S$. Following Theorem~1.1.11 in [10],
$$  F_1=\bigcup\limits_{s\in S_1} F_s\in\co\tau_1, \;\;\;
        F_2=\bigcup\limits_{s\in S_2} F_s\in\co\tau_2       $$
and since $(X,\tau_1,\tau_2)$ is $p$-normal,
$$  F_1\cap\Phi_2=\vnth=F_2\cap\Phi_1, \;\;\text{where}\;\;
        \Phi_i=\Big(X\setminus\bigcup\limits_{s\in S_j} U_s\Big)\in\co\tau_i,   $$
there are neighborhoods $U(\Phi_i)\in\tau_j$ such that
$$  F_1\cap\tau_2\cl U(\Phi_2)=\vnth=F_2\cap\tau_1\cl U(\Phi_1).    $$
Let $\boldsymbol{\cV}=\{\boldsymbol{\cV}_1,\boldsymbol{\cV}_2\}$,
where $  \boldsymbol{\cV}_i=\big\{V_s=U_s\setminus\tau_i\cl
            U(\Phi_i):\;s\in S_j\big\}$.
Then, it is clear that $F_s\sbs V_s$ for each $s\in S$ and
$\Phi_i\cup\big(\bigcup\limits_{s\in S_j} U_s\big)=X$.

It remains to prove only, that the $p$-open $d\wedge p$-disjoint
family
\allowdisplaybreaks
\begin{eqnarray*}
    &\ds \boldsymbol{\cV}=\{\boldsymbol{\cV}_1,\boldsymbol{\cV}_2\}= \\
    &\ds =\!\Big\{\big\{V_s\!=\!U_s\setminus\tau_1\cl U(\Phi_1):\;s\!\in\!S_2\big\},
        \big\{V_s\!=\!U_s\setminus\tau_2\cl U(\Phi_2):\;s\!\in\!S_1\big\}\Big\}
\end{eqnarray*}
is $\Wp$-discrete.

Indeed, let, first, $x\in\Phi_i$. Then $x\in U(\Phi_i)\in\tau_j$
and $U(\Phi_i)\cap V_s=\vnth$ for each $s\in S_j$. If $x\in
\bigcup\limits_{s\in S_j} U_s$, then $x\in U_{s_0}\in\tau_i$ for
some $s_0\in S_j$. Since
$\boldsymbol{\cU}=\{\boldsymbol{\cU}_1,\boldsymbol{\cU}_2\}$ is
$d\wedge p$-disjoint, we have $U_{s_0}\cap U_s=\vnth$ for each
$s\in S_i$ and by determination of the double family
$\boldsymbol{\cV}=\{\boldsymbol{\cV}_1,\boldsymbol{\cV}_2\}$,
$U_{s_0}\cap V_s=\vnth$ for each $s\in S_i$.

Thus, by Definition~2.24, the $p$-open family
$\boldsymbol{\cV}=\{\boldsymbol{\cV}_1,\boldsymbol{\cV}_2\}$ is
\linebreak $\Wp$-discrete.
\end{pf}

\begin{definition}{2.28}
A $\BS$ $(X,\tau_1,\tau_2)$ is said to be $(i,j)$-separately
normal if for each pair of disjoint sets $A,\,B\in\co\tau_i$ there
are disjoint sets $U,\,V\in\tau_j$ such that $A\sbs U$ and $B\sbs
U$.
\end{definition}

It is evident that every $p$-collectionwise normal $\BS$
$(X,\tau_1,\tau_2)$ is       \linebreak      $p$-separately
normal.

\begin{proposition}{2.29}
A $\BS$ $(X,\tau_1,\tau_2)$ is $(i,j)$-separately normal if and
only if for every set $A\in\co\tau_i$ and every neighborhood
$U(A)\in\tau_i$ there is a neighborhood $V(A)\in\tau_j$ such that
$\tau_j\cl V(A)\sbs U(A)$.
\end{proposition}

\begin{pf}
Let, first, $X$ be $(i,j)$-separately normal, $A\in\co\tau_i$ and
$U(A)\in\tau_i$. Then $A\cap B=\vnth$, where $B=X\setminus
U(A)\in\co\tau_i$ and so, there are $V(A),\,V(B)\in\tau_j$ such
that $ \tau_j\cl V(A)\cap V(B)=\vnth$. Hence
$$  \tau_j\cl V(A)\sbs X\setminus V(B)\sbs X\setminus B=U(A).   $$

Conversely, let the condition is satisfied, $A,\,B\in\co\tau_i$
and $A\cap B=\vnth$. Then $A\sbs X\setminus B=U(A)\in\tau_i$ and
hence, there is $V(A)\in\tau_j$ such that $\tau_j\cl V(A)\sbs
U(A)=X\setminus B$ and thus, $B\sbs X\setminus\tau_j\cl
V(A)=V(B)\in\tau_j$. Clearly, $V(A)\cap V(B)=\vnth$.
\end{pf} \vskip+0.2cm

For a $\BS$ $(X,\tau_1<\tau_2)$ take place the following
implications:
$$  
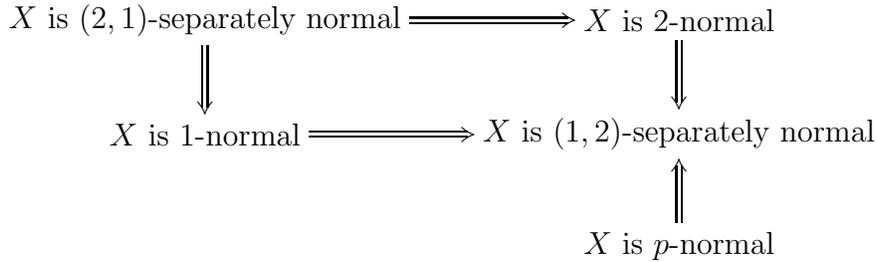
\xymatrix{
    \text{$X$ is $(2,1)$-separately normal} \ar@{=>}[r] \ar@{=>}[d]
        & \text{$X$ is $2$-normal} \ar@{=>}[d] \\
    \text{$X$ is $1$-normal} \ar@{=>}[r]
        & \text{$X$ is $(1,2)$-separately normal} \\
    & \text{$X$ is $p$-normal} \ar@{=>}[u] }        $$

Moreover, by (4) of Corollary~2.2.8 in [8], for a $\BS$
$(X,\tau_1<_C\tau_2)$ we have {\small
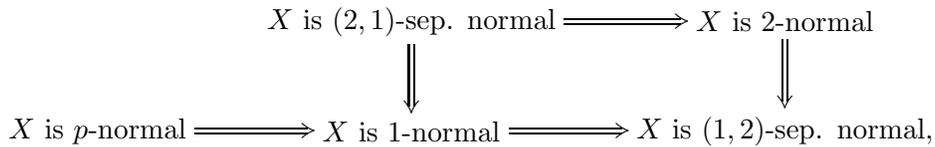
$$  \xymatrix{
    & \text{$X$ is $(2,1)$-sep. normal} \ar@{=>}[r] \ar@{=>}[d]
        & \text{$X$ is $2$-normal} \ar@{=>}[d] \\
    \text{$X$ is $p$-normal} \ar@{=>}[r]
        & \text{$X$ is $1$-normal} \ar@{=>}[r]
        & \text{$X$ is $(1,2)$-sep. normal}, }     $$
} \noindent and by (4) of Corollary~2.3.13 in [8], for a $\BS$
$(X,\tau_1<_N\tau_2)$ we have
$$  \xymatrix{
        \text{$X$ is $(2,1)$-separately normal}\ar@{=>}[r] &
                \text{$X$ is $2$-normal} \ar@{=>}[d] \\
        \text{$X$ is $1$-normal} \ar@{=>}[d] & \text{$X$ is $p$-normal} \ar@{=>}[l] \\
        \text{$X$ is $(1,2)$-separately normal}. & }        $$

\begin{definition}{2.30}
Let $(Y,\tau_1',\tau_2')$ be a $\BsS$ of a $\BS$
$(X,\tau_1,\tau_2)$. Then
\begin{enumerate}
\item[(1)] $Y$ is $(i,j)$-quasi separately normal in $X$ if for
each pair of disjoint sets $A,\,B\in\co\tau_i$ there are disjoint
sets $U,\,V\in\tau_j'$ such that $A\cap Y\sbs U$ and $B\cap Y\sbs
V$.

\item[(2)] $Y$ is $(i,j)$-separately normal in $X$ if for each
pair of disjoint sets $A,\,B\in\co\tau_i$ there are disjoint sets
$U,\,V\in\tau_j$ such that $A\cap Y\sbs U$ and $B\cap Y\sbs V$.

\item[(3)] $Y$ is $(i,j)$-strongly separately normal in $X$ if for
each pair of disjoint sets $A,\,B\in\co\tau_i'$ there are disjoint
sets $U,\,V\in\tau_j$ such that $A\sbs U$ and $B\sbs V$.

\item[(4)] $Y$ is $(i,j)$-separately supernormal in $X$ if for
each pair of disjoint sets $A\in\co\tau_i'$, $B\in\co\tau_i$ there
are disjoint sets $U,V\in\tau_j$, such that $A\sbs U$ and $B\sbs
V$.

\item[(5)] $Y$ is $(i,j)\dd\WS$-separately normal in $X$ if for
each pair of disjoint sets $A\in\co\tau_i'$, $B\in\co\tau_i$ there
are disjoint sets $U\in\tau_j'$, $V\in\tau_j$ such that $A\sbs U$
and $B\cap Y\sbs V$.

\item[(6)] $Y$ is $(i,j)\dd\WS$-separately supernormal in $X$ if
for each pair of disjoint sets $A\in\co\tau_i'$, $B\in\co\tau_i$
there are disjoint sets $U\in\tau_j'$, $V\in\tau_j$ such that
$A\sbs U$ and $B\sbs V$.
\end{enumerate}
\end{definition}

It is clear that if $(Y,\tau_1',\tau_2')$ is a $\BsS$ of a $\BS$
$(X,\tau_1,\tau_2)$, then {\small
$$  \xymatrix{
    \text{$X$\,is\,$(i,j)$-sep.\,norm.} \ar@{=>}[d]_{Y\in\co\tau_i} & \\
    \text{$Y$\,is\,$(i,j)$-sep.\,supernorm.\,in\,$X$} \ar@{=>}[r]\ar@{=>}[d]
        & \text{$Y$\,is\,$(i,j)$-strong.\,sep.\,norm.\,in\,$X$} \ar@{=>}[d] \\
    \text{$Y$\,is\,$(i,j)\dd\WS$-sep.\,supernorm.\,in\,$X$} \ar@{=>}[d]
        & \text{$Y$\,is\,$(i,j)$-sep.\,norm.\,in\,$X$} \ar@{=>}[d] \\
    \text{$Y$\,is\,$(i,j)\dd\WS$-sep.\,norm.\,in\,$X$}
        & \text{$Y$\,is\,$(i,j)$-quasi\,sep.\,norm.\,in\,$X$}. }   $$
}

By analogy with Proposition~2.9 one can prove

\begin{proposition}{2.31}
Let $(Y,\tau_1',\tau_2')$ be a $\BsS$ of a $\BS$
$(X,\tau_1,\tau_2)$. Then
\begin{enumerate}
\item[(1)] $Y$ is $(i,j)$-strongly separately normal in $X$ if and
only if for each set $F\in\co\tau_i'$ and any neighborhood
$U'(F)\in\tau_i'$ there is a neighborhood $U(F)\in\tau_j$ such
that $\tau_j\cl V(F)\cap Y\sbs U'(F)$.

\item[(2)] $Y$ is $(i,j)\dd\WS$-separately supernormal in $X$ if
and only if for each set $F\in\co\tau_i'$ and any neighborhood
$U(F)\in\tau_i$ there is a neighborhood $V'(F)\in\tau_j'$ such
that $\tau_j\cl V'(F)\sbs U(F)$.

\item[(3)] $Y$ is $(i,j)$-separately supernormal in $X$ if and
only if for each set $F\in\co\tau_i'$ and any neighborhood
$U(F)\in\tau_i$ there is a neighborhood $U(F)\in\tau_j$ such that
$\tau_j\cl V(F)\sbs U(F)$.
\end{enumerate}
\end{proposition}

\begin{definition}{2.32}
A $\BS$ $(X,\tau_1,\tau_2)$ is $(i,j)$-separately regular if for
each point $x\in X$ and each set $F\in\co\tau_i$,
$x\,\ol{\in}\,F$, there are disjoint sets $U,\,V\in \tau_j$ such
that $x\in U$ and $F\sbs V$.
\end{definition}

\begin{proposition}{2.33}
A $\BS$ $(X,\tau_1,\tau_2)$ is $(i,j)$-separately regular if and
only if for each point $x\in X$ and any neighborhood
$U(x)\in\tau_i$ there is a neighborhood $V(x)\in\tau_j$ such that
$\tau_j\cl V(x)\sbs U(x)$.
\end{proposition}

The proof is similar to the proof of Proposition~2.29.

Clearly, if $(X,\tau_1,\tau_2)$ is $i\dd\TT_1$, then
$(X,\tau_1,\tau_2)$ is $(i,j)$-separately normal implies that
$(X,\tau_1,\tau_2)$ is $(i,j)$-separately regular.

Take place the following simple but important

\begin{theorem}{2.34}
If a $\BS$ $(X,\tau_1,\tau_2)$ is $p$-separately regular, then
$(X,\tau_1)$ is Baire if and only if $(X,\tau_2)$ is Baire.
\end{theorem}

\begin{pf}
If $(X,\tau_1,\tau_2)$  is $p$-separately regular, then it follows
immediately from  Proposition~2.33 and (1) of Theorem~2.1.5 in [8]
that $\tau_1S\tau_2$. Hence, it remains to use Proposition~3.4 in
[25].
\end{pf}
\vskip+0.2cm

Therefore, all results from [8] and [25], which are connected with
\linebreak       $S$-related topologies, are also correct for
$p$-separately regular bitopological spaces.

Below we give the relative versions of bitopological separate
regularity.

\begin{definition}{2.35}
Let $(Y,\tau_1',\tau_2')$ be a $\BsS$ of a $\BS$
$(X,\tau_1,\tau_2)$. Then
\begin{enumerate}
\item[(1)] $Y$ is $(i,j)\dd\WS$-quasi separately regular in $X$ if
for $x\in Y$, $F\in\co\tau_i$ and $x\,\ol{\in}\,F$, there are
disjoint sets $U,\,V\in\tau_j'$ such that $x\in U$ and $F\cap
Y\sbs V$.

\item[(2)] $Y$ is $(i,j)$-separately regular in $X$ if for $x\in
Y$, $F\in\co\tau_i$ and $x\,\ol{\in}\,F$ there are disjoint sets
$U,\,V\in\tau_j$ such that $x\in U$ and $F\cap Y\sbs V$.

\item[(3)] $Y$ is $(i,j)$-strongly separately regular in $X$ if
for $x\in Y$,         \linebreak        $F\in\co\tau_i'$ and
$x\,\ol{\in}\,F$ there are disjoint sets $U,\,V\in\tau_j$ such
that $x\in U$ and $F\sbs V$.

\item[(4)] $Y$ is $(i,j)$-separately superregular in $X$ if for
$x\in Y$, $F\in\co\tau_i$ and $x\,\ol{\in}\,F$ there are disjoint
sets $U,\,V\in\tau_j$ such that $x\in U$ and $F\sbs V$.

\item[(5)] $Y$ is $(i,j)\dd\WS$-separately regular in $X$ if for
$x\in Y$, $F\in\co\tau_i$ and $x\,\ol{\in}\,F$ there are disjoint
sets $U\in\tau_j'$, $V\in\tau_j$ such that $x\in U$ and $F\cap
Y\sbs V$.

\item[(6)] $Y$ is $(i,j)\dd\WS$-separately superregular in $X$ if
for $x\in Y$,         \linebreak       $F\in\co\tau_i$ and
$x\,\ol{\in}\,F$ there are disjoint sets $U\in\tau_j'$,
$V\in\tau_j$ such that $x\in U$ and $F\sbs V$.

\item[(7)] $Y$ is $(i,j)$-free separately regular in $X$ if for
$x\in X$, $F\in\co\tau_i$ and $x\,\ol{\in}\,F$ there are disjoint
sets $U,\,V\in\tau_j$ such that $x\in U$ and $F\cap Y\sbs V$.
\end{enumerate}
\end{definition}

It is evident that if $(Y,\tau_1',\tau_2')$ is a $\BsS$ of a $\BS$
$(X,\tau_1,\tau_2)$, then the following implications hold: {\small
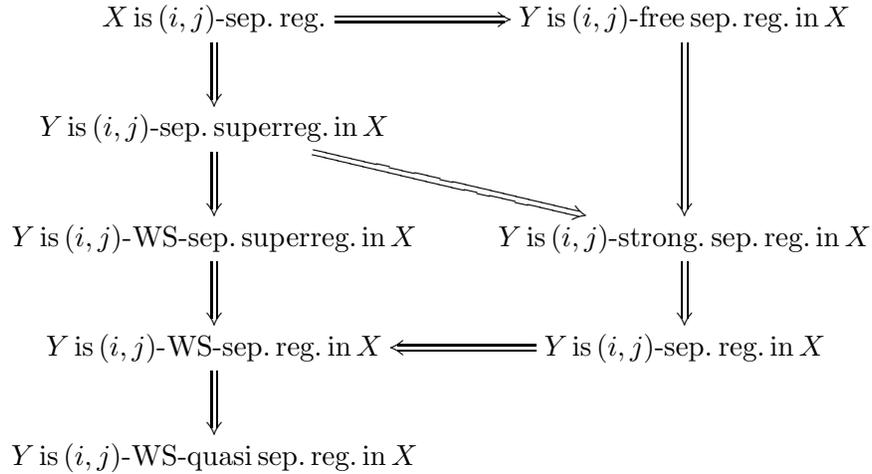
$$  \xymatrix{
    \text{$X$\,is\,$(i,j)$-sep.\,reg.} \ar@{=>}[d] \ar@{=>}[r]
        & \text{$Y$\,is\,$(i,j)$-free\,sep.\,reg.\,in\,$X$} \ar@{=>}[dd] \\
    \text{$Y$\,is\,$(i,j)$-sep.\,superreg.\,in\,$X$}\ar@{=>}[d] \ar@{=>}[dr] & \\
    \text{$Y$\,is\,$(i,j)\dd\WS$-sep.\,superreg.\,in\,$X$} \ar@{=>}[d] &
        \text{$Y$\,is\,$(i,j)$-strong.\,sep.\,reg.\,in\,$X$} \ar@{=>}[d] \\
    \text{$Y$\,is\,$(i,j)\dd\WS$-sep.\,reg.\,in\,$X$} \ar@{=>}[d]
        & \text{$Y$\,is\,$(i,j)$-sep.\,reg.\,in\,$X$} \ar@{=>}[l] \\
    \text{$Y$\,is\,$(i,j)\dd\WS$-quasi\,sep.\,reg.\,in\,$X$} & }      $$
}

Take place

\begin{proposition}{2.36}
Let $(Y,\tau_1',\tau_2')$ be a $\BsS$ of a $\BS$
$(X,\tau_1,\tau_2)$. Then
\begin{enumerate}
\item[(1)] $Y$ is $(i,j)$-strongly separately regular in $X$ if
and only if for each point $x\in Y$ and any neighborhood
$U'(x)\in\tau_i'$ there is a neighborhood $V(x)\in\tau_j$ such
that $\tau_j\cl V(x)\cap Y\sbs U'(x)$.

\item[(2)] $Y$ is $(i,j)\dd\WS$-separately superregular in $X$ if
and only if for each point $x\in Y$ and any neighborhood
$U(x)\in\tau_i$ there is a neighborhood $V'(x)\in\tau_j'$ such
that $\tau_j\cl V'(x)\sbs U(x)$.

\item[(3)] $Y$ is $(i,j)$-separately superregular in $X$ if and
only if for each point $x\in Y$ and any neighborhood
$U(x)\in\tau_i$ there is a neighborhood $V(x)\in\tau_j$ such that
$\tau_j\cl V(x)\sbs U(x)$.
\end{enumerate}
\end{proposition}

The proof is analogical to the proof of Proposition~2.3 and can be
omitted.

\begin{corollary}{2.37}
A $\TsS$ $(Y,\tau')$ of a $\TS$ $(X,\tau)$ is strongly separately
regular $($respectively, $\WS$-separately superregular, separately
su\-per\-re\-gu\-lar$)$ in $X$ if and only if for each point $x\in
Y$ and any neighborhood $U'(x)\in\tau'$ $($respectively,
$U(x)\in\tau)$ there is a neighborhood $V(x)\in\tau$
$($res\-pec\-ti\-ve\-ly, $V'(x)\in\tau'$, $V(x)\in\tau)$ such that
$\tau\cl V(x)\cap Y\sbs U'(x)$ $($respectively, $\tau\cl V'(x)\sbs
U(x)$, $\tau\cl V(x)\sbs U(X))$.
\end{corollary}

Hence, it follows immediately from the remark before
Definition~2.4 and Corollary~2.37 that if $(Y,\tau')$ is a $\TsS$
of a $\TS$, then $Y$ is strongly  regular (respectively,
$\WS$-superregular, superregular) in $X$ if and only if $Y$ is
strongly separately regular (respectively, $\WS$-separately
superregular, separately superregular) in $X$.

At the end of the section let us consider the bitopological and
topological versions of relative extremal disconnectedness and
relative connectedness.

\begin{proposition}{2.38}
A $\BS$ $(X,\tau_1,\tau_2)$ is $p$-connected if and only if for
each $A\sbs X$, $\vnth\neq A\neq X$, we have $(1,2)\dd\Fr
A\neq\vnth\neq(2,1)\dd\Fr A$.
\end{proposition}

\begin{pf}
Clearly, $(1,2)\dd\Fr A=\vnth$ if and only if
$A\in\tau_2\cap\co\tau_1$ and $(2,1)\dd\Fr A=\vnth$ if and only if
$A\in\tau_1\cap\co\tau_2$. Hence, it remains to use (7) of
Definition~1.1.
\end{pf}

\begin{definition}{2.39}
A $\BsS$ $(Y,\tau_1',\tau_2')$ of a $\BS$ $(X,\tau_1,\tau_2)$ is
$(i,j)$-ex\-tre\-mal\-ly disconnected in $X$ if $\tau_j\cl
U'\in\tau_i$ for each set $U'\in\tau_i'$.
\end{definition}

Clearly, if $Y\in\tau_i$ and $X$ is $(i,j)$-extremally
disconnected, then $Y$ is $(i,j)$-extremally disconnected in $X$
and so, by reasonings after Definition~0.1.18 in [8], if
$Y\in\tau_1\cap \tau_2$, then the following three conditions are
equivalent: $Y$ is $(1,2)$-extremally disconnected in $X$, $Y$ is
$(2,1)$-extremally disconnected in $X$ and $Y$ is $p$-extremally
disconnected in~$X$.

\begin{definition}{2.40}
A $\TsS$ $(Y,\tau')$ of a $\TS$ $(X,\tau)$ is extremally
disconnected in $X$ if $\tau\cl U'\in\tau$ for each set
$U'\in\tau'$.
\end{definition}

Clearly, if $X$ is extremally disconnected and $Y\in\tau$, then
$Y$ is extremally disconnected in $X$.

\begin{proposition}{2.41}
If $(Y,\tau_1',\tau_2')$ is a $(1,2)$-extremally disconnected or
$(2,1)$-extremally disconnected $\BsS$ of a $p$-connected $\BS$
$(X,\tau_1,\tau_2)$, then $Y$ is $p$-connected too.
\end{proposition}

\begin{pf}
Contrary: $Y$ is not $p$-connected. Then $Y=A\cup B$, where
$A\in(\tau_1'\cap\co\tau_2')\setminus\{\vnth\}$,
$B\in(\tau_2'\cap\co\tau_1')\setminus\{\vnth\}$ and $A\cap
B=\vnth$. Let, for example, $Y$ is $(2,1)$-extremally disconnected
in $X$. Then $\tau_1\cl B\in\tau_2$, $\tau_1\cl B\neq X$ as
$\tau_1\cl B\cap A=\vnth$ and so, by (7) of Definition~1.1, $X$ is
not $p$-connected.
\end{pf}

\begin{corollary}{2.42}
If $(Y,\tau')$ is an extremally disconnected $\TsS$ of a connected
$\TS$ $X$, then $Y$ is connected too.
\end{corollary}

However, we have the following elementary

\begin{example}{2.43}
Let $X=\{a,b,c,d,e\}$, $\tau=\big\{\vnth,\{a\},\{a,b\},\{a,b,c\},
        \linebreak      \{d,e\},\{a,d,e\},\{a,b,d,e\},X\big\}$ and $Y=\{a,b,c\}$. Then
$\tau'=\big\{\vnth,\{a\},         \linebreak      \{a,b\},Y\}$,
$(Y,\tau')$ is connected, $Y$ is extremally disconnected in $X$
and $X$ is disconnected.
\end{example}

\begin{proposition}{2.44}
If $(X,\tau_1,\tau_2)$ is a $p$-connected $p$-Tychonoff $\BS$
which contains more than one point, then $|X|\geq \bold{c}$, where
$\bold{c}$ is the cardinality of $\bR$.
\end{proposition}

\begin{pf}
Let $X$ be a $p$-Tychonoff $\BS$, $x_1,\,x_2\in X$, $x_1\neq x_2$.
Since $X$ is $\RR\dd\,p\,\dd\TT_1$, by (3) of Definition~1.1 and
Proposition~0.1.4 in [8], there is a $d$-continuous function
$f:(X,\tau_1,\tau_2)\to (I,\om_1',\om_2')$ such that $f(x_1)=0$
and $f(x_2)=1$. Since $X$ is $p$-connected, by (7) of
Definition~1.1, the function $f$ has the Darboux property and
hence, $I\sbs f(X)$. Thus $|X|\geq \bold{c}$.
\end{pf}

\begin{definition}{2.45}
A $\BsS$ $(Y,\tau_1',\tau_2')$ of a $\BS$ $(X,\tau_1,\tau_2)$ is
$p$-connected in $X$ if for each set $A\sbs X$, $\vnth\neq A\neq
X$, we have $  (1,2)\dd\Fr A\cap Y\neq\vnth\neq(2,1)\dd\Fr A\cap
Y.    $
\end{definition}

\begin{proposition}{2.46}
If $(Y,\tau_1',\tau_2')$ is a $\BsS$ of a $\BS$
$(X,\tau_1,\tau_2)$ and $Y$ is $p$-connected in $X$, then $X$ is
$p$-connected. Moreover, if $Y\in\tau_1\cap\tau_2$, then $Y$ is
$p$-connected too.
\end{proposition}

\begin{pf}
The first part follows directly from Definition~2.45 and
Proposition~2.38. Now, let $Y\in\tau_1\cap\tau_2$ and $A\sbs Y$,
$\vnth\neq A\neq Y$. Then $\vnth\neq A\neq X$ and since
$Y\in\tau_1\cap\tau_2$, $(i,j)\dd\Fr A\cap Y=(i,j)\dd\Fr_YA$.
Thus, it remains to apply once more Proposition~2.38.
\end{pf}

\begin{corollary}{2.47}
If $(Y,\tau')$ is a connected $\TsS$ of a $\TS$ $(X,\tau)$, i.e.,
if $\vnth\neq A\neq X$ implies that $\Fr A\cap Y\neq\vnth$, then
$X$ is connected. Moreover, if $Y\in\tau$, then $Y$ is connected
too.
\end{corollary}

\vskip+0.5cm
\section*{\textbf{3. $(i,j)$-Relative Regularities}}
\vskip+0.2cm

From Definition~2.2 it follows immediately

\begin{proposition}{3.1}
If
$(Z,\tau_1'',\tau_2'')\sbs(Y,\tau_1',\tau_2')\sbs(X,\tau_1,\tau_2)$,
then
\begin{enumerate}
\item[(1)] If $Y$ is $(i,j)$-quasi regular $($respectively,
$(i,j)$-regular, $(i,j)$-stron\-gly regular,
$(i,j)$-su\-per\-re\-gu\-lar, $(i,j)\dd\WS$-regular,
$(i,j)\dd\WS$-su\-per\-re\-gu\-lar, $(i,j)$-free regular$)$ in
$X$, then $Z$ is also $(i,j)$-quasi regular
$($res\-pec\-ti\-ve\-ly, $(i,j)$-regular, $(i,j)$-strongly
regular,       \linebreak       $(i,j)$-superregular,
$(i,j)\dd\WS$-re\-gu\-lar, $(i,j)\dd\WS$-su\-per\-re\-gu\-lar,
\linebreak $(i,j)$-free regular$)$ in $X$.

\item[(2)] If $Z$ is $(i,j)$-quasi regular $($respectively,
$(i,j)$-regular, $(i,j)$-stron\-gly regular,
$(i,j)$-su\-per\-re\-gu\-lar, $(i,j)\dd\WS$-regular,
$(i,j)\dd\WS$-su\-per\-re\-gu\-lar, $(i,j)$-free regular$)$ in
$X$, then $Z$ is also $(i,j)$-quasi regular
$($res\-pec\-ti\-ve\-ly, $(i,j)$-regular, $(i,j)$-strongly
regular, $(i,j)$-superregular, $(i,j)\dd\WS$-re\-gu\-lar,
$(i,j)\dd\WS$-superregular,          \linebreak  $(i,j)$-free
regular$)$ in $Y$.
\end{enumerate}
\end{proposition}

\begin{proposition}{3.2}
If $(Y,\tau_1'<\tau_2')$ is a $j$-dense $\BsS$ of a $\BS$
\linebreak    $(X,\tau_1<_S\tau_2)$, then $Y$ is $(i,j)$-regular
if and only if $Y$ is $(i,j)$-regular in $X$.
\end{proposition}

\begin{pf}
Let, first, $Y$ be $(i,j)$-regular in $X$, $x\in Y$,
$F\in\co\tau_i'$ and $x\,\ol{\in}\,F$. Since
$x\,\ol{\in}\,\tau_i\cl F\in\co\tau_i$, there are disjoint sets
$U\in\tau_i$, $V\in\tau_j$ such that $x\in U$, $\tau_i\cl F\cap
Y=F\sbs V$. It is clear, that $x\in U'=U\cap Y\in\tau_i'$, $F\sbs
V'=V\cap Y\in\tau_j'$ and so $Y$ is $(i,j)$-regular.

Conversely, let $Y$ be $(i,j)$-regular, $x\in Y$, $F\in\co\tau_i$
and $x\,\ol{\in}\,F$. Then there are disjoint sets
$U''\in\tau_i'$, $V'\in\tau_j'$ such that $x\in U''$ and $F\cap
Y\sbs V'$. Clearly $U''\cap\tau_i'\cl V'=\vnth$. Let
$U'\in\tau_i$, $V\in\tau_j$ such that $U'\cap Y=U''$ and $V\cap
Y=V'$. Then, by Lemma~2.16,
\begin{eqnarray*}
    \vnth &&\hskip-0.6cm =U'\cap\tau_i'\cl V'=
            U'\cap\tau_i'\cl(V\cap Y)=
        U'\cap\big(\tau_i\cl(V\cap Y)\cap Y\big)= \\
    &&\hskip-0.6cm =U'\cap(\tau_i\cl V\cap Y).
\end{eqnarray*}
Let $U=U'\setminus\tau_i\cl V=U'\setminus(\tau_i\cl V\setminus
Y)$. Then $x\in U\in\tau_i$, $F\cap Y\sbs V'\sbs V$ and $U\cap
V=\vnth$. Thus $Y$ is $(i,j)$-regular in $X$.
\end{pf}

\begin{corollary}{3.3}
If $(Y,\tau_1'<\tau_2')$ is a $2$-dense $\BsS$ of a $\BS$
$(X,\tau_1<\tau_2)$, then $Y$ is $(1,2)$-regular if and only if
$Y$ is $(1,2)$-regular in $X$.
\end{corollary}

\begin{pf}
Indeed, it suffices to note that by the first part of Lemma~2.16,
$\tau_1\cl U=\tau_1\cl(U\cap Y)$ for any set $U\in\tau_2$.
\end{pf}

\begin{proposition}{3.4}
If $(Y,\tau_1',\tau_2')$ is an $(i,j)$-superregular $\BsS$ in a
$\BS$ $(X,\tau_1,\tau_2)$, $A$ is an $i$-compact subset of $Y$,
$B\in\co\tau_i'$ and $A\cap B=\vnth$, then there are disjoint sets
$U\in\tau_i$, $V\in\tau_j$ such that $A\sbs U$ and $B\sbs V$.
\end{proposition}

\begin{pf}
Since $A,B\sbs Y$ and $A\cap B=\vnth$, we have $A\cap\tau_i\cl
B=\vnth$. Moreover, as $Y$ is $(i,j)$-superregular in $X$, for
each point $x\in A$ there are disjoint sets $U(x)\in\tau_i$ and
$U_x(\tau_i\cl B)\in\tau_j$. Since $A\sbs\bigcup\limits_{x\in A}
(U(x)\cap Y)$ and $A$ is $i$-compact in $Y$, there are
$x_1,x_2,\dots,x_n\in A$ such that $A\sbs\bigcup\limits_{k=1}^n
\big(U(x_k)\cap Y\big)$. Finally, if $U=\bigcup\limits_{k=1}^n
U(x_k)$ and $V=\bigcap\limits_{k=1}^n U_{x_k}(\tau_i\cl B)$, then
$U\in\tau_i$, $V\in\tau_j$, $A\sbs U$, $B\sbs V$ and $U\cap
V=\vnth$.
\end{pf}

\begin{corollary}{3.5}
If $(X,\tau_1,\tau_2)$ is $(i,j)$-regular, $(Y,\tau_1',\tau_2')$
is a $\BsS$ of $X$, $A$ is $i$-compact in $Y$, $B\in\co\tau_i'$
and $A\cap B=\vnth$, then there are disjoint sets $U\in\tau_i$,
$V\in\tau_j$ such that $A\sbs U$ and $B\sbs V$.
\end{corollary}

\begin{pf}
The condition is obvious since $Y$ is $(i,j)$-superregular in $X$.
\end{pf}

\begin{proposition}{3.6}
A $\BsS$ $(Y,\tau_1',\tau_2')$ is $(i,j)$-regular in a $\BS$
$(X,\tau_1,\tau_2)$ if and only if for each point $x\in Y$ and
each set $F\in\co\tau_i'$, $x\,\ol{\in}\,F$, there is a set
$U\in\tau_i$ such that $x\in U$ and $\tau_j\cl U\cap F=\vnth$.
\end{proposition}

\begin{pf}
Let, first, $Y$ be $(i,j)$-regular in $X$, $x\in Y$,
$F\in\co\tau_i'$ and $x\,\ol{\in}\,F$. Then
$x\,\ol{\in}\,\tau_i\cl F$ and hence, there are $U\in\tau_i$,
$V\in\tau_j$ such that $x\in U$, $\tau_i\cl F\cap Y=F\sbs V$ and
$U\cap V=\vnth$. Clearly,
$$  \tau_j\cl U\cap V=\vnth \;\;\text{and so}\;\;
            \tau_j\cl U\cap F=\vnth.        $$

Conversely, let the condition  be satisfied, $x\in Y$,
$F\in\co\tau_i$ and $x\,\ol{\in}\,F$. Then $x\,\ol{\in}\,F\cap
Y=F'\in\co\tau_i'$ and by condition, there is $U\in\tau_i$ such
that $x\in U$ and $\tau_j\cl U\cap F'=\vnth$. Let
$V=X\setminus\tau_j\cl U$. Then $V\in\tau_j$, $F'\sbs V$ and
$U\cap V=\vnth$. Thus $Y$ is $(i,j)$-regular in $X$.
\end{pf}

\begin{proposition}{3.7}
If $(Y,\tau_1',\tau_2')$ is $(i,j)\dd\WS$-superregular in $X$ and
$(i,j)$-extremally disconnected in $X$, then $(i,j)\dd\ind_xX=0$
for each point $x\in Y$.
\end{proposition}

\begin{pf}
Let $x\in Y$ be any point and $U(x)\in\tau_i$ be any neighborhood.
Since $Y$ is $(i,j)\dd\WS$-superregular in $X$, by (2) of
Proposition~2.3 there is $V'(x)\in\tau_i'$ such that $\tau_j\cl
V'(x)\sbs U(x)$. But $Y$ is $(i,j)$-extremally disconnected in $X$
and hence, $V(x)=\tau_j\cl V'(x)\in\tau_i$. Thus
$V(x)\in\tau_i\cap\co\tau_j$ and so $(i,j)\dd\ind_xX=0$.
\end{pf}

\begin{corollary}{3.8}
If $(Y,\tau')$ is $\WS$-superregular in $X$ and extremally
disconnected in $X$, then $\ind_xX=0$ for each point $x\in Y$.
\end{corollary}

\vskip+0.5cm
\section*{\textbf{4. $(i,j)$-Relative Normalities}}
\vskip+0.2cm

\begin{proposition}{4.1}
Let
$(Z,\tau_1'',\tau_2'')\sbs(Y,\tau_1',\tau_2')\sbs(X,\tau_1,\tau_2)$.
Then
\begin{enumerate}
\item[(1)] $Y$ is $p$-normal implies that $Z$ is $p$-quasi normal
in $X$, and $Y$ is $p$-normal in $X$ implies that $Z$ is
$p$-normal in $X$ too.

\item[(2)] If $Y\in\tau_1\cup\tau_2$, then $Z$ is $p$-strongly
normal in $Y$ implies that $Z$ is $p$-strongly normal in $X$.

\item[(3)] If $Y\in d\dd\,\cD(X)$, $\tau_1<_S\tau_2$ and $Z$ is
$p$-internally normal in $Y$, then $Z$ is $p$-internally normal in
$X$.
\end{enumerate}

Let $(Y,\tau_1''',\tau_2''')\sbs (Y_1,\tau_1'',\tau_2'')\sbs
                (X_1,\tau_1',\tau_2')\sbs (X,\tau_1,\tau_2)$.
Then
\begin{enumerate}
\item[(4)] $Y_1$ is $p$-quasi normal in $X_1$ implies that $Y$ is
$p$-quasi normal in~$X$.
\end{enumerate}
\end{proposition}

\begin{pf}
(1) Let, first, $A\in\co\tau_1$, $B\in\co\tau_2$ and $A\cap
B=\vnth$. Then $A\cap Y\in\co\tau_1'$, $B\cap Y\in\co\tau_2'$ are
disjoint and so, there are disjoint sets $U'\in\tau_2'$,
$V'\in\tau_1'$ such that $A\cap Y\sbs U'$ and $B\cap Y\sbs V'$. It
is evident that $A\cap Z\sbs U\in\tau_2''$ and $B\cap Z\sbs
V\in\tau_1''$, where $U=U''\cap Z$, $V=V''\cap Z$ and $ U''\cap
Y=U'$, $V''\cap Y=V'$ for $U''\in\tau_2$, $V''\in\tau_1''$.

The second condition is obvious.

(2) Now, let $Y\in\tau_1\cup\tau_2$, $A\in\co\tau_1''$,
$B\in\co\tau_2''$ and $A\cap B=\vnth$. Hence, by condition, there
are disjoint $U'\in\tau_2'$, $V'\in\tau_1'$ such that $A\sbs U'$
and $B\sbs V'$. Let, for example, $Y\in\tau_1$. Then $V'\in\tau_1$
and if $U\in\tau_2$, $U\cap Y=U'$, then we have $A\sbs U$, $B\sbs
V'$ and $U\cap V'=\vnth$. Thus, $Z$ is $p$-strongly normal in $X$.

(3) Let $A\in\co\tau_1$, $B\in\co\tau_2$, $A,B\sbs Z$ and $A\cap
B=\vnth$. Since $A\in\co\tau_1'$, $B\in\co\tau_2'$ and $Z$ is
$p$-internally normal in $Y$, there are disjoint sets
$U''\in\tau_2'$, $V''\in\tau_1'$ such that $A\sbs U''$ and $B\sbs
V''$.

Therefore, by the second part of Lemma~2.16 there are disjoint
sets $U\in\tau_2$, $V\in\tau_1$ such that $A\sbs U$ and $B\sbs V$
so that $Z$ is $p$-internally normal in $X$.

(4) Finally, let $A\in\co\tau_1$, $B\in\co\tau_2$ and $A\cap
B=\vnth$. Then
$$  A_1=A\cap X_1\in\co\tau_1', \;\;\; B_1=B\cap X_1\in\co\tau_2'  $$
and by condition, there are disjoint $U'\in\tau_2''$,
$V'\in\tau_1''$ such that $A_1\cap Y_1\sbs U'$ and $B_1\cap
Y_1\sbs V'$. It is clear, that
\begin{eqnarray*}
    A\cap Y &&\hskip-0.6cm =(A_1\cap Y_1)\cap Y\sbs U\in\tau_2', \\
    B\cap Y &&\hskip-0.6cm =(B_1\cap Y_1)\cap Y\sbs V\in\tau_1',
\end{eqnarray*}
where  $U=U'\cap Y$ and $V=V'\cap Y$. Thus $Y$ is $p$-quasi normal
in $X$.~\end{pf} \vskip+0.2cm

Our next purpose is to establish the conditions under which the
converse implications to the implications after Definition~2.11
hold.

\begin{lemma}{4.2}
If in a $\BS$ $(X,\tau_1,\tau_2)$ for every pair of disjoint sets
$A\in\co\tau_1$, $B\in\co\tau_2$ and every pair of sets
$U\in\tau_2$, $V\in\tau_1$, where $A\sbs U$ and $B\sbs V$, there
are sequences $\{U_n\}_{n=1}^\infty\sbs\tau_2$,
$\{V_n\}_{n=1}^\infty\sbs\tau_1$ such that
$$  A\sbs \bigcup\limits_{n=1}^\infty U_n, \;\;\;
            B\sbs \bigcup\limits_{n=1}^\infty V_n, \;\;\;
        \tau_1\cl U_n\sbs U \;\;\text{and}\;\; \tau_2\cl V_n\sbs V  $$
for each $n=\ol{1,\infty}$, then $X$ is $p$-normal.
\end{lemma}

\begin{pf}
Let $A\in\co\tau_1$, $B\in\co\tau_2$ and $A\cap B=\vnth$. If
$U=X\setminus B\in\tau_2$ and $V=X\setminus A\in\tau_1$, then
$A\sbs U$, $B\sbs V$ and, by condition, there are sequences
$\{U_n\}_{n=1}^\infty\sbs\tau_2$, $\{V_n\}_{n=1}^\infty\sbs\tau_1$
such that
$$  A\sbs \bigcup\limits_{n=1}^\infty U_n \;\;\text{and}\;\;
        B\cap\tau_1\cl U_n=\vnth, \;\;\text{for each}\;\;
                n=\ol{1,\infty},        $$
and
\begin{equation}
    B\sbs \bigcup\limits_{n=1}^\infty V_n \;\;\text{and}\;\;
        A\cap\tau_2\cl V_n=\vnth, \;\;\text{for each}\;\; n=\ol{1,\infty}.
                    \tag{$*$}
\end{equation}
Let us consider for each $k=\ol{1,\infty}$ the sets:
\begin{equation}
     G_k=U_k\setminus \bigcup\limits_{\ell\leq k} \tau_2\cl V_\ell, \;\;\;
        H_k=V_k\setminus \bigcup\limits_{\ell\leq k} \tau_1\cl U_\ell.
                    \tag{$**$}
\end{equation}
Then $G_k\in\tau_2$, $H_k\in\tau_1$ for each $k=\ol{1,\infty}$.
Moreover, by $(*)$,
$$  A\sbs \bigcup\limits_{k=1}^\infty G_k=G \;\;\text{and}\;\;
        B\sbs \bigcup\limits_{k=1}^\infty H_k=H.       $$
By $(**)$, $G_k\cap V_\ell=\vnth$ for $\ell\leq k$ and, hence,
$G_k\cap H_\ell=\vnth$ for $\ell\leq k$. Similarly,
$H_\ell\cap\,U_k=\vnth$ for $k\leq\ell$ and so $H_\ell\cap
G_k=\vnth$ for $k\leq\ell$. Thus $G_k\cap H_\ell=\vnth$ for each
$k=\ol{1,\infty}$, $\ell=\ol{1,\infty}$ and, hence, $G\cap
H=\vnth$.  Since $A\sbs G$ and $B\sbs H$, $X$ is $p$-normal.
\end{pf}

\begin{definition}{4.3}
We say that a $\BsS$ $(Y,\tau_1',\tau_2')$ is $(i,j)$-perfectly
located in a $\BS$ $(X,\tau_1,\tau_2)$ if for each set
$U\in\tau_j$, $Y\sbs U$, there is a countable family
$\{F_n\}_{n=1}^\infty\sbs\co\tau_i$ such that
$Y\sbs\bigcup\limits_{n=1}^\infty F_n\sbs U$.
\end{definition}

\begin{theorem}{4.4}
Let $(Y,\tau_1',\tau_2')$ be a $p$-quasi normal $\BsS$ of a $\BS$
$(X,\tau_1,\tau_2)$. Then the following conditions are satisfied:
\begin{enumerate}
\item[(1)] If $Y\in d\dd\,\cD(X)$ and $\tau_1<_S\tau_2$, then $Y$
is $p$-normal in $X$ and so $Y$ is $p$-internally normal in $X$.

\item[(2)] If $Y$ is $p$-perfectly located in $X$ and
$\tau_i\cl\tau_j\cl F\cap Y=F$ for each $F\in\co\tau_j'$, then $Y$
is $p$-normal.
\end{enumerate}
\end{theorem}

\begin{pf}
(1) Let $A\in\co\tau_1$, $B\in\co\tau_2$ and $A\cap B=\vnth$.
Since $Y$ is $p$-quasi normal in $X$, there are disjoint sets
$U''\in\tau_2'$, $V''\in\tau_1'$ such that $A\cap Y\sbs U''$ and
$B\cap Y\sbs V''$. Therefore, by the second part of Lemma~2.16,
there are $U\in\tau_2$, $V\in\tau_1$ such that $A\sbs U$, $B\sbs
V$ and $U\cap V=\vnth$. The rest is given by implications after
Definition~2.11.

(2) Let $A\in\co\tau_1'$, $B\in\co\tau_2'$ and $A\cap B=\vnth$.
Let us consider the sets
$$  P=\tau_1\cl\big(\tau_1\cl A\cap \tau_2\cl B\big) \;\;\text{and}\;\;
        Q=\tau_2\cl\big(\tau_1\cl A\cap \tau_2\cl B\big).       $$
Then $U=X\setminus P\in\tau_1$ and $V=X\setminus Q\in\tau_2$.
Moreover,
\begin{eqnarray*}
    P\cap Y  &&\hskip-0.6cm =
        \tau_1\cl\big(\tau_1\cl A\cap \tau_2\cl B\big)\cap Y\sbs
            \tau_1\cl A\cap\tau_1\cl\tau_2\cl B\cap Y= \\
    &&\hskip-0.6cm =A\cap B=\vnth
\end{eqnarray*}
and similarly $Q\cap Y=\vnth$. Hence, $Y\cap (P\cup Q)=\vnth$ so
that $Y\sbs X\setminus (P\cup Q)=U\cap V$. Since $Y$ is
$p$-perfectly located in $X$, there are sequence
$\{F_n\}_{n=1}^\infty\sbs\co\tau_1$ and
$\{\Phi_n\}_{n=1}^\infty\sbs\co\tau_2$ such that $Y\sbs
\bigcup\limits_{n=1}^\infty F_n\sbs V$ and $Y\sbs
\bigcup\limits_{n=1}^\infty \Phi_n\sbs U$. Let
\begin{eqnarray*}
    A_n &&\hskip-0.6cm =A\cap F_n=(A\cap Y)\cap F_n=
        A\cap (F_n\cap Y)\in\co\tau_1', \\
    B_n &&\hskip-0.6cm =B\cap\Phi_n=(B\cap Y)\cap \Phi_n=
        B\cap(\Phi_n\cap Y)\in\co\tau_2'
\end{eqnarray*}
for each $n=\ol{1,\infty}$. Then
\begin{eqnarray*}
    &\ds \tau_1\cl A\cap\tau_2\cl B_n= \\
    &\ds =\tau_1\cl A\cap\tau_2\cl(B\cap\Phi_n)\sbs
        \tau_1\cl A\cap\tau_2\cl B\cap\Phi_n\sbs P\cap\,U=\vnth,
\end{eqnarray*}
and similarly, $\tau_2\cl B\cap\tau_1\cl A_n=\vnth$ for each
$n=\ol{1,\infty}$. Since $Y$ is $p$-quasi normal in $X$ and
$\tau_1\cl A\cap\tau_2\cl B_n=\vnth$ for each $n=\ol{1,\infty}$,
for the sets $A$ and $B_n$ there is a sequence
$\{V_n\}_{n=1}^\infty\sbs\tau_1'$ such that $B_n\sbs V_n$ and
$A\cap\tau_2'\cl V_n=\vnth$ for each $n=\ol{1,\infty}$. Similarly,
since $\tau_2\cl B\cap\tau_1\cl A_n=\vnth$ for each
$n=\ol{1,\infty}$, for the sets $B$ and $A_n$, there is a sequence
$\{U_n\}_{n=1}^\infty\sbs\tau_2'$ such that $A_n\sbs U_n$ and
$B\cap\tau_1'\cl U_n=\vnth$ for each $n=\ol{1,\infty}$. It is also
obvious that
\begin{eqnarray*}
    A &&\hskip-0.6cm =A\cap\Big(\bigcup\limits_{n=1}^\infty F_n\Big)=
        \bigcup\limits_{n=1}^\infty (A\cap F_n)=\bigcup\limits_{n=1}^\infty A_n, \\
    B &&\hskip-0.6cm =B\cap\Big(\bigcup\limits_{n=1}^\infty \Phi_n\Big)=
        \bigcup\limits_{n=1}^\infty (B\cap \Phi_n)=\bigcup\limits_{n=1}^\infty B_n
\end{eqnarray*}
and $A\sbs \bigcup\limits_{n=1}^\infty U_n$, $B\sbs
\bigcup\limits_{n=1}^\infty V_n$. Moreover, $\tau_2'\cl V_n\sbs
V=Y\setminus A$ and $\tau_1'\cl U_n\sbs U=Y\setminus B$ for each
$n=\ol{1,\infty}$. Thus, by Lemma~4.2, $(Y,\tau_1',\tau_2')$ is
$p$-normal.
\end{pf}

\begin{corollary}{4.5}
For a $\BsS$ $(Y,\tau_1',\tau_2')$ of a $\BS$ $(X,\tau_1,\tau_2)$
we have:
\begin{enumerate}
\item[(1)] If $Y\in d\dd\,\cD(X)$ and $\tau_1<_S\tau_2$, then $Y$
is $p$-quasi normal in $X$ if and only if $Y$ is $p$-normal in
$X$.

\item[(2)] If $Y$ is $p$-perfectly located in $X$ and
$\tau_i\cl\tau_j\cl F\cap Y=F$ for each $F\in\co\tau_j'$, then $Y$
is $p$-quasi normal in $X$ if and only if $Y$ is $p$-normal.

\item[(3)] If $Y\in d\dd\,\cD(X)$, $\tau_1<_S\tau_2$ and $Y$ is
$p$-normal, then $Y$ is $p$-in\-ter\-nal\-ly normal in $X$.
\end{enumerate}
\end{corollary}

\begin{pf}
(1) and (2) follow respectively from (1) and (2) of Theorem~4.4
taking into account the implications after Definition~2.11.

(3) If $Y$ is $p$-normal, then by implications after
Definition~2.11, it is $p$-quasi normal in $X$. Hence, it remains
to use (1) of Theorem~4.4.
\end{pf}

\begin{theorem}{4.6}
A $\BsS$ $(Y,\tau_1',\tau_2')$ of a $\BS$ $(X,\tau_1,\tau_2)$ is
$p$-strongly normal in $X$ anyone of the following conditions is
satisfied:
\begin{enumerate}
\item[(1)] $\big(Y\in d\dd\,\cD(X)$, $\tau_1<_S\tau_2$ and $Y$ is
$p$-normal$\big)$ or $\big(Y\in\tau_1\cup\tau_2$ and $Y$ is
$p$-normal$\big)$.

\item[(2)] $\big(Y\in p\,\dd\,\Cl(X)$ and $Y$ is $p$-normal in
$X\big)$ or $\big(Y\in p\,\dd\dd\,\Cl(X)$ and $X$ is $p$-normal on
$Y\big)$.
\end{enumerate}
\end{theorem}

\begin{pf}
(1) Let, first, $Y\in d\dd\,\cD(X)$, $\tau_1<_S\tau_2$ and $Y$ is
$p$-normal. If $A\in\co\tau_1'$, $B\in\co\tau_2'$ and $A\cap
B=\vnth$, then there are disjoint $U''\in\tau_2'$, $V''\in\tau_1'$
such that $A\sbs U''$ and $B\sbs V''$.

Hence, by the second part of Lemma~2.16 there are $U\in\tau_2$,
$V\in\tau_1$ such that $A\sbs U$, $B\sbs V$ and $U\cap V=\vnth$.

Now, let for example, $Y\in\tau_2$, $A\in\co\tau_1'$,
$B\in\co\tau_2'$ and $A\cap B=\vnth$. Then there are
$U'\in\tau_2'$, $V'\in\tau_1'$ such that $A\sbs U'$, $B\sbs V'$
and $U'\cap V'=\vnth$. Clearly, $U'\in\tau_2$ and if $V\in\tau_1$,
$V\cap Y=V'$, then $U'$ and $V$ are the required disjoint
respectively $2$-open and $1$-open neighborhoods of $A$ and $B$.

(2) Let $A\in\co\tau_1'$, $B\in\co\tau_2'$ and $A\cap B=\vnth$.
Since $Y=\tau_1\cl Y\cap\tau_2\cl Y$ and $A,B\sbs Y$, we have
$\tau_1\cl A\cap\tau_2\cl B=\vnth$. Hence, by condition, there are
disjoint sets $U\in\tau_2$, $V\in\tau_1$ such that $\tau_1\cl
A\cap Y=A\sbs U$ and $\tau_2\cl B\cap Y=B\sbs V$.

Finally, if $A\!\in\!\co\tau_1'$, $B\!\in\!\co\tau_2'$ and $A\cap
B\!=~\!\!\vnth$, then $\tau_1\cl A\cap\tau_2\cl B=\vnth$. Since
$$  \tau_1\cl A=\tau_1\cl(\tau_1\cl A\cap Y) \;\;\text{and}\;\;
        \tau_2\cl B=\tau_2\cl(\tau_2\cl B\cap Y),       $$
the sets $\tau_1\cl A$ and $\tau_2\cl B$ are respectively
$1$-con\-cen\-tra\-ted and $2$-con\-cen\-tra\-ted on $Y$ and since
$X$ is $p$-normal on $Y$, there are disjoint sets $U\in\tau_2$,
$V\in\tau_1$ such that $A\sbs\tau_1\cl A\sbs U$ and $B\sbs
\tau_2\cl B\sbs V$.
\end{pf}

\begin{corollary}{4.7}
For a $\BsS$ $(Y,\tau_1',\tau_2')$ of a $\BS$ $(X,\tau_1,\tau_2)$
we have:
\begin{enumerate}
\item[(1)] If $Y\in d\dd\,\cD(X)$ and $\tau_1<_S\tau_2$ or if
$Y\in\tau_1\cup\tau_2$, then $Y$ is $p$-normal if and only if $Y$
is $p$-strongly normal in $X$.

\item[(2)] If $Y\in p\,\dd\,\Cl(X)$, then $Y$ is $p$-normal in $X$
if and only if $Y$ is $p$-strongly normal in $X$.
\end{enumerate}
\end{corollary}

\begin{pf}
(1) and (2) follow respectively from (1) and (2) of Theorem~4.6
taking into account the implications after Definition~2.11.
\end{pf}

\begin{corollary}{4.8}
If $(Y,\tau_1',\tau_2')$ is a $\BsS$ of a $p$-normal $\BS$
$(X,\tau_1,\tau_2)$ and $Y\in p\,\dd\,\Cl(X)$, then $Y$ is
$p$-strongly normal in $X$.
\end{corollary}

\begin{pf}
By implications after Definition~2.11, $Y$ is $p$-normal in $X$
and so it remains to use (2) of Corollary~4.7.

Corollary~4.7 is the reinforcement of Corollary~3.2.6 in [8].
\end{pf}
\vskip+0.2cm

The following simple statement is the reinforcement of the
hereditarily $p$-normality.

\begin{proposition}{4.9}
Every $\BsS$ $(Y,\tau_1',\tau_2')$ of a hereditarily $p$-normal
$\BS$ $(X,\tau_1,\tau_2)$ is $p$-strongly normal in $X$.
\end{proposition}

\begin{pf}
Let $A\in\co\tau_1'$, $B\in\co\tau_2'$ and $A\cap B=\vnth$. Then
$$  (\tau_1\cl A\cap B)\cup (A\cap \tau_2\cl B)=\vnth   $$
and by (5) of Definition~1.1, there are disjoint sets
$U\in\tau_2$, $V\in\tau_1$ such that $A\sbs U$ and $B\sbs V$. Thus
$Y$ is $p$-strongly normal in $X$.
\end{pf}
\vskip+0.2cm

In counterpart to the fact that not every $\BsS$ of a $p$-normal
$\BS$ is $p$-normal, we have the following simple

\begin{proposition}{4.10}
A $p$-normal $\BS$ is $p$-normal on its arbitrary $\BsS$.
\end{proposition}

\begin{theorem}{4.11}
If a $\BS$ $(X,\tau_1,\tau_2)$ is $p$-normal on a $\BsS$
$(Y,\tau_1',\tau_2')$, then
\begin{enumerate}
\item[(1)] $Y$ is $p$-normal in $X$ and so, $Y$ is $p$-internally
normal in $X$.

\item[(2)] $Y\in p\,\dd\,\Cl(X)$ implies that $Y$ is $p$-middly
normal.
\end{enumerate}
\end{theorem}

\begin{pf}
(1) Let $A\in\co\tau_1$, $B\in\co\tau_2$ and $A\cap B=\vnth$.
Clearly,
$$    \tau_1\cl(A\cap Y)\sbs\tau_1\cl A=A, \;\;\;
        \tau_2\cl(B\cap Y)\sbs\tau_2\cl B=B, \;\;\text{and} $$
\begin{eqnarray*}
    \tau_1\cl(A\cap Y) &&\hskip-0.6cm \!\sbs\!
        \tau_1\cl\big(\tau_1\cl(A\cap Y)\cap Y\big) \;\;\text{since}\;\;
            A\cap Y\!\!\sbs\!\tau_1\cl(A\cap Y)\cap Y, \\
    \tau_2\cl(B\cap Y) &&\hskip-0.6cm \!\sbs\!
        \tau_2\cl\big(\tau_2\cl(B\cap Y)\cap Y\big) \;\;\text{since}\;\;
            B\cap Y\!\!\sbs\!\tau_2\cl(B\cap Y)\cap Y.
\end{eqnarray*}
Hence $\tau_1\cl(A\cap Y)$ is $1$-concentrated on $Y$,
$\tau_2\cl(B\cap Y)$ is $2$-concentrated on $Y$ and
$\tau_1\cl(A\cap Y)\cap\tau_2\cl(B\cap Y)=\vnth$. Since $X$ is
$p$-normal on $Y$, there are disjoint sets $U\in\tau_2$, $V\in
\tau_1$ such that
$$  A\cap Y\sbs\tau_1\cl(A\cap Y)\sbs U, \;\;\;
        B\cap Y\sbs\tau_2\cl(B\cap Y)\sbs V.        $$
Thus $Y$ is $p$-normal in $X$. The rest follows from implications
after Definition~2.11.

(2) Let
$$  A'=\tau_1'\cl\tau_2'\nt A\in(1,2)\dd\,\cC\cD(Y), \;\;\;
        B'=\tau_2'\cl\tau_1'\nt B\in(2,1)\dd\,\cC\cD(Y)       $$
in $Y=\tau_1\cl Y\cap\tau_2\cl Y$ and $A'\cap B'=\vnth$. Let us
consider the disjoint sets $A=\tau_1\cl\tau_2'\nt A'$ and
$B=\tau_2\cl\tau_1'\nt B'$. Since $A\in\co\tau_1$,
$A=\tau_1\cl\tau_2'\nt A'$ for the set $\tau_2'\nt A'\sbs Y$ and
$B\in\co\tau_2$, $B=\tau_2\cl\tau_1'\nt B'$ for the set
$\tau_1'\nt B'\sbs Y$, by the first proposition in Remark~2.12,
$$  A\sbs\tau_1\cl(A\cap Y), \;\;\; B\sbs\tau_2\cl(A\cap Y) $$
and since $X$ is $p$-normal on $Y$, there are disjoint sets
$U\in\tau_2$, $V\in \tau_1$ such that $A\sbs U$ and $B\sbs V$.
Clearly,
$$  A'\sbs U'=U\cap Y\in\tau_2', \;\;\; B'\sbs V'=V\cap Y\in\tau_1' $$
and $U'\cap V'=\vnth$. Thus $Y$ is $p$-middly normal.
\end{pf}

\begin{corollary}{4.12}
If a $\TS$ $(X,\tau)$ is normal on a $\TsS$ $(Y,\tau')$, where
$Y\in\co\tau$, then $(Y,\tau')$ is $\bold{k}$--normal in the sense
of $[21]$.
\end{corollary}

\begin{corollary}{4.13}
Let $(Y,\tau_1',\tau_2')$ be a $\BsS$ of a $\BS$
$(X,\tau_1,\tau_2)$. Then
\begin{enumerate}
\item[(1)] For $Y\in p\,\dd\,\Cl(X)$ we have:
$$  \xymatrix{\text{$X$ is $p$-normal on $Y$} \ar@{=>}[r] &
            \text{$Y$ is $p$-normal in $X$} \ar@{<=>}[d] \\
    & \text{$Y$ is $p$-strongly normal in $X$}. }       $$

\item[(2)] For $Y\in\co\tau_1\cap\co\tau_2$ the following
conditions are equivalent: $X$ is $p$-normal on $Y$, $Y$ is
$p$-normal in $X$ and $Y$ is $p$-strongly normal in $X$.
\end{enumerate}
\end{corollary}

\begin{pf}
(1) The implication is given by (1) of Theorem~4.11, and the
equivalence is given by (2) of Corollary~4.7.

(2) Let $Y$ be $p$-normal in $X$, $A\in\co\tau_1$,
$B\in\co\tau_2$, $A\cap B=\vnth$, $A\sbs\tau_1\cl(A\cap Y)$ and
$B\sbs\tau_2\cl(B\cap Y)$. Then
$$  A\sbs\tau_1\cl A\cap\tau_1\cl Y=A\cap Y, \;\;\;
        B\sbs\tau_2\cl B\cap\tau_2\cl Y=B\cap Y,        $$
so that $A\sbs Y$, $B\sbs Y$ and so there are disjoint sets
$U\in\tau_2$, $V\in \tau_1$ such that $A=A\cap Y\sbs U$ and
$B=B\cap Y\sbs V$. Thus $X$ is $p$-normal on $Y$ and since
$\co\tau_1\cap\co\tau_2\sbs p\,\dd\,\Cl(X)$, it remains to use (1)
above.
\end{pf}

\begin{proposition}{4.14}
If $(Y,\tau_1',\tau_2')$ is a $\BsS$ of a $p$-regular $\BS$
$(X,\tau_1,\tau_2)$, $Y\in p\,\dd\,\Cl(X)$ and $(\tau_1\cl
Y\cup\tau_2\cl Y,\tau_1'',\tau_2'')$ is $p$-internally normal in
$X$, then $Y$ is $p$-Tychonoff.
\end{proposition}

\begin{pf}
Since $X$ is $p$-regular, $Y$ is $p$-regular too and by
Theorem~3.2 in [22], it suffices to prove only that $Y$ is
$p$-middly normal. Let
$$  A'=\tau_1'\cl\tau_2'\nt A'\in(1,2)\dd\,\cC\cD(Y)  $$
and
$$  B'=\tau_2'\cl\tau_1'\nt B'\in(2,1)\dd\,\cC\cD(Y).     $$
Clearly, $Y=\tau_1\cl Y\cap\tau_2\cl Y$ implies that the sets
$$  A=\tau_1\cl\tau_2'\nt A' \;\;\text{and}\;\; B=\tau_2\cl\tau_1'\nt B'    $$
are disjoint and $A,B\sbs\tau_1\cl Y\cup\tau_2\cl Y$. Since
$(\tau_1\cl Y\cup\tau_2\cl Y,\tau_1'',\tau_2'')$ is $p$-internally
normal in $X$, there are disjoint sets $U\in\tau_2$, $V\in \tau_1$
such that $A\sbs U$ and $B\sbs V$. Let $U'=U\cap Y$, $V'=V\cap Y$.
Then $A\cap Y=A'\sbs U'$, $B\cap Y=B'\sbs V'$ and so $Y$ is
$p$-middly normal.
\end{pf}

\begin{theorem}{4.15}
Let $(Y,\tau_1'',\tau_2'')\sbs (p\,\dd\cl
Y,\tau_1',\tau_2')\sbs(X,\tau_1,\tau_2)$. Then $X$ is $p$-normal
on $p\,\dd\cl Y\lra X$ $p$-normal on $Y\lra p\,\dd\cl Y$ is
$p$-normal on $Y$.
\end{theorem}

\begin{pf}
Let $A\in\co\tau_1$, $B\in\co\tau_2$, $A\cap B=\vnth$,
$A=\tau_1\cl(A\cap Y)$ and $B=\tau_2\cl(B\cap Y)$. Then
$$  A\sbs\tau_1\cl(A\cap p\,\dd\cl Y), \;\;\;
        B\sbs\tau_2\cl(B\cap p\,\dd\cl Y)       $$
and since $X$ is $p$-normal on $p\,\dd\cl Y$, there are disjoint
sets $U\in\tau_2$, $V\in\tau_1$ such that $A\sbs U$ and $B\sbs V$.
Thus, $X$ is $p$-normal on $Y$.

Now, let $A\in\co\tau_1'$, $B\in\co\tau_2'$, $A\cap B=\vnth$,
$A=\tau_1'\cl(A\cap Y)$ and $B=\tau_2'\cl(B\cap Y)$. Then
$\tau_1\cl A\cap\tau_2\cl B=\vnth$ since $p\,\dd\cl Y=\tau_1\cl
Y\cap \tau_2\cl Y$. Moreover,
$$  \tau_1\cl(\tau_1\cl A\cap Y)=\tau_1\cl A, \;\;\;
        \tau_2\cl(\tau_2\cl B\cap Y)=\tau_2\cl B,       $$
that is, $\tau_1\cl A$ is $1$-concentrated on $Y$ and $\tau_2\cl
B$ is $2$-concentrated on $Y$. Since $X$ is $p$-normal on $Y$,
there are disjoint sets $U\in\tau_2$, $V\in\tau_1$ such that
$\tau_1\cl A\sbs U$ and $\tau_2\cl B\sbs V$. It is evident that
$$  \tau_1\cl A\cap Y=A\sbs U\cap p\,\dd\cl Y=U'\in\tau_2'  $$
and
$$  \tau_2\cl B\cap Y=B\sbs V\cap p\,\dd\cl Y=V'\in\tau_1'      $$
so that $p\,\dd\cl Y$ is $p$-normal on $Y$.
\end{pf}

\begin{definition}{4.16}
We say that topologies $\tau_1$ and $\tau_2$ on a set $X$ are
$p$-strongly $S$-related (briefly $\tau_1S(p)\tau_2)$ if the
$S$-relation is hereditary with respect to $p$-closed subsets of
$X$.
\end{definition}

It is clear that $\tau_1S(p)\tau_2\lra\tau_1S(i)\tau_2$ in the
sense of Definition~2.1.9 in [8] $((\tau_1S(i)\tau_2)$ means that
the $S$-relation is hereditary with respect to $i$-closed subsets
of $X)$, and $\tau_1<_{S(p)}\tau_2\llra
(\tau_1S(p)\tau_2\wedge\tau_1\sbs \tau_2)$.

\begin{theorem}{4.17}
If a $p$-regular $\BS$ $(X,\tau_1<_{S(p)}\tau_2)$ is $p$-normal on
a $\BsS$ $(Y,\tau_1'',\tau_2'')$, then $(Y,\tau_1'',\tau_2'')$ is
$p$-Tychonoff.
\end{theorem}

\begin{pf}
By the second implication in Theorem~4.15, the $\BS$ $(p\,\dd\cl
Y,\tau_1',\tau_2')$ is also $p$-normal on $(Y,\tau_1'',\tau_2'')$.
Since $Y\in d\dd\,\cD(p\,\dd\cl Y)$, it is obvious that the $\BS$
$(p\,\dd\cl Y,\tau_1'<_S\tau_2')$ is $(d,p)$-densely normal.
Hence, by Proposition~2.18, $(p\,\dd\cl Y,\tau_1',\tau_2')$ is
$p$-middly normal. Since $X$ is $p$-regular, $p\,\dd\cl Y$ is also
$p$-regular, and so, by Theorem~3.2 in [22], $(p\,\dd\cl
Y,\tau_1',\tau_2')$ is $p$-Tychonoff. Thus,
$(Y,\tau_1'',\tau_2'')$ is $p$-Tychonoff.
\end{pf}

\begin{theorem}{4.18}
For a $\BsS$ $(Y,\tau_1',\tau_2')$ of a $\BS$ $(X,\tau_1,\tau_2)$,
where $Y=X_1^i\cap X_2^i\in d\dd\,\cD(X)$, the following
conditions are equivalent:
\begin{enumerate}
\item[(1)] $X$ is $p$-middly normal.

\item[(2)] $X$ is $p$-normal on $Y$.

\item[(3)] $X$ is $(d,p)$-densely normal.
\end{enumerate}
\end{theorem}

\begin{pf}
Since $Y\in d\dd\,\cD(X)$, $(2)\lra (3)$. Besides, by
proposition~2.18, $(3)\lra (1)$. Hence, it remains to prove only
that $(1)\lra (2)$. Let $A\in\co\tau_1$, $B\in\co\tau_2$, $A\cap
B=\vnth$, $A=\tau_1\cl(A\cap Y)$, and $B=\tau_2\cl(B\cap Y)$. But
if $x\in Y$, then $\{x\}\in\tau_1\cap\tau_2$ so that $A\cap
Y\in\tau_2$, $B\cap Y\in\tau_1$ and so
$$  A=\tau_1\cl(A\cap Y)\in (1,2)\dd\,\cC\cD(X), \;\;
        B=\tau_2\cl(B\cap Y)\in (2,1)\dd\,\cC\cD(X).        $$

Since $X$ is $p$-middly normal, there are disjoint sets
$U\in\tau_2$, $V\in\tau_1$ such that $A\sbs U$ and $B\sbs V$.
\end{pf}

\begin{theorem}{4.19}
A $\BsS$ $(Y,\tau_1',\tau_2')$ of a $\BS$ $(X,\tau_1,\tau_2)$ is
$p$-strongly normal in $X$ if and only if $(Y,\tau_1',\tau_2')$ is
$p$-normal and for each       \linebreak      $(1,2)$-l.u.s.c.
function $f:(Y,\tau_1',\tau_2')\to (I,\om)$ there is a
$Y\dd\,(1,2)$-l.u.s.c. function $f^{*}:(X,\tau_1,\tau_2)\to
(I,\om)$ which is an extension of $f$, that is, $f^{*}(x)=f(x)$
for each point $x\in Y$.
\end{theorem}

\begin{pf}
If $Y$ is $p$-strongly normal in $X$, then by implications after
Definition~2.11, $Y$ is $p$-normal. Let us consider the families
$$  \cS_i=\big\{E:\;E\sbs X\setminus Y\big\}\cup\tau_i        $$
as subbases on new topologies $\gm_1$ and $\gm_2$ on $X$. First,
let us prove that $(X,\gm_1,\gm_2)$ is $p$-normal. Clearly,
$Y\in\co\gm_1\cap\co\gm_2$ and $\gm_i\big|_Y=\tau_i\big|_Y$. If
$A\in\co\gm_1$, $B\in\co\gm_2$ and $A\cap B=\vnth$, then $A,B\sbs
Y$ and since $\gm_i\big|_Y=\tau_i\big|_Y$, we have
$A\in\co\tau_1'$, $B\in\co\tau_2'$. Since $Y$ is $p$-strongly
normal in $(X,\tau_1,\tau_2)$, there are disjoint sets
$U\in\tau_2$, $V\in\tau_1$ such that $A\sbs U$ and $B\sbs V$. But
$\tau_i\sbs\gm_i$ so that $U\in\gm_2$, $V\in\gm_1$ and so
$(X,\gm_1,\gm_2)$ is $p$-normal. Since $Y\in\co\gm_1\cap\co\gm_2$,
by Theorem~2.6 in [14], every $(1,2)$-l.u.s.c. function
$f:(Y,\gm_1',\gm_2')\to (I,\om)$ can be extended to a
$(1,2)$-l.u.s.c. function $f^{*}:(X,\gm_1,\gm_2)\to (I,\om)$. Let
$y\in Y$ be any point. Since the families
$\cB_i=\{U\in\tau_i:\;y\in U\}$ are bases of $(X,\gm_1,\gm_2)$  at
$y$, the function $f^{*}$ is $(1,2)$-l.u.s.c. at $y$ with respect
to the original topologies $\tau_1$ and $\tau_2$.

Conversely, let $A\in\co\tau_1'$, $B\in\co\tau_2'$ and $A\cap
B=\vnth$. Since $Y$ is   \linebreak        $p$-normal, by
Theorem~2.7 in [13], there exists a $(1,2)$-l.u.s.c. function
$f:(Y,\tau_1',\tau_2')\to (I,\om)$ such that $f(A)=0$ and
$f(B)=1$. By condition, $f$ can be extended to a
$Y\dd\,(1,2)$-l.u.s.c. function $f^{*}:(X,\tau_1,\tau_2)\to
(I,\om)$. Let us consider the disjoint sets
$$  U=\tau_2\nt\Big\{x\in X:\;f^{*}(x)<\frac{1}{2}\Big\}    $$
and
$$  V=\tau_1\nt\Big\{x\in X:\;f^{*}(x)>\frac{1}{2}\Big\}.   $$
Then $A\sbs U$, $B\sbs V$, that is, $Y$ is $p$-strongly normal in
$X$.
\end{pf}
\vskip+0.2cm

Recall that by Definition~2.4 in [14], a $\BsS$
$(Y,\tau_1',\tau_2')$ of a $\BS$ $(X,\tau_1,\tau_2)$ is
$\SC^{*}$-embedded in $X$ if every bounded real-valued
$(1,2)$-l.u.s.c. function on $Y$ can be extended to a
$(1,2)$-l.u.s.c. function on~$X$.

Taking into account the condition of Theorem~4.19, we introduce

\begin{definition}{4.20}
We say that a $\BsS$ $(Y,\tau_1',\tau_2')$ is weakly
$\SC^{*}$-em\-bed\-ded (briefly, $\WSC^{*}$-embedded) in
$(X,\tau_1,\tau_2)$ if every bounded real-valued $(1,2)$-l.u.s.c.
function on $Y$ can be extended to a real-valued function on $X$
which is $Y\dd\,(1,2)$-l.u.s.c.
\end{definition}

\begin{corollary}{4.21}
If $(X,\tau_1,\tau_2)$ is hereditarily $p$-normal, then every
$\BsS$ $(Y,\tau_1',\tau_2')$ of $X$ is $p$-strongly normal and
$\WSC^{*}$-embedded in $X$.
\end{corollary}

\begin{pf}
Follows directly from Proposition~4.9 and Theorem~4.19.
\end{pf}

\begin{theorem}{4.22}
Let $(Y,\tau_1',\tau_2')$ be a $p$-normal $\BsS$ of a $\BS$
$(X,\tau_1,\tau_2)$ such that $Y$ $\SC^{*}$-embedded in $X$. Then
$Y$ is $p$-normal on $Y$.
\end{theorem}

\begin{pf}
Let $A\in\co\tau_1$, $B\in\co\tau_2$, $A\cap B=\vnth$,
$A=\tau_1\cl(A\cap Y)$ and $B=\tau_2\cl(B\cap Y)$. Let $A_1=A\cap
Y$, $B_1=B\cap Y$. Then $A\in\co\tau_1'$, $B\in\co\tau_2'$ in $Y$
and since $Y$ is $p$-normal, by Theorem~2.7 in [13], there is a
$(1,2)$-l.u.s.c. function $f:(Y,\tau_1',\tau_2')\to (I,\om)$ such
that $f(A_1)=0$ and $f(B_1)=1$. Since $Y$ is $\SC^{*}$-embedded in
$X$, $f$ can be extended to a $(1,2)$-l.u.s.c. function
$f^{*}:(X,\tau_1,\tau_2)\to (I,\om)$. Since $A_1$ is $1$-dense in
$A$ and $B_1$ is $2$-dense in $B$, we have $f^{*}(A)=0$ and
$f^{*}(B)=1$. Indeed, let $x\in A\setminus A_1=\tau_1\cl
A_1\setminus A_1$ be a point such that $f^{*}(x)>0$. Since $f^{*}$
is $1$-l.s.c., there is $U(x)\in\tau_1$ such that $f^{*}(y)>0$ for
each point $g\in U(x)$. Since $f^{*}(A_1)=f(A_1)=0$, we have
$U(x)\cap A_1=\vnth$, which contradicts the fact $x\in\tau_1\cl
A_1$. Similarly, since $f^{*}$ is $2$-u.s.c. and $\tau_2\cl
B_1=B$, we have $f^{*}(B)=1$. Clearly,
$$  U=\Big\{x\in X:\;f^{*}(x)<\frac{1}{2}\Big\}\in\tau_2, \;\;\;
        V=\Big\{x\in X:\;f^{*}(x)>\frac{1}{2}\Big\}\in\tau_1,   $$
$A\sbs U$, $B\sbs V$ and $U\cap V=\vnth$. Thus $Y$ is $p$-normal
in $X$.
\end{pf}

\begin{corollary}{4.23}
If for a $\BS$ $(X,\tau_1<_S\tau_2)$ there exists a $d$-dense
$p$-normal $\BsS$ $(Y,\tau_1',\tau_2')$ which is
$\SC^{*}$-embedded in $X$, then $(X,\tau_1,\tau_2)$ is
$(d,p)$-densely normal and hence, $p$-middly normal.
\end{corollary}

\begin{pf}
By Theorem~4.22, $X$ is $p$-normal on $Y$. Hence, by
De\-fi\-ni\-ti\-on~2.17, $(X,\tau_1,\tau_2)$ is $(d,p)$-densely
normal and therefore, it remains to use Proposition~2.18.
\end{pf}
\vskip+0.2cm

The end of the section is devoted to study of relative $p$-strong
nor\-ma\-li\-ty and relative $\WS$-supernormality together with
their interrelations for both the topological and the
bitopological case (see also [9]).

First of all note that if $(X,\tau)$ is normal, then every closed
subspace of $(X,\tau)$ is $\WS$-supernormal in $X$, and if
$(Y,\tau')$ is $\WS$-supernormal in $X$, then $(Y,\tau')$ is
normal. But for the opposite implication we have the following
elementary

\begin{example}{4.24}
Let us consider the natural $\TS$ $(\bR,\om)$ and an open interval
$(a,b)\sbs \bR$, which is normal. If $A=[c,b)$ and $B=[b,d]$ are
sets, where $a<c<b<d$, then  $A\in \co\om'$ in $((a,b),\om')$,
$B\in \co\om$ and $A\cap B=\vnth$. But for any neighborhood
$V(B)\in \om$ we have $A\cap V(B)\neq\vnth$.
\end{example}

Evidently, if $(X,\tau_1,\tau_2)$ is $p\,$-normal, then for every
subset $Y\in \co\tau_i$ the $\BsS$ $(Y,\tau_1',\tau_2')$ is
$(i,j)\dd\WS$-supernormal in $X$ and so, if $(X,\tau_1,\tau_2)$ is
$p\,$-normal and $Y\in \co\tau_1\cap \co\tau_2$, then
$(Y,\tau_1',\tau_2')$ is $p\,\dd\WS$-supernormal in $X$ (clearly,
for a $\BS$ $(X,\tau_1<\tau_2)$ the last fact is correct for $Y\in
\co\tau_1)$.

\begin{example}{4.25}
\rm Let $(\bR,\om_1,\om_2)$ be the natural $\BS$. Then it is not
difficult to see that every set $Y\in \om_i\setminus \{\vnth\}$ is
$(i,j)\dd\WS$-supernormal in the $p\,$-normal $\BS$
$(\bR,\om_1,\om_2)$, but it is not $(j,i)\dd\WS$-supernormal in
$(\bR,\om_1,\om_2)$. Moreover, if $Y=[a,b]=(-\infty,b]\cap
[a,+\infty)\in p\,\dd\,\Cl(\bR)$, then it is
$p\,\dd\WS$-supernormal in $(\bR,\om_1,\om_2)$.
\end{example}

\begin{proposition}{4.26}
In a $\BS$ $(X,\tau_1<\tau_2)$ we have:
\begin{enumerate}
\item[$(1)$] If $Y\in \tau_2$ and $\tau_1<_N\tau_2$, then the
following implications hold:
$$  \xymatrix{
        \text{$Y$\,is $2\dd\WS$-supernorm.\,in\,$X$} \ar@{=>}[r] \ar@{=>}[d]
            & \text{$Y$\,is\,$(1,2)\dd\WS$-supernorm.\,in\,$X$} \ar@{=>}[d] \\
        \text{$Y$\,is $(2,1)\dd\WS$-supernorm.\,in\,$X$} \ar@{=>}[r]
            & \text{$Y$\,is $1\dd\WS$-supernorm.\,in\,$X$}\,. }   $$

\item[$(2)$] If $\tau_1<_C\tau_2$, then
\begin{eqnarray*}
    &\ds \ \hskip-2cm (Y,\tau_1',\tau_2') \;\mbox{is $p\,$-strongly normal in}\;
            X \lra \\
    &\ds\ \hskip+2cm \lra (Y,\tau_1',\tau_2') \;\mbox{is $1$-strongly
                    normal in}\; X
\end{eqnarray*}
and if $\tau_1<_N\tau_2$, then
$$  \xymatrix{
        Y\;\mbox{is $2$-strongly normal in} \ar@{=>}[r] \ar@{=>}[d]
            & Y\;\mbox{is $d$-strongly normal in}\;X \ar@{=>}[d] \\
        Y\;\mbox{is $p\,$-strongly normal in}\ar@{=>}[r]
            & Y\;\mbox{is $1$-strongly normal in} }     $$
\end{enumerate}
\end{proposition}

\begin{pf}
(1) First of all, let us note that if $Y\in \tau_2$, then by (2)
of Corollary~2.3.13 in [8], $\tau_1<_N\tau_2$ implies that
$\tau_1'<_N\tau_2'$.

For the upper horizontal implication let $A\in \co\tau_1'$, $B\in
\co\tau_2$ and $A\cap B=\vnth$. Since $A\in \co\tau_1'\sbs
\co\tau_2'$ and $(Y,\tau_1',\tau_2')$ is $2\dd\WS$-supernormal in
$X$, there are neighborhoods $U'(A)\in \tau_2'\sbs \tau_2$,
$U'(B)\in \tau_2$ such that $U'(A)\cap U'(B)=\vnth$ and so
$\tau_2\cl U'(A)\cap U'(B)=\vnth$. According to (2) of
Corollary~2.3.12 in [8],
$$  U'(A)\sbs \tau_2\nt\tau_2\cl U'(A)=\tau_1\nt\tau_2\cl U'(A)=
        U''(A)\in \tau_1\sbs\tau_2  $$
and for the set $U(A)=U''(A)\cap Y\in \tau_1'\sbs \tau_2'$ we have
$U(A)\cap \tau_1\cl U'(B)=\vnth$. Moreover, by Corollary~2.3.10 in
[8], $\tau_1<_N\tau_2$ implies that $\tau_1<_C\tau_2$ and hence,
(2) of Corollary~2.2.7 in [8] gives that
$$  U'(B)\sbs \tau_2\nt\tau_1\cl U'(B)=\tau_1\nt\tau_1\cl U'(B)=
        U(B)\in \tau_1.     $$
Now, it is clear that $U(A)\cap U(B)=\vnth$ and, hence,
$(Y,\tau_1',\tau_2')$ is $(2,1)\dd\WS$-supernormal in $X$.

For the lower horizontal implication let $A\in \co\tau_1'$, $B\in
\co\tau_1$ and $A\cap B=\vnth$. Since $A\in \co\tau_1'\sbs
\co\tau_2'$ and $(Y,\tau_1',\tau_2')$ is $(2,1)\dd\WS$-supernormal
in $X$, there are neighborhoods $U(A)\in \tau_1'$, $U'(B)\in
\tau_2$ such that $U(A)\cap U'(B)=\vnth$. Since $U(A)\in
\tau_1'\sbs\tau_2'\sbs \tau_2$, we have $U(A)\cap \tau_2\cl
U'(B)=\vnth$. But
$$  U'(B)\sbs \tau_2\nt\tau_2\cl U'(B)=\tau_1\nt\tau_2\cl U'(B)=
            U(B)\in \tau_1      $$
and $U(A)\cap U(B)=\vnth$. Thus $(Y,\tau_1',\tau_2')$ is
$1\dd\WS$-supernormal in $X$.

For the left vertical implication let $A\in \co\tau_2'$, $B\in
\co\tau_1$ and $A\cap B=\vnth$. Since $B\in \co\tau_1\sbs
\co\tau_2$ and $(Y,\tau_1',\tau_2')$ is $2\dd\WS$-supernormal in
$X$, there are neighborhoods $U'(A)\in \tau_2'$, $U(B)\in \tau_2$
such that $U'(A)\cap U(B)=\vnth$ and so $\tau_2\cl U'(A)\cap
U(B)=\vnth$. Furthermore, according to (2) of Corollary~2.3.12 in
[8],
$$  U'(A)\sbs \tau_2'\nt\tau_2'\cl U'(A)=\tau_1'\nt\tau_2'\cl U'(A)=
            U(A)\in \tau_1'     $$
and $U(A)\cap U(B)=\vnth$ so that $(Y,\tau_1',\tau_2')$ is
$(2,1)\dd\WS$-supernormal in~$X$.

Finally, for the right vertical implication let $A\in \co\tau_1'$,
$B\in \co\tau_1$ and $A\cap B=\vnth$. Since $B\in \co\tau_1\sbs
\co\tau_2$ and $(Y,\tau_1',\tau_2')$ is $(1,2)\dd\WS$-supernormal
in $X$, there are neighborhoods $U'(A)\in \tau_2'$, $U(B)\in
\tau_1$ such that $U'(A)\cap U(B)=\vnth$. Hence $\tau_1\cl
U'(A)\cap U(B)=\vnth$ and by (2) of Corollary~2.3.12 in [8],
$$  U'(A)\sbs \tau_2'\nt\tau_2'\cl U'(A)=\tau_1'\nt\tau_2'\cl U'(A)=
            U(A)\in \tau_1'.     $$
Clearly $U(A)\cap U(B)=\vnth$ and, thus, $(Y,\tau_1',\tau_2')$ is
$1\dd\WS$-supernormal in $X$.

(2) First, let $A,\,B\in \co\tau_1'$ and $A\cap B=\vnth$. Since
$B\in\co\tau_1'\sbs \co\tau_2'$ and $(Y,\tau_1',\tau_2')$ is
$p\,$-strongly normal in $X$, there are neighborhoods $U'(A)\in
\tau_2$, $U(B)\in \tau_1$ such that $U'(A)\cap U(B)=\vnth$.
Clearly, $\tau_1\cl U'(A)\cap U(B)=\vnth$ and according to (2) of
Corollary~2.2.7 in [8],
$$  U'(A)\sbs \tau_2\nt\tau_1\cl U'(A)=\tau_1\nt\tau_1\cl U'(A)=
            U(A)\in \tau_1,     $$
where $U(A)\cap U(B)=\vnth$. Thus $(Y,\tau_1',\tau_2')$ is
$1$-strongly normal in $X$.

Furthermore, let $A\in \co\tau_1'$, $B\in \co\tau_2'$ and $A\cap
B=\vnth$. Since $A\in \co\tau_1'\sbs \co\tau_2'$ and
$(Y,\tau_1',\tau_2')$ is $2$-strongly normal in $X$, there are
neighborhoods $U(A),\,U'(B)\in \tau_2$ such that $U(A)\cap
U'(B)=\vnth$ and so $U(A)\cap\tau_2\cl U'(B)=\vnth$. Hence, by (2)
of Corollary~2.3.12 in [8],
$$  U'(B)\sbs \tau_2\nt\tau_2\cl U'(B)=\tau_1\nt\tau_2\cl U'(B)=
            U(B)\in \tau_1,     $$
$U(A)\cap U(B)=\vnth$ and, hence, $(Y,\tau_1',\tau_2')$ is
$p\,$-strongly normal in $X$. The horizontal implications follow
from the first implication of (2) since $\tau_1<_N\tau_2$ implies
$\tau_1<_C\tau_2$.
\end{pf}

\begin{proposition}{4.27}
For a $\BsS$ $(Y,\tau_1',\tau_2')$ of a $\BS$ $(X,\tau_1,\tau_2)$
the following conditions are satisfied:
\begin{enumerate}
\item[$(1)$] If $Y\in \co\tau_i$ and $(Y,\tau_1',\tau_2')$ is
$p\,$-strongly normal in $X$, then $(Y,\tau_1',\tau_2')$ is
$(i,j)\dd\WS$-supernormal in $X$.

\item[$(2)$] If $Y\in \tau_j$ and $(Y,\tau_1',\tau_2')$ is
$(i,j)\dd\WS$-supernormal in $X$, then $(Y,\tau_1',\tau_2')$ is
$p\,$-strongly normal in $X$.
\end{enumerate}
\end{proposition}

\begin{pf}
(1) Let $Y\in \co\tau_i$, $A\in \co\tau_i'$, $B\in \co\tau_j$ and
$A\cap B=\vnth$. Then $B'=B\cap Y\in \co\tau_j'$ and since
$(Y,\tau_1',\tau_2')$ is $p\,$-strongly normal in $X$, there are
disjoint sets $U'\in \tau_j$, $V'\in \tau_i$ such that $A\sbs U'$
and $B'\sbs V'$. Evidently, if $U=U'\cap Y$, $V=V'\cup (X\setminus
Y)$, then $U\in \tau_j'$, $V\in \tau_i$, $A\sbs U$, $B\sbs V$ and
$U\cap V=\vnth$ so that $(Y,\tau_1',\tau_2')$ is
$(i,j)\dd\WS$-supernormal in~$X$.

(2) Let $Y\in \tau_j$, $A\in \co\tau_i'$, $B\in\co\tau_j'$ and
$A\cap B=\vnth$. Then $A\cap \tau_j\cl B=\vnth$ and since
$(Y,\tau_1',\tau_2')$ is $(i,j)\dd\WS$-supernormal in $X$, there
are sets $U\in \tau_j'$, $V\in \tau_i$ such that $A\sbs U$,
$\tau_j\cl B\sbs V$ and $U\cap V=\vnth$. Clearly $U\in \tau_j$,
$B\sbs V$ and, hence, $(Y,\tau_1',\tau_2')$ is $p\,$-strongly
normal in~$X$.
\end{pf}

\begin{corollary}{4.28}
For a $\BsS$ $(Y,\tau_1',\tau_2')$ of a $\BS$ $(X,\tau_1,\tau_2)$
the following conditions are satisfied:
\begin{enumerate}
\item[$(1)$] If $Y\in \co\tau_1\cap \co\tau_2$ and
$(Y,\tau_1',\tau_2')$ is $p\,$-strongly normal in $X$, then
$(Y,\tau_1',\tau_2')$ is $p\dd\WS$-supernormal in $X$.

\item[$(2)$] If $Y\in \tau_1\cap\tau_2$ and $(Y,\tau_1',\tau_2')$
is $p\,\dd\WS$-supernormal in $X$, then $(Y,\tau_1',\tau_2')$ is
$p\,$-strongly normal in $X$.

\item[$(3)$] If $Y\in \tau_1\cap\tau_2\cap\co\tau_1\cap\co\tau_2$,
then the following conditions are equivalent:
$(Y,\tau_1',\tau_2')$ is $(1,2)\dd\WS$-supernormal in $X$,
$(Y,\tau_1',\tau_2')$ is $(2,1)\dd\WS$-supernormal in $X$,
$(Y,\tau_1',\tau_2')$ is $p\,\dd\WS$-supernormal in $X$ and
$(Y,\tau_1',\tau_2')$ is $p\,$-strongly normal in $X$.
\end{enumerate}

Therefore, the equivalences remain valid for a $\BS$
$(X,\tau_1<\tau_2)$ if $Y\in \tau_1\cap\co\tau_1$.
\end{corollary}

\begin{pf}
(1) and (2) follow directly from (1) and (2) of Proposition~4.27,
respectively.

(3) By (1) and (2) of Proposition~4.27, if
$Y\in\tau_2\cap\co\tau_1$, then $(Y,\tau_1',\tau_2')$ is
$(1,2)\dd\WS$-supernormal in $X$ if and only if
$(Y,\tau_1',\tau_2')$ is $p\,$-strongly normal in $X$ and if
$Y\in\tau_1\cap\co\tau_2$, then $(Y,\tau_1',\tau_2')$ is
$(2,1)\dd\WS$-supernormal in $X$ if and only if
$(Y,\tau_1',\tau_2')$ is $p\,$-strongly normal in $X$.

The rest is obvious.
\end{pf}

\begin{corollary}{4.29}
For a $\TsS$ $(Y,\tau')$ of a $\TS$ $(X,\tau)$ the following
conditions are satisfied:
\begin{enumerate}
\item[$(1)$] If $Y\in \co\tau$ and $(Y,\tau')$ is strongly normal
in $X$, then $(Y,\tau')$ is $\WS$-supernormal in $X$, and if $Y\in
\tau$ and $(Y,\tau')$ is $\WS$-supernormal in $X$, then
$(Y,\tau')$ is strongly normal in $X$.

\item[$(2)$] If $Y\in \tau\cap \co\tau$, then $(Y,\tau')$ is
$\WS$-supernormal in $X$ if and only if $(Y,\tau')$ is strongly
normal in $X$.
\end{enumerate}
\end{corollary}

In connection with Corollary~4.29 and so, with the bitopological
case, we consider it necessary to give the following elementary

\begin{example}{4.30}
\rm Let $X=\{a,b,c,d,e\}$,
$\tau=\{\vnth,\{c\},\{a,b,c\},\{c,d,e\},      \linebreak
 X\}$ and $Y\in \{a,b,d\}$. Then $(Y,\tau')$ is $\WS$-supernormal
in $X$, but it is not strongly normal in $X$.

Now, if
$\tau\!=\!\{\vnth,\!\{a\},\!\{b\},\{a,b\},\{d,e\},\{a,d,e\},\{b,d,e\},\{a,b,d,e\},X\}$
and $Y=\{a,d,e\}$, then $(Y,\tau')$ is strongly normal in $X$, but
it is not $\WS$-supernormal in $X$.
\end{example}

\vskip+0.5cm
\section*{\textbf{5. $(i,j)$-Relative Real Normalities}}
\vskip+0.2cm

\begin{proposition}{5.1}
Let $(Z,\tau_1'',\tau_2'')\sbs (Y,\tau_1',\tau_2')\sbs
(X,\tau_1,\tau_2)$. Then the following conditions are satisfied:
\begin{enumerate}
\item[(1)] If $Y$ is $p$-normal or $Y$ is $p$-weakly realnormal in
$X$, then $Z$ is $p$-weakly realnormal in $X$.

\item[(2)] If $Y$ is $p$-realnormal in $X$, then $Z$ is also
$p$-realnormal in $X$.
\end{enumerate}

Moreover, if $(Y,\tau_1''',\tau_2''')\sbs
(Y_1,\tau_1'',\tau_2'')\sbs (X_1,\tau_1',\tau_2')\sbs
(X,\tau_1,\tau_2)$, then
\begin{enumerate}
\item[(3)] $Y_1$ is $p$-weakly realnormal in $X_1$ implies that
$Y$ is $p$-weakly realnormal in $X$ too.
\end{enumerate}
\end{proposition}

\begin{pf}
(1) First, let $Y$ be $p$-normal, $A\in\co\tau_1$, $B\in\co\tau_2$
and $A\cap B=\vnth$. Then $A\cap Y\in\co\tau_1'$, $B\cap
Y\in\co\tau_2'$ and by Theorem~2.7 in [13], there is a
$(1,2)$-l.u.s.c. function $f:(Y,\tau_1',\tau_2')\to (I,\om)$ such
that $f(A\cap Y)=0$ and $f(B\cap Z)=1$. It is clear that $f$ is
$Z\dd\,(1,2)$-l.u.s.c., $f(A\cap Z)\sbs\{0\}$ and $f(B\cap
Y)\sbs\{1\}$, so that $Z$ is $p$-weakly realnormal in $X$.

Now, if $Y$ is $p$-weakly realnormal in $X$, then evidently $Z$ is
also $p$-weakly realnormal in $X$.

(2) The condition is obvious.

(3) Let $A\in\co\tau_1$, $B\in\co\tau_2$ and $A\cap B=\vnth$.
Since $Y_1$ is $p$-weakly realnormal in $X_1$, for the disjoint
sets $A\cap X_1\in\co\tau_1'$ and $B\cap X_1\in\co\tau_2'$ there
is a $(1,2)$-l.u.s.c. function $f:(Y_1,\tau_1'',\tau_2'')\to
(I,\om)$ such that $f(A\cap Y_1)\sbs\{0\}$ and $f(B\cap
Y_1)\sbs\{1\}$. Clearly,
$\varphi=f\big|_Y:(Y,\tau_1''',\tau_2''')\to (I,\om)$ is
$(1,2)$-l.u.s.c., $\vf(A\cap Y)\sbs\{0\}$, $\vf(B\cap Y)\sbs\{1\}$
and so $Y$is $p$-weakly realnormal in $X$.
\end{pf}

\begin{corollary}{5.2}
If a $\BsS$ $(Y,\tau_1',\tau_2')$ of a $\BS$ $(X,\tau_1,\tau_2)$
is $p$-normal, then $Y$ is $p$-weakly realnormal in $X$.
\end{corollary}

\begin{pf}
Follows directly from (1) of Proposition~5.1 for $Y=Z$.
\end{pf}

\begin{proposition}{5.3}
If a $\BsS$ $(Y,\tau_1',\tau_2')$ is $p$-realnormal in a $\BS$
\linebreak           $(X,\tau_1,\tau_2)$, then $Y$ is $p$-normal
in $X$.
\end{proposition}

\begin{pf}
Let $A\in\co\tau_1$, $B\in\co\tau_2$ and $A\cap B=\vnth$. Then, by
condition, there is a $Y\dd\,(1,2)$-l.u.s.c. function
$f:(X,\tau_1,\tau_2)\to (I,\om)$ such that $f(A)=\{0\}$ and
$f(B)=\{1\}$. Let $U=\tau_2\nt\{x\in X:\;f(x)<1/2\}$ and
$V=\tau_1\nt\{x\in X:\;f(x)>1/2\}$. Then $A\cap Y\sbs U$ and
$B\cap Y\sbs V$.~\end{pf}

\begin{theorem}{5.4}
A $\BsS$ $(Y,\tau_1',\tau_2')$ of a $\BS$ $(X,\tau_1,\tau_2)$ is
$p$-strongly realnormal in $X$ if and only if $Y$ is $p$-srtongly
normal in $X$.
\end{theorem}

\begin{pf}
First, let $Y$ be $p$-strongly realnormal in $X$,
$A\in\co\tau_1'$, \linebreak       $B\in\co\tau_2'$ and $A\cap
B=\vnth$. Then there is a $Y\dd\,(1,2)$-l.u.s.c. function
$f:(X,\tau_1,\tau_2)\to (I,\om)$ such that $f(A)=\{0\}$ and
$f(B)=\{1\}$. Let $U=\tau_2\nt\{x\in X:\;f(x)<1/2\}$ and
$V=\tau_1\nt\{x\in X:\;f(x)>1/2\}$. Then $A\sbs U$, $B\sbs V$ and
$U\cap V=\vnth$, so that $Y$ is $p$-srtongly normal in $X$.

Conversely, let $A\in\co\tau_1'$, $B\in\co\tau_2'$ and $A\cap
B=\vnth$. Then, by Theorem~4.19, there is a $Y\dd\,(1,2)$-l.u.s.c.
function $f^{*}:(X,\tau_1,\tau_2)\to (I,\om)$ such that
$f^{*}(A)=\{0\}$ and $f^{*}(B)=\{1\}$, that is, $Y$ is
$p$-srtongly normal in $X$.
\end{pf}

\begin{corollary}{5.5}
For a $p$-closed $\BsS$ $(Y,\tau_1',\tau_2')$ of a $\BS$
$(X,\tau_1,\tau_2)$ we have:
$$  \xymatrix{
        \text{$X$ is $p$-normal on $Y$} \ar@{=>}[r] &
            \text{$Y$ is $p$-normal in $X$} \ar@{<=>}[d] \\
        \text{$Y$ is $p$-strongly realnormal in $X$} \ar@{<=>}[r] \ar@{=>}[d] &
            \text{$Y$ is $p$-strongly normal in $X$} \\
        \text{$Y$ is $p$-realnormal in $X$} \ar@{=>}[r] &
            \text{$Y$ is $p$-weakly realnormal in $X$.} }        $$
\end{corollary}

\begin{pf}
The implications and equivalences are immediate consequences of
Theorem~5.4, (1) of Corollary~4.13 and the implications after
Proposition~2.15.
\end{pf}

\begin{theorem}{5.6}
For a $\BsS$ $(Y,\tau_1',\tau_2')$ of a $\BS$ $(X,\tau_1,\tau_2)$,
where      \linebreak         $Y\in\tau_1\cap \tau_2$, we have:
\begin{enumerate}
\item[(1)] $Y$ is $p$-weakly realnormal in $X$ if and only if $Y$
is $p$-realnormal in $X$.

\item[(2)] $Y$ is $p$-normal if and only if $Y$ is $p$-strongly
realnormal in $X$.
\end{enumerate}
\end{theorem}

\begin{pf}
(1) Clearly, it suffices to prove only the implication from the
left to the right. Let $A\in\co\tau_1$, $B\in\co\tau_2$ and $A\cap
B=\vnth$. Then, by (1) of Definition~2.14, there is a
$(1,2)$-l.u.s.c. function $f:(Y,\tau_1',\tau_2')\to (I,\om)$ such
that $f(A\cap Y)\sbs\{0\}$ and $f(B\cap Y)\sbs\{1\}$. Let us
define $f^{*}:(X,\tau_1,\tau_2)\to (I,\om)$ as follows:
$f^{*}(x)=f(x)$ if $x\in Y$ and $f^{*}(X\setminus Y)=\{1\}$. Then,
taking into account (3) of Proposition~2.15, $f^{*}$ is the
required $Y\dd\,(1,2)$-l.u.s.c. function on $X$.

(2) By Theorem~5.4 and the implications after Definition~2.11, it
suffices to prove only the implication from the left to the right.
let $A\in\co\tau_1'$, $B\in\co\tau_2'$ and $A\cap B=\vnth$. Then,
by (1) of Theorem~2.7 in [13], there is a $(1,2)$-l.u.s.c.
function $f:(Y,\tau_1',\tau_2')\to (I,\om)$ such that $f(A)=\{0\}$
and $f(B)=\{1\}$. Let $f^{*}:(X,\tau_1,\tau_2)\to (I,\om)$ be any
function such that $f^{*}\big|_Y=f$. Then $Y\in\tau_1\cap\tau_2$
implies that $f^{*}$ is $Y\dd\,(1,2)$-l.u.s.c., and so $Y$ is
$p$-strongly normal in $X$.
\end{pf}

\begin{proposition}{5.7}
Any $\BsS$ $(Y,\tau_1',\tau_2')$ of a $p$-normal $\BS$
$(X,\tau_1,\tau_2)$ is $p$-realnormal in $X$.
\end{proposition}

\begin{pf}
Let $A\in\co\tau_1$, $B\in\co\tau_2$ and $A\cap B=\vnth$. Then, by
Theorem~2.7 in [13], there is a $(1,2)$-l.u.s.c. function
$f:(X,\tau_1,\tau_2)\to (I,\om)$ such that $f(A)=0$ and $f(B)=1$.
It is clear, that for any $Y\sbs X$, the function $f$ is
$Y\dd\,(1,2)$-l.u.s.c. and, hence, it remains to use (2) of
De\-fi\-ni\-ti\-on~2.14.~\end{pf}

\begin{proposition}{5.8}
If $(Y,\tau_1',\tau_2')$ is $p$-weakly realnormal in a
$\RR\dd\,p\,\dd\TT_1$ $\BS$ $(X,\tau_1,\tau_2)$, then
$(Y,\tau_1',\tau_2')$ is $p$-Tychonoff and so, $p$-almost
completely regular.
\end{proposition}

\begin{pf}
Let $(Y,\tau_1',\tau_2')$ be $p$-weakly realnormal in $X$, $x\in
Y$, $F\in\co\tau_i'$ and $x\,\ol{\in}\,F$. Then
$x\,\ol{\in}\,\tau_i\cl F$ and since $X$ is $\RR\dd\,p\,\dd\TT_1$.
By (1) of Definition~2.14, there is a $Y\dd\,(i,j)$-l.u.s.c.
function $f:(Y,\tau_1',\tau_2')\to (I,\om)$ such that
$$  f(\tau_i\cl F\cap Y)=f(F)\sbs\{0\} \;\;\text{and}\;\;
                f(x)\sbs\{1\}.    $$
Thus $Y$ is $p$-Tychonoff, and it remains to use the implications
after Definition~2.6.
\end{pf}

\begin{theorem}{5.9}
If $(X,\tau_1<_S\tau_2)$ is $p$-normal on $(Y,\tau_1',\tau_2')$,
where $Y\in d\dd\,\cD(X)$, then $Y$ is $p$-realnormal in $X$.
\end{theorem}

\begin{pf}
We will use the condition (2) of Proposition~2.15. Let $A,B\sbs
Y$, $A\neq\vnth\neq B$ and $\tau_1\cl A\cap\tau_2\cl B=\vnth$.
Since $X$ is $p$-normal on $Y$, taking into account Remark~2.12,
there are disjoint sets $U\in\tau_2$, $V\in\tau_1$ such that
$\tau_1\cl A\sbs U$ and $\tau_2\cl B\sbs V$. Clearly, $\tau_1\cl
U\cap V=\vnth=U\cap\tau_2\cl V$, and so, $\tau_1\cl U\cap\tau_2\cl
B=\vnth$. Now, by Remark~1.3.4 in [8], $\tau_1\cl
U\in(1,2)\dd\,\cC\cD(X)$ and by remark after Lemma~2.16,
$\tau_1\cl U$ is $1$-concentrated on $Y$ since $Y\in
2\dd\,\cD(X)$. $X$ is $p$-normal on $Y$ implies that there are
disjoint sets $E\in\tau_2$, $W\in\tau_1$ such that $\tau_1\cl
U\sbs E$ and $\tau_2\cl B\sbs W$. It is evident that
$E\cap\tau_2\cl W=\vnth$. Therefore, we have $\tau_1\cl A\sbs
\tau_1\cl U\sbs E$, $\tau_2\cl B\sbs \tau_2\cl W $ and $\tau_1\cl
U\cap\tau_2\cl W=\vnth$, where $\tau_1\cl
U\in(1,2)\dd\,\cC\cD(X)$, $\tau_2\cl W\in (2,1)\dd\,\cC\cD(X)$.
Since $X$ is $(d,p)$-densely normal, by Proposition~2.18, $X$ is
$p$-middly normal and so, by Theorem~2.2 in [22], there is a
$(1,2)$-l.u.s.c. function $f:(X,\tau_1,\tau_2)\to (I,\om)$ such
that $f(\tau_1\cl U)=0$ and $f(\tau_2\cl W)=1$. Clearly, $f$ is
$Y\dd(1,2)$-l.u.s.c. $f(\tau_1\cl A)\sbs\{0\}$, $f(\tau_2\cl
B)\sbs\{1\}$ and by (2) of Proposition~2.15, $Y$ is $p$-realnormal
in $X$.
\end{pf}

\vskip+0.5cm
\section*{\textbf{6. $(i,j)$-Relative Properties of Compactness Type}}
\vskip+0.2cm

\begin{definition}{6.1}
A $\BsS$ $(Y,\tau_1',\tau_2')$ of a $\BS$ $(X,\tau_1,\tau_2)$ is
$p$-compact      \linebreak       $(i$-compact) in $X$ if every
$p$-open $(i$-open) covering $\boldsymbol{\cU}$ of $X$ contains a
finite subcovering $\boldsymbol{\cU}'$ such that
$Y\sbs\cup\boldsymbol{\cU}'$.
\end{definition}

\begin{proposition}{6.2}
If a $\BsS$ $(Y,\tau_1',\tau_2')$ is $(i,j)\dd\sS\dd\TT_2$ in
$(X,\tau_1,\tau_2)$ and $p$-compact in $X$, then $Y$ is
$(i,j)$-superregular in $X$.
\end{proposition}

\begin{pf}
If $z\in Y$, $z\,\ol{\in}\,F$ and $F\in\co\tau_i$, then for each
point $x\in F$ there are disjoint neighborhoods $U_x(z)\in\tau_i$,
$U_z(x)\in\tau_j$ and hence, the family
$\boldsymbol{\cU}=\big\{U_z(x):\;x\in F\big\}\cup\{X\setminus F\}$
is a $p$-open covering of $X$. Since $Y$ is $p$-compact in $X$,
there are points $x_1,x_2,\dots,x_n\in F$ such that
$Y\sbs\bigcup\limits_{k=1}^n U_z(x_k)\cup(X\setminus F)$. Then
$F\sbs\bigcup\limits_{k=1}^n U_z(x_k)=U\in \tau_j$,
$z\in\bigcap\limits_{k=1}^n U_{x_k}(z)=V\in \tau_i$ and $U\cap
V=\vnth$. Thus $Y$ is $(i,j)$-superregular in $X$.
\end{pf}

\begin{corollary}{6.3}
If $(Y,\tau_1',\tau_2')$ is $p\,\dd\sS\dd\TT_2$ in
$(X,\tau_1,\tau_2)$ and $p$-compact in $X$, then $Y$ is
$p$-superregular in $X$.
\end{corollary}

\begin{proposition}{6.4}
If $(Y,\tau_1',\tau_2')$ is $(1,2)\dd\sS\dd\TT_2$ or
$(2,1)\dd\sS\dd\TT_2$ in $X$ and $p$-compact in $X$, then $Y$ is
$p$-strongly normal in $X$.
\end{proposition}

\begin{pf}
Let $A\in\co\tau_1'$, $B\in\co\tau_2'$, $A\cap B=\vnth$ and, for
example, $Y$ is $(1,2)\dd\sS\dd\TT_2$ in $X$. Clearly, $\tau_1\cl
A\cap B=\vnth$, and by Proposition~6.2, $Y$ is
$(1,2)$-superregular in $X$. Hence, $Y$ is $(1,2)$-strongly
regular in $X$ and so for each point $x\in\tau_1\cl A$ there are
disjoint neighborhoods $U_B(x)\in\tau_2$, $U_x(B)\in\tau_1$. It is
clear that the family
$$  \boldsymbol{\cU}=\big\{U_B(x):\;x\in\tau_1\cl A\big\}\cup\{X\setminus\tau_1\cl A\} $$
is a $p$-open covering of $X$ and since $Y$ is $p$-compact in $X$,
there are points $x_1,x_2,\dots,x_n\in\tau_1\cl A$ such that
$Y\sbs\bigcup\limits_{k=1}^n U_B(x_k)\cup(X\setminus\tau_1\cl A)$.
It is obvious that
$$  A=\tau_1\cl A\cap Y\sbs \bigcup\limits_{k=1}^n U_B(x_k)=U\in\tau_2, \;\;\;
        B\sbs\bigcap\limits_{k=1}^n U_{x_k}(B)=V\in\tau_1    $$
and $U\cap V=\vnth$. Thus $Y$ is $p$-strongly normal in $X$.
\end{pf}

\begin{definition}{6.5}
A $\BsS$ $(Y,\tau_1',\tau_2')$ of a $\BS$ $(X,\tau_1,\tau_2)$ is
$(i,j)$-internally compact in $X$ if every $i$-closed subset of
$X$, contained in $Y$, is $j$-com\-pact.
\end{definition}

\begin{proposition}{6.6}
If $(Y,\tau_1',\tau_2')$ is a $p$-internally compact $\BsS$ in a
$p\,\dd\TT_2$ $\BS$ $(X,\tau_1,\tau_2)$, then $Y$ is
$p$-internally normal in $X$.
\end{proposition}

\begin{pf}
Let $A\in\co\tau_1$, $B\in\co\tau_2$, $A,B\sbs Y$ and $A\cap
B=\vnth$. Then, by condition, $A$ is $2$-compact and $B$ is
$1$-compact. If $x\in A$ is any point, then for each point $y\in
B$ there are disjoint neighborhoods $U_x(y)\in\tau_1$ and
$U_y(x)\in\tau_2$. Clearly, $B\sbs\cup\,U_x(y)$ and since $B$ is
$1$-compact, there are points $y_1,y_2,\dots,y_m\in B$ such that
$B\sbs\bigcup\limits_{k=1}^m U_x(y_k)=U_x(B)\in\tau_1$. If
$U(x)=\bigcap\limits_{k=1}^m U_{y_k}(x)\in\tau_2$, then
$U(x)\cap\,U_x(B)=\vnth$. Hence, for each $x\in A$ there are
disjoint sets $U(x)\in\tau_2$ and $U_x(B)\in\tau_1$. Clearly,
$A\sbs\bigcup\limits_{x\in A} U(x)$ and since $A$ is $2$-compact,
there are points $x_1,x_2,\dots,x_n$ such that
$A\sbs\bigcup\limits_{k=1}^n U(x_k)=U(A)\in\tau_2$. Let
$U(B)=\bigcap\limits_{k=1}^n U_{x_k}(B)\in\tau_1$. Then
$U(A)\cap\,U(B)=\vnth$ and thus $Y$ is $p$-internally normal in
$X$.
\end{pf}

\begin{definition}{6.7}
A $\BS$ $(X,\tau_1,\tau_2)$ is $p$-Lindel\"{o}ff
$(i$-Lindel\"{o}ff) if every $p$-open $(i$-open) covering of $X$
contains a countable subcovering and hence, a $\BsS$
$(Y,\tau_1',\tau_2')$ of a $\BS$ $(X,\tau_1,\tau_2)$ is
$p$-Lindel\"{o}ff $(i$-Lindel\"{o}ff) in $X$ if every $p$-open
$(i$-open) covering $\boldsymbol{\cU}$ of $X$ contains a countable
subcovering $\boldsymbol{\boldsymbol{\cU}}'$ such that
$Y\sbs\cup\boldsymbol{\cU}'$.
\end{definition}

\begin{lemma}{6.8}
Let $\boldsymbol{\xi}=\big\{U_n:\;n=\ol{1,\infty}\big\}\sbs\tau_1$
and
$\boldsymbol{\eta}=\big\{V_n:\;n=\ol{1,\infty}\big\}\sbs\tau_2$ be
any two countable families in a $\BS$ $(X,\tau_1,\tau_2)$. Then
there are disjoint sets $U\in\tau_1$ and $V\in\tau_2$ such that
$\cup\,\boldsymbol{\xi}\setminus\cup\,\big\{\tau_1\cl V_n:\;
            n=\ol{1,\infty}\big\}\sbs U$
and $\cup\,\boldsymbol{\eta}\setminus\cup\,\big\{\tau_2\cl U_n:\;
            n=\ol{1,\infty}\big\}\sbs V$.
\end{lemma}

\begin{pf}
Obviously we may assume that $U_n\sbs U_{n+1}$ and $V_n\sbs
V_{n+1}$ for each $n=\ol{1,\infty}$. Now, let
$G_k=U_k\setminus\tau_1\cl V_k\in\tau_1$ and
$W_k=V_k\setminus\tau_2\cl U_k\in\tau_2$ for each
$k=\ol{1,\infty}$. It follows that the sets
$U=\bigcup\limits_{k=1}^\infty G_k\in\tau_1$ and
$V=\bigcup\limits_{l=1}^\infty W_l\in\tau_2$ are disjoint,
$\cup\,\boldsymbol{\xi}\setminus\cup\,\big\{\tau_1\cl
V_n:\;n=\ol{1,\infty}\big\}\sbs U$ and
$\cup\,\boldsymbol{\eta}\setminus\cup\,\big\{\tau_2\cl U_n:\;
            n=\ol{1,\infty}\big\}\sbs V$.
\end{pf}

\begin{lemma}{6.9}
If a $\BS$ $(X,\tau_1,\tau_2)$ is $p$-Lindel\"{o}ff, then every
$i$-closed subset of $X$ is $j$-Lindel\"{o}ff.
\end{lemma}

\begin{pf}
Let $\boldsymbol{\cU}'=\{U\}\sbs\tau_j$ be any $j$-open covering
of a set $F\in\co\tau_i$. Then
$\boldsymbol{\cU}=\boldsymbol{\cU}'\cup\{X\setminus F\}$ is a
$p$-open covering of $X$ and since $X$ is \linebreak
   $p$-Lindel\"{o}ff, by Definition~6.7 there is a subfamily $\boldsymbol{\cV}\sbs
\boldsymbol{\cU}$ such that $|\boldsymbol{\cV}|\leq\aleph_0$ and
$X=\cup\,\boldsymbol{\cV}$. If
$\boldsymbol{\cV}'=\boldsymbol{\cV}\setminus\{X\setminus F\}$,
then $\boldsymbol{\cV}'\sbs\boldsymbol{\cU}'$,
$|\boldsymbol{\cV}'|\leq\aleph_0$ and
$F\sbs\cup\,\boldsymbol{\cV}'$, so that $F$ is $j$-Lindel\"{o}ff.
\end{pf}

\begin{proposition}{6.10}
If a $\BsS$ $(Y,\tau_1',\tau_2')$ of a $\BS$ $(X,\tau_1,\tau_2)$
is $p$-regular in $X$ and $p$-Lindel\"{o}ff, then $Y$ is
$p$-strongly normal in $X$.
\end{proposition}

\begin{pf}
Let $A\in\co\tau_1'$, $B\in\co\tau_2'$ and $A\cap B=\vnth$. Since
$Y$ is $p$-regular in $X$, for each point $a\in A$ and each point
$b\in B$ there are neighborhoods $U(a)\in\tau_2$ and
$U(b)\in\tau_1$ such that
$$  \big(\tau_1\cl U(a)\cap B\big)\cup \big(\tau_2\cl U(b)\cap A\big)=\vnth.  $$
Since $Y$ is $p$-Lindel\"{o}ff, $A\in\co\tau_1'$ and
$B\in\co\tau_2'$, by Lemma~6.9, there are countable sets $A_1\sbs
A$ and $B_1\sbs B$ such that
\begin{equation}
    A\sbs\cup\,\big\{U(a):\;a\in A_1\big\} \;\;\text{and}\;\;
        B\sbs\cup\,\big\{U(b):\;b\in B_2\big\}.     \tag{$*$}
\end{equation}
Thus, by Lemma~6.8 there are disjoint sets $U\in\tau_2$,
$V\in\tau_1$ such that
$$  \bigcup\limits_{a\in A_1} U(a)\setminus
        \cup\,\big\{\tau_2\cl U(b):\;b\in B_1\big\}\sbs U     $$
and
$$  \bigcup\limits_{b\in B_1} U(b)\setminus
        \cup\,\big\{\tau_1\cl U(a):\;a\in A_1\big\}\sbs V.     $$
Clearly
$$  \Big(A\cap \Big(\bigcup\limits_{b\in B_1} \tau_2\cl U(b)\Big)\Big)\cup
        \Big(B\cap\Big(\bigcup\limits_{a\in A_1} \tau_1\cl
                U(a)\Big)\Big)=\vnth    $$
and by $(*)$, $A\sbs U$, $B\sbs V$. Thus $Y$ is $p$-storngly
normal in $X$.
\end{pf}

\begin{corollary}{6.11}
If a $\BsS$ $(Y,\tau_1',\tau_2')$ of a $\BS$ $(X,\tau_1,\tau_2)$
is $p$-Lindel\"{o}ff and $X$ is $p$-regular $($respectively, $Y$
is $p$-superregular in $X$, $Y$ is $p$-stron\-gly regular in $X$,
$Y$ is $p$-free regular in $X)$, then $Y$ is $p$-strongly normal
in $X$.
\end{corollary}

\begin{pf}
Follows directly from implications after Definition~2.2.
\end{pf}

\begin{corollary}{6.12}
Under the hypotheses of Proposition~$6.10$, $Y$ is $p$-strongly
realnormal in $X$.
\end{corollary}

\begin{pf}
Follows directly from Theorem~5.4.
\end{pf}

\begin{corollary}{6.13}
Under the hypotheses of Proposition~$6.10$, $Y$ is $p$-normal, $Y$
is $p$-quasi normal in $X$ and $Y$ is $p$-normal in $X$.
\end{corollary}

\begin{pf}
The conditions are immediate consequences of implications after
De\-fi\-ni\-ti\-on~2.11.~\end{pf}

\begin{corollary}{6.14}
Under the hypotheses of Corollary~$6.11$, $Y$ is     \linebreak
$p$-nor\-mal, $Y$ is $p$-quasi normal in $X$ and $Y$ is $p$-normal
in $X$.
\end{corollary}

\begin{proposition}{6.15}
Let $\gm_i\sbs\tau_i$ for a $\BS$ $(X,\tau_1,\tau_2)$. If
$(Y,\tau_1',\tau_2')$ is $p$-compact in $(X,\tau_1,\tau_2)$ and
$(Y,\gm_1',\gm_2')$ is $p\,\dd\sS\dd\TT_2$ in $(X,\gm_1,\gm_2)$
$($in particular, if $(X,\gm_1,\gm_2)$ is $p\,\dd\TT_2)$, then
$\gm_i'=\tau_i'$.
\end{proposition}

\begin{pf}
Contrary: $\gm_i'\neq\tau_i'$. Then there is a set
$F\in\co\tau_i'$ and a point $y\in Y\setminus F$ such that
$y\in\gm_i'\cl F$ since $\gm_i\sbs\tau_i$. By condition
$(Y,\gm_1',\gm_2')$ is $p\,\dd\sS\dd\TT_2$ in $(X,\gm_1,\gm_2)$
and by (18) of Definition~2.1 for each point $x\in
X\setminus\{y\}$ there are disjoint sets $U(x)\in\gm_2$
$(U(x)\in\gm_1)$ and $U_x(y)\in\gm_1$ $(U_x(y)\in\gm_2)$. Let
$\Phi=\tau_i\cl F$. Then
$$  \big\{U(x):\;U(x)\in\gm_j,\;x\in X\setminus\{y\}\big\}\cup
        \{X\setminus\Phi\}          $$
is a $p$-open covering of $(X,\tau_1,\tau_2)$ as
$\gm_i\sbs\tau_i$. Since $(Y,\tau_1',\tau_2')$ is $p$-compact in
$(X,\tau_1,\tau_2)$, there is a finite set $A\sbs X\setminus\{y\}$
such that $Y\sbs(X\setminus\Phi)\cup\big(\bigcup\limits_{x\in A}
U(x)\big)$. Hence for the neighborhood $U'(y)\in\gm_i$, where
$U'(y)=\bigcap\limits_{x\in
A}\big\{U_x(y):\;U_x(y)\in\gm_i\big\}$, we have $U'(y)\sbs
X\setminus\Phi$ and so $U'(y)\cap F=\vnth$. If $U(y)=U'(y)\cap Y$,
then $U(y)\in\gm_i'$ and $U(y)\cap F=\vnth$ so that
$y\,\ol{\in}\,\gm_i'\cl F$. A contradiction.
\end{pf}

\begin{proposition}{6.16}
Let $f:(X,\tau_1,\tau_2)\to (X_1,\gm_1,\gm_2)$ be a
$d$-con\-ti\-nu\-o\-us function and $(Y,\!\tau_1',\!\tau_2')$ be
$p$-Lindel\"{o}ff in $(X,\!\tau_1,\!\tau_2)$. Then
$(f(Y),\!\gm_1',\!\gm_2')$ is $p$-Lindel\"{o}ff in
$(X_1,\gm_1,\gm_2)$.
\end{proposition}

\begin{pf}
Let $\boldsymbol{\cU}=\{U_\al:\;\al\in D\}$ be any $p$-open
covering of $X_1$. Then the family
$\boldsymbol{\cV}=f^{-1}(\boldsymbol{\cU})$ is a $p$-open covering
of $X$ and by condition there is a countable subfamily
$\boldsymbol{\cV}'\sbs \boldsymbol{\cV}$ such that
$Y\sbs\cup\boldsymbol{\cV}'$. Clearly, for each
$V_\alpha\in\boldsymbol{\cV}$ there is
$U_\alpha\in\boldsymbol{\cU}$ such that
$V_\alpha=f^{-1}(U_\alpha)$. Hence
$$  f(Y)\sbs f\big(\cup \boldsymbol{\cV}'\big)=
            \cup\,\big\{f(V_\al):\;V_\al\in\boldsymbol{\cV}'\big\}=
        \cup\big\{U_\al:\;\al\in\boldsymbol{\cU}'\big\},     $$
where $\boldsymbol{\cU}'\sbs\boldsymbol{\cU}$ and
$\boldsymbol{\cU}'$ is countable.
\end{pf}

\begin{lemma}{6.17}
If a $\BsS$ $(Y,\tau_1',\tau_2')$ of a $\BS$ $(X,\tau_1,\tau_2)$
is $(j,i)$-su\-per\-re\-gu\-lar in $X$ and $Y$ is
$j$-Lindel\"{o}ff, then $Y$ is is $(i,j)$-perfectly located in
$X$.
\end{lemma}

\begin{pf}
Let $Y\sbs U$, where $U\in\tau_j$. Then $Y\cap F=\vnth$, where
$F=X\setminus U\in\co\tau_j$. Since $Y$ is $(j,i)$-superregular in
$X$, for each point $x\in Y$ there are disjoint neighborhoods
$U(x)\in\tau_j$ and $U_x(F)\in\tau_i$. Clearly, $\{U(x):\;x\in
Y\}$ is a $j$-open covering of $Y$ and by condition there is a
countable set of point $\{x_k\}_{k=1}^\infty\sbs Y$ such that
$Y\sbs\bigcup\limits_{k=1}^\infty U(x_k)$. Moreover, since
$U(x_k)\cap\,U_{x_k}(F)=\vnth$, we have $\tau_i\cl
U(x_k)\cap\,U_{x_k}(F)=\vnth$ and so
$$  \bigcup\limits_{k=1}^\infty\tau_i\cl U(x_k)\cap
        \Big(\bigcap\limits_{k=1}^\infty U_{x_k}(F)\Big)=\vnth.    $$
Therefore
$$  Y\sbs \bigcup\limits_{k=1}^\infty \tau_i\cl U(x_k)=
            \bigcup\limits_{k=1}^\infty\Phi_k\sbs
        X\setminus \bigcap\limits_{k=1}^\infty U_{x_k}(F)\sbs
            X\setminus F=U,     $$
where $\Phi_k\in\co\tau_i$ for each $k=\ol{1,\infty}$. Thus $Y$ is
$(i,j)$-perfectly located in $X$.
\end{pf}

\begin{proposition}{6.18}
Let $(Y,\tau_1',\tau_2')$ be a $(j,i)$-superregular $\BsS$ in a
$\BS$ $(X,\tau_1,\tau_2)$. Then the following conditions are
satisfied:
\begin{enumerate}
\item[(1)] If $Y$ is $j$-Lindel\"{o}ff, then $Y$ is
$j$-Lindel\"{o}ff in $X$ and $Y$ is $(i,j)$-perfectly located in
$X$.

\item[(2)] If $Y$ is $p$-Lindel\"{o}ff in $X$ and $Y$ is
$(j,i)$-perfectly located in $X$, then $Y$ is $i$-Lindel\"{o}ff.
\end{enumerate}
\end{proposition}

\begin{pf}
(1) By Lemma~6.17, it suffices to prove only that if $Y$ is
\linebreak       $j$-Lindel\"{o}ff, then $Y$ is $j$-Lindel\"{o}ff
in $X$. Indeed, if $\boldsymbol{\cU}=\big\{U_\al:\;\al\in D\big\}$
is any $i$-open covering of $X$, then
$\boldsymbol{\cU}_Y=\big\{U_\al\cap Y:\;\al\in D\big\}$ is the
$j$-open covering of $Y$, consisting of $j$-open in $Y$ sets.
Hence, by condition, there is a countable subfamily
$\boldsymbol{\cU}_Y'=\big\{U_{\al_k}\cap
Y:\;k=\ol{1,\infty}\big\}$ such that
$Y=\bigcup\limits_{k=1}^\infty (U_{\al_k}\cap Y)$. Clearly, $Y\sbs
\bigcup\limits_{k=1}^\infty U_{\al_k}$ and so $Y$ is
$j$-Lindel\"{o}ff in $X$.

(2) Now, let $Y$ be $p$-Lindel\"{o}ff in $X$ and $Y$ be
$(j,i)$-perfectly located in $X$. Suppose that
$\boldsymbol{\cU}'=\big\{U_\al':\;\al\in D\big\}\sbs\tau_i'$ is
any $i$-open covering of $Y$. Let $U_\al\in\tau_i$, $U_\al\cap
Y=U_\al'$ for each $\al\in D$ and $W=\bigcup\limits_{\al\in D}
U_\al$. Then $W\in\tau_i$ and $Y\sbs W$. Since $Y$ is
$(j,i)$-perfectly located in $X$, there is a countable family
$\big\{F_n:\;n=\ol{1,\infty}\big\}\sbs\co\tau_j$ such that
$Y\sbs\bigcup\limits_{n=1}^\infty F_n\sbs W$. Hence $X\setminus
W\sbs \bigcap\limits_{n=1}^\infty (X\setminus F_n)$, i.e. $
X\setminus W\sbs X\setminus F_n$ for each $n=\ol{1,\infty}$.
Hence, for each $n=\ol{1,\infty}$, the family
$\boldsymbol{\cU}_n=\big\{U_\al:\;\al\in D\big\}\cup\{X\setminus
F_n\}$ is a $p$-open covering of $X$, and since $Y$ is
$p$-Lindel\"{o}ff in $X$ and $(Y\cap F_n)\cap (X\setminus
F_n)=\vnth$ for each $n=\ol{1,\infty}$, there is a countable
subfamily $\boldsymbol{\cU}_n'=\big\{U_\al'=U_\al\cap
Y:\;U_\al\in\boldsymbol{\cU}_n\big\}\sbs\boldsymbol{\cU}'$ such
that $Y\cap
F_n\sbs\cup\,\big\{U_\al':\;U_\al'\in\boldsymbol{\cU}_n'\big\}$.
Let $\boldsymbol{\cU}_c'=\bigcup\limits_{n=1}^\infty
\boldsymbol{\cU}_n'$. Then
$$  \bigcup\limits_{n=1}^\infty (F_n\cap Y)=Y=
        \cup\,\big\{U_\al':U_\al'\in\boldsymbol{\cU}_c'\big\},  $$ where
$\boldsymbol{\cU}_e'$ is also countable and
$\boldsymbol{\cU}_c'\sbs \boldsymbol{\cU}'$, i.e., $Y$ is
$i$-Lindel\"{o}ff.
\end{pf}

\begin{theorem}{6.19}
Let $f:(X,\tau_1,\tau_2)\to (X_1,\gm_1,\gm_2)$ be a $d$-closed and
$d$-continuous function and
$(Y,\gm_1',\gm_2')\sbs(X_1,\gm_1,\gm_2)$ be $p$-Lindel\"{o}ff in
$X_1$. Moreover, let $f^{-1}(x)$ be $p$-Lindel\"{o}ff in $X$ for
each $x\in X_1$. Then $(f^{-1}(Y),\tau_1',\tau_2')$ is
$p$-Lindel\"{o}ff in $X$.
\end{theorem}

\begin{pf}
Let $\boldsymbol{\cU}=\{U_\al:\;\al\in D\}$ be any $p$-open
covering of a $\BS$ $(X,\tau_1,\tau_2)$. For each point $x\in X_1$
there exists a countable subfamily
$\boldsymbol{\cU}_x\sbs\boldsymbol{\cU}$ such that $f^{-1}(x)\sbs
U_x=\cup\,\boldsymbol{\cU}_x=U_x^1\cup\,U_x^2$, where
$U_x^i\in\tau_i$. Since $f$ is $d$-closed, and so, $i$-closed,
there is an $i$-open neighborhood $W_x^i\in\gm_i$ such that
$f^{-1}(W_x^i)\sbs U_x^i$. Since $(Y,\gm_1',\gm_2')$ is
$p$-Lindel\"{o}ff in $X_1$, there is a countable subset $X'\sbs
X_1$ such that $Y\sbs W=\cup\,\big\{W_x:x\in X'\big\}$, where
$W_x=W_x^1\cup W_x^2$ for each $x\in X'$. Therefore
\begin{eqnarray*}
    &\ds f^{-1}(Y)\sbs f^{-1}(W)\sbs \\
    &\ds \sbs \cup\,\big\{U_x=U_x^1\cup\,U_x^2:\;x\in X'\big\}=
        \cup\,\big\{\cup\,\boldsymbol{\cU}_x:\;x\in X'\big\}.
\end{eqnarray*}
Hence, the subfamily
$\boldsymbol{\cU}'=\cup\{\boldsymbol{\cU}_x:x\in X'\}$ of the
family $\boldsymbol{\cU}$ covers $f^{-1}(Y)$ and since each
$\boldsymbol{\cU}_x$ is countable, $\boldsymbol{\cU}'$ is also
countable and so, $f^{-1}(Y)$ is $p$-Lindel\"{o}ff in $X$.
\end{pf}

\begin{theorem}{6.20}
For a $\BS$ $(X,\tau_1,\tau_2)$ the following conditions are
satisfied:
\begin{enumerate}
\item[(1)] If $(X,\tau_1,\tau_2)$ is $p$-regular,
$(Z,\tau_1'',\tau_2'')\sbs (Y,\tau_1',\tau_2')\sbs
            (X,\tau_1,\tau_2)       $
and $Y$ is $p$-Lindel\"{o}ff, then $Z$ is $p$-realnormal in $X$.

\item[(2)] If $  (Y,\tau_1''',\tau_2''')\sbs
(Y_1,\tau_1'',\tau_2'')\sbs
        (X_1,\tau_1',\tau_2')\sbs (X,\tau_1,\tau_2)$,
$X_1\in d\dd\,\cD(X)$ and $Y_1$ is $p$-Lindel\"{o}ff in $X_1$,
then $Y$ is $p$-Lindel\"{o}ff in~$X$.

\item[(3)] If $(Y,\tau_1''',\tau_2''')\sbs (p\,\dd\cl
Y,\tau_1'',\tau_2'')\sbs
        (X_1,\tau_1',\tau_2')\sbs (X,\tau_1,\tau_2)$
and $X_1\in d\dd\,\cD(X)$, then $Y$ is $p$-Lindel\"{o}ff in $X$ if
and only if $Y$ is $p$-Lindel\"{o}ff in $X_1$.

\item[(4)] If $(Y,\tau_1',\tau_2')\sbs (X,\tau_1,\tau_2)$, $Y$ is
$p$-Lindel\"{o}ff in $X$, $X_1\in p\,\dd\,\Cl(X)$ and $Y_1\sbs
Y\cap X_1$, then $Y_1$ is $p$-Lindel\"{o}ff in $X_1$.
\end{enumerate}
\end{theorem}

\begin{pf}
(1) It is evident that $Y$ is $p$-regular in $X$ and by
Proposition~6.10, $Y$ is $p$-strongly normal in $X$. Therefore, by
Theorem~5.4, $Y$ is $p$-stron\-gly realnormal in $X$. Let
$A\in\co\tau_1$, $B\in\co\tau_2$ and $A\cap B=\vnth$. Then $A\cap
Y\in\co\tau_1'$, $B\cap Y\in\co\tau_2'$ and by (3) of
Definition~2.14, there is a $Y\dd\,(1,2)$-l.u.s.c. function
$f:(X,\tau_1,\tau_2)\to (I,\om)$ such that $f(A\cap Y)=\{0\}$ and
$f(B\cap Y)=\{1\}$. Clearly, $f$ is $Z\dd\,(1,2)$-l.u.s.c. and
$f(A\cap Z)\sbs\{0\}$, $f(B\cap Z)\sbs\{1\}$. Thus, by (3) of
Proposition~2.15, $Z$ is $p$-realnormal in $X$.

(2) Let $\boldsymbol{\cU}=\{U_\al:\;\al\in D\}$ be a $p$-open
covering of $X$. Then $\boldsymbol{\cU}_{X_1}=\big\{U_\al\cap
X_1:\;\al\in D\big\}$ is a $p$-open covering of $X_1$ since
$X_1\in d\dd\,\cD(X)$. Hence, by condition, there is a countable
subfamily $\{U_{\al_k}\cap
X_1\}_{k=1}^\infty\sbs\boldsymbol{\cU}_{X_1}$ such that
$Y_1\sbs\bigcup\limits_{k=1}^\infty (U_{\al_k}\cap X_1)$. Now, it
is obvious that $Y\sbs\bigcup\limits_{k=1}^\infty U_{\al_k}$ and
thus $Y$ is $p$-Lindel\"{o}ff in $X$.

(3) Let, first, $Y$ be $p$-Lindel\"{o}ff in $X_1$ and
$\boldsymbol{\cU}=\{U_\al:\al\in D\}$ be a $p$-open covering of
$X$. Then $\boldsymbol{\cU}_{X_1}=\big\{U_\al\cap X_1:\;\al\in
D\big\}$ is a $p$-open covering of $X_1$ since $X_1\in
d\dd\,\cD(X)$. Hence, by condition, there is a countable subfamily
$\{U_{\al_k}\cap X_1\}_{k=1}^\infty\sbs\boldsymbol{\cU}_{X_1}$
such that $Y\sbs\bigcup\limits_{k=1}^\infty (U_{\al_k}\cap X_1)$.
Hence, it is obvious that $Y\sbs\bigcup\limits_{k=1}^\infty
U_{\al_k}$ and so $Y$ is $p$-Lindel\"{o}ff in $X$.

Conversely, let $Y$ be $p$-Lindel\"{o}ff in $X$ and
$\boldsymbol{\cU}_{X_1}=\big\{U_\al':\;\al\in D\big\}$ be a
$p$-open covering of $X_1$. Since $p\,\dd\cl Y=\tau_1\cl
Y\cap\tau_2\cl Y\sbs X_1$, we have
$$  X\setminus X_1\sbs X\setminus p\,\dd\cl Y=
        (X\setminus\tau_1\cl Y)\cup(X\setminus\tau_2\cl Y). $$
Therefore
$$  \boldsymbol{\cU}=\boldsymbol{\cU}_1\cup\,\big\{X\setminus\tau_1\cl Y,
        X\setminus\tau_2\cl Y\big\},        $$
where $\boldsymbol{\cU}_1\cap X_1=\boldsymbol{\cU}_{X_1}$, is a
$p$-open covering of $X$. Since $Y$ is $p$-Lindel\"{o}ff in $X$
and
$$  Y\cap\big((X\setminus\tau_1\cl Y)\cup
        (X\setminus\tau_2\cl Y)\big)=\vnth,     $$
there is a countable subfamily
$\{U_{\al_k}\}_{k=1}^\infty\sbs\boldsymbol{\cU}_1$ such that
$Y\sbs\bigcup\limits_{k=1}^\infty U_{\al_k}$. But $Y\sbs X_1$ and
so $Y\sbs \bigcup\limits_{k=1}^\infty U_{\al_k}'$. Thus, $Y$ is
$p$-Lindel\"{o}ff in $X_1$.

(4) Let $\boldsymbol{\cU}_{X_1}=\{U_\al':\;\al\in D\}$ be a
$p$-open covering of $X_1$ and
$$  \boldsymbol{\cU}=\big\{U_\al:\;U_\al\cap X_1=U_\al',\;\al\in D\big\}\cup
        \big\{X\setminus\tau_1\cl X_1,X\setminus\tau_2\cl X_1\big\}.     $$
Then $\boldsymbol{\cU}$ is a $p$-open covering of $X$ and since
$Y$ is $p$-Lindel\"{o}ff in $X$, there is a countable subfamily
$\{U_{\al_k}\}_{k=1}^\infty\sbs\boldsymbol{\cU}$ such that
$Y_1\sbs Y\sbs\bigcup\limits_{k=1}^\infty U_{\al_k}$. Since
$Y_1\sbs X_1$, we have
$$  Y_1\cap\big((X\setminus\tau_1\cl X_1)\cup
            (X\setminus\tau_2\cl X_1)\big)=\vnth,        $$
so that $\{U_{\al_k}\cap
X_1\}_{k=1}^\infty=\{U_{\al_k}'\}_{k=1}^\infty\sbs\boldsymbol{\cU}_{X_1}$.
Thus $Y_1$ is $p$-Lindel\"{o}ff in~$X_1$.~\end{pf}

\begin{corollary}{6.21}
For a $\BS$ $(X,\tau_1,\tau_2)$ the following conditions are
satisfied:
\begin{enumerate}
\item[(1)] If $Y\in d\dd\,\cD(X)$, then $Y$ is $p$-Lindel\"{o}ff
in $X$ if and only if $Y$ is $p$-Lindel\"{o}ff.

\item[(2)] If $(Y,\tau_1'',\tau_2'')\sbs (Z,\tau_1',\tau_2')\sbs
            (X,\tau_1,\tau_2)   $
and $Z$ is $p$-Lindel\"{o}ff, then $Y$ is $p$-Lindel\"{o}ff in
$X$.
\end{enumerate}
\end{corollary}

\begin{pf}
(1) It suffices to suppose in (3) of Theorem~6.20 that
$Y=p\,\dd\cl Y=X_1$.

(2) Suppose in (2) of Theorem~6.20 that $Y_1=X_1=Z$.
\end{pf}

\begin{remark}{6.22}
It is clear that Theorem~6.20 remains correct if
 \linebreak   $p$-Lindel\"{o}ffnees is changed by $p$-compactness.
\end{remark}

\begin{definition}{6.23}
Let $(Y,\tau_1',\tau_2')$ be a $\BsS$ of a $\BS$
$(X,\tau_1,\tau_2)$. Then
\begin{enumerate}
\item[(1)] $Y$ is $(i,j)$-quasi paracompact in $X$ if for each
$i$-open covering $\boldsymbol{\cV}$ of $X$ one can find a
covering $\boldsymbol{\cU}$ of $Y$ by $i$-open in $Y$ sets which
refines $\boldsymbol{\cV}$ and is $j$-locally finite at each point
of $Y$.

\item[(2)] $Y$ is $(i,j)$-paracompact in $X$ if for each $i$-open
covering $\boldsymbol{\cV}$ of $X$ there exists a family
$\boldsymbol{\cU}\sbs\tau_i$ which refines $\boldsymbol{\cV}$,
$Y\sbs\cup\,\boldsymbol{\cU}$ and $\boldsymbol{\cU}$ is
$j$-locally finite at each point of $Y$.

\item[(3)] $Y$ is $(i,j)$-strongly paracompact in $X$ if for each
$i$-open covering $\boldsymbol{\cV}$ of $X$ there exists an
$i$-open covering $\boldsymbol{\cU}$ of $X$ which refines
$\boldsymbol{\cV}$ and is $j$-locally finite at each point of $Y$.
\end{enumerate}
\end{definition}

It is clear that if $(X,\tau_1,\tau_2)$ is
$(i,j)\dd\RR\RR$-paracompact in the sense of (6) of
Definition~1.1, then any $\BsS$ $Y$ of $X$ is $(i,j)$-strongly
paracompact in $X$, and so $Y$ is $(i,j)$-paracompact in $X$ and
$(i,j)$-quasi paracompact in $X$.

It is clear that if $Y\in\tau_i$, then $Y$ is $(i,j)$-quasi
paracompact in $X$ if and only if $Y$ is $(i,j)$-paracompact in
$X$.

\begin{remark}{6.24}
For a $\BS$ $(Y,\tau_1'<\tau_2')$ of a $\BS$ $(X,\tau_1<\tau_2)$
we have the following implications: {\small
$$  \xymatrix{
        \text{$Y$\,is\,$1$-strongly\,paracomact\,in\,$X$} \ar@{=>}[r] \ar@{=>}[d]
            & \text{$Y$\,is\,$(1,2)$-strongly\,paracompact\,in\,$X$} \ar@{=>}[d] \\
        \text{$Y$\,is\,$1$-paracompact\,in\,$X$} \ar@{=>}[r] \ar@{=>}[d]
            & \text{$Y$\,is\,$(1,2)$-paracompact\,in\,$X$} \ar@{=>}[d] \\
        \text{$Y$\,is\,$1$-quasi\,paracompact\,in\,$X$} \ar@{=>}[r]
            & \text{$Y$\,is\,$(1,2)$-quasi\,paracompact\,in\,$X$}\,, }
            $$ }
and {\small
$$  \xymatrix{
        \text{$Y$\,is\,$(2,1)$-strongly\,paracompact\,in\,$X$} \ar@{=>}[r] \ar@{=>}[d]
            & \text{$Y$\,is\,$2$-strongly\,paracompact\,in\,$X$} \ar@{=>}[d]  \\
        \text{$Y$\,is\,$(2,1)$-paracompact\,in\,$X$} \ar@{=>}[r] \ar@{=>}[d]
            & \text{$Y$\,is\,$2$-paracompact\,in\,$X$} \ar@{=>}[d]  \\
        \text{$Y$\,is\,$(2,1)$-quasi\,paracompact\,in\,$X$} \ar@{=>}[r]
            & \text{$Y$\,is\,$2$-quasi\,paracompact\,in\,$X$}\,. } $$ }
\end{remark}

\begin{theorem}{6.25}
If $(Y,\tau_1',\tau_2')$ is a $\BsS$ of an $(i,j)$-regular $\BS$
\linebreak           $(X,\tau_1,\tau_2)$, $\tau_jC\tau_i$,
$\tau_iN\tau_j$ and $Y$ is $j$-Lindel\"{o}ff in $X$, then $Y$ is
$(i,j)$-paracompact in $X$.
\end{theorem}

\begin{pf}
Let $\boldsymbol{\cV}=\{U\}$ be any $i$-open covering of $X$. Then
for each point $x\in X$ there is a neighborhood $U(x)\in
\boldsymbol{\cV}$ such that $x\in U(x)$. Since
$x\,\ol{\in}\,X\setminus U(x)\in\co\tau_i$ and $X$ is
$(i,j)$-regular, there is $V(x)\in\tau_i$ such that $\tau_j\cl
V(x)\sbs U(x)$ and so, $\tau_j\cl V(x)\cap(Y\setminus
U(x))=\vnth$. Since $\tau_jC\tau_i$, by Theorem~2 in [27],
$\tau_j\cl V(x)$ is a $j$-neighborhood of $x$, i.e., there is
$E(x)\in\tau_j$ such that $E(x)\sbs\tau_j\cl V(x)$. Clearly,
$\boldsymbol{\cU}=\{E(x):\;x\in X\}$ is a $j$-open covering of
$X$, which refines $\boldsymbol{\cV}$. Since $Y$ is
$j$-Lindel\"{o}ff in $X$, there is a countable subfamily
$\boldsymbol{\cU}'=\{E(x_k)\}_{k=1}^\infty\sbs\boldsymbol{\cU}$
such that $Y\sbs\cup\,\boldsymbol{\cU}'$. For each
$k=\ol{1,\infty}$ let $F_k=\bigcup\limits_{l<k} \tau_i\cl E(x_l)$
and $F_0=\vnth$. Then $F_k\in\co\tau_i$ and if $W_k=U(x)\setminus
F_k$, then $W_k\in\tau_i$ for each $k=\ol{1,\infty}$. Moreover,
the family $\boldsymbol{\cW}=\{W_k\}_{k=1}^\infty$ refines
$\boldsymbol{\cV}$.

Hence, it remains to prove that $Y\sbs\cup\boldsymbol{\cW}$ and
$\boldsymbol{\cW}$ is $j$-locally finite at each point $y\in Y$.
Let $y\in Y$ be any point. Since $Y\sbs\cup\,\boldsymbol{\cU}'$
and $\boldsymbol{\cU}'$ refines $\boldsymbol{\cV}$, let $m$ be the
smallest natural number $k\in N^{+}$ such that $y\in U(x_k)$, that
is, $y\in U(x_m)$ and $y\ol{\in}\,U(x_l)$, i.e. $y\in (Y\setminus
U(x_l))$ if $l<m$. Since $y\ol{\in}\,U(x_l)$ and $\tau_j\cl
E(x_l)\sbs \tau_j\cl V(x_l)\sbs U(x_l)$, we obtain, that
$y\ol{\in}\,\tau_j\cl E(x_l)$ for $l<m$. Moreover, since
$\tau_iN\tau_j$ and $E(x)\in\tau_j$, we have $\tau_i\cl
E(x_l)\sbsq \tau_j\cl E(x_l)$ and so $y\ol{\in}\,\tau_j\cl E(x_l)$
for $l<m$. Thus $y\ol{\in}\,F_m=\bigcup\limits_{l<m} \tau_j\cl
E(x_l)$. Therefore, we obtain $y\in U(x_m)\setminus F_m=
W_m\in\boldsymbol{\cW}$ and thus $Y\sbs\cup\boldsymbol{\cW}$.

Finally, let $y\in Y$ be any point. Since
$Y\sbs\cup\,\boldsymbol{\cU}'$, there is $k\in N^{+}$ such that
$y\in E(x_k)$. Then $E(x_k)\cap W_l=\vnth$ for $l>k$. It follows
that $E(x)$ is a $j$-open neighborhood of $y$, which intersects
only finitely many elements of the family $\boldsymbol{\cW}$.
\end{pf}

\begin{definition}{6.26}
We will say that a $\BsS$ $(Y,\tau_1',\tau_2')$ of a $\BS$
$(X,\tau_1,\tau_2)$ is $(i,j)$-strongly quasi Lindel\"{o}ff in $X$
if for every covering $\boldsymbol{\cV}\sbs\tau_i$ of $Y$ there
exists a countable subfamily $\boldsymbol{\cU}\sbs
\boldsymbol{\cV}$ such that $Y\sbs\tau_j\cl\cup\boldsymbol{\cU}$.
\end{definition}

It is clear that for a $\BsS$ $(Y,\tau_1'<\tau_2')$ of a $\BS$
$(X,\tau_1<\tau_2)$, we have: $Y$ is $(2,1)$-strongly quasi
Lindel\"{o}ff in $X$ implies that $Y$ is $1$-strongly quasi
Lindel\"{o}ff in $X$ and $Y$ is $2$-strongly quasi Lindel\"{o}ff
in $X$ implies that $Y$ is $(2,1)$-strongly quasi Lindel\"{o}ff in
$X$.

\begin{theorem}{6.27}
Let $(Y,\tau_1',\tau_2')$ be a $\BsS$ of a $\BS$
$(X,\tau_1,\tau_2)$, $Y\in i\dd\,\cD(X)$, $Y$ be $(j,i)$-strongly
quasi Lindel\"{o}ff in $X$ and $Y$ be $(i,j)$-strongly paracompact
in $X$. Then $X$ is $i$-Lindel\"{o}ff.
\end{theorem}

\begin{pf}
Let $\boldsymbol{\cV}$ be an $i$-open covering of $X$. Since $Y$
is $(i,j)$-strongly paracompact in $X$, there is an $i$-open
covering $\boldsymbol{\cU}$ of $X$ which refines
$\boldsymbol{\cV}$ and such that $\boldsymbol{\cU}$ is $j$-locally
finite at each point $y\in Y$. Let $\boldsymbol{\cW}$ be the
family of all $j$-open subsets of $X$ which intersect only
finitely many elements of $\boldsymbol{\cU}$. Clearly,
$Y\sbs\cup\boldsymbol{\cW}$ and since $Y$ is $(j,i)$-strongly
quasi Lindel\"{o}ff in $X$, we can find a countable subfamily
$\boldsymbol{\cT}\sbs\boldsymbol{\cW}$ such that
$Y\sbs\tau_i\cl\cup\boldsymbol{\cT}$. Furthermore, every set
$U\in\boldsymbol{\cU}$, $U\neq\vnth$ intersects at least one
element of $\boldsymbol{\cT}$. Indeed, if $U\cap V=\vnth$ for each
$V\in\boldsymbol{\cT}$, then $U\cap(\cup\boldsymbol{\cT})=\vnth$
and since $U\in\tau_i$, we have
$U\cap\tau_i\cl(\cup\,\boldsymbol{\cT})=\vnth$. But
$Y\sbs\tau_i\cl\cup\boldsymbol{\cT}$ and so $U\cap Y=\vnth$ which
is impossible since $Y\in i\dd\,\cD(X)$. On the other hand,
$\boldsymbol{\cT}$ is countable and every element of
$\boldsymbol{\cT}$ intersects only finitely many elements of
$\boldsymbol{\cU}$. It follows that $\boldsymbol{\cU}$ is
countable and thus $Y$ is $i$-Lin\-de\-loff.
\end{pf}

\begin{theorem}{6.28}
If $(Y,\tau_1'<\tau_2')$ is a $\BsS$ of a $\BS$
$(X,\tau_1<\tau_2)$, $Y\in 2\dd\,\cD(X)$ and $Y$ is $(2,1)$-quasi
paracompact in $X$, then $Y$ is $(2,1)$-paracompact in $X$.
\end{theorem}

\begin{pf}
Let $\boldsymbol{\cV}$ be any $2$-open covering of $X$,
$\boldsymbol{\cU}''=\big\{V_\al'':\;\al\in D\big\}\sbs\tau_2'$
refines $\boldsymbol{\cV}$, $Y=\bigcup\limits_\al V_\al''$ and
$\boldsymbol{\cU}''$ is $1$-locally finite at each point of $Y$.
Let $\boldsymbol{\cU}'=\big\{V_\al':\;\al\in D\big\}\sbs\tau_2$
such that $V_\al'\cap Y=V_\al''$ for each $\al\in D$. Since
$\boldsymbol{\cU}''$ refines $\boldsymbol{\cV}$, for each
$V_\al''$ there is $U_\al\in\boldsymbol{\cV}$ such that
$V_\al''\sbs U_\al$. If
$\boldsymbol{\cU}=\big\{W_\al:\;W_\al=V_\al'\cap\,U_\al,\;\al\in
D\big\}$, then $\boldsymbol{\cU}$ refines $\boldsymbol{\cV}$ and
$Y=\bigcup\limits_\al V_\al''\sbs\bigcup\limits_\al W_\al$. Hence,
it remains to prove only that $\boldsymbol{\cU}$ is $1$-locally
finite at each point $x\in Y$. Contrary: there is a point $x_0\in
Y$ such that for each neighborhood $U(x_0)\in\tau_1$ we have
$|\{W_\al:\;W_\al\in\boldsymbol{\cU},\;W_\al\cap
U(x_0)\neq\vnth\}|\geq \aleph_0$. Since $\boldsymbol{\cU}''$ is
$1$-locally finite at $x_0$, there is $V(x_0)\in\tau_1$ such that
$V(x_0)$ intersects only finitely many elements of
$\boldsymbol{\cU}''$. Hence, there is a finite set
$\{\al_k\}_{k=1}^n$ such that $V(x_0)\cap V_{\al_k}''\neq\vnth$
and $V(x_0)\cap V_\al''=\vnth$ if $\al\neq\al_k$. Let
$\al\neq\al_k$ and $W_\al\cap V(x_0)\neq\vnth$. The inclusion
$\tau_1\sbs\tau_2$ implies that $W_\al\cap V(x_0)\in\tau_2$ and
since $Y\in 2\dd\,\cD(X)$, we have $W_\al\cap V(x_0)\cap
Y\neq\vnth$. On the other hand,
$$  W_\al\cap V(x_0)\cap Y=
        V(x_0)\cap (V_\al'\cap\,U_\al\cap Y)=V(x_0)\cap V_\al''=\vnth $$
as $\al\neq \al_k$. A contradiction. Hence,
$\boldsymbol{\cU}\sbs\tau_2$, $\boldsymbol{\cU}$ refines
$\boldsymbol{\cV}$, $Y\sbs\cup\boldsymbol{\cU}$ and
$\boldsymbol{\cU}$ is $1$-locally finite at each point $x\in Y$.
Thus $Y$ is $(2,1)$-paracompact in $X$.
\end{pf}

\begin{theorem}{6.29}
If $(Y,\tau_1'<\tau_2')$ is $(2,1)$-paracompact in
$(X,\tau_1<\tau_2)$ and $(2,1)$-free regular in $X$, then $Y$ is
$p$-normal in $X$.
\end{theorem}

\begin{pf}
Let $A\in\co\tau_1$, $B\in\co\tau_2$ and $A\cap B=\vnth$. Since
$Y$ is    \linebreak    $(2,1)$-free regular in $X$, for each
point $a\in A$ there are disjoint neighborhoods $U(a)\in\tau_2
\;\;\text{and}\;\; U_a(B\cap Y)\in\tau_1$. Clearly, the family
$$  \boldsymbol{\cV}=\big\{U(a):\;a\in A\big\}\cup\{X\setminus A\}  $$
is a $2$-open covering of $X$ as $\tau_1\sbs\tau_2$. Since $Y$ is
$(2,1)$-paracompact in $X$, there is a family
$\boldsymbol{\cU}\sbs\tau_2$ such that $\boldsymbol{\cU}$ refines
$\boldsymbol{\cV}$, $Y\sbs\cup\boldsymbol{\cU}$ and
$\boldsymbol{\cU}$ is $1$-locally finite at each point $x\in Y$.
Let $\boldsymbol{\cW}=\big\{V\in\boldsymbol{\cU}:\;V\cap
A\neq\vnth\big\}$ and $U=\cup\,\boldsymbol{\cW}\in\tau_2$. Then
$A\cap Y\sbs U$. Clearly, $\boldsymbol{\cW}$ is also $1$-locally
finite at each point $x\in Y$. Hence, if $x\in Y$, then
$x\in\tau_1\cl U$ if and only if $x\in\tau_1\cl V$ for some
$V\in\boldsymbol{\cW}$. It is also obvious that $\boldsymbol{\cW}$
refines $\boldsymbol{\cV}$. Moreover, since
$\boldsymbol{\cW}=\big\{V\in\boldsymbol{\cU}:\;V\cap
A\neq\vnth\big\}$, for each $V\in\boldsymbol{\cW}$ there is
$U(a)\in\boldsymbol{\cV}$ such that $V\sbs U(a)$. Since
$U(a)\cap\,U_a(B\cap Y)=\vnth$, $V\sbs U(a)$ and $U_a(B\cap
Y)\in\tau_1$, we have $\tau_1\cl V\cap\,U_a(B\cap Y)=\vnth$. Hence
$\tau_1\cl U\cap (B\cap Y)=\vnth$. Let $W=X\setminus\tau_1\cl
U\in\tau_1$. Then $B\cap Y\sbs W$, $A\cap Y\sbs U$ and $U\cap
W=\vnth$. Thus $Y$ is $p$-normal in $X$.
\end{pf}

\begin{corollary}{6.30}
If $(Y,\tau_1'<\tau_2')$ is $(2,1)$-paracompact in a
$(2,1)$-re\-gu\-lar $\BS$ $(X,\tau_1<\tau_2)$, then $Y$ is
$p$-normal in $X$.
\end{corollary}

\begin{proposition}{6.31}
Let $(X,\tau_1,\tau_2)$ be a $\BS$. Then the following conditions
are satisfied:
\begin{enumerate}
\item[(1)] If $(Z,\tau_1'',\tau_2'')\sbs (Y,\tau_1',\tau_2')\sbs
            (X,\tau_1,\tau_2)$
and $Y$ is $(i,j)\dd\RR\RR$-pa\-ra\-com\-pact, then $Z$ is
$(i,j)$-quasi paracompact in $X$.

\item[(2)] If $(Y,\tau_1''',\tau_2''')\sbs
(Y_1,\tau_1'',\tau_2'')\sbs
        (X_1,\tau_1',\tau_2')\sbs (X,\tau_1,\tau_2)$
and $Y_1$ is $(i,j)$-quasi paracompact in $X_1$, then $Y$ is
$(i,j)$-quasi paracompact in $X$.
\end{enumerate}
\end{proposition}

\begin{pf}
(1) Let $\boldsymbol{\cV}$ be any $i$-open covering of $X$. Then
$\boldsymbol{\cV}\cap Y=\big\{U\cap Y:\;U\in
\boldsymbol{\cV}\big\}$ is an $i$-open covering of $Y$ and, hence,
there is a $i$-open covering $\boldsymbol{\cU}$ of $Y$ such that
$\boldsymbol{\cU}$ refines $\boldsymbol{\cV}$ and
$\boldsymbol{\cU}$ is $j$-locally finite at each point of $Y$.
Clearly, $\boldsymbol{\cU}\cap Z$ is an $i$-open covering of $Z$,
$\boldsymbol{\cU}\cap Z$ refines $\boldsymbol{\cV}$,
$\boldsymbol{\cU}\cap Z$ is $j$-locally finite at each point of
$Z$ and thus $Z$ is $(i,j)$-quasi paracompact in $X$.

(2) Let $\boldsymbol{\cV}$ be any $i$-open covering of $X$. Then
$\boldsymbol{\cV}\cap X_1=\big\{U\cap
X_1:\;U\in\boldsymbol{\cV}\big\}$ is an $i$-open covering of
$X_1$. Since $Y_1$ is $(i,j)$-quasi paracompact in $X_1$, one can
find a covering $\boldsymbol{\cU}$ of $Y_1$ by $i$-open in $Y_1$
sets which refines $\boldsymbol{\cV}\cap X_1$ and so --
$\boldsymbol{\cV}$, and which is $j$-locally finite at each point
of $Y_1$. It is clear, that $\boldsymbol{\cU}\cap Y=\{V\cap
Y:\;V\in\boldsymbol{\cU}\}$ is a covering of $Y$ by $i$-open in
$Y$ sets, $j$-locally finite at each point of $Y$ and since
$\boldsymbol{\cU}$ refines $\boldsymbol{\cV}\cap X_1$,
$\boldsymbol{\cU}\cap Y$ refines $\boldsymbol{\cV}$, that is, $Y$
is $(i,j)$-quasi paracompact in $X$.
\end{pf}

\begin{theorem}{6.32}
Let $(Y,\tau_1'<\tau_2')$ be a $\BsS$ of a $\BS$
$(X,\tau_1<\tau_2)$. Then  the following conditions are satisfied:
\begin{enumerate}
\item[(1)] If $X$ is $1\dd\TT_2$ and $Y$ is $(2,1)$-quasi
paracompact in $X$, then $Y$ is $(1,2)$-regular.

\item[(2)] If $X$ is $1\dd\TT_2$ and $Y$ is $(2,1)$-paracompact in
$X$, then $Y$ is $(1,2)$-regular in $X$.

\item[(3)] If $Y$ is $(2,1)$-quasi paracompact in $X$ and
$(2,1)$-free regular in $X$, then $Y$ is $p$-quasi normal in $X$.
\end{enumerate}
\end{theorem}

\begin{pf}
(1) Let $x\in Y$, $F\in\co\tau_1'$ and $x\,\ol{\in}\,F$. Then
$x\,\ol{\in}\,\tau_1\cl F\in\co\tau_1$. Since $X$ is $1\dd\TT_2$,
it is $p\,\dd\TT_2$ and so for each point $y\in\tau_1\cl F$ there
are disjoint neighborhoods $U_y(x)\in\tau_1$ and
$U_x(y)\in\tau_2$. Clearly,
$$  \boldsymbol{\cV}=\big\{U_x(y):\;y\in\tau_1\cl F\big\}\cup
        \big\{X\setminus\tau_1\cl F\big\}       $$
is a $2$-open covering of $X$ as $\tau_1\sbs\tau_2$. Since $X$ is
$(2,1)$-quasi paracompact in $X$, there is a family
$\boldsymbol{\cU}=\{V\}\sbs\tau_2'$ such that
$Y=\cup\boldsymbol{\cU}$, $\boldsymbol{\cU}$ refines
$\boldsymbol{\cV}$ and $\boldsymbol{\cU}$ is $1$-locally finite at
each point of $Y$. Let
$\boldsymbol{\cW}=\big\{V'\in\boldsymbol{\cU}:\;V'\cap
F\neq\vnth\big\}$ and $U=\cup\boldsymbol{\cW}$. Then it is clear
that $\boldsymbol{\cW}$ refines $\boldsymbol{\cV}$ and
$\boldsymbol{\cW}$ is $1$-locally finite at each point $y\in Y$,
so that $t\in\tau_1\cl U$  if and only if there is
$V'\in\boldsymbol{\cW}$ such that $t\in\tau_1\cl V'$. Since
$\boldsymbol{\cW}=\{V'\in\boldsymbol{\cU}:\;V'\cap F\neq\vnth\}$,
for each $V'\in\boldsymbol{\cW}$ there is
$U_x(y)\in\boldsymbol{\cV}$ such that $V'\sbs U_x(y)$. Since
$U_x(y)\cap\,U_y(x)=\vnth$ for such $U_x(y)$ and the corresponding
$U_y(x)\in\tau_1$, we have $x\,\ol{\in}\,\tau_1\cl V'$ for each
$V'\in\boldsymbol{\cW}$. Since $\boldsymbol{\cW}$ is $1$-locally
finite at each point of $Y$, we have $x\,\ol{\in}\,\tau_1\cl U$.
Clearly, $x\in X\setminus\tau_1\cl U=W'\in\tau_1$ and so $W'\cap
Y=W\in\tau_1'$, $F\sbs U\in\tau_2'$, where $U\cap W=\vnth$. Thus
$Y$ is $(1,2)$-regular.

(2) Let $x\in Y$, $F\in\co\tau_1$ and $x\,\ol{\in}\,F$. Then for
each point $y\in F$, there are disjoint neighborhoods
$U_y(x)\in\tau_1$ and $U_x(y)\in\tau_2$. The family
$$  \boldsymbol{\cV}=\big\{U_x(y):\;y\in\tau_1\cl F\big\}\cup \{X\setminus F\}  $$
is a $2$-open covering of $X$ as $\tau_1\sbs\tau_2$. Since $Y$ is
$(2,1)$-paracompact in $X$, there is a family
$\boldsymbol{\cU}=\{V\}\sbs\tau_2$, $\boldsymbol{\cU}$ refines
$\boldsymbol{\cV}$, $Y\sbs\cup\boldsymbol{\cU}$ and
$\boldsymbol{\cU}$ is $1$-locally finite at each point of $Y$. Let
$\boldsymbol{\cW}=\big\{V\in\boldsymbol{\cU}:\;V\cap
F\neq\vnth\big\}$ and $U=\cup\boldsymbol{\cW}\in\tau_2$. Clearly,
$\boldsymbol{\cW}$ refines $\boldsymbol{\cV}$ and
$\boldsymbol{\cW}$ is $1$-locally finite at each point of $Y$.
Hence for each $t\in Y$, $t\in\tau_1\cl U$ if and only if
$t\in\tau_1\cl V$ for some $V\in\boldsymbol{\cW}$. For each
$V\in\boldsymbol{\cW}$ there is $U_x(y)\in\boldsymbol{\cV}$ such
that $V\sbs U_x(y)$. Since $U_x(y)\cap\,U_y(x)=\vnth$ and each
$U_y(x)\in\tau_1$, we have $x\ol{\in}\,\tau_1\cl V$ for each
$V\in\boldsymbol{\cW}$. Since $\boldsymbol{\cW}$ is $1$-locally
finite, $x\ol{\in}\,\tau_1\cl U$. Clearly, $x\in X\setminus
\tau_1\cl U=W\in\tau_1$, $F\cap Y\sbs U$ and $W\cap\,U=\vnth$.
Thus $Y$ is $(1,2)$-regular in $X$.

(3) Let $A\in\co\tau_1$, $B\in\co\tau_2$ and $A\cap B=\vnth$.
Since $Y$ is $(2,1)$-free regular in $X$, for each point $a\in A$
there are disjoint neighborhoods $U(a)\in\tau_2$ and $U_a(B\cap
Y)\in\tau_1$. Clearly, $\boldsymbol{\cV}=\big\{U(a):\;a\in
A\big\}\cup\{X\setminus A\}$ is a $2$-open covering of $X$ as
$\tau_1\sbs\tau_2$. Since $Y$ is $(2,1)$-quasi paracompact in $X$,
there is a family $\boldsymbol{\cU}=\{V'\}\sbs\tau_2'$,
$\boldsymbol{\cU}$ refines $\boldsymbol{\cV}$,
$Y=\cup\boldsymbol{\cU}$ and $\boldsymbol{\cU}$ is $1$-locally
finite at each point of $Y$. Let
$\boldsymbol{\cW}=\big\{V'\in\boldsymbol{\cU}:\;V'\cap
A\neq\vnth\big\}$ and $U=\cup\boldsymbol{\cW}\in\tau_2'$. Then
$A\cap Y\sbs U$ and $\boldsymbol{\cW}$ is $1$-locally finite at
each point of $Y$. Hence, if $x\in Y$, then
\begin{equation}
    x\in \tau_1\cl U \;\;\text{if and only if}\;\;
        x\in\tau_1\cl V' \;\;\text{for some}\;\;
                    V'\in\boldsymbol{\cW}       \tag{$*$}
\end{equation}
Since $\boldsymbol{\cW}\sbs\boldsymbol{\cU}$, $\boldsymbol{\cW}$
refines $\boldsymbol{\cV}$. Moreover, since
$\boldsymbol{\cW}=\{V'\in\boldsymbol{\cU}:\;V'\cap A\neq\vnth\}$,
for each $V'\in\boldsymbol{\cW}$ there is $a\in A$ such that
$V'\sbs U(a)$. Since $U(a)\cap\,U_a(B\cap Y)=\vnth$ for such
$U(a)$ and $V'\sbs U(a)$ we have $V'\cap\,U_a(B\cap Y)=\vnth$.
Hence
$$  \tau_1\cl V'\cap\,U_a(B\cap Y)=\vnth\;\;\text{for each}\;\; V'\in\boldsymbol{\cW}   $$
since $U_a(B\cap Y)\in\tau_1$. By $(*)$, $\tau_1\cl U\cap(B\cap
Y)=\vnth$. Let
$$  W'=X\setminus \tau_1\cl U\in\tau_1 \;\;\text{and}\;\;
        W=W'\cap Y\in\tau_1'.       $$
Then
$$  A\cap Y\sbs U\in\tau_2', \;\;\; B\cap Y\sbs W\in\tau_1' $$
and $U\cap W=\vnth$. Thus $Y$ is $p$-quasi normal in $X$.
\end{pf}

\begin{definition}{6.33}
We will say that a $\BS$ $(X,\tau_1,\tau_2)$ has an $(i,j)$-strong
paracompactness property if for every set $F\in\co\tau_i$,
$\vnth\neq F\neq X$, and every family of sets
$\boldsymbol{\cV}=\big\{U_\al:\;\al\in A\big\}\sbs\tau_j$ such
that $F\sbs\cup\boldsymbol{\cV}$ there is an $i$-locally finite
family $\boldsymbol{\cU}=\big\{V_\bt:\;\bt\in B\big\}\sbs\tau_j$
such that $F\sbs\cup\boldsymbol{\cU}$ and $\boldsymbol{\cU}$
refines $\boldsymbol{\cV}$.
\end{definition}

\begin{proposition}{6.34}
Every $p$-Hausdorff $\BS$ having a $p$-strong paracompactness
property, is $p$-regular.
\end{proposition}

\begin{pf}
Let $a\in X$, $F\in\co\tau_i$ and $a\,\ol{\in}\,F$. Then for each
point $x\in F$ there are disjoint sets $U_x(a)\in\tau_i$ and
$U_a(x)\in\tau_j$. Hence,
$$  U_x(a)\cap\tau_i\cl U_a(x)=\vnth \;\;\text{for each}\;\;
        x\in F.     $$
Clearly, $\boldsymbol{\cV}=\big\{U_a(x):\;x\in F\big\}\sbs\tau_j$
and by condition there is a family
$\boldsymbol{\cU}=\{V_\al:\;\al\in A\}\sbs\tau_j$ such that
$F\sbs\cup\boldsymbol{\cU}$, $\boldsymbol{\cU}$ refines
$\boldsymbol{\cV}$, and $\boldsymbol{\cU}$ is $i$-locally finite.
Hence, for each $\al\in A$ there is $x\in F$ such that
$$  V_\al\sbs U_a(x)\sbs X\setminus U_x(a)=\tau_i\cl(X\setminus U_x(a))\sbs
            X\setminus\{a\}.  $$
Therefore, $\tau_i\cl V_\al\sbs\tau_i\cl U_a(x)\sbs
X\setminus\{a\}$, for each $\al\in A$. Since $\boldsymbol{\cU}$ is
$i$-locally finite, $\tau_i\cl\bigcup\limits_{\al\in A} V_\al=
        \bigcup\limits_{\al\in A} \tau_i\cl V_\al  $
and hence,
$$  \tau_i\cl\bigcup\limits_{\al\in A} V_\al=
        \bigcup\limits_{\al\in A} \tau_i\cl V_\al\sbs X\setminus\{a\}.  $$
Let $U(a)=X\setminus \tau_i\cl\bigcup\limits_{\al\in A}
V_\al\in\tau_i  $ and $U(F)=\bigcup\limits_{\al\in A}
V_\al\in\tau_j$. Then $U(a)\cap\,U(F)=\vnth$ and thus $X$ is
$p$-regular.
\end{pf}

\begin{proposition}{6.35}
Every $(1,2)$-regular $\BS$ having the $(2,1)$-strong
paracompactness property or  every $(2,1)$-regular $\BS$ having
the $(1,2)$-strong paracompactness property, is $p$-normal.
\end{proposition}

\begin{pf}
For example, let us prove the first case. Let $A\in\co\tau_1$,
$B\in\co\tau_2$ and $A\cap B=\vnth$. Since $X$ is $(1,2)$-regular,
for each $x\in B$ there are disjoint neighborhoods
$U_x(A)\in\tau_2$ and $U_A(x)\in\tau_1$. Hence
$$  U_x(A)\cap\tau_2\cl U_A(x)=\vnth.       $$
Clearly, $\boldsymbol{\cV}=\big\{U_A(x):\;x\in B\big\}\sbs\tau_1$
is a $1$-open covering of $B$ and since $X$ has a $(2,1)$-strong
paracompactness property, there is family
$\boldsymbol{\cU}=\{V_\al:\;\al\in D\}\sbs\tau_1$ such that
$\boldsymbol{\cU}$ is $2$-locally finite, $\boldsymbol{\cU}$
refines $\boldsymbol{\cV}$ and $B\sbs\cup\boldsymbol{\cU}$.
Clearly, for each $\al\in D$ there is $x\in B$ such that
$$  V_\al\sbs U_A(x)\sbs X\setminus U_x(A)=
        \tau_2\cl(X\setminus U_x(A))\sbs X\setminus A.  $$
Hence
$$  \tau_2\cl V_\al\sbs\tau_2\cl U_A(x)\sbs
            X\setminus A \;\;\text{for each}\;\; \al\in D.   $$
Since $\boldsymbol{\cU}$ is $2$-locally finite, we have
$$  \tau_2\cl \bigcup\limits_{\al\in D} V_\al=
        \bigcup\limits_{\al\in D} \tau_2\cl V_\al     $$
and, therefore,
$$  \tau_2\cl \bigcup\limits_{\al\in D} V_\al=
        \bigcup\limits_{\al\in D} \tau_2\cl V_\al\sbs X\setminus A.  $$
Let $U(A)=X\setminus \tau_2\cl \bigcup\limits_{\al\in D} V_\al $
and $  U(B)=\bigcup\limits_{\al\in D} V_\al$. Then
$U(A)\in\tau_2$, $U(B)\in\tau_1$ and $U(A)\cap\,U(B)=\vnth$. Thus,
$X$ is $p$-normal.
\end{pf}

\begin{definition}{6.36}
A $\BsS$ $(Y,\tau_1',\tau_2')$ of a $\BS$ $(X,\tau_1,\tau_2)$ is
$(i,j)$-strongly pseudocompact in $X$ if every family
$\boldsymbol{\cV}=\big\{U_\al:\;\al\in D\big\}\sbs\tau_i$ such
that $U_\al\cap Y\neq\vnth$ for each $\al\in D$ and
$\boldsymbol{\cV}$ is $j$-locally finite at each point of $Y$, is
finite.
\end{definition}

For a $\BsS$ $(Y,\tau_1'<\tau_2')$ of a $\BS$ $(X,\tau_1<\tau_2)$,
we have: $Y$ is $(1,2)$-strongly pseudocompact in $X$ implies that
$Y$ is $1$-strongly pseudocompact in $X$ and $Y$ is $2$-strongly
pseudocompact in $X$ implies that $Y$ is $(2,1)$-strongly
pseudocompact in $X$.

\begin{proposition}{6.37}
If $(Y,\tau_1',\tau_2')$ is a $\BsS$ of a $\BS$
$(X,\tau_1,\tau_2)$, $Y$ is $(i,j)$-strongly pseudocompact in $X$
and $(i,j)$-paracompact in $X$, then $Y$ is $i$-compcat in $X$.
\end{proposition}

\begin{pf}
Let $\boldsymbol{\cV}$ be an arbitrary $i$-open covering of $X$.
Since $Y$ is $(i,j)$-paracompact in $X$, there is a family
$\boldsymbol{\cU}\sbs\tau_i$ such that $\boldsymbol{\cU}$ refines
$\boldsymbol{\cV}$, $Y\sbs\cup\boldsymbol{\cU}$ and
$\boldsymbol{\cU}$ is $j$-locally finite at each point of $Y$. We
may also assume that $U\cap Y\neq\vnth$ for each set
$U\in\boldsymbol{\cU}$. Since $Y$ is $(i,j)$-strongly
pseudocompact in $X$, $\boldsymbol{\cU}$ is finite and so, $Y$ is
$i$-compact in $X$.
\end{pf}

\begin{corollary}{6.38}
If $(Y,\tau_1',\tau_2')$ is a $\BsS$ of an $(i,j)$-regular $\BS$
\linebreak        $(X,\tau_1,\tau_2)$, $\tau_jC\tau_i$,
$\tau_iN\tau_j$, $Y$ is $j$-Lindel\"{o}ff in $X$ and $Y$ is
$(i,j)$-strongly pseudocompact in $X$, then $X$ is $i$-compact.
\end{corollary}

\begin{pf}
By Theorem~6.25, $Y$ is $(i,j)$-paracompact in $X$ and it remains
to use Proposition~6.37.
\end{pf}

\begin{corollary}{6.39}
If $(Y,\tau_1',\tau_2')$ is $(i,j)$-strongly pseudocompact in $X$
and $(i,j)$-paracompact in $X$, where $X$ is $i$-regular, then
$\tau_i\cl Y$ is     \linebreak        $i$-compact.
\end{corollary}

\begin{pf}
By Proposition~6.37, $Y$ is $i$-compact in $X$. Since $X$ is
$i$-regular, by [19], $\tau_i\cl Y$ is $i$-compact.
\end{pf}

\begin{corollary}{6.40}
If $(X,\tau_1,\tau_2)$ is $i$-regular, $(i,j)$-regular,
$\tau_jC\tau_i$, $\tau_iN\tau_j$, $Y$ is $(i,j)$-strongly
pseudocompact in $X$ and $Y$ is $j$-Lindel\"{o}ff in $X$, then
$\tau_i\cl Y$ is $i$-compact
\end{corollary}

\begin{pf}
Indeed, by Theorem~6.25, $Y$ is $(i,j)$-paracompact in $X$ and
hence, it remains, to use Corollary~6.39.
\end{pf}

\begin{corollary}{6.41}
Let $(Y,\tau_1'<\tau_2')$ is $2$-Lindel\"{o}ff in
$(X,\tau_1<_N\tau_2)$, $Y$ is $(1,2)$-strongly pseudocompact in
$X$ and $X$ is $1$-regular, then $\tau_1\cl Y$ is $1$-compact.
\end{corollary}

\begin{pf}
Indeed, since $\tau_1\sbs\tau_2$, we have $\tau_2C\tau_1$ and $X$
is $1$-regular implies that $X$ is $(1,2)$-regular. Hence, by
Theorem~6.25 for $i=1$, $j=2$, we have that $Y$ is
$(1,2)$-paracompact in $X$. Therefore, by Corollary~6.39,
$\tau_1\cl Y$ is $1$-compact.
\end{pf}

\begin{theorem}{6.42}
For a $\BS$ $(X,\tau_1,\tau_2)$ the following conditions are
satisfied:
\begin{enumerate}
\item[(1)] If $(Y,\tau_1'',\tau_2'')$ is $(i,j)$-strongly
pseudocompact in $X$, then        \linebreak $(\tau_i\cl
Y,\tau_1',\tau_2')$ is $(i,j)$-strongly pseudocompact in $X$ too.

\item[(2)] If $(Z,\tau_1'',\tau_2'')\sbs (Y,\tau_1',\tau_2')\sbs
            (X,\tau_1,\tau_2)$
and $Z$ is $(i,j)$-strongly pseudocompact in $Y$, then $Z$ is
$(i,j)$-strongly pseudocompact in $X$ as well. Moreover, if $Y\in
i\dd\,\cD(X)$, then $Z$ is $(i,j)$-strongly pseudocompact in $Y$
if and only if $Z$ is $(i,j)$-strongly pseudocompact in $X$.

\item[(3)] If $(Z,\tau_1'',\tau_2'')$ is $(i,j)$-strongly
pseudocompact in $X$, then there is a $\BsS$ $(Y,\tau_1',\tau_2')$
of $X$ such that $Y\in i\dd\,\cD(X)$, $Z\sbs Y$, $Z\in\co\tau_i'$
and $Z$ is $(i,j)$-strongly pseudocompact in $Y$.
\end{enumerate}
\end{theorem}

\begin{pf}
(1) Let $\boldsymbol{\cV}=\{U_\al:\;\al\in A\}\sbs\tau_i$,
$\boldsymbol{\cV}$ be $j$-locally finite at each point of
$\tau_i\cl Y$ and $U_\al\cap\tau_i\cl Y\neq\vnth$ for each $\al\in
A$. Clearly, $\boldsymbol{\cV}$ is also $j$-locally finite at each
point of $Y$ and since $Y\in\tau_i$, $U_\al\cap\tau_i\cl
Y\neq\vnth$ implies that $U_\al\cap Y\neq\vnth$ for each $\al\in
A$. But $Y$ is $(i,j)$-strongly pseudocompact in $X$ and so
$\boldsymbol{\cV}$ is finite.

(2) First, let $Z$ be $(i,j)$-strongly pseudocompact in $Y$,
$\boldsymbol{\cV}=\{U_\al:\;\al\in A\}\sbs\tau_i$,
$\boldsymbol{\cV}$ be $j$-locally finite at each point of $Z$ and
$U_\al\cap Z\neq\vnth$ for each $\al\in A$. Let us consider the
family
$$  \boldsymbol{\cV}'=\big\{U_\al'=U_\al\cap Y:\;\al\in A\big\}\sbs\tau_i'.  $$
Since $\boldsymbol{\cV}$ is $j$-locally finite at each point of
$Z$, $\boldsymbol{\cV}'$ is also $j$-locally finite at each point
of $Z$. Moreover, $Z\sbs Y$ gives
$$  U_\al'\cap Z=(U_\al\cap Y)\cap Z=U_\al\cap Z\neq\vnth
        \;\;\text{for each} \;\; \al\in A.  $$
Since $Z$ is $(i,j)$-strongly pseudocompact in $Y$,
$\boldsymbol{\cV}'$ is finite and thus, $\boldsymbol{\cV}$ is
finite a well. Hence, $Z$ is $(i,j)$-pseudocompact in $X$.

Now, let $Y\in i\dd\,\cD(X)$ and $Z$ be $(i,j)$-strongly
pseudocompact in $X$. Let $\boldsymbol{\cU}'=\big\{U_\al':\;\al\in
A\big\}\sbs\tau_i'$, $\boldsymbol{\cU}'$ be $j$-locally finite at
each point of $Z$ and $U_\al'\cap Z\neq\vnth$ for each $\al\in A$.
Let $\boldsymbol{\cU}=\big\{U_\al:\;U_\al\cap
Y=U_\al'\big\}\sbs\tau_i$. Clearly,
$$  U_\al\cap Z=(U_\al\cap Y)\cap Z=U_\al'\cap Z\neq\vnth
        \;\;\text{for each}\;\; \al\in A.       $$
Let us show that $\boldsymbol{\cU}$ is $j$-locally finite at each
point of $Z$. Contrary: there is $z\in Z$ such that
$\boldsymbol{\cU}$ is not $j$-locally finite at $z$. Hence, for
each neighborhood $U(z)\in\tau_j$ we have
$$  \Big|\big\{U_\al:\;U_\al\in\boldsymbol{\cU},\;
        U_\al\cap\,U(z)\neq\vnth\big\}\Big|\geq\aleph_0.  $$
Since $Y\in i\dd\,\cD(X)$, for each $U\in\tau_i\setminus\{\vnth\}$
we have $U\cap Y\neq\vnth$ and so, if
$$  \boldsymbol{\cA}=\big\{U_\al\in\boldsymbol{\cU}:\;U_\al\cap\,U(z)\neq\vnth\big\},  $$
then $U_\al\cap Y\neq\vnth$ for each $U_\al\in\boldsymbol{\cA}$.
Since $Z\sbs Y$, we have $U_\al'\cap\,U(z)\neq\vnth$ for each
$U_\al'\in\boldsymbol{\cA}'$, where
$\boldsymbol{\cA}'=\{U_\al'=U_\al\cap
Y:\;U_\al\in\boldsymbol{\cA}\}$. Hence
$|\boldsymbol{\cA}'|\geq\aleph_0$ which is impossible. Therefore,
$\boldsymbol{\cU}$ is $j$-locally finite at each point of $Z$ and
since $Z$ is $(i,j)$-strongly pseudocompact in $X$,
$\boldsymbol{\cU}$ is finite. Hence $\boldsymbol{\cU}'$ is also
finite and thus, $Z$ is $(i,j)$-strongly pseudocompact in $Y$.

(3) Let $Y=X\setminus\big(\tau_i\cl Z\setminus Z\big)$. Then
$$  \tau_i\cl Y=\tau_i\big((X\setminus\tau_i\cl Z)\cup Z\big)=X $$
and $Z\sbs Y$. Hence, by the second part of (2), $Z$ is
$(i,j)$-strongly pseudocompact in $Y$. Finally,
$$  \tau_i\cl Z\cap Y=
        \tau_i\cl Z\cap\big((X\setminus\tau_i\cl Z)\cup Z\big)=Z $$
and so $Z\in\co\tau_i'$.
\end{pf}
\vskip+0.2cm

The next five sections are devoted to introducing and studying of
the relative bitopological inductive and covering dimension
functions, their interrelations and also, to the new inductive and
relative inductive dimension functions, the so called separately
inductive dimension functions. Note also here that the relative
bitopological inductive dimension functions will be characterized
in two different ways: as in terms of neighborhoods so by means of
partitions. Relative topological dimension functions have been
considered in [5], [23], [24], [26].

\vskip+0.5cm
\section*{\textbf{7. $(i,j)$-Small Relative Inductive Dimension Functions}}
\vskip+0.2cm

\begin{definition}{7.1}
Let $(Y,\tau_1',\tau_2')$ be a $\BsS$ of a $\BS$
$(X,\tau_1,\tau_2)$ and $n$ denote a nonnegative integer. We say
that:
\begin{enumerate}
\item[(1)] $(i,j)\dd\ind(Y,X)=-1$ if and only if $Y=\vnth$.

\item[(2)] $(i,j)\dd\ind(Y,X)\leq n$ if for each point $y\in Y$
and any neighborhood $U(y)\in\tau_i$ there is a neighborhood
$V(y)\in\tau_i$ such that $\tau_j\cl V(y)\sbs U(y)$ and
$$  (i,j)\dd\ind\big((j,i)\dd\Fr_XV(y)\cap Y,X\big)\leq n-1. $$

\item[(3)] $(i,j)\dd\ind(Y,X)=n$ if $(i,j)\dd\ind(Y,X)\leq n$ and
the inequality $(i,j)\dd\ind(Y,X)\leq n-1$ does not hold.

\item[(4)] $(i,j)\dd\ind(Y,X)=\infty$ if the inequality
$(i,j)\dd\ind(Y,X)\leq n$ does not hold for any $n$.
\end{enumerate}
\end{definition}

As a rule,
$$  p\,\dd\ind(Y,X)\leq n\llra
        \big((1,2)\dd\ind(Y,X)\leq n\wedge
                (2,1)\dd\ind(Y,X)\leq n\big).       $$

Recall also here that the small relative inductive dimension was
introduced by V.~V.~Filippov.

From the our point of view this fact is motivated by the
well-known topological lemma, following which for a $\TsS$
$(Y,\tau')$ of a hereditarily normal $\TS$ $(X,\tau)$ we have
$\ind Y\leq n$ if and only if for each point $x\in Y$ and any
neighborhood $U(x)\in\tau$ there is a neighborhood $V(x)\in\tau$
such that $\tau\cl V(x)\sbs U(x)$ and $\ind(\Fr_XV(x)\cap Y)\leq
n-1$, and is confirmed by Proposition~1 in [23] which says that
for any $\TsS$ $(Y,\tau')$ of a hereditarily normal $\TS$
$(X,\tau)$ take place the equality $\ind Y=\ind (Y,X)$.

\begin{remark}{7.2}
Let us show that for any subset $Y\sbs X$ and any nonnegative
integer $n$ we have: $(i,j)\dd\ind(Y,X)\leq n$ if and only if for
each point $x\in Y$ and any neighborhood $U(x)\in\tau_i$ there is
a neighborhood $V(x)\in\tau_i$ such that
\begin{equation}
    \tau_j\cl V(x)\sbs U(x) \;\;\text{and}\;\;
        (i,j)\dd\ind\big((j,i)\dd\Fr_X V(x)\cap Y\big)\leq n-1
                    \tag{$*$}
\end{equation}
in the sense of (1) of Corollary~3.1.6 in [8].
\end{remark}

We shall prove this equivalence by induction under a nonnegative
integer $n=(i,j)\dd\ind(Y,X)$.

Indeed, by (1) and (2) of Definition~7.1 and (1) of
Definition~3.1.3 in [8], for $n=0$ we have: $(i,j)\dd\ind(Y,X)\leq
0\llra$ $\big($for each point $x\in Y$ and any neighborhood
$U(x)\in\tau_i$ there is a neighborhood $V(x)\in\tau_i$ such that
$\tau_j\cl V(x)\sbs U(x)$ and $(i,j)\dd\ind((j,i)\dd\Fr_X V(x)\cap
Y,X)=-1\big)\llra$ $\big($for each point $x\in Y$ and any
neighborhood $U(x)\in\tau_i$ there is a neighborhood
$V(x)\in\tau_i$ such that $\tau_j\cl V(x)\sbs U(x)$ and
$(j,i)\dd\Fr_X V(x)\cap Y=\vnth\big)\llra$ $\big($for each point
$x\in Y$ and any neighborhood $U(x)\in\tau_i$ there  is a
neighborhood $V(x)\in\tau_i$ such that $\tau_j\cl V(x)\sbs U(x)$
and $(i,j)\dd\ind((j,i)\dd\Fr_X V(x)\cap Y)=-1\big)$.

Now, suppose that for any subset $Y\sbs X$ the equivalences:
\linebreak         $(i,j)\dd\ind(Y,X)\leq k\llra$ $\big($for each
point $x\in Y$ and any neighborhood $U(x)\in\tau_i$ there is a
neighborhood $V(x)\in\tau_i$ such that $\tau_j\cl V(x)\sbs U(x)$
and $(i,j)\dd\ind((j,i)\dd\Fr_X V(x)\cap Y,X)\leq k-1\big)\llra$
$\big($for each point $x\in Y$ and any neighborhood
$U(x)\in\tau_i$ there  is a neighborhood $V(x)\in\tau_i$ such that
$\tau_j\cl V(x)\sbs U(x)$ and
\begin{equation}
    (i,j)\dd\ind\big((j,i)\dd\Fr_X V(x)\cap Y\big)\leq k-1\big)
                \tag{$**$}
\end{equation}
are proved for $n\leq k-1$ and prove they for $n=k$.

Once more applying (2) of Definition~7.1 we obtain:
$(i,j)\dd\ind(Y,X)=k\llra$ $\big($for each point $x\in Y$ and any
neighborhood $U(x)\in\tau_i$ there  is a neighborhood
$V(x)\in\tau_i$ such that $\tau_j\cl V(x)\sbs U(x)$ and
     \linebreak $(i,j)\dd\ind((j,i)\dd\Fr_X V(x)\cap Y,X)\leq
k-1\big)$.

Since $(**)$ holds for any subset $Y\sbs X$, by inductive
hypothesis the equivalence $(*)$ is proved.

Hence, it should be noted, that the definitions of relative
inductive functions $(i,j)\dd\ind (Y,X)$ given in Remark~7.2, are
generalizations of the topological lemma, mentioned before
Remark~7.2.

\begin{theorem}{7.3}
If $(Y,\tau_1',\tau_2')$ is a $\BsS$ of a $\BS$
$(X,\tau_1,\tau_2)$, then the following conditions are satisfied:
\begin{enumerate}
\item[(1)] $(i,j)\dd\ind Y\leq (i,j)\dd\ind(Y,X)\leq (i,j)\dd\ind
X$.

\item[(2)] $(i,j)\dd\ind (X,X)=(i,j)\dd\ind X$.

\item[(3)] If $Z\sbs Y$, then
$$  (i,j)\dd\ind(Z,Y)\leq (i,j)\dd\ind(Z,X)\leq (i,j)\dd\ind(Y,X). $$
\end{enumerate}

Moreover, if $\tau_1\sbs\tau_2$, then
\begin{enumerate}
\item[(4)] $(1,2)\dd\ind Y\leq (1,2)\dd\ind(Y,X)\leq (1,2)\dd\ind
X\leq 1\dd\ind X$ and
$$  2\dd\ind Y\leq (2,1)\dd\ind Y\leq (2,1)\dd\ind(Y,X)\leq
            (2,1)\dd\ind X.     $$
\end{enumerate}
\end{theorem}

\begin{pf}
(1) First, let us prove the left inequality by induction under a
nonnegative integer $n=(i,j)\dd\ind(Y,X)$. Suppose that the
inequality is proved for $n\leq k-1$ and prove it for $n=k$. Let
$x\in Y$ be any point and $U'(x)\in\tau_i'$ be any neighborhood.
If $U(x)\cap Y=U'(x)$, then by $(i,j)\dd\ind(Y,X)=k$ there exists
$V(x)\in\tau_i$ such that $\tau_j\cl V(x)\sbs U(x)$ and
$$  (i,j)\dd\ind\big((j,i)\dd\Fr_X V(x)\cap Y,X\big)\leq k-1.   $$
Hence, by inductive assumption, $(i,j)\dd\ind\big((j,i)\dd\Fr_X
V(x)\cap Y\big)\leq k-1$. Moreover, $\tau_j'\cl V'(x)\sbs U'(x)$
for $V'(x)=V(x)\cap Y$ and since
$$  (j,i)\dd\Fr_Y V'(x)\sbsq (j,i)\dd\Fr_X V(x)\cap Y,  $$
by the monotonicity, i.e., by (2) of Proposition~3.1.4 in [8],
$$  (i,j)\dd\ind(j,i)\dd\Fr_Y V'(x)\leq k-1.        $$
Therefore, by (1) of Corollary~3.1.6 in [8] $(i,j)\dd\ind Y\leq k$
and so $(i,j)\dd\ind Y\leq (i,j)\dd\ind(Y,X)$.

Now, let us prove the right inequality by induction under $n=
      \linebreak     (i,j)\dd\ind X$. Suppose that the inequality is proved for $n\leq
k-1$ and prove it for $n=k$. Let $x\in Y\sbs X$ be any point and
$U(x)\in\tau_i$ be any neighborhood. Then by (1) of
Corollary~3.1.6 in [8] there is a neighborhood $V(x)\in\tau_i$
such that $\tau_j\cl V(x)\sbs U(x)$ and $(i,j)\dd\ind(j,i)\dd\Fr_X
V(x) \linebreak       \leq k-1$. Hence, by inductive assumption
$$  (i,j)\dd\ind\big((j,i)\dd\Fr_XV(x)\cap Y,(j,i)\dd\Fr_XV(x)\big)\leq k-1  $$
and by the left inequality
$$  (i,j)\dd\ind\big((j,i)\dd\Fr_X V(x)\cap Y\big)\leq k-1.  $$
Thus, by Remark~7.2, $(i,j)\dd\ind(Y,X)\leq k$ and so
$(i,j)\dd\ind(Y,X)\leq (i,j)\dd\ind X$.

(2) The equality follows directly from (1) for $Y=X$.

(3) First, let us prove the right inequality, i.e.,
$$  (i,j)\dd\ind(Z,X)\leq (i,j)\dd\ind(Y,X)         $$
(\textbf{the monotonicity property of the $(i,j)$-small relative
inductive dimensions}).

Let $n=(i,j)\dd\ind(Y,X)$, the inequality be proved for $n\leq
k-1$ and prove it for $n=k$. If $x\in Z$ is any point and
$U(x)\in\tau_i$ is any neighborhood, then $x\in Z\sbs Y$ and
$(i,j)\dd\ind(Y,X)=k$ imply that there is a neighborhood
$V(x)\in\tau_i$ such that $\tau_j\cl V(x)\sbsq U(x)$ and
$$  (i,j)\dd\ind\big((j,i)\dd\Fr_X V(x)\cap Y,X\big)\leq k-1.   $$
Since $Z\sbs Y$, we have
$$  (j,i)\dd\Fr_X V(x)\cap Z\sbsq (j,i)\dd\Fr_X V(x)\cap Y  $$
and by inductive assumption,
\begin{eqnarray*}
    &\ds (i,j)\dd\ind\big((j,i)\dd\Fr_X V(x)\cap Z,X\big)\leq \\
    &\ds \leq (i,j)\dd\ind\big((j,i)\dd\Fr_X V(x)\cap Y,X\big)\leq k-1.
\end{eqnarray*}
Thus $(i,j)\dd\ind(Z,X)\leq k$ and so $(i,j)\dd\ind(Z,X)\leq
(i,j)\dd\ind(Y,X)$.

Now, let us prove the left inequality by induction under $n=
          \linebreak         (i,j)\dd\ind(Z,X)$. Let the inequality be proved for $n\leq k-1$
and prove it for $n=k$. Let $x\in Z$ be any point and
$U'(x)\in\tau_i'$ be any neighborhood. Since $x\in Z\sbs Y$ and
$(i,j)\dd\ind(Z,X)=k$, for $U(x)\in\tau_i$, where $U(x)\cap
Y=U'(x)$, there is $V(x)\in\tau_i$ such that $\tau_j\cl V(x)\sbsq
U(x)$ and
$$  (i,j)\dd\ind\big((j,i)\dd\Fr_X V(x)\cap Z,X\big)\leq k-1.  $$
If $V'(x)=V(x)\cap Y$, then $\tau_j'\cl V'(x)\sbs U'(x)$. Since
$$  (j,i)\dd\Fr_Y V'(x)\sbsq (j,i)\dd\Fr_X V(x)\cap Y   $$
and $Z\sbs Y$, we have
$$  (j,i)\dd\Fr_Y V'(x)\cap Z\sbsq (j,i)\dd\Fr_X V(x)\cap Z.    $$
Hence, by the right inequality,
\begin{align}
    (i,j)\dd\ind\big((j,i)\dd\Fr_Y V'(x)\cap Z,X\big) & \leq
                            \nonumber \\
    \leq (i,j)\dd\ind\big((j,i)\dd\Fr_X V(x)\cap Z,X\big) & \leq k-1.
                \tag{$*$}
\end{align}

Therefore, by the inductive hypothesis and $(*)$,
\begin{eqnarray*}
    (i,j)\dd\ind\big((j,i)\dd\Fr_Y V'(x)\cap Z,Y\big) &&\hskip-0.6cm \leq \\
    \leq (i,j)\dd\ind\big((j,i)\dd\Fr_Y V'(x)\cap Z,X\big) &&\hskip-0.6cm \leq k-1.
\end{eqnarray*}
Thus, $(i,j)\dd\ind(Z,Y)\leq k$ and so
$(i,j)\dd\ind(Z,Y)\leq(i,j)\dd\ind(Z,X)$.

(4) Finally, for the case $\tau_1\sbs\tau_2$ it suffices to use
(1) above and (1) of Theorem~3.1.36 in [8].
\end{pf}

\begin{theorem}{7.4}
For an $i$-open $\BsS$ $(Y,\tau_1',\tau_2')$ of an $(i,j)$-regular
$\BS$ $(X,\tau_1,\tau_2)$ the following conditions are satisfied:
\begin{enumerate}
\item[(1)] If $p\,\dd\cl Z\sbsq Y$, then
$(i,j)\dd\ind(Z,X)=(i,j)\dd\ind(Z,Y)$.

\item[(2)] $(i,j)\dd\ind(Y,X)=(i,j)\dd\ind Y$.
\end{enumerate}
\end{theorem}

\begin{pf}
(1) By (3) of Theorem~7.3, it suffices to prove the inequality
$(i,j)\dd\ind(Z,X)\leq (i,j)\dd\ind(Z,Y)$. We shall use the
induction under $n=(i,j)\dd\ind(Z,Y)$. Suppose that the inequality
is proved for $n\leq k-1$ and prove it for $n=k$. Let $x\in Z$ be
any point and $U(x)\in\tau_i$ be any neighborhood. As
$Y\in\tau_i$, one can assume that $U(x)\sbs Y$. Since $X$ is
$(i,j)$-regular and $(i,j)\dd\ind(Z,Y)=k$, there is
$V(x)\in\tau_i$ such that $\tau_j\cl V(x)\sbs U(x)$ and
$$  (i,j)\dd\ind\big((j,i)\dd\Fr_Y V(x)\cap Z,Y\big)\leq k-1.  $$
Clearly,
$$  p\,\dd\cl\big((j,i)\dd\Fr_Y V(x)\cap Z\big)\sbsq p\,\dd\cl Z\sbsq Y  $$
and so, one can use the inductive assumption under the set
\linebreak        $(j,i)\dd\Fr_Y V(x)\cap Z$, i.e.
\begin{eqnarray*}
    (i,j)\dd\ind\big((j,i)\dd\Fr_Y V(x)\cap Z,X\big) &&\hskip-0.6cm \leq \\
    \leq (i,j)\dd\ind\big((j,i)\dd\Fr_Y V(x)\cap Z,Y\big) &&\hskip-0.6cm \leq k-1.
\end{eqnarray*}
On the other hand, since $\tau_j\cl V(x)\sbs U(x)\sbs Y$ and
hence,
$$  (j,i)\dd\Fr_Y V(x)=(j,i)\dd\Fr_X V(x),      $$
we have
$$  (i,j)\dd\ind\big((j,i)\dd\Fr_X V(x)\cap Z,X\big)\leq k-1    $$
so that $(i,j)\dd\ind(Z,X)\leq k$. Thus $(i,j)\dd\ind(Z,X)\leq
(i,j)\dd\ind(Z,Y)$.

(2) By (1) of Theorem~7.3, it suffices to prove the inequality
\linebreak          $(i,j)\dd\ind(Y,X)\leq (i,j)\dd\ind Y$, where
$(i,j)\dd\ind Y=n$. Let the inequality be proved for $n\leq k-1$
and prove it for $n=k$. Let $x\in Y$ be any point and
$U(x)\in\tau_i$ be any neighborhood. But $Y\in\tau_i$ and one can
assume that $U(x)\sbs Y$. Since $X$ is $(i,j)$-regular and
$(i,j)\dd\ind Y=k$, there is $V(x)\in\tau_i$ such that $\tau_j\cl
V(x)\sbs U(x)\sbs Y$ and
$$  (i,j)\dd\ind(j,i)\dd\Fr_YV(x)=
        (i,j)\dd\ind(j,i)\dd\Fr_XV(x)\leq k-1.        $$
But $(i,j)\dd\ind Y=(i,j)\dd\ind(Y,Y)$ and one can suppose that
$$  (i,j)\dd\ind\big((j,i)\dd\Fr_X V(x)\cap Y,Y\big)\leq k-1,   $$
where $(j,i)\dd\Fr_X V(x)\cap Y\sbs Y$ and so
$$  (i,j)\dd\ind\big((j,i)\dd\Fr_X V(x),Y\big)\leq k-1. $$
Since
$$  p\,\dd\cl(j,i)\dd\Fr_X V(x)=(j,i)\dd\Fr_X V(x)\sbs Y, $$
one can consider in (1) the set $(j,i)\dd\Fr_X V(x)$ instead of
$Z$, and so
\begin{eqnarray*}
    &\ds (i,j)\dd\ind\big((j,i)\dd\Fr_X V(x)\cap Y,X\big)=
        (i,j)\dd\ind\big((j,i)\dd\Fr_X V(x),X\big)= \\
    &\ds =(i,j)\dd\ind\big((j,i)\dd\Fr_X V(x),Y\big)\leq k-1.
\end{eqnarray*}
Therefore $(i,j)\dd\ind(Y,X)\!\leq\!k$ and thus
$(i,j)\dd\ind(Y,X)\!\leq\!(i,j)\dd\ind Y$.~\end{pf}

\begin{proposition}{7.5}
If $(Y,\tau_1',\tau_2')$ is a $\BsS$ of a $\BS$
$(X,\tau_1,\tau_2)$ and $Y=Y_1\cup Y_2$, then
$$  (i,j)\dd\ind(Y,X)\leq
        (i,j)\dd\ind(Y_1,X)+(i,j)\dd\ind(Y_2,X)+1.  $$
\end{proposition}

\begin{pf}
We shall use the induction under
$$  n=(i,j)\dd\ind(Y_1,X)+(i,j)\dd\ind(Y_2,X).  $$
Let the inequality be proved for $n\leq k-1$ and prove it for
$n=k$. Without loss of generality let $x\in Y_1\sbs Y$. Then for
any neighborhood $U(x)\in\tau_i$ there is $V(x)\in\tau_i$ such
that $\tau_j\cl V(x)\sbs U(x)$ and
$$  (i,j)\dd\ind\big((j,i)\dd\Fr_X V(x)\cap Y_1,X\big)\leq
        (i,j)\dd\ind(Y_1,X)-1.      $$
Since
$$  (j,i)\dd\Fr_X V(x)\cap Y=
            \big((j,i)\dd\Fr_X V(x)\cap Y_1\big)\cup
        \big((j,i)\dd\Fr_X V(x)\cap Y_2\big)        $$
and
\begin{eqnarray*}
    &\ds (i,j)\dd\ind\big((j,i)\dd\Fr_X V(x)\cap Y_1,X\big)+ \\
    &\ds +(i,j)\dd\ind\big((j,i)\dd\Fr_X V(x)\cap Y_2,X\big)\leq \\
    &\ds \leq (i,j)\dd\ind(Y_1,X)-1+(i,j)\dd\ind(Y_2,X)=k-1,
\end{eqnarray*}
one can use the induction under the set $(j,i)\dd\Fr_X V(x)\cap
Y$. Hence,
$$  (i,j)\dd\ind\big((j,i)\dd\Fr_X V(x)\cap Y,X\big)\leq
            k-1+1=k     $$
and so $(i,j)\dd\ind(Y,X)\leq k+1$.
\end{pf}

\begin{proposition}{7.6}
For any $\BsS$ $(Y,\tau_1',\tau_2')$ of a hereditarily $p$-normal
$\BS$ $(X,\tau_1,\tau_2)$ we have $(i,j)\dd\ind(Y,X)=(i,j)\dd\ind
Y$.
\end{proposition}

\begin{pf}
First, let us prove that if $U'\in\tau_i'$, then there is
$U\in\tau_i$ such that $U\cap Y=U'$ and $(j,i)\dd\Fr_X U\cap
Y=(j,i)\dd\Fr_Y U'$. Indeed, for the sets $U'$ and
$Y\setminus\tau_j'\cl U'$ we have:
$$  \big(\tau_j\cl U'\cap(Y\setminus\tau_j'\cl U')\big)\cup
        \big(U'\cap\tau_i\cl(Y\setminus\tau_j'\cl U')\big)=\vnth  $$
and by (5) of Definition~1.1, there are $U\in\tau_i$, $V\in\tau_j$
such that $U'\sbs U$, $Y\setminus\tau_j'\cl U'\sbs V$ and $U\cap
V=\vnth$. Hence $\tau_j\cl U\cap(Y\setminus\tau_j'\cl U')=\vnth$,
i.e., $\tau_j\cl U\cap Y=\tau_j'\cl U'$ and so
$$  (j,i)\dd\Fr_X U\cap Y=(\tau_j\cl U\setminus U)\cap Y=
        \tau_j'\cl U'\setminus U'=(j,i)\dd\Fr_Y U'.     $$

Evidently by (1) of Theorem~7.3, it suffices to prove that
\linebreak         $(i,j)\dd\ind(Y,X)\leq(i,j)\dd\ind Y$. Let
$(i,j)\dd\ind Y=n$, the inequality be proved for $n\leq k-1$ and
prove it for $n=k$. Let $x\in Y$ and $U(x)\in\tau_i$ be any
neighborhood. Since $(i,j)\dd\ind Y=n$, there is $V'(x)\in\tau_i'$
such that $\tau_j'\cl V'(x)\sbs U(x)\cap Y$ and
$(i,j)\dd\ind\big((j,i)\dd\Fr_Y V'(x)\big)\leq k-1$. But, by the
first part of the proof, there is $V(x)\in\tau_i$ such that
$V(x)\cap Y=V'(x)$, $(j,i)\dd\Fr_YV'(x)=(j,i)\dd\Fr_XV(x)\cap Y$
and so
$$  (i,j)\dd\ind\big((j,i)\dd\Fr_X V(x)\cap Y\big)\leq k-1. $$
Hence, by the inductive assumption,
\begin{eqnarray*}
    &\ds (i,j)\dd\ind\big((j,i)\dd\Fr_X V(x)\cap Y,X\big)\leq \\
    &\ds \leq (i,j)\dd\ind\big((j,i)\dd\Fr_X V(x)\cap Y\big)\leq k-1
\end{eqnarray*}
so that $(i,j)\dd\ind(Y,X)\leq k$ and thus $(i,j)\dd\ind(Y,X)\leq
(i,j)\dd\ind Y$.
\end{pf}

\begin{proposition}{7.7}
If $(Y,\tau_1',\tau_2')$ is a $p$-perfectly normal $\BsS$ of a
      \linebreak       $p$-normal $\BS$ $(X,\tau_1,\tau_2)$ and
$Y\in\co\tau_1\cap\co\tau_2$, then $(i,j)\dd\ind(Y,X)=(i,j)\dd\ind
Y$.
\end{proposition}

\begin{pf}
By Proposition~7.6 it suffices to prove that if $U'\in\tau_i'$,
then there is $U\in\tau_i$ such that $U\cap Y=U'$ and
$(i,j)\dd\Fr_XU\cap Y=(j,i)\dd\Fr_YU'$. Since
$(Y,\tau_1',\tau_2')$ is $p$-perfectly normal and
$Y\in\co\tau_1\cap\co\tau_2$, we have
\begin{eqnarray*}
    &\ds U'\in j\dd\FF_\sg(Y)\sbs j\dd\FF_\sg(X), \;\;\;
        Y\setminus\tau_j'\cl U'\in i\dd\FF_\sg(Y)\sbs i\dd\FF_\sg(X), \\
    &\ds \tau_j'\cl U'=\tau_j\cl U' \;\;\text{and}\;\;
        \tau_i'\cl(Y\setminus\tau_j'\cl U')=
                \tau_i\cl(Y\setminus\tau_j'\cl U').
\end{eqnarray*}
Hence
\begin{equation}
    \big(\tau_j\cl U'\cap(Y\setminus\tau_j'\cl U')\big)\cup
        \big(U'\cap \tau_i\cl(Y\setminus\tau_j\cl U')\big)=\vnth.
                    \tag{$*$}
\end{equation}
Since $(X,\tau_1,\tau_2)$ is $p$-normal, $U'\in j\dd\FF_\sg(X)$,
$Y\setminus\tau_j'\cl U'\in i\dd\FF_\sg(X)$ and $(*)$ is
satisfied, by Lemma~2.21, there are $U\in\tau_i$, $V\in\tau_j$
such that $U'\sbs U$, $Y\setminus\tau_j'\cl U'\sbs V$ and $U\cap
V=\vnth$. Clearly, $\tau_j\cl U\cap(Y\setminus\tau_j'\cl
U')=\vnth$ and hence, $\tau_j\cl U\cap Y=\tau_j'\cl U'$. Thus,
$$  (j,i)\dd\Fr_X U\cap Y=(\tau_j\cl U\setminus U)\cap Y=
        \tau_j'\cl U'\setminus U'=(j,i)\dd\Fr_Y U'.     $$

Therefore, it remains to use the second part of the proof of
Pro\-po\-si\-ti\-on~7.6.~\end{pf}

\begin{theorem}{7.8}
If $f:(X,\tau_1<\tau_2)\to (X_1,\gm_1<\gm_2)$ is a $d$-continuous
surjection and $(Y,\tau_1'<\tau_2')$ is a $\BsS$ of $X$ such that
the restriction $f\big|_Y:(Y,\tau_1'<\tau_2')\to
(Y_1=f(Y),\gm_1'<\gm_2')$ is a $d$-homeomorphism, then
$(1,2)\dd\ind(Y,X)\leq (1,2)\dd\ind(Y_1,X_1)$. Moreover, if
$(X_1,\gm_1,\gm_2)$ is hereditarily $p$-normal, then
$$  (1,2)\dd\ind(Y,X)=(1,2)\dd\ind(Y_1,X_1)=(1,2)\dd\ind Y.     $$
\end{theorem}

\begin{pf}
First, let us prove the inequality $(1,2)\dd\ind(Y,X)\leq
\linebreak        (1,2)\dd\ind(Y_1,X_1)$ by induction under
$n=(1,2)\dd\ind(Y_1,X_1)$. Suppose that the inequality is proved
for $n\leq k-1$ and prove it for $n=k$. Let $x\in Y$ be any point
and $U(x)\in\tau_1$ be any neighborhood. Since $f\big|_Y$ is a
$d$-homeomorphism, we have $f(U(x)\cap Y)\in\gm_1'$. Furthermore,
since $(1,2)\dd\ind(Y_1,X_1)=k$, for a set $U\in\gm_1$, where
$U\cap Y_1=f(U(x)\cap Y)$, there is $W\in\gm_1$ such that
$\gamma_2\cl W\sbs U$ and
$$  (1,2)\dd\ind\big((2,1)\dd\Fr_{X_1}W\cap Y_1,X_1\big)\leq k-1.  $$
Let $V(x)=f^{-1}(W)\cap U(x)\in\tau_1$ and let us prove that
$(2,1)\dd\Fr_XV(x)\cap Y\sbsq f^{-1}((2,1)\dd\Fr_{X_1}W)\cap Y$.

Indeed,
\begin{eqnarray*}
    &\ds (2,1)\dd\Fr_XV(x)=(2,1)\dd\Fr_X\big(f^{-1}(W)\cap U(x)\big)= \\
    &\ds =\big((2,1)\dd\Fr_Xf^{-1}(W)\cap(2,1)\dd\Fr_XU(x)\big)\cup \\
    &\ds \cup \big((2,1)\dd\Fr_Xf^{-1}(W)\cap U(x)\big)\cup
        \big((2,1)\dd\Fr_XU(x)\cap f^{-1}(W)\big),
\end{eqnarray*}
as $(2,1)\dd\Fr_XV(x)=2\dd\Fr_XV(x)$ for
$V(x)\in\tau_1\sbs\tau_2$. Since $W\cap Y_1\sbsq U\cap
Y_1=f(U(x)\cap Y)$, it is evident that $\big((2,1)\dd\Fr_XU(x)\cap
f^{-1}(W)\big)\cap Y=\vnth$ so that
\begin{eqnarray*}
    &\ds \big((2,1)\dd\Fr_X f^{-1}(W)\cap(2,1)\dd\Fr_XU(x)\big)\cup \\
    &\ds \cup \big((2,1)\dd\Fr_Xf^{-1}(W)\cap U(x)\big)\sbsq \\
    &\ds \sbsq (2,1)\dd\Fr_Xf^{-1}(W)\sbsq
    f^{-1}\big((2,1)\dd\Fr_{X_1}W\big).
\end{eqnarray*}
Therefore
$$  (2,1)\dd\Fr_XV(x)\cap Y\!\sbsq\!f^{-1}\big((2,1)\dd\Fr_{X_1}W\big)\cap Y\!=\!
        f^{-1}\big((2,1)\dd\Fr_{X_1}W\cap Y_1\big)  $$
and so
$$  f\big((2,1)\dd\Fr_XV(x)\cap Y\big)\!\sbs\!
        f\big((2,1)\dd\Fr_XV(x)\big)\cap Y_1\!\sbs\!
            (2,1)\dd\Fr_{X_1}W\cap Y_1.         $$
Hence, by inductive assumption and (3) of Theorem~7.3,
\allowdisplaybreaks
\begin{eqnarray*}
    &\ds (1,2)\dd\ind\big((2,1)\dd\Fr_XV(x)\cap Y,X\big)\leq \\
    &\ds \leq (1,2)\dd\ind\big(f\big((2,1)\dd\Fr_XV(x)\cap Y\big),X_1\big)\leq \\
    &\ds \leq (1,2)\dd\ind\big((2,1)\dd\Fr_{X_1}W\cap Y_1,X_1\big)\leq k-1,
\end{eqnarray*}
so that, $(1,2)\dd\ind(Y,X)\!\leq\!k$. Thus
$(1,2)\dd\ind(Y,X)\!\leq\!(1,2)\dd\ind(Y_1,X_1)$.

For the second part, first of all note that
$f\big|_Y:(Y,\tau_1'<\tau_2')\to (Y_1=f(Y),\gm_1'<\gm_2')$ is a
$d$-homeomorphism implies that $(1,2)\dd\ind Y=(1,2)\dd\ind Y_1$.
On the other hand, since $(X_1,\gm_1,\gm_2)$ is hereditarily
\linebreak       $p$-normal, by Proposition~7.6,
$(1,2)\dd\ind(Y_1,X_1)=(1,2)\dd\ind Y_1$. Hence, by (1) of
Theorem~7.3, it suffices to prove only that
$(1,2)\dd\ind(Y,X)\leq(1,2)\dd\ind Y$. Let $(1,2)\dd\ind Y=n$, the
inequality be proved for $n\leq k-1$ and prove it for $n=k$. If
$x\in Y$ is any point and $U(x)\in\tau_1$ is any neighborhood,
then there is $V'(x)\in\tau_1'$ such that $\tau_2'\cl V'(x)\sbsq
U(x)\cap Y$ and $(1,2)\dd\ind\big((2,1)\dd\Fr_YV'(x)\big)\leq
k-1$. Since $f\big|_Y$ is a $d$-homeomorphism, we have
$f(V'(x))\in\gm_1'$ and since $(X_1,\gm_1,\gm_2)$ is hereditarily
$p$-normal, by the proof of first part of Proposition~7.6, there
is $W\in\gm_1$ such that
$$  W\cap Y_1=f(V'(x)) \;\;\text{and}\;\;
        (2,1)\dd\Fr_{X_1}W\cap Y_1=(2,1)\dd\Fr_{Y_1}f(V'(x)).  $$
Let $V(x)=f^{-1}(W)\cap U(x)$. Then $V(x)\in\tau_1$ and $V(x)\cap
Y=V'(x)$. Let us prove that $(2,1)\dd\Fr_XV(x)\cap
Y=(2,1)\dd\Fr_YV'(x)$. First, let us show that
$$  (2,1)\dd\Fr_Xf^{-1}(W)\sbsq
         f^{-1}\big((2,1)\dd\Fr_{X_1}W\big).    $$
Indeed, let
$$  x\in (2,1)\dd\Fr_Xf^{-1}(W)=
        \tau_2\cl f^{-1}(W)\setminus f^{-1}(W).     $$
Then $f(x)\in f(\tau_2\cl f^{-1}(W))$ and $f(x)\,\ol{\in}\,W$.
Since $f$ is $d$-continuous, it is $2$-continuous and so
$$  f\big(\tau_2\cl f^{-1}(W)\big)\sbs
        \gm_2\cl\big(ff^{-1}(W)\big)=\gm_2\cl W.    $$
Therefore,
$$  f(x)\in\gm_2\cl W\setminus W=(2,1)\dd\Fr_{X_1}W $$
so that $x\in f^{-1}((2,1)\dd\Fr_{X_1}W)$.

Hence
\begin{eqnarray*}
    &\ds (2,1)\dd\Fr_Xf^{-1}(W)\cap Y\sbsq
            f^{-1}\big((2,1)\dd\Fr_{X_1}W\big)\cap Y\sbsq \\
    &\ds \sbsq f^{-1}\big((2,1)\dd\Fr_{X_1}W\cap Y_1\big)\cap Y= \\
    &\ds =f^{-1}\big((2,1)\dd\Fr_{Y_1}f(V'(x))\cap Y=
            (2,1)\dd\Fr_YV'(x).
\end{eqnarray*}
On the other hand,
\begin{eqnarray*}
    &\ds (2,1)\dd\Fr_XV(x)=
        \big((2,1)\dd\Fr_XU(x)\cap (2,1)\dd\Fr_Xf^{-1}(W)\big)\cup \\
     &\ds \cup \big((2,1)\dd\Fr_XU(x)\cap f^{-1}(W)\big)\cup
        \big((2,1)\dd\Fr_Xf^{-1}(W)\cap U(x)\big)
\end{eqnarray*}
and so $(2,1)\dd\Fr_XV(x)\cap Y=(2,1)\dd\Fr_YV'(x)$. Therefore, it
remains to use the inductive assumption with respect to
$(2,1)\dd\Fr_YV'(x)$, i.e.,
\begin{eqnarray*}
    &\ds (1,2)\dd\ind\big((2,1)\dd\Fr_XV(x)\cap Y,X\big)\!=\!
        (1,2)\dd\ind\big((2,1)\dd\Fr_YV'(x),X\big)\!\leq \\
    &\ds \leq (1,2)\dd\ind\big((2,1)\dd\Fr_YV'(x)\big)\leq k-1.
\end{eqnarray*}
Hence $(1,2)\dd\ind(Y,X)\leq k$ and thus,
$(1,2)\dd\ind(Y,X)\leq(1,2)\dd\ind Y$.~\end{pf}

\begin{proposition}{7.9}
If $(X,\tau_1<_S\tau_2)$ is $(i,j)$-regular, $Y\in 2\dd\,\cD(X)$
and $Z\sbs Y\sbs X$, then $(i,j)\dd\ind(Z,Y)=(i,j)\dd\ind(Z,X)$.
\end{proposition}

\begin{pf}
By (3) of Theorem~7.3, it suffices to prove that
$$  (i,j)\dd\ind(Z,X)\leq (i,j)\dd\ind(Z,Y) $$
by induction under $n=(i,j)\dd\ind(Z,Y)$. Suppose that the
inequality is proved for $n\leq k-1$ and prove it for $n=k$. Let
$x\in Z$ be any point and $U(x)\in\tau_i$ be any neighborhood.
Since $(i,j)\dd\ind(Z,Y)=k$, there is $V'(x)\in\tau_i'$ in
$(Y,\tau_1',\tau_2')$ such that $\tau_j'\cl V'(x)\sbs U(x)\cap Y$
and $(i,j)\dd\ind\big((j,i)\dd\Fr_YV'(x)\cap Z,Y\big)\leq k-1$.
Furthermore, by inductive hypothesis, we have
\begin{eqnarray*}
    &\ds (i,j)\dd\ind\big((j,i)\dd\Fr_YV'(x)\cap Z,X\big)\leq \\
    &\ds \leq (i,j)\dd\ind\big((j,i)\dd\Fr_YV'(x)\cap Z,Y\big)\leq k-1.
\end{eqnarray*}
Let $W(x)\in\tau_i$ and $W(x)\cap Y=V'(x)$. Since $V'(x)\sbs
U(x)$, for $V(x)=W(x)\cap U(x)\in\tau_i$ we have $V(x)\sbs U(x)$
and $V(x)\cap Y=V'(x)$. By $Y\in 2\dd\,\cD(X)$, $V(x)\cap Y=V'(x)$
and lemma 2.16,
$$  \tau_j\cl V(x)=\tau_j\cl\big(V(x)\cap Y\big)=\tau_j\cl V'(x). $$
Therefore, since $Z\sbs Y$, we have
\begin{eqnarray*}
    &\ds (j,i)\dd\Fr_XV(x)\cap Z=
        \big(\tau_j\cl V(x)\setminus V(x)\big)\cap Z= \\
    &\ds =\big(\big(\tau_j\cl V(x)\setminus V(x)\big)\cap Y\big)\cap Z=
        \big((\tau_j\cl V'(x)\setminus V(x))\cap Y\big)\cap Z= \\
    &\ds =\big(\tau_j'\cl V'(x)\setminus V'(x)\big)\cap Z=
            (j,i)\dd\Fr_YV'(x)\cap Z.
\end{eqnarray*}
Moreover, since $X$ is $(i,j)$-regular, one can assume that
$\tau_j\cl V(x)\sbs U(x)$. Hence, for each $x\in Z$ and any
neighborhood $U(x)\in\tau_i$ there is $V(x)\in\tau_i$ such that
$\tau_j\cl V(x)\sbs U(x)$ and
\begin{eqnarray*}
    (i,j)\dd\ind\big((j,i)\dd\Fr_XV(x)\cap Z,X\big) &&\hskip-0.6cm = \\
    =(i,j)\dd\ind\big((j,i)\dd\Fr_YV'(x)\cap Z,X\big) &&\hskip-0.6cm \leq \\
    \leq (i,j)\dd\ind\big((j,i)\dd\Fr_YV'(x)\cap Z,Y\big) &&\hskip-0.6cm \leq k-1.
\end{eqnarray*}
Therefore, $(i,j)\dd\ind(Z,X)\leq k$ and thus,
$(i,j)\dd\ind(Z,X)\leq     \linebreak
(i,j)\dd\ind(Z,Y)$.
\end{pf}

\begin{proposition}{7.10}
If $(i,j)\dd\ind(Y,X)$ is finite, then $(Y,\tau_1',\tau_2')$ is
$(i,j)$-superregular in $X$ and so, $Y$ is
$(i,j)\dd\WS$-superregular in $X$, $Y$ is $(i,j)$-strongly regular
in $X$, $Y$ is $(i,j)$-regular in $X$, $Y$ is $(i,j)\dd\WS$-quasi
regular in $X$, $Y$ is $(i,j)\dd\WS$-regular in $X$ and $Y$ is
$(i,j)$-regular.
\end{proposition}

\begin{pf}
Let $x\in Y$, $F\in\co\tau_i$ and $x\,\ol{\in}\,F$. Then $x\in
X\setminus F=U(x)\in\tau_i$ and since $(i,j)\dd\ind(Y,X)<\infty$,
there is $V(x)\in\tau_i$ such that $\tau_j\cl V(x)\sbs U(x)$.
Hence
$$  F=X\setminus U(x)\sbs X\setminus \tau_j\cl V(x)=
        \tau_j\nt\big(X\setminus V(x)\big)=U(F)\in\tau_j   $$
and $U(x)\cap U(F)=\vnth$. Thus $Y$ is $(i,j)$-superregular in
$X$.

The rest follows from implications after Definition~2.2.
\end{pf}

\begin{corollary}{7.11}
If $\ind(Y,X)$ is finite, then $(Y,\tau')$ is superregular in $X$
and so, $Y$ is $\WS$-superregular in $X$, $Y$ is strongly regular
in $X$, $Y$ is regular in $X$, $Y$ is $\WS$-quasi regular in $X$,
$Y$ is $\WS$-regular in $X$ and $Y$ is regular.
\end{corollary}

\begin{proposition}{7.12}
If $(X,\tau_1,\tau_2)$ is $p$-normal, $Y\in\co\tau_1\cap\co\tau_2$
and $(i,j)\dd\ind Y=0$, then $(i,j)\dd\ind(Y,X)=0$.
\end{proposition}

\begin{pf}
We shall prove that for each point $x\in Y$ and any neighborhood
$U(x)\in\tau_i$ there is $V(x)\in\tau_i$ such that $\tau_j\cl
V(x)\sbs U(x)$ and $(\tau_j\cl V(x)\setminus V(x))\cap Y=\vnth$.
Since $(i,j)\dd\ind Y=0$, for $U(x)\in\tau_i$ there is
$V\in\tau_i'\cap\co\tau_j'$ in $(Y,\tau_1',\tau_2')$ such that
$x\in V\sbs U(x)\cap Y$. Since $Y\in\co\tau_1\cap\co\tau_2$ and
$(X,\tau_1,\tau_2)$ is $p$-normal, for $V\in\co\tau_j$ and
$Y\setminus V\in\co\tau_i$ there are $U(V)\in\tau_i$,
$U(Y\setminus V)\in\tau_j$ such that $U(V)\cap U(Y\setminus
V)=\vnth$. Moreover, since $V\sbs U(x)$, by (4) of Definition~1.1,
there is $E(V)\in\tau_i$ such that $V\sbs E(V)\sbs\tau_j\cl
E(V)\sbs U(x)$. Let $W(V)=E(V)\cap U(V)$. Then
$$  V\sbs W(V)\sbs\tau_j\cl W(V)\sbs
        X\setminus(Y\setminus V)=(X\setminus Y)\cup V.    $$
It is evident that $\tau_j\cl W(V)\sbs U(x)$ and $\tau_j\cl
W(V)\cap Y=\tau_j\cl V=V$. Hence
\[  (j,i)\dd\Fr_XW(V)\cap Y=
        \big(\tau_j\cl W(V)\setminus W(V)\big)\cap Y=
                V\setminus V=\vnth.    \]
\vskip-0.75cm
\end{pf}
\vskip+0.2cm

\begin{proposition}{7.13}
If $(Y,\tau_1',\tau_2')$ is $(i,j)\dd\WS$-superregular in $X$ and
$Y$ is $(i,j)$-extremally disconnected in $X$, then
$(i,j)\dd\ind(Y,X)=0$.
\end{proposition}

\begin{pf}
Indeed, by Proposition~3.7, $(i,j)\dd\ind_xX=0$ for each point
$x\in Y$, i.e., for any $U(x)\in\tau_i$ there is
$V(x)\in\tau_i\cap\co\tau_j$ such that $V(x)\sbs U(x)$. Clearly,
$(j,i)\dd\Fr_XV(x)=\vnth$ implies that $(j,i)\dd\Fr_XV(x)\cap
Y=\vnth$ and thus $(i,j)\dd\ind(Y,X)=0$.
\end{pf}

\begin{corollary}{7.14}
If $(Y,\tau')$ is $\WS$-superregular in $X$ and $Y$ is extremally
disconnected in $X$, then $\ind(Y,X)=0$.
\end{corollary}

\begin{definition}{7.15}
For a $\BsS$ $(Y,\tau_1',\tau_2')$ of a $\BS$ $(X,\tau_1,\tau_2)$
a pair $(x,A)$, where $x\in Y$, $A\in\co\tau_i$ and
$x\ol{\in}\,A$, is said to be a relatively $i$-re\-gu\-lar pair. A
relative partition, corresponding to the pair $(x,A)$, is a
$p$-closed set $T\sbs X$ such that $X\setminus T$ is not
$p$-connected and
\begin{eqnarray*}
    &\ds (X\setminus T)\cap Y=Y\setminus T\sbs H_1\cup H_2, \;\;\;
            x\in H_i'=H_i\cap Y\in\tau_i', \\
    &\ds A\sbs H_j\in\tau_j \;\;\text{and}\;\; H_1\cap H_2=\vnth.
\end{eqnarray*}
\end{definition}

\begin{remark}{7.16}
Let $(Y,\tau_1',\tau_2')$ be a $\BsS$ of a $\BS$
$(X,\tau_1,\tau_2)$ and $(x,A)$ be a relatively $i$-regular pair.
Then
\begin{enumerate}
\item[(1)] If there is a neighborhood $U(x)\in\tau_i$
$(U(A)\in\tau_j)$ such that $\tau_j\cl U(x)\sbs X\setminus A$
$(\tau_i\cl U(A)\sbs X\setminus\{x\})$, then $(j,i)\dd\Fr_XV(x)$
$((i,j)\dd\Fr_XU(A))$ is a relative partition, corresponding to
$(x,A)$.

\item[(2)] If $T$ is a relative partition, corresponding to
$(x,A)$, then          \linebreak       $(j,i)\dd\Fr_XH_i\cap
Y\sbs T\cap Y$ $((i,j)\dd\Fr_XH_j\cap Y\sbs T\cap Y)$.
\end{enumerate}
\end{remark}

Indeed, suppose, for example, in (1) that $U(x)\in\tau_i$ and
$\tau_j\cl U(x)\sbs X\setminus A$. Then
\begin{eqnarray*}
    &\ds X\setminus(j,i)\dd\Fr_XU(x)=
        \big(X\setminus \tau_j\cl U(x)\big)\cup U(x)=
            H_j\cup H_i, \\
    &\ds H_i\in\tau_i\setminus\{\vnth\}, \;\;\;
        H_1\cap H_2=\vnth, \\
    &\ds x\in H_i'\sbs H_i=U(x), \;\;\;
            A\sbs H_j=X\setminus\tau_j\cl U(x).
\end{eqnarray*}
Hence, taking into account Definition~7.15, $(j,i)\dd\Fr_XU(x)$ is
a relative partition, corresponding to $(x,A)$ since
$$  \big(X\setminus(j,i)\dd\Fr_XU(x)\big)\cap Y=
        Y\setminus(j,i)\dd\Fr_XU(x)\sbs H_1\cup H_2.    $$

(2) Since $T$ is a relative partition, corresponding to $(x,A)$,
we have
$$  (X\setminus T)\cap Y=Y\setminus T\sbs H_1\cup H_2,   $$
where $H_i\in\tau_i\setminus\{\vnth\}$, $H_1\cap H_2=\vnth$, $x\in
H_i'=H_i\cap Y$ and $A\sbs H_j$. Hence
\begin{eqnarray*}
    &\ds (j,i)\dd\Fr_XH_i\cap(Y\setminus T)= \\
    &\ds =\!\big(\tau_j\cl H_i\cap(X\!\setminus\!H_i)\big)\cap(Y\!\setminus\!T)\!\sbs\!
        (\tau_j\cl H_i\setminus\!H_i)\cap(H_i\!\cup\!H_j)\!=\!\vnth
\end{eqnarray*}
since $H_j\in\tau_j$. Thus $(j,i)\dd\Fr_XH_i\cap Y\sbs T\cap Y$.

\begin{theorem}{7.17}
Let $(Y,\tau_1',\tau_2')$ be a $\BsS$ of a $\BS$
$(X,\tau_1,\tau_2)$. Then $(i,j)\dd\ind(Y,X)\leq n$ if and only if
to every relatively $i$-regular pair $(x,A)$, $x\in Y$,
$A\in\co\tau_i$ and $x\ol{\in}\,A$, there corresponds a relative
partition $T$ such that $(i,j)\dd\ind(T\cap Y,X)\leq n-1$.
\end{theorem}

\begin{pf}
Let $(i,j)\dd\ind(Y,X)\leq n$ and $(x,A)$ be any relatively
$i$-regular pair. Then $x\in U(x)=X\setminus A\in\tau_i$. Since
$(i,j)\dd\ind(Y,X)\leq n$, by (2) of Definition~7.1, there is
$V(x)\in \tau_i$ such that $\tau_j\cl V(x)\sbs U(x)$ and
$$  (i,j)\dd\ind\big((j,i)\dd\Fr_XV(x)\cap Y,X\big)\leq n-1.    $$

By (1) of Remark~7.16, $(j,i)\dd\Fr_XV(x)$ is a relative
partition, corresponding to $(x,A)$ and so, the first part is
proved.

Conversely, let us suppose that the condition is satisfied and
prove that $(i,j)\dd\ind(Y,X)\leq n$. If $x\in Y$ is any point and
$U(x)\in\tau_i$ is any neighborhood, then the pair
$(x,A=X\setminus U(x))$ is a relatively $i$-regular pair. Hence,
by condition, there is a relative partition $T$ for $(x,A)$ such
that $(i,j)\dd\ind(T\cap Y,X)\leq n-1$. But $Y\setminus T\sbs
H_1\cup H_2$, where $x\in H_i'\sbs H_i$, $A\sbs H_j$ and
$\tau_j\cl H_i\cap H_j=\vnth$, so that $\tau_j\cl H_i\sbs
X\setminus H_j\sbs X\setminus A=U(x)$. Let $H_i=V(x)$. Then
$\tau_j\cl V(x)\sbs U(x)$ and by (2) of Remark~7.16,
$$  (j,i)\dd\Fr_XH_i\cap Y=(j,i)\dd\Fr_XV(x)\cap Y\sbs T\cap Y. $$
Since $(i,j)\dd\ind(T\cap Y,X)\leq n-1$ and $(j,i)\dd\Fr_XV(x)\cap
Y\sbs T\cap Y$, by (3) of Theorem~7.3,
$$  (i,j)\dd\ind\big((j,i)\dd\Fr_XV(x)\cap Y,X\big)\leq n-1 $$
and thus $(i,j)\dd\ind(Y,X)\leq n$.
\end{pf}

\begin{definition}{7.18}
Let $(Y,\tau')$ be a $\TsS$ of a $\TS$ $(X,\tau)$. Then a pair
$(x,A)$, where $x\in Y$, $A\in\co\tau$ and $x\ol{\in}\,A$, is said
to be a relatively regular pair. A relative partition,
corresponding to the pair $(x,A)$, is a closed set $T\sbs X$ such
that $(X\setminus T)\cap Y=Y\setminus T\sbs H_1\cup H_2$, $x\in
H_i'=H_i\cap Y\in\tau'$, $A\sbs H_j\in\tau$ and $H_1\cap
H_2=\vnth$.
\end{definition}

Now, if $(Y,\tau')$ is a $\TsS$ of a $\TS$ $(X,\tau)$ and $(x,A)$
is a relatively regular pair, then the topological version of
Remark~7.14 gives:
\begin{enumerate}
\item[(1)] If there is a neighborhood $U(x)\in\tau$
$(U(A)\in\tau)$ such that $\tau\cl U(x)\sbs X\setminus A$
$(\tau\cl U(A)\sbs X\setminus\{x\})$, then $\Fr_XU(x)$ \linebreak
$(\Fr_XU(A))$ is a relative partition, corresponding to $(x,A)$.

\item[(2)] If $T$ is a relative partition, corresponding to
$(x,A)$, then         \linebreak       $\Fr_X H_i\cap Y\sbs T\cap
Y$.
\end{enumerate}

By analogy with Theorem~7.17, we obtain the following new
characterization of the small relative inductive dimension.

\begin{theorem}{7.19}
Let $(Y,\tau')$ be a $\TsS$ of a $\TS$ $(X,\tau)$. Then \linebreak
$\ind(Y,X)\leq n$ in the sense of Filippov if and only if to every
relatively regular pair $(x,A)$, $x\in Y$, $A\in\co\tau$ and
$x\ol{\in}\,A$, there corresponds a relative partition $T$ such
that $\ind(T\cap Y,X)\leq n-1$.
\end{theorem}

\vskip+0.5cm
\section*{\textbf{8. $(i,j)$-Large Relative Inductive Dimension Functions}}
\vskip+0.2cm

\begin{definition}{8.1}
Let $(Y,\tau_1',\tau_2')$ be a $\BsS$ of a $\BS$
$(X,\tau_1,\tau_2)$ and $n$ denote a nonnegative integer. We say
that:
\begin{enumerate}
\item[(1)] $(i,j)\dd\Ind(Y,X)=-1$ if and only if $Y=\vnth$.

\item[(2)] $(i,j)\dd\Ind(Y,X)\leq n$ if for any set
$F\in\co\tau_j'$ and any neighborhood $U(F)\in\tau_i$ there is a
neighborhood $V(F)\in\tau_i$ such that $\tau_j\cl V(F)\sbs U(F)$
and $(i,j)\dd\ind\big((j,i)\dd\Fr_XV(F)\cap Y,X\big)\leq n-1$.

\item[(3)] $(i,j)\dd\Ind(Y,X)\!=\!n$ if
$(i,j)\dd\Ind(Y,X)\!\leq\!n$ and $(i,j)\dd\ind(Y,X)\!\leq n-1$
does not hold.

\item[(4)] $(i,j)\dd\Ind(Y,X)=\infty$ if the inequality
$(i,j)\dd\Ind(Y,X)\leq n$ does not hold for any $n$.
\end{enumerate}
\end{definition}

As a rule
$$  p\,\dd\Ind(Y,X)\leq n\llra \big((1,2)\dd\Ind(Y,X)\leq n\wedge
        (2,1)\dd\Ind(Y,X)\leq n\big).       $$

It follows immediately from (2) of Definition~8.1 and (3) of
Proposition~2.9 that if $(i,j)\dd\Ind(Y,X)$ is finite, then $Y$ is
$(j,i)$-supernormal in $X$ and hence, $Y$ is $p$-strongly normal
in $X$, $Y$ is $(j,i)\dd\WS$-supernormal in $X$, $Y$ is $p$-normal
in $X$, $Y$ is $(j,i)\dd\WS$-normal in $X$, $Y$ is $p$-quasi
normal in $X$, $Y$ is $p$-internally normal in $X$, $Y$ is
$p$-normal and $Y$ is $p$-normal in $X$ from inside.

\begin{remark}{8.2}
By analogy with Remark~7.2 one can prove that $(i,j)\Ind(Y,X)\leq
n$ if and only if for any set $F\in\co\tau_j'$ and any
neighborhood $U(F)\in\tau_i$ there is a neighborhood
$V(F)\in\tau_i$ such that $\tau_j\cl V(F)\sbs U(F)$ and
$(i,j)\dd\Ind\big((j,i)\dd\Fr_X V(F)\cap Y\big)\leq n-1$.
\end{remark}

\begin{definition}{8.3}
Let $(Y,\tau')$ be a $\TsS$ of a $\TS$ $(X,\tau)$ and $n$ denote a
nonnegative integer. We say that:
\begin{enumerate}
\item[(1)] $\Ind(Y,X)=-1$ if and only if $Y=\vnth$.

\item[(2)] $\Ind(Y,X)\leq n$ if for any set $F\in\co\tau'$ and any
neighborhood $U(F)\in\tau$ there is a neighborhood $V(F)\in\tau$
such that $\tau\cl V(F)\sbs U(F)$ and $\Ind(\Fr_XV(F)\cap Y,X)\leq
n-1(\llra$ for any set $F\in\co\tau'$ and any neighborhood
$U(F)\in\tau$ there is a neighborhood $V(F)\in\tau$ such that
$\tau\cl V(F)\sbs U(F)$ and $\Ind(\Fr_XV(F)\cap Y)\leq n-1)$.

\item[(3)] $\Ind(Y,X)=n$ if $\Ind(Y,X)\leq n$ and $\Ind(Y,X)\leq
n-1$ does not hold.

\item[(4)] $\Ind(Y,X)=\infty$ if $\Ind(Y,X)\leq n$ does not hold
for any $n$.
\end{enumerate}
\end{definition}

Hence, if $\Ind(Y,X)$ is finite, then $(Y,\tau')$ is supernormal
in $X$, $Y$ is strongly normal in $X$, $Y$ is $\WS$-supernormal is
$X$, $Y$ is normal in $X$, $Y$ is $\WS$-normal in $X$, $Y$ is
quasi normal in $X$, $Y$ is internally normal in $X$, $Y$ is
normal and $Y$ is normal in $X$ from inside.

\begin{proposition}{8.4}
If $(Y,\tau_1',\tau_2')$ is a $j\dd\TT_1$ $\BsS$ of a $\BS$
$(X,\tau_1,\tau_2)$, then $(i,j)\dd\ind(Y,X)\leq
(i,j)\dd\Ind(Y,X)$.

Hence, if $(Y,\tau')$ is a $T_1$ $\TsS$ of a $\TS$ $(X,\tau)$,
then $\ind(Y,X)\leq\Ind(Y,X)$.
\end{proposition}

\begin{theorem}{8.5}
For a $\BsS$ $(Y,\tau_1',\tau_2')$ of a $\BS$ $(X,\tau_1,\tau_2)$
the following conditions are satisfied:
\begin{enumerate}
\item[(1)] $(i,j)\dd\Ind Y\leq (i,j)\dd\Ind(Y,X)$.

\item[(2)] If $Y\in\co\tau_j$, then $(i,j)\dd\Ind(Y,X)\leq
(i,j)\dd\Ind X$.

\item[(3)] $(i,j)\dd\Ind(X,X)=(i,j)\dd\Ind X$.
\end{enumerate}
\end{theorem}

\begin{pf}
(1) We shall use the induction under a nonnegative integer
$n=(i,j)\dd\Ind(Y,X)$. Let the inequality be proved for $n\leq
k-1$ and prove it for $n=k$. Let $F\in\co\tau_j'$ and
$U'(F)\in\tau_i'$. If $U(F)\in\tau_i$, $U(F)\cap Y=U'(F)$, then
$(i,j)\dd\Ind(Y,X)=k$ implies that there is $V(F)\in\tau_i$ such
that $\tau_j\cl V(F)\sbs U(F)$ and
$(i,j)\dd\Ind\big((j,i)\dd\Fr_XV(F)\cap Y,X\big)\leq k-1$. Hence,
by inductive assumption,
$$  (i,j)\dd\Ind\big((j,i)\dd\Fr_XV(F)\cap Y\big)\leq k-1.  $$
Let $V'(F)=V(F)\cap Y$. Then
\begin{eqnarray*}
    &\ds (j,i)\dd\Fr_YV'(F)=\tau_j'\cl V'(F)\setminus V'(F)= \\
    &\ds =\big(\tau_j\cl V'(F)\cap Y\big)\setminus
                \big(V(F)\cap Y\big)\sbsq
        \big(\tau_j\cl V(F)\setminus V(F)\big)\cap Y= \\
    &\ds =(j,i)\dd\Fr_XV(F)\cap Y
\end{eqnarray*}
and since $(j,i)\dd\Fr_YV'(F)$ is $p$-closed in $Y$, it is also
$p$-closed in       \linebreak  $(j,i)\dd\Fr_XV(F)\cap Y$.
Therefore, by (2) of Proposition~3.2.7 in [8],
$$  (i,j)\dd\Ind\big((j,i)\dd\Fr_YV'(F)\big)\leq
        (i,j)\dd\Ind\big((j,i)\dd\Fr_XV(F)\cap Y\big)\leq k-1. $$
Hence, by Corollary~3.2.9 in [8], $(i,j)\dd\Ind Y\leq k$ and so
$(i,j)\dd\Ind Y\leq (i,j)\dd\Ind(Y,X)$.

(2) The induction will be used under a nonnegative integer
$n=(i,j)\dd\Ind X$. Let the inequality be proved for $n\leq k-1$
and prove it for $n=k$. Let $F\in\co\tau_j'\sbs\co\tau_j$ and
$U(F)\in\tau_i$. Since $(i,j)\dd\Ind X=k$, by Corollary~3.2.9 in
[8] there is $V(F)\in\tau_i$ such that $\tau_j\cl V(F)\sbs U(F)$
and $(i,j)\dd\Ind(j,i)\dd\Fr_XV(F)\leq k-1$. Clearly,
$Y\in\co\tau_j\sbs p\,\dd\,\Cl(X)$ implies that the set
$(j,i)\dd\Fr_XV(F)\cap Y$ is $p$-closed in $X$ and so, by the
second part of Remark~2.12, it is $p$-closed in
$(j,i)\dd\Fr_XV(F)$ too. Hence, by inductive assumption,
$$  (i,j)\dd\Ind\big((j,i)\dd\Fr_XV(F)\cap Y,(j,i)\dd\Fr_XV(F)\big)\leq k-1. $$
and by (1), $(i,j)\dd\Ind\big((j,i)\dd\Fr_XV(F)\cap Y\big)\leq
k-1$. Therefore, by Remark~8.2, $(i,j)\dd\Ind(Y,X)\leq k$ and thus
$(i,j)\dd\Ind(Y,X)\leq (i,j)\dd\Ind X$.

(3) If $Y=X$, then (1) and (2) imply that
$$  (i,j)\dd\Ind X\leq (i,j)\dd\Ind(X,X)\leq (i,j)\dd\Ind X $$
and so $(i,j)\dd\Ind(X,X)=(i,j)\dd\Ind X$.
\end{pf}

\begin{corollary}{8.6}
Let $(Y,\tau')$ be a $\TsS$ of a $\TS$ $(X,\tau)$. Then
\begin{enumerate}
\item[(1)] $\Ind Y\leq\Ind(Y,X)$.

\item[(2)] If $Y\in\co\tau$, then $\Ind(Y,X)\leq\Ind X$.

\item[(3)] $\Ind(X,X)=\Ind X$.
\end{enumerate}
\end{corollary}

\begin{proposition}{8.7}
For a $\BsS$ $(Y,\tau_1',\tau_2')$ of a hereditarily $p$-normal
$\BS$ $(X,\tau_1,\tau_2)$ we have $(i,j)\dd\Ind(Y,X)=(i,j)\dd\Ind
Y$.
\end{proposition}

\begin{pf}
By (1) of Theorem~8.5, it suffices to prove only that
$$  (i,j)\dd\Ind(Y,X)\leq (i,j)\dd\Ind Y.   $$
First of all, recall that by the proof of the first part of
Proposition~7.6, for any set $U'\in\tau_i'$ there is a set
$U\in\tau_i$ such that $U\cap Y=U'$ and $(j,i)\dd\Fr_XU\cap
Y=(j,i)\dd\Fr_YU'$.

Now, we shall use the induction under a nonnegative integer $n=
     \linebreak       (i,j)\dd\Ind Y$. Let the inequality be proved for $n\leq k-1$ and
prove it for $n=k$. Let $F\in\co\tau_j'$ and $U(F)\in\tau_i$ be an
arbitrary neighborhood. Since $(i,j)\dd\Ind Y=k$, there is
$V'(F)\in\tau_i'$ such that $\tau_j'\cl V'(F)\sbs U(F)\cap Y$ and
$(i,j)\dd\Ind\big((j,i)\dd\Fr_YV'(F)\big)\leq k-1$. For the set
$V'(F)$ there is a set $V(F)\in\tau_i$ such that $V(F)\cap
Y=V'(F)$ and
$$  (i,j)\dd\Ind\big((j,i)\dd\Fr_XV(F)\cap Y\big)\leq k-1. $$
Hence, by inductive assumption,
\begin{eqnarray*}
    &\ds (i,j)\dd\Ind\big((j,i)\dd\Fr_XV(F)\cap Y,X\big)\leq \\
    &\ds \leq (i,j)\dd\Ind\big((j,i)\dd\Fr_XV(F)\cap Y\big)\leq k-1.
\end{eqnarray*}
Therefore, $(i,j)\dd\Ind(Y,X)\!\leq\!k$ and so,
$(i,j)\dd\Ind(Y,X)\!\leq\!(i,j)\dd\Ind Y$.~\end{pf}

\begin{corollary}{8.8}
If $(Y,\tau')$ is a $\TsS$ of a hereditarily normal $\TS$
$(X,\tau)$, then $\Ind(Y,X)=\Ind Y$.
\end{corollary}

Take place the following important \textbf{monotonicity property
of the $(i,j)$-large relative inductive dimension functions}.

\begin{theorem}{8.9}
If $(Z,\tau_1'',\tau_2'')\sbs (Y,\tau_1',\tau_2')\sbs
(X,\tau_1,\tau_2)$, the $\BS$ $(X,\tau_1,\tau_2)$ is hereditarily
$p$-normal and $Z\in p\,\dd\,\Cl(Y)$, then
$$  (i,j)\dd\Ind(Z,X)\leq (i,j)\dd\Ind(Y,X).        $$
\end{theorem}

\begin{pf}
Since $(X,\tau_1,\tau_2)$ is hereditarily $p$-normal, by
Proposition~8.7, $(i,j)\dd\Ind(Y,X)=(i,j)\dd\Ind Y$ and
$(i,j)\dd\Ind(Z,X)=(i,j)\dd\Ind Z$. On the other hand, since $Z\in
p\,\dd\,\Cl(Y)$, by (2) of Proposition~3.2.7 in [8],
$$  (i,j)\dd\Ind(Z,X)=(i,j)\dd\Ind Z\leq
        (i,j)\dd\Ind Y=(i,j)\dd\Ind(Y,X).       $$
\vskip-1cm
\end{pf}
\vskip+0.2cm

\begin{corollary}{8.10}
If $(Z,\tau'')\sbs (Y,\tau')\sbs (X,\tau)$, the $\TS$ $(X,\tau)$
is hereditarily normal and $Z\in\co\tau'$, then $\Ind(Z,X)\leq
\Ind(Y,X)$.
\end{corollary}

\begin{proposition}{8.11}
If $(Y,\!\tau_1',\!\tau_2')$ is $(j,\!i)\dd\WS$-supernormal in
$(X,\!\tau_1,\!\tau_2)$ and $Y$ is $(i,j)$-extremally disconnected
in $X$, then $(i,j)\dd\Ind(Y,X)=0$.
\end{proposition}

\begin{pf}
Let $F\in\co\tau_j'$ and $U(F)\in\tau_i$. Then by (2) of
Proposition~2.9, there is $V'(F)\in\tau_i'$ such that $\tau_j\cl
V'(F)\sbs U(F)$. Since $Y$ is $(i,j)$-extremally disconnected in
$X$, $\tau_j\cl V'(F)=V(F)\in\tau_i\cap\co\tau_j$ and so
$(j,i)\dd\Fr_XV(F)=\vnth$. Hence $(j,i)\dd\Fr_XV(F)\cap Y=\vnth$
and thus $(i,j)\dd\Ind(Y,X)=0$.
\end{pf}

\begin{corollary}{8.12}
If $(Y,\tau')$ is $\WS$-supernormal in $(X,\tau)$ and $Y$ is
extremally disconnected in $X$, then $\Ind(Y,X)=0$.
\end{corollary}

\begin{definition}{8.13}
A $\BsS$ $(Y,\tau_1',\tau_2')$ of a $\BS$ $(X,\tau_1,\tau_2)$ is
said to be hereditarily $(i,j)\dd\WS$-supernormal in $X$ if every
set $Z\in p\,\dd\,\Cl(Y)$ is $(i,j)\dd\WS$-supernormal in $X$.
\end{definition}

It is evident that every hereditarily $p\,\dd\WS$-supernormal
$\BsS$ $Y$ of $X$ is $p\,\dd\WS$-supernormal in $X$ and, moreover,
by the second part of Remark~2.12, every $p$-closed subset of a
hereditarily $p\,\dd\WS$-supernormal $\BsS$ $Y$ of $X$ is also
hereditarily $p\,\dd\WS$-supernormal in $X$.

\begin{theorem}{8.14}
If $(Z,\tau_1'',\tau_2'')\sbs (Y,\tau_1',\tau_2')\sbs
(X,\tau_1,\tau_2)$, where $Y$ is hereditarily
$(j,i)\dd\WS$-supernormal in $X$ and $Z\in p\,\dd\,\Cl(Y)$, then
$$  (i,j)\dd\Ind(Z,Y)\leq (i,j)\dd\Ind(Z,X)\leq (i,j)\dd\Ind(Y,X).  $$
\end{theorem}

\begin{pf}
First, let us prove the right inequality, i.e., the monotonicity
property of the $(i,j)$-large relative inductive dimensions, by
induction under a nonnegative integer $n=(i,j)\dd\Ind(Y,X)$. Let
this inequality be proved for $n\leq k-1$ and prove it for $n=k$.
If $F\in\co\tau_j''$ and $U(F)\in\tau_i$ are any sets, then
$F\cap\Phi=\vnth$, where $\Phi=X\setminus U(F)\in\co\tau_i$. Since
$Z\in p\,\dd\,\Cl(Y)$, by condition $Z$ is
$(j,i)\dd\WS$-supernormal in $X$ and so, there are
$V(F)\in\tau_i''$, $V(\Phi)\in\tau_j$ such that $\tau_j\cl
V(F)\cap V(\Phi)=\vnth$. Hence $\tau_j\cl F\cap Y=\tau_j'\cl F\sbs
U(F)$ as
$$  \tau_j\cl F\sbs\tau_j\cl V(F)\sbs X\setminus V(\Phi)\sbs
                    X\setminus\Phi=U(F).        $$
Since $(i,j)\dd\Ind(Y,X)=k$, there is $W(\tau_j'\cl
F)=W(F)\in\tau_i$ such that $\tau_j\cl W(F)\sbs U(F)$ and
$$  (i,j)\dd\Ind\big((j,i)\dd\Fr_XW(F)\cap Y,X\big)\leq k-1.    $$
Clearly, we have
$$  (j,i)\dd\Fr_XW(F)\cap Z\in
        p\,\dd\,\Cl\big((j,i)\dd\Fr_XW(F)\cap Y\big),  $$
where by remark before Theorem~8.14, $(j,i)\dd\Fr_XW(F)\cap Y\in
p\,\dd\,\Cl(Y)$ implies that $(j,i)\dd\Fr_XW(F)\cap Y$ is also
hereditarily $(j,i)\dd\WS$-su\-per\-nor\-mal in $X$. Therefore, by
inductive hypothesis,
\begin{eqnarray*}
    (i,j)\dd\Ind\big((j,i)\dd\Fr_XW(F)\cap Z,X\big) &&\hskip-0.6cm \leq \\
    \leq (i,j)\dd\Ind\big((j,i)\dd\Fr_XW(F)\cap Y,X\big) &&\hskip-0.6cm \leq k-1
\end{eqnarray*}
and hence, $(i,j)\dd\Ind(Z,X)\!\leq\!k$. Thus
$(i,j)\dd\Ind(Z,X)\!\leq\!(i,j)\dd\Ind(Y,X)$.

Now, we also use the induction under a nonnegative integer
$n=(i,j)\dd\Ind(Z,X)$. Let the left inequality be proved for
$n\leq k-1$ and prove it for $n=k$. Let $F\in\co\tau_i''$ and
$U'(F)\in\tau_i'$ be any sets. Since $(i,j)\dd\Ind(Z,X)=k$, for
$F$ and $U(F)\in\tau_i$, where $U(F)\cap Y=U'(F)$, there is
$V(F)\in\tau_i$ such that $\tau_j\cl V(F)\sbs U(F)$ and
$$  (i,j)\dd\Ind\big((j,i)\dd\Fr_XV(F)\cap Z,X\big)\leq k-1.    $$
Let $V'(F)=V(F)\cap Y$. Then $\tau_j'\cl V'(F)\sbs U'(F)$ and
since
$$  (j,i)\dd\Fr_YV'(F)\sbsq (j,i)\dd\Fr_XV(F)\cap Y, \;\;\;
                Z\sbs Y,  $$
we have
$$  (j,i)\dd\Fr_YV'(F)\cap Z\sbsq (j,i)\dd\Fr_XV(F)\cap Z.  $$
Since by the second part of Remark~2.12,
$$  (j,i)\dd\Fr_XV(F)\cap Z\in p\,\dd\,\Cl(Z)\sbs p\,\dd\,\Cl(Y),   $$
the $\BsS$ $(j,i)\dd\Fr_XV(F)\cap Z$ is hereditarily
$p\,\dd\WS$-supernormal in $X$. Moreover,
$$  (j,i)\dd\Fr_YV'(F)\cap Z\in
        p\,\dd\,\Cl\big((j,i)\dd\Fr_XV(F)\cap Z\big)  $$
and by the right inequality,
\begin{eqnarray*}
    (i,j)\dd\Ind\big((j,i)\dd\Fr_YV'(F)\cap Z,X\big) &&\hskip-0.6cm \leq \\
    \leq (i,j)\dd\Ind\big((j,i)\dd\Fr_XV(F)\cap Z,X\big) &&\hskip-0.6cm \leq k-1.
\end{eqnarray*}
Hence, by inductive assumption,
\begin{eqnarray*}
    (i,j)\dd\Ind\big((j,i)\dd\Fr_YV'(F)\cap Z,Y\big) &&\hskip-0.6cm \leq \\
    \leq (i,j)\dd\Ind\big((j,i)\dd\Fr_YV'(F)\cap Z,X\big) &&\hskip-0.6cm \leq k-1
\end{eqnarray*}
and so, $(i,j)\dd\Ind(Z,Y)\!\leq\!k$. Thus,
$(i,j)\dd\Ind(Z,Y)\!\leq\!(i,j)\dd\Ind(Z,X)$.~\end{pf}
\vskip+0.2cm

It is clear, that for the topological case every $\WS$-supernormal
$\TsS$ $(Y,\tau')$ of a $\TS$ $(X,\tau)$ is hereditarily
$\WS$-supernormal in $X$ so that a $\TsS$ $(Y,\tau')$ is
$\WS$-supernormal in $X$ if and only if $(Y,\tau')$ is
hereditarily $\WS$-supernormal in $X$. Hence, take place

\begin{corollary}{8.15}
If $(Z,\tau'')\sbs (Y,\tau')\sbs (X,\tau)$, where $Y$ is
$\WS$-supernormal in $X$ and $Z\in\co\tau'$, then $\Ind(Z,Y)\leq
\Ind(Z,X)\leq \Ind(Y,X)$.
\end{corollary}

\begin{theorem}{8.16}
If $(Z,\tau_1''<\tau_2'')\sbs (Y,\tau_1'<\tau_2')\sbs
(X,\tau_1<\tau_2)$ and $Z\in\co\tau_2'$, then
$$  (1,2)\dd\Ind(Z,Y)\leq (1,2)\dd\Ind(Z,X)\leq (1,2)\dd\Ind(Y,X). $$
\end{theorem}

\begin{pf}
First, let us prove the right inequality by induction under a
nonnegative integer $n=(1,2)\dd\Ind(Y,X)$. Let the inequality be
proved for $n\leq k-1$ and prove it for $n=k$. Let
$F\in\co\tau_2''$ and $U(F)\in\tau_1$. Then $F\in\co\tau_2'$ and
since $(1,2)\dd\Ind(Y,X)=k$, there is $V(F)\in\tau_1$ such that
$\tau_2\cl V(F)\sbs U(F)$ and
$$  (1,2)\dd\Ind\big((2,1)\dd\Fr_XV(F)\cap Y,X\big)\leq k-1.    $$
Clearly,
$$  (2,1)\dd\Fr_XV(F)\cap Z\sbsq (2,1)\dd\Fr_XV(F)\cap Y,   $$
and $(2,1)\dd\Fr_XV(F)\cap Z\!\!\in\!\!\co\tau_2'''$ in the $\BsS$
$((2,1)\dd\Fr_XV(F)\cap Y,\tau_1'''\!<~\!\!\!\tau_2''')$. Hence,
by inductive assumption,
\begin{eqnarray*}
    (1,2)\dd\Ind\big((2,1)\dd\Fr_XV(F)\cap Z,X\big) &&\hskip-0.6cm \leq \\
    \leq (1,2)\dd\Ind\big((2,1)\dd\Fr_XV(F)\cap Y,X\big) &&\hskip-0.6cm \leq k-1.
\end{eqnarray*}
Thus, $(1,2)\dd\Ind(Z,X)\leq k$ and so $(1,2)\dd\Ind(Z,X)\leq
(1,2)\dd\Ind(Y,X)$.

Now, let us prove the left inequality by induction under a
nonnegative integer $n=(1,2)\dd\Ind(Z,X)$. Let the inequality be
proved for $n\leq k-1$ and prove it for $n=k$. Let
$F\in\co\tau_2''$ and $U'(F)\in\tau_1'$ be any sets. Since $F\sbsq
Z\sbsq Y$ and $(1,2)\dd\Ind(Z,X)=k$, for $U(F)\in\tau_1$, where
$U(F)\cap Y=U'(F)$, there is $V(F)\in\tau_1$ such that $\tau_2\cl
V(F)\sbs U(F)$ and
$$  (1,2)\dd\Ind\big((2,1)\dd\Fr_XV(F)\cap Z,X\big)\leq k-1.  $$
Let $V'(F)=V(F)\cap Y$. Then $\tau_2'\cl V'(F)\sbs U'(F)$. Since
$$  (2,1)\dd\Fr_YV'(F)\sbsq (2,1)\dd\Fr_XV(F)\cap Y     $$
and $Z\sbs Y$, we have
$$  (2,1)\dd\Fr_YV'(F)\cap Z\sbsq (2,1)\dd\Fr_XV(F)\cap Z.   $$
Since $(2,1)\dd\Fr_YV'(F)\cap Z$ is $2$-closed in
$(2,1)\dd\Fr_XV(F)\cap Z$, by the right inequality we have:
\begin{eqnarray*}
    (1,2)\dd\Ind\big((2,1)\dd\Fr_YV'(F)\cap Z,X\big) &&\hskip-0.6cm \leq \\
    \leq (1,2)\dd\Ind\big((2,1)\dd\Fr_XV(F)\cap Z,X\big) &&\hskip-0.6cm \leq k-1.
\end{eqnarray*}
Therefore, by inductive hypothesis,
\begin{eqnarray*}
    (1,2)\dd\Ind\big((2,1)\dd\Fr_YV'(F)\cap Z,Y\big) &&\hskip-0.6cm \leq \\
    \leq (1,2)\dd\Ind\big((2,1)\dd\Fr_YV'(F)\cap Z,X\big) &&\hskip-0.6cm \leq k-1.
\end{eqnarray*}
Thus $(1,2)\dd\Ind(Z,Y)\leq k$ and so
$(1,2)\dd\Ind(Z,Y)\leq(1,2)\dd\Ind(Z,X)$.~\end{pf}

\begin{definition}{8.17}
Let $(Y,\tau_1',\tau_2')$ be a $\BsS$ of a $\BS$
$(X,\tau_1,\tau_2)$. Then a pair $(A,B)$, where $A\in\co\tau_j'$,
$B\in\co\tau_i$ and $A\cap B=\vnth$, is said to be a relatively
$p$-normal pair. A relative partition, corresponding to the pair
$(A,B)$, is a $p$-closed set $T\sbs X$ such that $X\setminus T$ is
not $p$-connected, $(X\setminus T)\cap Y=Y\setminus T\sbs H_1\cup
H_2$, $A\sbs H_i'=H_i\cap Y\in\tau_i'$, $B\sbs H_j\in\tau_j$ and
$H_1\cap H_2=\vnth$.
\end{definition}

\begin{remark}{8.18}
Let $(Y,\tau_1',\tau_2')$ be a $\BsS$ of a $\BS$
$(X,\tau_1,\tau_2)$ and $(A,B)$ be a relatively $p$-normal pair.
Then
\begin{enumerate}
\item[(1)] If there is a neighborhood $U(A)\in\tau_i$
$(U(B)\in\tau_j)$ such that $\tau_j\cl U(A)\sbs X\setminus B$
$(\tau_i\cl U(B)\sbs X\setminus A)$, then $(j,i)\dd\Fr_XU(A)$
$((i,j)\dd\Fr_XU(B))$ is a relative partition, corresponding to
$(A,B)$.

\item[(2)] If $T$ is a relative partition, corresponding to
$(A,B)$, then       \linebreak       $(j,i)\dd\Fr_XH_i\cap Y\sbs
T\cap Y$ $((i,j)\dd\Fr_XH_j\cap Y\sbs T\cap Y)$.
\end{enumerate}
\end{remark}

The proof is similar to the proof of Remark~7.16.

\begin{theorem}{8.19}
If $(Y,\tau_1',\tau_2')$ is a $\BsS$ of a hereditarily $p$-normal
$\BS$ $(X,\tau_1,\tau_2)$, then $(i,j)\dd\Ind(Y,X)\leq n$ if and
only if to any relatively $p$-normal pair $(A,B)$,
$A\in\co\tau_j'$, $B\in\co\tau_i$ and $A\cap B=\vnth$, there
corresponds a relative partition $T$ such that $(i,j)\dd\Ind(T\cap
Y,X)\leq n-1$.
\end{theorem}

\begin{pf}
Let, first, $(i,j)\dd\Ind(Y,X)\leq n$ and $(A,B)$ be any
relatively          \linebreak         $p$-normal pair. Then
$A\!\sbs\!U(A)\!=\!X\!\setminus\!B\!\in\!\tau_i$ and since
$(i,j)\dd\Ind(Y,X)\!\leq~\!\!\!n$, by (2) of Definition~8.1, there
is $V(A)\in\tau_i$ such that
$$  \tau_j\cl V(A)\sbs U(A) \;\;\text{and}\;\;
        (i,j)\dd\Ind\big((j,i)\dd\Fr_XV(A)\cap Y,X\big)\leq n-1. $$
But, by (1) of Remark~8.18, $(j,i)\dd\Fr_XV(A)$ is a relative
partition, corresponding to $(A,B)$, and thus, the first part is
proved.

Conversely, let us suppose that the condition is satisfied,
$A\in\co\tau_j'$ be any set and $U(A)\in\tau_i$ be any
neighborhood. Then $(A,B=X\setminus U(A))$ is a relatively
$p$-normal pair and by condition, there is a relative partition
$T$, corresponding to $(A,B)$, such that $(i,j)\dd\Ind(T\cap
Y,X)\leq n-1$. Since by Definition~8.17, $Y\setminus T\sbs H_1\cup
H_2$, where $A\sbs H_i\in\tau_i$, $B\sbs H_j\in\tau_j$ and
$H_1\cap H_2=\vnth$, we have $\tau_j\cl H_i\cap H_j=\vnth$ and so
$$  \tau_j\cl H_i\sbs X\setminus H_j\sbs X\setminus B=U(A). $$
If $H_i=V(A)\in\tau_i$, then $\tau_j\cl V(A)\sbs U(A)$. Moreover,
by (2) of Remark~8.18, $(j,i)\dd\Fr_XV(A)\cap Y\sbs T\cap Y$ and
since $(j,i)\dd\Fr_XV(A)\in p\,\dd\,\Cl(X)$, we have
$(j,i)\dd\Fr_XV(A)\cap Y\in p\,\dd\,\Cl(Y)$. By the second part of
Remark~2.12, $(j,i)\dd\Fr_XV(A)\cap Y\sbs T\cap Y$ implies
$$  (j,i)\dd\Fr_XV(A)\cap Y\in p\,\dd\,\Cl(T\cap Y)   $$
and by Theorem~8.9,
$$  (i,j)\dd\Ind\big((j,i)\dd\Fr_XV(A)\cap Y,X\big)\leq
        (i,j)\dd\Ind(T\cap Y,X)\leq n-1     $$
so that $(i,j)\dd\Ind(Y,X)\leq n$.
\end{pf}

\begin{definition}{8.20}
Let $(Y,\tau')$ be a $\TsS$ of a $\TS$ $(X,\tau)$. Then a pair
$(A,B)$, where $A\in\co\tau'$, $B\in\co\tau$ and $A\cap B=\vnth$,
is said to be a relatively normal pair. A relative partition,
corresponding to $(A,B)$, is a set $T\in\co\tau$ such that
$X\setminus T$ is not connected, $(X\setminus T)\cap Y=Y\setminus
T\sbs H_1\cup H_2$, $A\sbs H_1'=H_1\cap Y\in\tau'$, $B\sbs
H_2\in\tau$ and $H_1\cap H_2=\vnth$.
\end{definition}

\begin{remark}{8.21}
Let $(Y,\tau')$ be a $\TsS$ of a $\TS$ $(X,\tau)$ and $(A,B)$ be a
relatively normal pair. Then
\begin{enumerate}
\item[(1)] If there is a neighborhood $U(A)\in\tau$
$(U(B)\in\tau)$ such that $\tau\cl U(A)\!\sbs\!X\setminus B$
$(\tau\cl U(B)\!\sbs\!X\setminus A)$, then $\Fr_XU(A)$
$(\Fr_XU(B))$ is a relative partition, corresponding to $(A,B)$.

\item[(2)] If $T$ is a partition, corresponding to $(A,B)$, then
$\Fr H_i\cap Y\sbs T\cap Y$.
\end{enumerate}
\end{remark}

\begin{proposition}{8.22}
If $(Y,\tau')$ is a $\TsS$ of a $\TS$ $(X,\tau)$, where $X$ is
hereditarily normal, then $\Ind(Y,X)\leq n$ if and only if to
every relatively normal pair $(A,B)$ there corresponds a relative
partition $T$ such that $\Ind(T\cap Y,X)\leq n-1$.
\end{proposition}

\begin{proposition}{8.23}
If $(Y,\tau_1'<\tau_2')$ is a $\BsS$ of a $\BS$
$(X,\tau_1<\tau_2)$, then $(1,2)\dd\Ind(Y,X)\leq n$ if and only if
to every relatively $p$-normal pair $(A,B)$, $A\in\co\tau_2'$,
$B\in\co\tau_1$ and $A\cap B=\vnth$, there corresponds a relative
partition $T$ such that $(1,2)\dd\Ind(T\cap Y,X)\leq n-1$.
\end{proposition}

\begin{pf}
The proof of the first part is given by the proof of the first
part of Theorem~8.19 for $i=1$, $j=2$.

Conversely, if the condition is satisfied, then by analogy with
the proof of the second part of Theorem~8.19, for $A\in\co\tau_2'$
and $U(A)\in\tau_1$ there is $V(A)\in\tau_1$ such that $\tau_2\cl
V(A)\sbs U(A)$ and $(2,1)\dd\Fr_XV(A)\cap Y\sbs T\cap Y$. Since
$(2,1)\dd\Fr_XV(A)\cap Y\in\co\tau_2'$, we have
$(2,1)\dd\Fr_XV(A)\cap Y\in\co\tau_2''$ in $(T\cap
Y,\tau_1'',\tau_2'')$. Hence, by Proposition~8.16,
$$  (1,2)\dd\Ind\big((2,1)\dd\Fr_XV(A)\cap Y,X\big)\leq
        (1,2)\dd\Ind(T\cap Y,X)\leq n-1    $$
and thus $(1,2)\dd\Ind(Y,X)\leq n$.
\end{pf}

\begin{theorem}{8.24}
If $(Y,\tau_1',\tau_2')$ is a hereditarily
$(j,i)\dd\WS$-supernormal $\BsS$ of a $\BS$ $(X,\tau_1,\tau_2)$,
then $(i,j)\dd\Ind(Y,X)\leq n$ if and only if to any relatively
$p$-normal pair $(A,B)$, $A\in\co\tau_j'$, $B\in\co\tau_i$ and
$A\cap B=\vnth$, there corresponds a relative partition $T\sbs X$
such that $(i,j)\dd\Ind(T\cap Y,X)\leq n-1$.
\end{theorem}

\begin{pf}
The proof of the first part is similar to the proof of the first
part of Theorem~8.19.

For the proof of the second part, let us recall that by the proof
of the second part of Theorem~8.19, for $A\in\co\tau_j'$ and
$U(A)\in\tau_i$ there is $V(A)\in\tau_i$ such that $\tau_j\cl
V(A)\sbs U(A)$ and $(j,i)\dd\Fr_XV(A)\cap Y\sbs T\cap Y$.

Clearly, $(j,i)\dd\Fr_XV(A)\cap Y\in p\,\dd\,\Cl(Y)$ and by the
second part of Remark~2.12,
$$  (j,i)\dd\Fr_XV(A)\cap Y\in p\,\dd\,\Cl(T\cap Y).  $$
Hence, by Definition~8.13, $(j,i)\dd\Fr_XV(A)\cap Y$ is
$(j,i)\dd\WS$-supernormal in $X$ and Theorem~8.14 gives
$$  (i,j)\dd\Ind\big((j,i)\dd\Fr_XV(A)\cap Y,X\big)\leq
        (i,j)\dd\Ind(T\cap Y,X)\leq n-1.     $$
Thus $(i,j)\dd\Ind(Y,X)\leq n$.
\end{pf}

\begin{corollary}{8.25}
If $(Y,\tau')$ is a $\WS$-supernormal $\TsS$ of a $\TS$
$(X,\tau)$, then $\Ind(Y,X)\leq n$ if and only if to any
relatively normal pair $(A,B)$ there corresponds a relative
partition $T\sbs X$ such that $\Ind(T\cap Y,X)\leq n-1$.
\end{corollary}

\begin{pf}
It suffices to recall that by the remark before Corollary~8.15,
every $\WS$-supernormal $\TsS$ $(Y,\tau')$ of a $\TS$ $(X,\tau)$
is hereditarily $\WS$-supernormal in $X$.
\end{pf}
\vskip+0.2cm

At the end of the section let us recall once more that
Theorem~8.19 says that for an arbitrary $\BsS$
$(Y,\tau_1',\tau_2')$ of a hereditarily $p$-normal $\BS$
$(X,\tau_1,\tau_2)$ its relative inductive dimension functions
$(i,j)\dd\Ind(Y,X)$ can be characterized in two equivalent ways:
as by neighborhoods so by partitions. The similar result is given
by Theorem~8.24 but for hereditarily $(j,i)\dd\WS$-supernormal
$\BsS$ $(Y,\tau_1',\tau_2')$ in a $\BS$ $(X,\tau_1,\tau_2)$. In
both cases the characterizations are conditioned by the
monotonicity property of the relative inductive dimension
functions $(i,j)\dd\Ind(Y,X)$. i.e., by the inequality
$(i,j)\dd\Ind(Z,X)\leq (i,j)\dd\Ind(Y,X)$ for $Z\in
p\,\dd\,\Cl(Y)$ (Theorems~8.9 and 8.14).

Moreover, for the topological case the requirement of hereditary
$\WS$-supernormality of a $\TsS$ $(Y,\tau')$ in a $\TS$ $(X,\tau)$
can be weakened to the requirement of $\WS$-supernormality of $Y$
in $X$ (Corollaries~8.15 and 8.25).

\vskip+0.5cm
\section*{\textbf{9. $(i,j)$-Small Separate Inductive Dimension Functions
and $(i,j)$-Small Separate Relative Inductive Dimension
Functions}} \vskip+0.2cm

\begin{definition}{9.1}
Let $(X,\tau_1,\tau_2)$ be a $\BS$ and $n$ denote a nonnegative
integer. Then
\begin{enumerate}
\item[(1)] $(i,j)\dd\,\sind X=-1$ if and only if $X=\vnth$.

\item[(2)] $(i,j)\dd\,\sind X\leq n$ if for each point $x\in X$
and any neighborhood $U(x)\in\tau_i$ there is a neighborhood
$V(x)\in\tau_j$ such that $\tau_j\cl V(x)\sbs U(x)$ and
$(i,j)\dd\,\sind(j\dd\Fr_XV(x))\leq n-1$.

\item[(3)] $(i,j)\dd\,\sind X=n$ if $(i,j)\dd\,\sind X\leq n$ and
$(i,j)\dd\,\sind X\leq n-1$ does not hold.

\item[(4)] $(i,j)\dd\,\sind X=\infty$ if $(i,j)\dd\,\sind X\leq n$
does not hold for any $n$.
\end{enumerate}
\end{definition}

As usual
$$  p\,\dd\,\sind X\leq n\llra
        \big((1,2)\dd\,\sind X\leq n\wedge (2,1)\dd\,\sind X\leq n\big). $$

It follows immediately from Proposition~2.33 that if
$(i,j)\dd\,\sind X$ is finite, then $(X,\tau_1,\tau_2)$ is
$(i,j)$-separately regular.

\begin{proposition}{9.2}
If $(Y,\tau_1',\tau_2')$ is an arbitrary $\BsS$ of a $\BS$
$(X,\tau_1,\tau_2)$, then $(i,j)\dd\,\sind Y\leq (i,j)\dd\,\sind
X$.
\end{proposition}

\begin{pf}
We shall use the induction under a nonnegative integer
$n=(i,j)\dd\,\sind X$. Let the inequality be proved for $n\leq
k-1$ and prove it for $n=k$. If $x\in Y$ and $U'(x)\in\tau_i'$,
then there is $U(x)\in\tau_i$ such that $U(x)\cap Y=U'(x)$. Since
$(i,j)\dd\,\sind X=k$, there is $V(x)\in\tau_j$ such that
$\tau_j\cl V(x)\sbs U(x)$ and $(i,j)\dd\,\sind(j\dd\Fr_XV(x))\leq
k-1$. Now, for $V'(x)=V(x)\cap Y$,
$$  \tau_j'\cl V'(x)=\tau_j\cl V'(x)\cap Y\sbs U(x)\cap Y=U'(x) $$
and $j\dd\Fr_YV'(x)\sbs j\dd\Fr_XV(x)$. Then, by inductive
assumption,
$$  (i,j)\dd\,\sind(j\dd\Fr_YV'(x))\leq
        (i,j)\dd\,\sind(j\dd\Fr_XV(x))\leq k-1,    $$
and hence, $(i,j)\dd\,\sind Y\leq k$ so that $(i,j)\dd\,\sind
Y\leq(i,j)\dd\,\sind X$.
\end{pf}

\begin{theorem}{9.3}
The following inequalities hold for a $\BS$ $(X,\tau_1<\tau_2)$:
\begin{eqnarray*}
    &\ds (1,2)\dd\,\sind X\leq i\dd\ind X\leq (2,1)\dd\,\sind X, \\
    &\ds (1,2)\dd\ind X\leq (2,1)\dd\,\sind X \;\;\text{and}\;\;
        (1,2)\dd\,\sind X\leq (2,1)\dd\ind X.
\end{eqnarray*}
\end{theorem}

\begin{pf}
For $i=1$ first let us prove the right upper inequality by
induction under a nonnegative integer $n=(2,1)\dd\,\sind X$. Let
the inequality be proved for $n\leq k-1$ and prove it for $n=k$.
If $x\in X$ is any point and $U(x)\in\tau_1$ is any neighborhood,
then $U(x)\in\tau_2$ and since $(2,1)\dd\,\sind X=k$, there is
$V(x)\in\tau_1$ such that
$$  \tau_1\cl V(x)\sbs U(x) \;\;\text{and}\;\;
        (2,1)\dd\,\sind\big(1\dd\Fr_XV(x)\big)\leq k-1. $$
Hence, by inductive hypothesis, $1\dd\ind(1\dd\Fr_XV(x))\leq k-1$
and so $1\dd\ind X\leq k$. Thus $1\dd\ind X\leq (2,1)\dd\,\sind
X$.

Let the left upper inequality be proved for a nonnegative integer
$n\leq k-1$, where $n=1\dd\ind X$, and prove it for $n=k$. If
$x\in X$ is any point and $U(x)\in\tau_1$ is any neighborhood,
then by $k=1\dd\ind X$ there is $V(x)\in\tau_1$ such that
$\tau_1\cl V(x)\sbs U(x)$ and $1\dd\ind(1\dd\Fr_XV(x))\leq k-1$.
Now, by monotonicity of the small inductive dimension
      \linebreak       $1\dd\ind(2\dd\Fr_XV(x))\leq k-1$ as $2\dd\Fr_XV(x)\sbs
1\dd\Fr_XV(x)$. Hence, by inductive assumption
$(1,2)\dd\,\sind(2\dd\Fr_XV(x))\leq k-1$. Clear\-ly,
$V(x)\in\tau_2$ and $\tau_2\cl V(x)\sbs U(x)$ so that
$(1,2)\dd\,\sind X\leq k$. Thus $(1,2)\dd\,\sind X\leq 1\dd\,\sind
X$.

For $i=2$ first let us prove the left upper inequality by
induction under a nonnegative integer $n=2\dd\ind X$. Let the
inequality be proved for $n\leq k-1$ and prove it for $n=k$. Let
$x\in X$ be any point and $U(x)\in\tau_1$ be any neighborhood.
Since $U(x)\in\tau_1\sbs\tau_2$ and $2\dd\ind X=k$, there is
$V(x)\in\tau_2$ such that $\tau_2\cl V(x)\sbs U(x)$ and
$2\dd\ind(2\dd\Fr_XV(x))\leq k-1$. Then, by inductive assumption,
$(1,2)\dd\,\sind(2\dd\Fr_XV(x))\leq k-1$ and so $(1,2)\dd\,\sind
X\leq k$. Hence $(1,2)\dd\,\sind X\leq 2\dd\ind X$.

Now, let $(2,1)\dd\,\sind X=n$ and the right upper inequality be
proved for $n\leq k-1$. Let us prove this inequality for $n=k$. If
$x\in X$ is any point and $U(x)\in\tau_2$ is any neighborhood,
then $(2,1)\dd\,\sind X=k$ implies that there is
$V(x)\in\tau_1\sbs\tau_2$ such that
$$  \tau_1\cl V(x)\sbs U(x) \;\;\text{and}\;\;
        (2,1)\dd\,\sind\big(1\dd\Fr_XV(x)\big)\leq k-1. $$
Hence, by inductive hypothesis, $2\dd\ind(1\dd\Fr_XV(x))\leq k-1$.
On the other hand, $\tau_2\cl V(x)\sbs U(x)$ and so
$2\dd\Fr_XV(x)\sbs 1\dd\Fr_XV(x)$. Now, by monotonicity,
$2\dd\ind(2\dd\Fr_XV(x))\leq k-1$. Thus $2\dd\ind X\leq k$ and so
$2\dd\ind X\leq (2,1)\dd\,\sind X$.

Now, let us consider the lower inequalities. Let $(2,1)\dd\sind
X=n$ and let us prove that $(1,2)\dd\ind X\leq (2,1)\dd\sind X$.
Let the inequality be proved for $n\leq k-1$ and prove it for
$n=k$. If $x\in X$ is any point and $U(x)\in\tau_1$ is any
neighborhood, then $U(x)\in\tau_2$ and by $(2,1)\dd\,\sind X=k$
there is $V(x)\in\tau_1$ such that $\tau_1\cl V(x)\sbs U(x)$ and
$(2,1)\dd\,\sind(1\dd\Fr_XV(x))\leq k-1$. Hence, by inductive
assumption, $(1,2)\dd\ind(1\dd\Fr_XV(x))\leq k-1$. Evidently,
$\tau_2\cl V(x)\sbs U(x)$, $(2,1)\dd\Fr_XV(x)\sbs 1\dd\Fr_XV(x)$
and by monotonicity of the small bitopological dimension,
$(1,2)\dd\ind\big((2,1)\dd\Fr_XV(x)\big)\leq k-1$. Thus,
$(1,2)\dd\ind X\leq k$ and so $(1,2)\dd\ind X\leq (2,1)\dd\,\sind
X$.

Finally, let us prove the inequality $(1,2)\dd\,\sind X\leq
(2,1)\dd\ind X$ by induction under $n=(2,1)\dd\ind X$. Let the
inequality be proved for $n\leq k-1$ and prove it for $n=k$. If
$x\in X$ is any point and $U(x)\in\tau_1$ is any neighborhood,
then $U(x)\in\tau_2$ and by $(2,1)\dd\ind X=k$ there is
$V(x)\in\tau_2$ such that
$$  \tau_1\cl V(x)\sbs U(x) \;\;\text{and}\;\;
        (2,1)\dd\ind\big((1,2)\dd\Fr_XV(x)\big)\leq k-1. $$
Hence, by inductive assumption,
$(1,2)\dd\,\sind((1,2)\dd\Fr_XV(x))\!\leq~\!\!\!k\!-~\!\!\!1$.
Cle\-ar\-ly, $\tau_2\cl V(x)\sbs U(x)$, $2\dd\Fr_XV(x)\sbs
(1,2)\dd\Fr_XV(x)$ and by Proposition~9.2,
$(1,2)\dd\,\sind(2\dd\Fr_XV(x))\leq k-1$. Hence, $(1,2)\dd\,\sind
X\leq k$ and so $(1,2)\dd\,\sind X\leq (2,1)\dd\ind X$.
\end{pf}

\begin{corollary}{9.4}
For a $\BS$ $(X,\tau_1,\tau_2)$ we have:
\begin{enumerate}
\item[(1)] If $\tau_1<_C\tau_2$, then $(1,2)\dd\,\sind X\leq
(1,2)\dd\ind X\leq (2,1)\dd\,\sind X$.

\item[(2)] If $\tau_1<_N\tau_2$, then $(1,2)\dd\,\sind X\leq
(i,j)\dd\ind X\leq (2,1)\dd\,\sind X$.
\end{enumerate}
\end{corollary}

\begin{pf}
Follows directly from Theorem~9.3 taking into account (5) and (6)
of Theorem~3.1.36 in [8].
\end{pf}

\begin{proposition}{9.5}
If a $\BS$ $(X,\tau_1,\tau_2)$ is $(i,j)$-separately regular and
$j$-extremally disconnected, then $(i,j)\dd\,\sind X=0$.
\end{proposition}

\begin{pf}
Let $x\in X$ be any point and $U(x)\in\tau_i$ be any neighborhood.
Then, by proposition~2.33, there is a neighborhood $V(x)\in\tau_j$
such that $\tau_j\cl V(x)\sbs U(x)$. Since $X$ is $j$-extremally
disconnected, the set $W(x)=\tau_j\cl V(x)\in\tau_j\cap\co\tau_j$
and hence, $j\dd\Fr_XW(x)\!=\!\vnth$. Thus, by (2) of
Definition~9.1, $(i,j)\dd\,\sind X=0$.
\end{pf}

\begin{definition}{9.6}
Let $(x,A)$ be an $i$-regular pair in a $\BS$ $(X,\tau_1,\tau_2)$,
that is, $x\in X$, $A\in\co\tau_i$ and $x\ol{\in}\,A$. Then we say
that a $j$-closed set $T\sbs X$ is a separate partition,
corresponding to the pair $(x,A)$, if $X\setminus T$ is not
$j$-connected so that $X\setminus T=H_1\cup H_2$, where
$H_1,\,H_2\in\tau_j\setminus\{\vnth\}$, $x\in H_1$, $A\sbs H_2$
and $H_1\cap H_2=\vnth$.
\end{definition}

\begin{remark}{9.7}
One can easily to verify  that in a $\BS$ $(X,\tau_1,\tau_2)$ the
following conditions are satisfied for an $i$-regular pair
$(x,A)$:
\begin{enumerate}
\item[(1)] If there exists a $j$-open neighborhood $U(x)$
$(j$-open neighborhood $U(A))$ such that $\tau_j\cl U(x)\sbs
X\setminus A$ $(\tau_j\cl U(A)\sbs X\setminus\{x\})$, then the set
$j\dd\Fr_XU(x)$ $(j\dd\Fr_XU(A))$ is a separate partition,
corresponding to $(x,A)$.

\item[(2)] If $T$ is a separate partition, corresponding to
$(x,A)$, then     \linebreak        $j\dd\Fr_XH_i\sbs T$.
\end{enumerate}
\end{remark}

Indeed, (1) $X\setminus j\dd\Fr_XU(x)=
        \big(X\setminus \tau_j\cl U(x)\big)\cup U(x)$,
where $x\in U(x)=H_1\in\tau_j$, $A\sbs X\setminus\tau_j\cl
U(x)=H_2\in\tau_j$ and $H_1\cap H_2=\vnth$. Similarly for $U(A)$.

(2) Since $(X\setminus T,\tau_1',\tau_2')$ is not $j$-connected,
we have $H_i\in\tau_j'\cap\co\tau_j'$ and hence
\begin{eqnarray*}
    &\ds j\dd\Fr_XH_i\cap (X\setminus T)= \\
    &\ds =\big(\tau_j\cl H_i\cap (X\setminus T)\big)\cap
        \big((X\setminus H_i)\cap (X\setminus T)\big)= \\
    &\ds =\tau_j'\cl H_i\cap H_j=H_i\cap H_j=\vnth
\end{eqnarray*}
so that $j\dd\Fr_XH_i\sbs T$.

\begin{proposition}{9.8}
Let $(X,\tau_1,\tau_2)$ be a $\BS$ and $n$ denote a nonnegative
integer. Then $(i,j)\dd\,\sind X\leq n$ if and only if to every
$i$-regular pair $(x,A)$ there corresponds a separate partition
$T$ such that $(i,j)\dd\,\sind T\leq n-1$.
\end{proposition}

\begin{pf}
First, let $(i,j)\dd\,\sind X\leq n$ and $(x,A)$ be an $i$-regular
pair. Then $U(x)=X\setminus A\in\tau_i$ is an $i$-open
neighborhood of $x$ and since $(i,j)\dd\,\sind X\leq n$, there is
a $j$-open neighborhood $V(x)$ such that
$$  \tau_j\cl V(x)\sbs U(x) \;\;\text{and}\;\;
        (i,j)\dd\,\sind(j\dd\Fr_XV(x)\leq n-1.  $$
It follows from (1) of Remark~9.7 that  $j\dd\Fr_XV(x)$ is a
separate partition, corresponding to $(x,A)$.

Conversely, let the condition be satisfied, $x\in X$ and
$U(x)\in\tau_i$. Then for an $i$-regular pair $(x,A=X\setminus
U(x))$ there is a separate partition $T$ such that $X\setminus
T=H_1\cup H_2$, $H_1,\,H_2\in\tau_j\setminus\{\vnth\}$, $x\in
H_1$, $A\sbs H_2$, $H_1\cap H_2=\vnth$ and $(i,j)\dd\,\sind T\leq
n-1$. It is evident that $\tau_j\cl H_1\cap H_2=\vnth$ so that
$\tau_j\cl H_1\sbs X\setminus A=U(x)$. Let $H_1=V(x)$. Then
$\tau_j\cl V(x)\sbs U(x)$ and by (2) Remark~9.7 together with
Proposition~9.2, $(i,j)\dd\,\sind(j\dd\Fr_XV(x))\leq n-1$.
\end{pf}

\begin{definition}{9.9}
Let $(Y,\tau_1',\tau_2')$ be a $\BsS$ of a $\BS$
$(X,\tau_1,\tau_2)$ and $n$ denote a nonnegative integer. Then
\begin{enumerate}
\item[(1)] $(i,j)\dd\,\sind(Y,X)=-1$ if and only if $X=\vnth$.

\item[(2)] $(i,j)\dd\,\sind(Y,X)\leq n$ if for each point $x\in Y$
and any neighborhood $U(x)\in\tau_i$ there is a neighborhood
$V(x)\in\tau_j$ such that $\tau_j\cl V(x)\sbs U(x)$ and
$(i,j)\dd\,\sind(j\dd\Fr_XV(x)\cap Y,X)\leq n-1(\llra$ for each
point $x\in Y$ and any neighborhood $U(x)\in\tau_i$ there is a
neighborhood $V(x)\in\tau_j$ such that $\tau_j\cl V(x)\sbs U(x)$
and $(i,j)\dd\,\sind(j\dd\Fr_XV(x)\cap Y)\leq n-1)$.

\item[(3)] $(i,j)\dd\,\sind(Y,X)=n$ if $(i,j)\dd\,\sind(Y,X)\leq
n$ and             \linebreak       $(i,j)\dd\,\sind(Y,X)\leq n-1$
does not hold.

\item[(4)] $(i,j)\dd\,\sind(Y,X)=\infty$ if
$(i,j)\dd\,\sind(Y,X)\leq n$ does not hold for any $n$.
\end{enumerate}
\end{definition}

As usual,
$$  p\,\dd\,\sind(Y,X)\!\leq\!n\!\llra\!
        \big((1,2)\dd\,\sind(Y,X)\!\leq\!n\wedge(2,1)\dd\,\sind(Y,X)\!\leq\!n\big). $$

It follows immediately from (2) of Definition~9.9 and (3) of
Proposition~2.36 that  if $(i,j)\dd\,\sind (Y,X)$ is finite, then
$Y$ is $(i,j)$-separately superregular in $X$. Hence, by
implications after Definition~2.35, $Y$ is $(i,j)$-strongly
separately regular in $X$, $Y$ is $(i,j)\dd\WS$-separately
superregular in $X$, $Y$ is $(i,j)$-separately regular in $X$, $Y$
is $(i,j)\dd\WS$-separately regular in $X$ and $Y$ is
$(i,j)\dd\WS$-quasi separately regular in~$X$.

\begin{theorem}{9.10}
If $(Y,\tau_1',\tau_2')$ is a $\BsS$ of a $\BS$
$(X,\tau_1,\tau_2)$, then the following conditions are satisfied:
\begin{enumerate}
\item[(1)] $(i,j)\dd\,\sind Y\leq
(i,j)\dd\,\sind(Y,X)\leq(i,j)\dd\,\sind X$.

\item[(2)] $(i,j)\dd\,\sind(X,X)=(i,j)\dd\,\sind X$

\item[(3)] If $(Z,\tau_1'',\tau_2'')\sbsq (Y,\tau_1',\tau_2')$,
then
$$  (i,j)\dd\,\sind(Z,Y)\leq (i,j)\dd\,\sind(Z,X)\leq
        (i,j)\dd\,\sind(Y,X).     $$
\end{enumerate}
\end{theorem}

\begin{pf}
(1) First, let us prove the left inequality by induction under a
nonnegative integer $n=(i,j)\dd\,\sind(Y,X)$. Let the inequality
be proved for $n\leq k-1$ and prove it for $n=k$. Let $x\in Y$ be
any point and $U'(x)\in\tau_i'$ be any neighborhood. If
$U(x)\in\tau_i$, $U(x)\cap Y=U'(x)$, then $(i,j)\dd\,\sind(Y,X)=k$
implies that there is $V(x)\in\tau_j$ such that
$$  \tau_j\cl V(x)\sbs U(x) \;\;\text{and}\;\;
        (i,j)\dd\,\sind\big(j\dd\Fr_XV(x)\cap Y,X\big)\leq k-1. $$
Therefore, by inductive assumption,
$(i,j)\dd\,\sind(j\dd\Fr_XV(x)\cap Y)\leq k-1$. Moreover,
$\tau_j'\cl V'(x)\sbs U'(x)$, where $V'(x)=V(x)\cap Y$ and since
$j\dd\Fr_YV'(x)\sbsq j\dd\Fr_XV(x)\cap Y$, by Proposition~9.2,
$$  (i,j)\dd\,\sind(j\dd\Fr_YV'(x))\leq k-1.        $$
Hence $(i,j)\dd\,\sind Y\leq k$ and so $(i,j)\dd\,\sind Y\leq
(i,j)\dd\,\sind(Y,X)$.

Now, let us prove the right inequality by induction under a
nonnegative integer $n=(i,j)\dd\,\sind X$. Let the inequality be
proved for $n\leq k-1$ and prove it for $n=k$. Let $x\in Y$ be any
point and $U(x)\in\tau_i$ be any neighborhood. Since
$(i,j)\dd\,\sind X=k$, there is $V(x)\in\tau_j$ such that
$\tau_j\cl V(x)\sbs U(x)$ and
$(i,j)\dd\,\sind\big(j\dd\Fr_XV(x)\big)\leq k-1$. Hence, by
inductive assumption, $(i,j)\dd\,\sind(j\dd\Fr_XV(x)\cap
Y,j\dd\Fr_XV(x))\leq k-1$ and by the left inequality,
$(i,j)\dd\,\sind(j\dd\Fr_XV(x)\cap Y)\leq k-1$. Hence, by
condition (2) in brackets of Definition~9.9,
$(i,j)\dd\,\sind(Y,X)\leq k$ and so $(i,j)\dd\,\sind(Y,X)\leq
(i,j)\dd\,\sind X$.

(2) Follows from (1) for $Y=X$.

(3) First, let us prove the right inequality, i.e., \textbf{the
monotonicity property of $\boldsymbol{(i,j)\dd\,\sind(Y,X)}$} by
induction under a nonnegative integer $n=(i,j)\dd\,\sind(Y,X)$.
Let the inequality be proved for $n\leq k-1$ and prove it for
$n=k$. Let $x\in Z$ be any point and $U(x)\in\tau_i$ be any
neighborhood. Since $x\in Z\sbs Y$ and $(i,j)\dd\,\sind(Y,X)=k$
there is $V(x)\in\tau_j$ such that
$$  \tau_j\cl V(x)\sbs U(x) \;\;\text{and}\;\;
        (i,j)\dd\,\sind\big(j\dd\Fr_XV(x)\cap Y,X\big)\leq k-1. $$
But $Z\sbs Y$, $j\dd\Fr_XV(x)\cap Z\sbsq j\dd\Fr_XV(x)\cap Y$ and
by inductive assumption
\begin{eqnarray*}
    &\ds (i,j)\dd\,\sind\big(j\dd\Fr_XV(x)\cap Z,X\big)\leq \\
    &\ds \leq (i,j)\dd\,\sind\big(j\dd\Fr_XV(x)\cap Y,X\big)\leq k-1.
\end{eqnarray*}
Thus $(i,j)\dd\,\sind(Z,X)\leq k$ and so
$(i,j)\dd\,\sind(Z,X)\leq(i,j)\dd\,\sind(Y,X)$.

Finally, let us prove the left inequality by induction under a
nonnegative integer $n=(i,j)\dd\,\sind(Z,X)$. Let the inequality
be proved for $n\leq k-1$ and prove it for $n=k$. Let $x\in Z$ be
any point and $U'(x)\in\tau_i'$ be any neighborhood. Since
$(i,j)\dd\,\sind(Z,X)=k$, for $U(x)\in\tau_i$, $U(x)\cap Y=U'(x)$,
there is $V(x)\in\tau_j$ such that
$$  \tau_j\cl V(x)\sbs U(x) \;\;\text{and}\;\;
        (i,j)\dd\,\sind\big(j\dd\Fr_XV(x)\cap Z,X\big)\leq k-1. $$
If $V'(x)=V(x)\cap Y$, then $\tau_j'\cl V'(x)\sbs U'(x)$. For
$j\dd\Fr_YV'(x)\sbsq j\dd\Fr_XV(x)\cap Y$ and $Z\sbs Y$, we have
$$  j\dd\Fr_YV'(x)\cap Z\sbsq j\dd\Fr_XV(x)\cap Z.      $$
Hence, by the right inequality,
\begin{eqnarray*}
    (i,j)\dd\,\sind\big(j\dd\Fr_YV'(x)\cap Z,X\big) &&\hskip-0.6cm \leq \\
    \leq (i,j)\dd\,\sind\big(j\dd\Fr_XV(x)\cap Z,X\big) &&\hskip-0.6cm \leq k-1.
\end{eqnarray*}
Therefore, by inductive assumption,
\begin{eqnarray*}
    (i,j)\dd\,\sind\big(j\dd\Fr_YV'(x)\cap Z,Y\big) &&\hskip-0.6cm \leq \\
    \leq (i,j)\dd\,\sind\big(j\dd\Fr_YV'(x)\cap Z,X\big) &&\hskip-0.6cm \leq k-1.
\end{eqnarray*}
Thus $(i,j)\dd\,\sind(Z,Y)\!\leq\!k$ and so
$(i,j)\dd\,\sind(Z,Y)\!\leq\!(i,j)\dd\,\sind(Z,X)$.~\end{pf}

\begin{theorem}{9.11}
For an $i$-open $\BsS$ $(Y,\tau_1',\tau_2')$ of a $\BS$
$(X,\tau_1,\tau_2)$ the following conditions are satisfied:
\begin{enumerate}
\item[(1)] If $j\dd\cl Z\sbs Y$, then
$(i,j)\dd\,\sind(Z,X)=(i,j)\dd\,\sind(Z,Y)$.

\item[(2)] $(i,j)\dd\,\sind(Y,X)=(i,j)\dd\,\sind Y$.
\end{enumerate}
\end{theorem}

\begin{pf}
(1) By (3) of Theorem~9.10 it suffices to prove only that
\linebreak         $(i,j)\dd\,\sind(Z,X)\leq
(i,j)\dd\,\sind(Z,Y)$. Let us prove the inequality by induction
under a nonnegative integer $n=(i,j)\dd\,\sind(Z,Y)$. Let the
inequality be proved for $n\leq k-1$ and prove it for $n=k$. Let
$x\in Z$ be any point and $U(x)\in\tau_i$ be any neighborhood.
Since $Y\in\tau_i$, one can assume that $U(x)\sbs Y$ and
$(i,j)\dd\,\sind(Z,Y)=n$ implies that there is $V(x)\in\tau_j$
such that
$$  \tau_j\cl V(x)\sbs U(x) \;\;\text{and}\;\;
        (i,j)\dd\,\sind\big(j\dd\Fr_YV(x)\cap Z,Y\big)\leq k-1. $$
Clearly, $j\dd\Fr_YV(x)=j\dd\Fr_XV(x)$ and since
$$  j\dd\cl\big(j\dd\Fr_YV(x)\cap Z\big)\sbsq j\dd\cl Z\sbsq Y, $$
one can use the inductive assumption under the set
$j\dd\Fr_YV(x)\cap Z$. Hence
\begin{eqnarray*}
    (i,j)\dd\,\sind\big(j\dd\Fr_XV(x)\cap Z,X\big) &&\hskip-0.6cm \leq \\
    \leq (i,j)\dd\,\sind\big(j\dd\Fr_YV(x)\cap Z,Y\big) &&\hskip-0.6cm \leq k-1
\end{eqnarray*}
and so $(i,j)\dd\,\sind(Z,X)\!\leq\!k$. Thus
$(i,j)\dd\,\sind(Z,X)\!\leq\!(i,j)\dd\,\sind(Z,Y)$.

(2) By (1) of Theorem~9.10, it suffices to prove only that
$$  (i,j)\dd\,\sind(Y,X)\leq (i,j)\dd\,\sind Y.     $$
We shall use the induction under a nonnegative integer
$(i,j)\dd\,\sind Y=n$. Let the inequality be proved for $n\leq
k-1$ and prove it for $n=k$. If $x\in Y$ is any point and
$U(x)\in\tau_i$ is any neighborhood, then one can assume that
$U(x)\sbs Y$. Since $(i,j)\dd\,\sind Y=k$ and
$j\dd\Fr_YV(x)=j\dd\Fr_XV(x)\sbs Y$, where $V(x)\in\tau_j$ and
$\tau_j\cl V(x)\sbs U(x)$, we have
$$  (i,j)\dd\,\sind\big(j\dd\Fr_YV(x)\big)=
        (i,j)\dd\,\sind\big(j\dd\Fr_XV(x)\big)\leq k-1.  $$

By (2) of Theorem~9.10, $(i,j)\dd\,\sind Y=(i,j)\dd\,\sind(Y,Y)=k$
and therefore, by (2) of Definition~9.9, we have
$(i,j)\dd\,\sind(j\dd\Fr_YV(x)\cap
Y,Y)=(i,j)\dd\,\sind(j\dd\Fr_YV(x),Y)\leq k-1$. But
$j\dd\cl\big(j\dd\Fr_XV(x)\big)=j\dd\Fr_XV(x)\sbs Y$ and if we
consider in (1) the set $j\dd\Fr_XV(x)$ instead of $Z$, we obtain
\begin{eqnarray*}
    &\ds (i,j)\dd\,\sind\big(j\dd\Fr_XV(x)\cap Y,X\big)=
        (i,j)\dd\,\sind\big(j\dd\Fr_YV(x),X\big)= \\
    &\ds =(i,j)\dd\,\sind\big(j\dd\Fr_YV(x),Y\big)\leq k-1.
\end{eqnarray*}
Hence, $(i,j)\dd\,\sind(Y,X)\leq k$ and so
$(i,j)\dd\,\sind(Y,X)\leq (i,j)\dd\,\sind Y$.~\end{pf}

\begin{proposition}{9.12}
For a $\BsS$ $(Y,\tau_1',\tau_2')$ of a $\BS$ $(X,\tau_1,\tau_2)$
the equ\-a\-li\-ty $(i,j)\dd\,\sind(Y,X)=(i,j)\dd\,\sind Y$ holds
if anyone of the following two conditions is satisfied:
\begin{enumerate}
\item[(1)] $(X,\tau_1,\tau_2)$ is hereditarily $j$-normal.

\item[(2)] $(X,\tau_1,\tau_2)$ is $j$-normal,
$(Y,\tau_1',\tau_2')$ is $j$-perfectly normal $\BsS$ of $X$ and
$Y\in\co\tau_j$.
\end{enumerate}
\end{proposition}

\begin{pf}
(1) Since $(X,\tau_1,\tau_2)$ is hereditarily $j$-normal, by the
proof of the first part of Proposition~1 in [19], if
$U'\in\tau_j'$, then there is $U\in\tau_j$ such that $U\cap Y=U'$
and $j\dd\Fr_XU\cap Y=j\dd\Fr_YU'$.

Now, following (1) of Theorem~9.10 it suffices to prove only that
\linebreak        $(i,j)\dd\,\sind(Y,X)\leq(i,j)\dd\,\sind Y$ be
induction under a nonnegative integer $n=(i,j)\dd\,\sind Y$. Let
the inequality be proved for $n\leq k-1$ and prove it for $n=k$.
Since $(i,j)\dd\,\sind Y=k$, there is $V'(x)\in\tau_j'$ such that
$\tau_j'\cl V'(x)\sbs U(x)\cap Y$ and
$(i,j)\dd\,\sind\big(j\dd\Fr_YV'(x)\big)\leq k-1$, where $x\in Y$
is any point and $U(x)\in\tau_i$ is any neighborhood. If
$V(x)\in\tau_j$, $V(x)\cap Y=V'(x)$ and $j\dd\Fr_XV(x)\cap
Y=j\dd\Fr_YV'(x)$, then
$$  (i,j)\dd\,\sind\big(j\dd\Fr_XV(x)\cap Y\big)\leq k-1. $$
Hence, by inductive hypothesis,
$$  (i,j)\dd\,\sind\big(j\dd\Fr_XV(x)\cap Y,X\big)\!\leq\!
        (i,j)\dd\,\sind\big(j\dd\Fr_XV(x)\cap Y\big)\!\leq\!k\!-\!1   $$
and so $(i,j)\dd\,\sind(Y,X)\leq k$. Thus
$(i,j)\dd\,\sind(Y,X)\leq(i,j)\dd\,\sind Y$.

(2) By (1), it suffices to prove only that if $U'\in\tau_j'$, then
there is $U\in\tau_j$ such that $U\cap Y=U'$ and $j\dd\Fr_XU\cap
Y=j\dd\Fr_YU'$. But this fact is proved in the first part of the
proof of Proposition~2 in [23].
\end{pf}

\begin{proposition}{9.13}
Let $(Y,\tau_1',\tau_2')$ be a $\BsS$ of a $\BS$
$(X,\tau_1,\tau_2)$ and $Y=Y_1\cup Y_2$. Then
$$  (i,j)\dd\,\sind(Y,X)\leq
        (i,j)\dd\,\sind(Y_1,X)+(i,j)\dd\,\sind(Y_2,X)+1.  $$
\end{proposition}

\begin{pf}
We shall use the induction under a nonnegative integer
$n=(i,j)\dd\,\sind(Y_1,X)+(i,j)\dd\,\sind(Y_2,X)$. Let the
inequality be proved for $n\leq k-1$ and prove it for $n=k$. If
$x\in Y_1\sbs Y$ is any point, then for any $U(x)\in\tau_i$ there
is $V(x)\in\tau_j$ such that $\tau_j\cl V(x)\sbs U(x)$ and
$$  (i,j)\dd\,\sind\big(j\dd\Fr_XV(x)\cap Y_1,X\big)\leq
            (i,j)\dd\,\sind(Y_1,X)-1.         $$
Since
$$  j\dd\Fr_XV(x)\cap Y=\big(j\dd\Fr_XV(x)\cap Y_1\big)\cup
            \big(j\dd\Fr_XV(x)\cap Y_2\big)     $$
and
\begin{eqnarray*}
    &\ds (i,j)\dd\,\sind\big(j\dd\Fr_XV(x)\cap Y_1,X\big)\!+\!
        (i,j)\dd\,\sind\big(j\dd\Fr_XV(x)\cap Y_2,X\big)\!\leq \\
    &\ds \leq (i,j)\dd\,\sind(Y_1,X)-1+(i,j)\dd\,\sind(Y_2,X)=k-1,
\end{eqnarray*}
one can use the inductive assumption under the set
$j\dd\Fr_XV(x)\cap Y$. Hence
$$  (i,j)\dd\,\sind\big(j\dd\Fr_XV(x)\cap Y,X\big)\leq k-1+1=k  $$
and thus $(i,j)\dd\,\sind(Y,X)\leq k+1$.
\end{pf}

\begin{proposition}{9.14}
If $(X,\tau_1<\tau_2)$ is $2$-normal and $p$-normal,    \linebreak
$Y\in\co\tau_2$ and $(1,2)\dd\,\sind Y=0$, then
$(1,2)\dd\,\sind(Y,X)=0$.
\end{proposition}

\begin{pf}
We shall prove that for each $x\in Y$ and any $U(x)\in\tau_1$
there is $V(x)\in\tau_2$ such that $\tau_2\cl V(x)\sbs U(x)$ and
$(\tau_2\cl V(x)\setminus V(x))\cap Y=\vnth$.

Since $(1,2)\dd\,\sind Y=0$, there is $V\in\tau_2'\cap\co\tau_2'$
such that $x\in V\sbs U(x)\cap Y$. Clearly, $Y\in\co\tau_2$
implies that $V\in\co\tau_2$, $Y\setminus V\in\co\tau_2$ and since
$(X,\tau_1,\tau_2)$ is $2$-normal, there are $U(V)\in\tau_2$,
$U(Y\setminus V)\in\tau_2$ such that $U(V)\cap U(Y\setminus
V)=\vnth$. Moreover, since $X$ is $p$-normal and $V\sbs U(x)$, by
(4) of Definition~1.1, there is $E(V)\in\tau_1$ such that $V\sbs
E(V)\sbs\tau_2\cl E(V)\sbs U(x)$. If $W(V)=E(V)\cap U(V)$, then
$W(V)\in\tau_2$ as $\tau_1\sbs\tau_2$, and
$$  V\sbs W(V)\sbs \tau_2\cl W(V)\sbs X\setminus(Y\setminus V)=
        (X\setminus Y)\cup V        $$
since $\tau_2\cl W(V)\sbs\tau_2\cl U(V)$ and $\tau_2\cl U(V)\cap
(Y\setminus V)=\vnth$. It is clear that $\tau_2\cl W(V)\sbs U(x)$
and $\tau_2\cl W(V)\cap Y=\tau_2\cl V=V$. Let $W(V)=V(x)$. Then
\begin{eqnarray*}
    &\ds 2\dd\Fr_XV(x)\cap Y=
            \big(\tau_2\cl V(x)\setminus V(x)\big)\cap Y= \\
    &\ds =\big(\tau_2\cl W(V)\setminus W(V)\big)\cap Y=V\setminus V=\vnth.
\end{eqnarray*}
\vskip-0.7cm
\end{pf}
\vskip+0.2cm

\begin{theorem}{9.15}
If $(X,\tau_1<\tau_2)$ is $(1,2)$-separately regular, \linebreak
$Y\!\in\!2\dd\,\cD(X)$ and $Z\!\sbsq\!Y\!\sbsq\!X$, then
$(1,2)\dd\,\sind(Z,Y)\!=\!(1,2)\dd\,\sind(Z,X)$.
\end{theorem}

\begin{pf}
By (3) of Theorem~9.10 it suffices to prove only that \linebreak
$(1,2)\dd\,\sind(Z,X)\leq (1,2)\dd\,\sind(Z,Y)$. We shall use the
induction under a nonnegative integer $n=(1,2)\dd\,\sind(Z,Y)$.
Let the inequality be proved for $n\leq k-1$ and prove it for
$n=k$. If $x\in Z$ is any point and $U(x)\in\tau_1$ is any
neighborhood, then $k=(1,2)\dd\,\sind(Z,Y)$ implies that there is
$V'(x)\in\tau_2'$ such that $\tau_2'\cl V'(x)\sbs U(x)\cap Y$ and
$$  (1,2)\dd\,\sind\big(2\dd\Fr_YV'(x)\cap Z,Y\big)\leq k-1.    $$

Moreover, by inductive assumption, we have
\begin{eqnarray*}
    &\ds (1,2)\dd\,\sind\big(2\dd\Fr_YV'(x)\cap Z,X\big)\leq \\
    &\ds \leq (1,2)\dd\,\sind\big(2\dd\Fr_YV'(x)\cap Z,Y\big)\leq k-1.
\end{eqnarray*}

Let $W(x)\in\tau_2$ and $W(x)\cap Y=V'(x)$. Since $V'(x)\sbs U(x)$
and $\tau_1\sbs\tau_2$, for $V(x)=W(x)\cap U(x)\in\tau_2$ we have
$V(x)\sbs U(x)$ and $V(x)\cap Y=V'(x)$. Since $Y\in 2\dd\,\cD(X)$,
by Lemma~2.16 we have
$$  \tau_2\cl V(x)=\tau_2\cl(V(x)\cap Y)=\tau_2\cl V'(x).   $$
Therefore, $Z\sbs Y$ implies that
\begin{eqnarray*}
    &\ds 2\dd\Fr_XV(x)\cap Z=
        \big(\tau_2\cl V(x)\setminus V(x)\big)\cap Z= \\
    &\ds =\big(\big(\tau_2\cl V(x)\setminus V(x)\big)\cap Y\big)\cap Z=
        \big((\tau_2\cl V'(x)\setminus V(x))\cap Y\big)\cap Z= \\
    &\ds =\big(\tau_2'\cl V'(x)\setminus V'(x)\big)\cap Z=
            2\dd\Fr_YV'(x)\cap Z.
\end{eqnarray*}
Since $(X,\tau_1,\tau_2)$ is $(1,2)$-separately regular, one can
assume that        \linebreak     $\tau_2\cl V(x)\sbs U(x)$.
Hence, for each $x\in Z$ and any $U(x)\in\tau_1$ there is
$V(x)\in\tau_2$ such that $\tau_2\cl V(x)\sbs U(x)$ and
\begin{eqnarray*}
    &\ds (1,2)\dd\,\sind\big(2\dd\Fr_XV(x)\cap Z,X\big)\!=\!
        (1,2)\dd\,\sind\big(2\dd\Fr_YV'(x)\cap Z,X\big)\!\leq \\
    &\ds \leq (1,2)\dd\,\sind\big(2\dd\Fr_YV'(x)\cap Z,Y\big)\leq k-1.
\end{eqnarray*}
Therefore, $(1,2)\dd\,\sind(Z,X)\leq k$ and thus
$(1,2)\dd\,\sind(Z,X)\leq        \linebreak
(1,2)\dd\,\sind(Z,Y)$.
\end{pf}

\begin{theorem}{9.16}
If $f:(X,\tau_1<\tau_2)\to (X_1,\gm_1<\gm_2)$ is a
$d$-con\-ti\-nu\-ous surjection and $(Y,\tau_1'<\tau_2')$ is a
$\BsS$ of $X$ such that the restriction
$f\big|_Y:(Y,\tau_1'<\tau_2')\to (Y_1=f(Y),\gm_1'<\gm_2')$ is a
$d$-ho\-me\-o\-mor\-phism, then $(1,2)\dd\,\sind(Y,X)\leq
(1,2)\dd\,\sind(Y_1,X_1)$. Moreover, if $(X_1,\gm_1<\gm_2)$ is
hereditarily $2$-normal, then
$$  (1,2)\dd\,\sind(Y,X)=(1,2)\dd\,\sind(Y_1,X_1)=(1,2)\sind Y. $$
\end{theorem}

\begin{pf}
First, let us prove the inequality
$$  (1,2)\dd\,\sind(Y,X)\leq (1,2)\dd\,\sind(Y_1,X_1)   $$
by induction under a nonnegative integer
$n=(1,2)\dd\,\sind(Y_1,X_1)$. Suppose that the inequality is
proved for $n\leq k-1$ and prove it for $n=k$. Let $x\in Y$ be any
point and $U(x)\in\tau_1$ be any  neighborhood. Since $f\big|_Y$
is a $d$-homeomorphism, $f(U(x)\cap Y)\in\gm_1'$.

But $(1,2)\dd\,\sind(Y_1,X_1)=k$ and so for a set $U\in\gm_1$,
where $U\cap Y_1=f(U(x)\cap Y)$, there is a set $W\in\gm_2$ such
that $\gm_2\cl W\sbs U$ and
$$  (1,2)\dd\,\sind\big(2\dd\Fr_{X_1}W\cap Y_1,X_1\big)\leq k-1.  $$
On the other hand, the set $V(x)=f^{-1}(W)\cap U(x)\in\tau_2$ and
it suffices to prove that
$$  (1,2)\dd\,\sind\big(2\dd\Fr_XV(x)\cap Y,X\big)\leq k-1.   $$
By the inequality (6) in the proof of Lemma~1 in [24], we have
$$  2\dd\Fr_XV(x)\cap Y\sbsq f^{-1}(2\dd\Fr_{X_1}W)\cap Y=
        f^{-1}(2\dd\Fr_{X_1}W\cap Y_1).     $$
Hence, by inductive assumption and (3) of Theorem~9.10,
$$  (1,2)\dd\,\sind\big(2\dd\Fr_XV(x)\cap Y,X\big)\leq k-1    $$
so that $(1,2)\dd\,\sind(Y,X)\leq k$ and thus
$$  (1,2)\dd\,\sind(Y,X)\leq (1,2)\dd\,\sind(Y_1,X_1).  $$

For the second part, first of all note that since $f\big|_Y$ is a
$d$-ho\-me\-o\-mor\-phism, we have $(1,2)\dd\,\sind
Y=(1,2)\dd\,\sind Y_1$. Since $(X_1,\gm_1,\gm_2)$ is hereditarily
$2$-normal, by (1) of Proposition~9.12, for a $\BsS$ $Y_1\sbs X_1$
we have $(1,2)\dd\,\sind(Y_1,X_1)=(1,2)\dd\,\sind Y_1$. Therefore
it remains to prove the equality
$(1,2)\dd\,\sind(Y,X)=(1,2)\dd\,\sind Y$, i.e. by (1) of
Theorem~9.10, it suffices to prove only that
$(1,2)\dd\,\sind(Y,X)\leq (1,2)\dd\,\sind Y$. We shall use the
induction under a nonnegative integer $n=(1,2)\dd\,\sind Y$. Let
the inequality be proved for $n\leq k-1$ and prove it for $n=k$.
If $x\in Y$ and $U(x)\in\tau_1$, then there is $V'(x)\in\tau_2'$
such that
$$  \tau_2'\cl V'(x)\sbs U(x)\cap Y \;\;\text{and}\;\;
        (1,2)\dd\,\sind\big(2\dd\Fr_YV'(x)\big)\leq k-1.  $$
Since $f\big|_Y$ is a $d$-homeomorphism, we have
$f(V'(x))\in\gm_2'$. But the $\BS$ $(X_1,\gm_1,\gm_2)$ is
hereditarily $2$-normal and by the proof of the first part of
Proposition~1 in [23], there is $W\in\gm_2$ such that
$$  W\cap Y_1=f(V'(x)), \;\;\;
        2\dd\Fr_{X_1}W\cap Y_1=2\dd\Fr_{Y_1}f(V'(x)).   $$
Let $V(x)=f^{-1}(W)\cap U(x)$. It is clear that $V(x)\in\tau_2$
and $V(x)\cap Y=V'(x)$. By the proof of Lemma~2 in [24],
$2\dd\Fr_XV(x)\cap Y=2\dd\Fr_YV'(x)$. Therefore, it remains to use
the inductive assumption with respect to the set $2\dd\Fr_YV'(x)$,
i.e.,
\begin{eqnarray*}
    &\ds (1,2)\dd\,\sind\big(2\dd\Fr_XV(x)\cap Y,X\big)=
        (1,2)\dd\,\sind\big(2\dd\Fr_YV'(x),X\big)\leq \\
    &\ds \leq (1,2)\dd\,\sind(2\dd\Fr_YV'(x))\leq k-1.
\end{eqnarray*}
Hence $(1,2)\dd\,\sind(Y,X)\!\leq\!k$ and thus
$(1,2)\dd\,\sind(Y,X)\!\leq\!(1,2)\dd\,\sind Y$.~\end{pf}

\begin{definition}{9.17}
For a $\BsS$ $(Y,\tau_1',\tau_2')$ of a $\BS$ $(X,\tau_1,\tau_2)$
and a relatively $i$-regular pair $(x,A)$, $x\in Y$,
$A\in\co\tau_i$, $x\ol{\in}\,A$, a separate relative partition,
corresponding to $(x,A)$, is a $j$-closed set $T\sbs X$ such that
$X\setminus T$  is not $j$-connected and
$$  (X\setminus T)\cap Y=Y\setminus T\sbs H_1\cup H_2,  $$
where $x\in H_1'=H_1\cap Y\in\tau_j'$, $A\sbs H_2\in\tau_j$ and
$H_1\cap H_2=\vnth$.
\end{definition}

\begin{remark}{9.18}
Let $(Y,\tau_1',\tau_2')$ be a $\BsS$ of a $\BS$
$(X,\tau_1,\tau_2)$ and $(x,A)$ be a relatively $i$-regular pair.
Then
\begin{enumerate}
\item[(1)] If there is a neighborhood $U(x)\in\tau_j$
$(U(A)\in\tau_j)$ such that $\tau_j\cl U(x)\sbs X\setminus A$
$(\tau_j\cl U(A)\sbs X\setminus\{x\})$, then the set
$j\dd\Fr_XU(x)$ $(j\dd\Fr_XU(A))$ is a separate relative
partition, corresponding to $(x,A)$.

\item[(2)] If $T$ is a separate relative partition, corresponding
to $(x,A)$, then $j\dd\Fr_XH_i\sbs T$.
\end{enumerate}
\end{remark}

The proof is identical to the proof of Remark~9.7.

\begin{proposition}{9.19}
Let $(Y,\tau_1',\tau_2')$ be a $\BsS$ of a $\BS$
$(X,\tau_1,\tau_2)$ and $n$ denote a nonnegative integer. Then
$(i,j)\dd\,\sind(Y,X)\leq n$ if and only if to every relatively
$i$-regular pair $(x,A)$ there corresponds a separate relative
partition $T$ such that $(i,j)\dd\,\sind(T\cap Y,X)\leq n-1$.
\end{proposition}

\begin{pf}
First, let $(i,j)\dd\,\sind(Y,X)\leq n$ and $(x,A)$ be a
relatively    \linebreak        $i$-regular pair. Then $x\in
U(x)=X\setminus A\in\tau_i$ and by (2) of Definition~9.9 there is
a neighborhood $V(x)\in\tau_j$ such that
$$  \tau_j\cl V(x)\sbs U(x) \;\;\text{and}\;\;
        (i,j)\dd\,\sind\big(j\dd\Fr_XV(x)\cap Y,X\big)\leq k-1.   $$
But by (1) of Remark~9.18, $j\dd\Fr_XV(x)$ is a separate relative
partition, corresponding to $(x,A)$.

Conversely, let us suppose that the condition is satisfied and let
us prove that $(i,j)\dd\,\sind(Y,X)\leq n$. Let $x\in Y$ be any
point and $U(x)\in\tau_i$ be any neighborhood. Then the pair
$(x,A=X\setminus U(x))$ is a relatively $i$-regular pair and by
condition, there is a separate relative partition $T$ for $(x,A)$
such that $(i,j)\dd\,\sind(T\cap Y,X)\leq n-1$. But $Y\setminus
T\sbs H_1\cup H_2$, where $x\in H_1'\sbs H_1\in\tau_j$, $A\sbs
H_2\in\tau_j$, and $H_1\cap H_2=\vnth$. Clearly,
$$  \tau_j\cl H_1\sbs X\setminus H_2\sbs X\setminus A=U(x). $$
Let $H_1=V(x)$. Then $\tau_j\cl V(x)\sbs U(x)$ and by (2) of
Remark~9.18,
$$  j\dd\Fr_XV(x)\cap Y=j\dd\Fr_XH_1\cap Y\sbs T\cap Y. $$
Since $(i,j)\dd\,\sind(T\cap Y,X)\leq n-1$, by (3) of
Theorem~9.10,       \linebreak
       $(i,j)\dd\,\sind(j\dd\Fr_XV(x)\cap Y,X)\leq n-1$. Thus, it remains
to use (2) of Definition~9.9.
\end{pf}

\begin{proposition}{9.20}
If a $\BsS$ $(Y,\tau_1',\tau_2')$ of a $\BS$ $(X,\tau_1,\tau_2)$
is $(i,j)\dd\WS$-separately su\-per\-re\-gu\-lar in $X$ and
$j$-extremally disconnected in $X$, then $(i,j)\dd\,\sind(Y,X)=0$.
\end{proposition}

\begin{pf}
Let, $x\!\in\!Y$ and $U(x)\!\in\!\tau_i$ be any neighborhood.
Then, by (2) of Proposition~2.36, there is $V'(x)\!\in\!\tau_j'$
such that $\tau_j\cl V'(x)\!\sbs\!U(x)$, where, by topological
version of Definition~2.39, $\tau_j\cl
V'(x)\!\in\!\tau_j\cap\co\tau_j$. Hence $j\dd\Fr_XV'(x)\!=\!\vnth$
and by (2) of Definition~9.9,
$(i,j)\dd\,\sind(Y,X)\!=~\!\!\!0$.~\end{pf}

\vskip+0.5cm
\section*{\textbf{10. $(i,j)$-Large Separate Inductive Dimension Functions
and $(i,j)$-Large Separate Relative Inductive Dimension
Functions}} \vskip+0.2cm

\begin{definition}{10.1}
Let $(X,\tau_1,\tau_2)$ be a $\BS$ and $n$ denote a nonnegative
integer. Then
\begin{enumerate}
\item[(1)] $(i,j)\dd\,\sInd X=-1$ if and only if $X=\vnth$.

\item[(2)] $(i,j)\dd\,\sInd X\leq n$ if for each set
$F\in\co\tau_i$ and any neighborhood $U(F)\in\tau_i$ there is a
neighborhood $V(F)\in\tau_j$ such that $\tau_j\cl V(F)\sbs U(F)$
and $(i,j)\dd\,\sInd(j\dd\Fr_XV(F))\leq n-1$.

\item[(3)] $(i,j)\dd\,\sInd X=n$ if $(i,j)\dd\,\sInd X\leq n$ and
$(i,j)\dd\,\sInd X\leq n-1$ does not hold.

\item[(4)] $(i,j)\dd\,\sInd X=\infty$ if $(i,j)\dd\,\sInd X\leq n$
does not hold for any $n$.
\end{enumerate}
\end{definition}

As usual
$$  p\,\dd\,\sInd X\leq n\llra
        \big((1,2)\dd\,\sInd X\leq n\wedge (2,1)\dd\,\sInd X\leq n\big). $$

It follows immediately from Proposition~2.29 that if
$(i,j)\dd\,\sInd X$ is finite, then $(X,\tau_1,\tau_2)$ is
$(i,j)$-separately normal.

\begin{proposition}{10.2}
If $(Y,\tau_1'\!<\!\tau_2')$ is a $2$-closed $\BsS$ of a $\BS$
$(X,\tau_1\!<~\!\!\!\tau_2)$, then $(2,1)\dd\,\sInd Y\leq
(2,1)\dd\,\sInd X$.

Moreover, if $(Y,\tau_1'\!<_C\!\tau_2')$ is a $1$-closed $\BsS$ of
a $\BS$ $(X,\tau_1\!<~\!\!\!\tau_2)$, then $(1,2)\dd\,\sInd Y\leq
(2,1)\dd\,\sInd X$.
\end{proposition}

\begin{pf}
 First, let $Y\in\co\tau_2$ and $n=(2,1)\dd\,\sInd X$.
Let the inequality be proved for $n\leq k-1$ and prove it for
$n=k$. Let $F\in\co\tau_2'$ and $U'(F)\in\tau_2'$. Then
$F\in\co\tau_2$ and if $U(F)\in\tau_2$, $U(F)\cap Y=U'(F)$, then
$(2,1)\dd\,\sInd X=k$ implies that there is $V(F)\in\tau_1$ such
that $\tau_1\cl V(F)\sbs U(F)$ and
$(2,1)\dd\,\sInd(1\dd\Fr_XV(F))\leq k-1$. Let $V'(F)=V(F)\cap Y$.
Then
$$  \tau_1'\cl V'(F)=\tau_1\cl V'(F)\cap Y\sbs U(F)\cap Y=U'(F) $$
and $1\dd\Fr_YV'(F)\sbsq 1\dd\Fr_XV(F)$. Since $1\dd\Fr_YV'(F)$ is
$2$-closed in \linebreak $(1\dd\Fr_XV(F),\tau_1''<\tau_2'')$, by
inductive hypothesis,
$$  (2,1)\dd\,\sInd\big(1\dd\Fr_YV'(F)\big)\leq
        (2,1)\dd\,\sInd\big(1\dd\Fr_XV(F)\big)\leq k-1,  $$
and hence, $(2,1)\dd\,\sInd Y\leq k$ so that $(2,1)\dd\,\sInd
Y\leq(2,1)\dd\,\sInd X$.

Now, let us prove that $(1,2)\dd\,\sInd Y\leq (2,1)\dd\,\sInd X$
for $Y\in\co\tau_1$ and $n=(2,1)\dd\,\sInd X$. Let the inequality
be proved for $n\leq k-1$ and prove it for $n=k$. Let
$F\in\co\tau_1'$ and $U'(F)\in\tau_1'$. Then $Y\in\co\tau_1$
implies that $F\in\co\tau_1$. Let $U(F)\in\tau_1$ and $U(F)\cap
Y=U'(F)$. Since $F\in\co\tau_1\sbs\co\tau_2$,
$U(F)\in\tau_1\sbs\tau_2$ and $(2,1)\dd\,\sInd X\leq k$, there is
a neighborhood $V(F)\in\tau_1$ such that $\tau_1\cl V(F)\sbs U(F)$
and $(2,1)\dd\,\sInd(1\dd\Fr_XV(F))\leq k-1$. Since
$\tau_1'<_C\tau_2'$, by (2) of Corollary~2.2.7 in~[8],
$$  2\dd\Fr_YV'(F)=1\dd\Fr_YV'(F)\sbsq 1\dd\Fr_XV(F),       $$
where $V'(F)=V(F)\cap Y$ and as $1\dd\Fr_YV'(F)$ is $1$-closed in
\linebreak        $(1\dd\Fr_XV(F),\tau_1''<\tau_2'')$, by
inductive assumption,
$$  (1,2)\dd\,\sInd\big(2\dd\Fr_YV'(F)\big)\leq
        (2,1)\dd\,\sInd\big(1\dd\Fr_XV(F)\big)\leq k-1.   $$
Therefore, $(1,2)\dd\,\sInd Y\leq (2,1)\dd\,\sInd X$.
\end{pf}

\begin{proposition}{10.3}
For an $i\dd\TT_1$ $\BS$ $(X,\tau_1,\tau_2)$ we have
$$  (i,j)\dd\,\sind X\leq (i,j)\dd\,\sInd X.        $$
\end{proposition}

\begin{theorem}{10.4}
The following inequalities hold for a $\BS$
$(X,\tau_1\!<~\!\!\tau_2)$: $(1,2)\dd\,\sInd X\leq 2\dd\,\Ind X$
and $1\dd\,\Ind X\leq (2,1)\dd\,\sInd X$.
\end{theorem}

\begin{pf}
First, let us prove the right inequality by induction under a
nonnegative integer $n=(2,1)\dd\,\sInd X$. Let the inequality be
proved for $n\leq k-1$ and prove it for $n=k$. Let $F\in\co\tau_1$
and $U(F)\in\tau_1$ be any neighborhood. Then $F\in\co\tau_2$,
$U(F)\in\tau_2$ and since $(2,1)\dd\,\sInd X=k$, there is
$V(F)\in\tau_1$ such that
$$  \tau_1\cl V(F)\sbs U(F) \;\;\text{and}\;\;
        (2,1)\dd\,\sInd\big(1\dd\Fr_XV(F)\big)\leq k-1. $$
Hence, by inductive assumption
$$  1\dd\Ind\big(1\dd\Fr_XV(F)\big)\leq
        (2,1)\dd\,\sInd\big(1\dd\Fr_XV(F)\big)\leq k-1. $$
Thus $1\dd\,\Ind X\leq k$ and so $1\dd\,\Ind X\leq (2,1)\dd\,\sInd
X$.

Now, let $2\dd\Ind X=n$, the left inequality be proved for
$n\!\leq\!k\!-~\!\!\!1$ and prove it for $n=k$. Let
$F\in\co\tau_1$ and $U(F)\in\tau_1$ be any neighborhood. Since
$\tau_1\sbs\tau_2$ and $2\dd\Ind X=k$, there is $V(F)\in\tau_2$
such that $\tau_2\cl V(F)\sbs U(F)$ and
$2\dd\,\Ind(2\dd\Fr_XV(F))\leq k-1$. Hence, by inductive
assumption
$$  (1,2)\dd\,\sInd\big(2\dd\Fr_XV(F)\big)\leq
        2\dd\,\Ind\big(2\dd\Fr_XV(F)\big)\leq k-1,  $$
and so $(1,2)\dd\,\sInd X\leq k$. Thus $(1,2)\dd\,\sInd X\leq
2\dd\,\Ind X$.
\end{pf}

\begin{corollary}{10.5}
For a $\BS$ $(X,\tau_1,\tau_2)$ we have
\begin{enumerate}
\item[(1)] If $\tau_1<_C\tau_2$, then $1\dd\Ind X\leq
\min((1,2)\dd\Ind X,(2,1)\dd\,\sInd X)$.

\item[(2)] If $\tau_1<_N\tau_2$, then $1\dd\Ind X\leq
\min((1,2)\dd\Ind X,(2,1)\dd\,\sInd X)$ and $\max((2,1)\dd\Ind
X,(1,2)\dd\,\sInd X)\leq 2\dd\Ind X$.
\end{enumerate}
\end{corollary}

\begin{pf}
Follows from Theorem~10.4 taking into account (1) and (2) of
Theorem~3.2.38 in [6].
\end{pf}

\begin{proposition}{10.6}
If a $\BS$ $(X,\tau_1,\tau_2)$ is $(i,j)$-separately normal and
$j$-extremally disconnected, then $(i,j)\dd\,\sInd X=0$.
\end{proposition}

\begin{pf}
Let $F\in\co\tau_i$ and $U(F)\in\tau_i$ be any neighborhood. Then,
by Proposition~2.29, there is a neighborhood $V(F)\in\tau_j$ such
that $\tau_j\cl V(F)\sbs U(F)$. Since $X$ is $j$-extremally
disconnected, $W(F)=\tau_j\cl V(F)\in\tau_j\cap\co\tau_j$ and,
hence, $j\dd\Fr_XW(F)=\vnth$. Thus, by (2) of Definition~10.1,
$(i,j)\dd\,\sInd X=0$.
\end{pf}

\begin{definition}{10.7}
Let $(A,B)$ be an $i$-normal pair in a $\BS$ $(X,\tau_1,\tau_2)$,
i.e., $A,\,B\in\co\tau_i$ and $A\cap B=\vnth$. Then we say that a
$j$-closed set $T$ is a separate partition, corresponding to the
$i$-normal pair $(A,B)$, if $X\setminus T$ is not $j$-connected so
that $X\setminus T=H_1\cup H_2$, where
$H_1,\,H_2\in\tau_j\setminus\{\vnth\}$, $A\sbs H_1$, $B\sbs H_2$
and $H_1\cap H_2=\vnth$.
\end{definition}

\begin{remark}{10.8}
One can easily to verify that in a $\BS$ $(X,\tau_1,\tau_2)$ the
following conditions are satisfied for an $i$-normal pair $(A,B)$:
\begin{enumerate}
\item[(1)] If there exists a $j$-open neighborhood $U(A)$
$(j$-open neighborhood $U(B))$ such that $\tau_j\cl U(A)\sbs
X\setminus B$ $(\tau_j\cl U(B)\sbs X\setminus A)$, then the set
$j\dd\Fr_XU(A)$ $(j\dd\Fr_XU(B))$ is a separate partition,
corresponding to $(A,B)$.

\item[(2)] If $T$ is a separate partition, corresponding to
$(A,B)$, then       \linebreak      $j\dd\Fr_XH_i\sbs T$.
\end{enumerate}
\end{remark}

The proof is similar to the proof of Remark~9.7 and can be
omitted.

\begin{proposition}{10.9}
Let $(Y,\tau_1',\tau_2')$ be a $\BsS$ of a $\BS$
$(X,\tau_1,\tau_2)$ and $n$ denote a nonnegative integer. If
$(i,j)\dd\,\sInd X\leq n$, then to every $i$-normal pair $(A,B)$
there corresponds a separate partition $T$ such that
$(i,j)\dd\,\sInd T\leq n-1$.
\end{proposition}

\begin{pf}
Let $(i,j)\dd\,\sInd X\leq n$ and $(A,B)$ be an $i$-normal pair.
Then $U(A)=X\setminus B\!\in\!\tau_i$ is an $i$-open neighborhood
and since $(i,j)\dd\,\sInd X\!\leq~\!\!\!n$, there is
$V(A)\in\tau_j$ such that $\tau_j\cl V(A)\sbs U(A)$ and
$$  (i,j)\dd\,\sInd\big(j\dd\Fr_XV(A)\big)\leq n-1.     $$
Thus, it remains to use (1) of Remark~10.8.
\end{pf}

\begin{corollary}{10.10}
Let $(X,\tau_1<\tau_2)$ be a $\BS$ and $n$ denote a nonnegative
integer. Then $(2,1)\dd\,\sInd X\leq n$ if and only if to every
\linebreak          $2$-nor\-mal pair $(A,B)$ there corresponds a
separate partition $T$ such that $(2,1)\dd\,\sInd T\leq n-1$.
\end{corollary}

\begin{pf}
The first part is proved by Proposition~10.9.

Conversely, let the condition is satisfied, $F\in\co\tau_2$ and
$U(F)\in\tau_2$. Then for a $2$-normal pair $(F,X\setminus U(F))$
there is a separate partition $T$ such that $X\setminus T=H_1\cup
H_2$, $H_1,\,H_2\in\tau_1\setminus\{\vnth\}$, $F\sbs H_1$,
$X\setminus U(F)\sbs H_2$, $H_1\cap H_2=\vnth$ and
$(2,1)\dd\,\sInd T\leq n-1$. Clearly,
$$  \tau_1\cl H_1\sbs X\setminus(X\setminus U(F))=U(F). $$
Let $H_1=V(F)$. Then $\tau_1\cl V(F)\sbs U(F)$ and by (2) of
Remark~10.8 together with the first part of Proposition~10.2,
$(2,1)\dd\,\sInd(1\dd\Fr_XV(F))\!\leq n-1$. Thus $(2,1)\dd\,\sInd
X\leq n$.
\end{pf}

\begin{definition}{10.11}
Let $(Y,\tau_1',\tau_2')$ be a $\BsS$ of a $\BS$
$(X,\tau_1,\tau_2)$ and $n$ denote a nonnegative integer. Then
\begin{enumerate}
\item[(1)] $(i,j)\dd\,\sInd(Y,X)=-1$ if and only if $Y=\vnth$.

\item[(2)] $(i,j)\dd\,\sInd(Y,X)\leq n$ if for any set
$F\in\co\tau_i'$ and any neighborhood $U(F)\in\tau_i$ there is a
neighborhood $V(F)\in\tau_j$ such that $\tau_j\cl V(F)\sbs U(F)$
and $(i,j)\dd\,\sInd\big((j\dd\Fr_XV(F)\cap Y,X\big)\leq
n-1(\llra$ for any set $F\in\co\tau_i'$ and any neighborhood
$U(F)\in\tau_i$ there is a neighborhood $V(F)\in\tau_j$ such that
$\tau_j\cl V(F)\sbs U(F)$ and $(i,j)\dd\,\sInd(j\dd\Fr_XV(F)\cap
Y)\leq n-1)$.

\item[(3)] $(i,j)\dd\,\sInd(Y,X)=n$ if $(i,j)\dd\,\sInd(Y,X)\leq
n$ and         \linebreak      $(i,j)\dd\,\sInd(Y,X)\leq n-1$ does
not hold.

\item[(4)] $(i,j)\dd\,\sInd(Y,X)=\infty$ if
$(i,j)\dd\,\sInd(Y,X)\leq n$ does not hold for any $n$.
\end{enumerate}
\end{definition}

As usual
$$  p\,\dd\,\sInd(Y,X)\!\leq\!n\!\llra\!
        \big((1,2)\dd\,\sInd(Y,X)\!\leq\!n\wedge (2,1)\dd\,\sInd(Y,X)\!\leq\!n\big). $$

It follows immediately from (2) of Definition~10.11 and (3) of
Proposition~2.31 that if $(i,j)\dd\,\sInd(Y,X)$ is finite, then
$Y$ is $(i,j)$-separately supernormal in $X$. Hence, by
implications after Definition~2.30, $Y$ is $(i,j)$-strongly
separately normal in $X$, $Y$ is $(i,j)\dd\WS$-separately
supernormal in $X$, $Y$ is $(i,j)$-separately normal in $X$, $Y$
is $(i,j)\dd\WS$-separately normal in $X$ and $Y$ is $(i,j)$-quasi
separately normal in $X$.

\begin{proposition}{10.12}
For an $i\dd\TT_1$ $\BsS$ $(Y,\tau_1',\tau_2')$ of a $\BS$
$(X,\tau_1,\tau_2)$ we have $(i,j)\dd\,\sind(Y,X)\leq
(i,j)\dd\,\sInd(Y,X)$.
\end{proposition}

\begin{proposition}{10.13}
If $(Z,\tau_1''<\tau_2'')\sbs (Y,\tau_1'<\tau_2')\sbs
(X,\tau_1<\tau_2)$ and $Z\in\co\tau_2'$, then
$(2,1)\dd\,\sind(Z,X)\leq (2,1)\dd\,\sInd(Y,X)$.
\end{proposition}

\begin{pf}
Let us prove \textbf{the monotonicity property of
$\boldsymbol{(2,\!1)\dd\,\sInd(\!Y,\!X\!)}$} by induction under a
nonnegative integer $n=(2,1)\dd\,\sInd(Y,X)$. Let the inequality
be proved for $n\leq k-1$ and prove it for $n=k$. Let
$F\in\co\tau_2''$ and $U(F)\in\tau_2$. Then $F\in\co\tau_2'$ and
since $(2,1)\dd\,\sInd(Y,X)=k$, there is $V(F)\in\tau_1$ such that
$$  \tau_1\cl V(F)\sbs U(F) \;\;\text{and}\;\;
        (2,1)\dd\,\sInd\big(1\dd\Fr_XV(F)\cap Y,X\big)\leq k-1. $$
Since $1\dd\Fr_XV(F)\cap Z$ is $2$-closed in $(1\dd\Fr_X(V(F)\cap
Y,\tau_1''<\tau_2'')$, by inductive assumption,
\begin{eqnarray*}
    (2,1)\dd\,\sInd\big(1\dd\Fr_XV(F)\cap Z,X\big) &&\hskip-0.6cm \leq \\
    \leq (2,1)\dd\,\sInd\big(1\dd\Fr_XV(F)\cap Y,X\big) &&\hskip-0.6cm \leq k-1
            \;\; \text{and so}
\end{eqnarray*}
$(2,1)\dd\,\sInd(Z,X)\!\leq\!k$. Thus
$(2,1)\dd\,\sInd(Z,X)\!\leq\!(2,1)\dd\,\sInd(Y,X)$.~\end{pf}

\begin{theorem}{10.14}
For a $\BsS$ $(Y,\tau_1'<\tau_2')$ of a $\BS$ $(X,\tau_1<\tau_2)$
the following inequalities hold:
\begin{enumerate}
\item[(1)] $(1,2)\dd\,\sInd(Y,X)\!\leq\!2\dd\Ind(Y,X)$,
$i\dd\Ind(Y,X)\!\leq\!(2,1)\dd\,\sInd(Y,X)$ and hence,
$$  (1,2)\dd\,\sInd(Y,X)\leq 2\dd\Ind(Y,X)\leq (2,1)\dd\Ind(Y,X). $$
\end{enumerate}

Moreover, if $\tau_1<_C\tau_2$, then in addition, we have:
\begin{enumerate}
\item[(2)] $(1,2)\dd\,\sInd(Y,X)\leq 1\dd\Ind(Y,X)$ and hence
$$  (1,2)\dd\,\sInd(Y,X)\leq i\dd\Ind(Y,X)\leq (2,1)\dd\,\sInd(Y,X).  $$
\end{enumerate}
\end{theorem}

\begin{pf}
(1) First, let us prove the left inequality by induction under a
nonnegative integer $n=2\dd\Ind(Y,X)$. Let the inequality be
proved for $n\leq k-1$ and prove it for $n=k$. Let
$F\in\co\tau_1'$ and $U(F)\in\tau_1$ be any neighborhood. Since
$\tau_1\sbs\tau_2$, we have $F\in\co\tau_2'$, $U(F)\in\tau_2$ and
by $2\dd\Ind(Y,X)=k$ there is $V(F)\in\tau_2$ such that
$$  \tau_2\cl V(F)\sbs U(F) \;\;\text{and}\;\;
        2\dd\Ind\big(2\dd\Fr_XV(F)\cap Y,X\big)\leq k-1.    $$
Hence, by inductive assumption,
$$  (1,2)\dd\,\sInd\big(2\dd\Fr_XV(F)\cap Y,X\big)\!\leq\!
        2\dd\Ind\big(2\dd\Fr_XV(F)\cap Y,X\big)\!\leq\!k\!-\!1    $$
so that $(1,2)\dd\,\sInd(Y,X)\leq k$ and thus
$(1,2)\dd\,\sInd(Y,X)\leq 2\dd\Ind(Y,X)$.

Now, let $n=(2,1)\dd\,\sInd(Y,X)$, the right inequality
$1\dd\Ind(Y,X)\leq (2,1)\dd\,\sInd(Y,X)$ be proved for $n\leq k-1$
and prove it for $n=k$. Let $F\in\co\tau_1'$, $U(F)\in\tau_1$.
Then $F\in\co\tau_2'$, $U(F)\in\tau_2$ and since
$(2,1)\dd\,\sInd(Y,X)=k$, there is $V(F)\in\tau_1$ such that
$\tau_1\cl V(F)\sbs U(F)$ and
$(2,1)\dd\,\sInd\big(1\dd\Fr_XV(F)\cap Y,X\big)\leq k-1$. Hence,
by inductive hypothesis
$$  1\dd\Ind\big(1\dd\Fr_XV(F)\cap Y,X\big)\!\leq\!
        (2,1)\dd\,\sInd\big(1\dd\Fr_XV(F)\cap Y,X\big)\!\leq\!k\!-\!1.   $$
Therefore $1\dd\Ind(Y,X)\leq k$ and so $1\dd\Ind(Y,X)\leq
(2,1)\sInd(Y,X)$.

Furthermore, for a nonnegative integer $n=(2,1)\dd\,\sInd(Y,X)$
let us prove the inequality $2\dd\Ind(Y,X)\leq (2,1)\sInd(Y,X)$.
Let the inequality be proved for $n\leq k-1$ and prove it for
$n=k$. Let $F\in\co\tau_2'$, $U(F)\in\tau_2$. Then
$(2,1)\dd\,\sInd(Y,X)=k$ implies that there is $V(F)\in\tau_1$
such that
$$  \tau_1\cl V(F)\sbs U(F) \;\;\text{and}\;\;
        (2,1)\dd\,\sInd\big(1\dd\Fr_XV(F)\cap Y,X\big)\leq k-1.  $$
Since $2\dd\Fr_XV(F)\sbs 1\dd\Fr_XV(F)$, the set
$2\dd\Fr_XV(F)\cap Y$ is $2$-closed in $1\dd\Fr_XV(F)\cap Y$ and
by Proposition~10.13, $(2,1)\dd\,\sInd(2\dd\Fr_XV(F)\cap Y,X)\leq
k-1$. Hence, by inductive assumption, $2\dd\Ind(2\dd\Fr_XV(F)\cap
Y,X)\leq k-1$. Since $\tau_2\cl V(F)\sbs U(F)$, we have
$2\dd\Ind(Y,X)\leq k$ and thus $2\dd\Ind(Y,X)\leq
(2,1)\dd\,\sInd(Y,X)$.

(2) Let $\tau_1<_C\tau_2$ and use the induction under
$n=1\dd\Ind(Y,X)$. Suppose that the inequality is proved for
$n\leq k-1$ and prove it for $n=k$. If $F\in\co\tau_1'$,
$U(F)\in\tau_1$, then $1\dd\Ind(Y,X)=k$ implies that there is
$V(F)\in\tau_1$ such that
$$  \tau_1\cl V(F)\sbs U(F) \;\;\text{and}\;\;
        1\dd\Ind\big(1\dd\Fr_XV(F)\cap Y,X\big)\leq k-1.    $$
Since $V(F)\in\tau_1$, by (2) of Corollary~2.2.7 in [8],
$1\dd\Fr_XV(F)=     \linebreak       2\dd\Fr_XV(F)$ and so
$1\dd\Ind\big(2\dd\Fr_XV(F)\cap Y,X\big)\leq k-1$. Therefore, by
inductive assumption $(1,2)\dd\,\sInd\big(2\dd\Fr_XV(F)\cap
Y,X\big)\leq k-1$. Since $V(F)\in\tau_1\sbs\tau_2$ and $\tau_2\cl
V(F)\sbs\tau_1\cl V(F)\sbs U(F)$, we have
$(1,2)\dd\,\sInd(Y,X)\leq k$. Thus $(1,2)\dd\,\sInd(Y,X)\leq
1\dd\Ind(Y,X)$.
\end{pf}

\begin{corollary}{10.15}
For a $\BsS$ $(Y,\tau_1'<\tau_2')$ of a $\BS$ $(X,\tau_1<\tau_2)$
we have $i\dd\Ind Y\leq (2,1)\dd\,\sInd(Y,X)$.
\end{corollary}

\begin{pf}
Follows directly from (1) of Corollary~8.6 and (1) of
Theorem~10.14.
\end{pf}

\begin{theorem}{10.16}
Let $(Y,\tau_1'\!<\!\tau_2')$ be a $\BsS$ of a $\BS$
$(X,\tau_1\!<~\!\!\!\tau_2)$. Then $(2,1)\dd\,\sInd Y\leq
(2,1)\dd\,\sInd(Y,X)$. Moreover, if $Y\in\co\tau_2$, then
$(2,1)\dd\,\sInd(Y,X)\leq (2,1)\dd\,\sInd X$ and thus
$(2,1)\dd\,\sInd(X,X)=          \linebreak     (2,1)\dd\,\sInd X$.
\end{theorem}

\begin{pf}
We shall prove the inequality $(2,1)\dd\,\sInd Y\leq
(2,1)\dd\,\sInd(Y,X)$ by induction under a nonnegative integer
$n=(2,1)\dd\,\sInd(Y,X)$. Let the inequality be proved for $n\leq
k-1$ and prove it for $n=k$. Let $F\in\co\tau_2'$ and
$U'(F)\in\tau_2'$ be any neighborhood. Then there is
$U(F)\in\tau_2$ such that $U(F)\cap Y=U'(F)$. Since
$(2,1)\dd\,\sInd(Y,X)=k$, there is $V(F)\in\tau_1$ such that
$$  \tau_1\cl V(F)\sbs U(F) \;\;\text{and}\;\;
        (2,1)\dd\,\sInd\big(1\dd\Fr_XV(F)\cap Y,X\big)\leq k-1.  $$
Hence, by inductive assumption, $(2,1)\dd\,\sInd(1\dd\Fr_XV(F)\cap
Y)\leq k-1$. Since $1\dd\Fr_YV'(F)$ is $1$-closed and, hence,
$2$-closed in $1\dd\Fr_XV(F)\cap Y$, where $V'(F)=V(F)\cap Y$, by
Proposition~10.2,
$$  (2,1)\dd\,\sInd\big(1\dd\Fr_YV'(F)\big)\leq
        (2,1)\dd\,\sInd\big(1\dd\Fr_XV(F)\cap Y\big)\leq k-1.  $$
Thus $(2,1)\dd\,\sInd Y\leq k$ and so $(2,1)\dd\,\sInd Y\leq
(2,1)\dd\,\sInd(Y,X)$.

Now, let $Y\in\co\tau_2$, $(2,1)\dd\,\sInd X=n$, the inequality be
proved for $n\leq k-1$ and prove it for $n=k$. Let
$F\in\co\tau_2'\sbs\co\tau_2$ and $U(F)\in\tau_2$. Since
$(2,1)\dd\,\sInd X=k$, there is $V(F)\in\tau_1$ such that
$$  \tau_1\cl V(F)\sbs U(F) \;\;\text{and}\;\;
        (2,1)\dd\,\sInd\big(1\dd\Fr_XV(F)\big)\leq k-1.  $$
Therefore, by inductive assumption,
$$  (2,1)\dd\,\sind\big(1\dd\Fr_XV(F)\cap Y,1\dd\Fr_XV(F)\big)\leq k-1  $$
and by the first inequality of this theorem,
$(2,1)\dd\,\sInd\big(1\dd\Fr_XV(F)\cap Y\big)\leq k-1$. Hence, by
the condition in brackets in (2) of Definition~10.11,
$(2,1)\dd\,\sInd(Y,X)\leq k$ and thus $(2,1)\dd\,\sInd(Y,X)\leq
\linebreak        (2,1)\dd\,\sInd X$.

Finally, it is evident that $(2,1)\dd\,\sInd(X,X)=(2,1)\dd\,\sInd
X$.
\end{pf}
\vskip+0.2cm

Recall that by Definition~2.3.15 in [8], a topology $\tau_1$ is
$i$-strongly near a topology $\tau_2$ on a set $X$ (briefly,
$\tau_1N(i)\tau_2)$ if the $N$ relation is hereditary with respect
to $i$-closed subsets of $X$.

\begin{theorem}{10.17}
If $(Y,\tau_1'<\tau_2')$ is a $1$-closed $\BsS$ of a $\BS$
\linebreak       $(X,\tau_1<_{N(1)}\tau_2)$, then $(1,2)\dd\,\sInd
Y\leq (1,2)\dd\,\sInd(Y,X)\leq (1,2)\dd\,\sInd X$ and thus
$(1,2)\dd\,\sInd(X,X)=(1,2)\dd\,\sInd X$.
\end{theorem}

\begin{pf}
First, let us prove the right inequality. As usual, we shall use
the induction under $n=(1,2)\dd\,\sInd X$. Let the inequality be
proved for $n\leq k-1$ and prove it for $n=k$. Let
$F\in\co\tau_1'$ and $U(F)\in\tau_1$. Then $F\in\co\tau_1$ and
since $(1,2)\dd\,\sInd X=k$, there is $V(F)\in\tau_2$ such that
$$  \tau_2\cl V(F)\sbs U(F) \;\;\text{and}\;\;
        (1,2)\dd\,\sInd\big(2\dd\Fr_XV(F)\big)\leq k-1.  $$
Since $\tau_1<_{N(1)}\tau_2$, by (2) of Corollary~2.3.12 in [8],
$2\dd\Fr_XV(F)=         \linebreak       1\dd\Fr_XV(F)$. Clearly
$1\dd\Fr_XV(F)\cap Y\in\co\tau_1$ and by Corollary~2.3.10 in [8],
$\tau_1<_{N(1)}\tau_2$ implies that $\tau_1<_{C(1)}\tau_2$. Hence,
$1\dd\Fr_XV(F)\cap Y\in\co\tau_1$ implies that
$\tau_1''<_C\tau_2''$ in $(1\dd\Fr_XV(F)\cap Y,\tau_1''<\tau_2'')$
and since $1\dd\Fr_XV(F)\cap Y$ is $1$-closed in $2\dd\Fr_XV(F)$,
by the second part of Proposition~10.2,
\begin{eqnarray*}
    &\ds (1,2)\dd\,\sInd\big(1\dd\Fr_XV(F)\cap Y\big)=
        (1,2)\dd\,\sInd\big(2\dd\Fr_XV(F)\cap Y\big)\leq \\
    &\ds \leq (1,2)\dd\,\sInd\big(2\dd\Fr_XV(F)\big)\leq k-1.
\end{eqnarray*}
Therefore, by the condition in brackets in (2) of
Definition~10.11, $(1,2)\dd\,\sInd(Y,X)\leq k$ and so,
$(1,2)\dd\,\sInd(Y,X)\leq (2,1)\dd\,\sInd X$.

Now, we shall prove the left inequality by induction under a
nonnegative integer $n=(1,2)\dd\,\sInd(Y,X)$. Let the inequality
be proved for $n\leq k-1$ and prove it for $n=k$. Suppose that
$F\in\co\tau_1'$, $U'(F)\in\tau_1'$ and $U(F)\in\tau_1$ such that
$U(F)\cap Y=U'(F)$. Since $(1,2)\dd\,\sInd(Y,X)=k$, there is
$V(F)\in\tau_2$ such that
$$  \tau_2\cl V(F)\sbs U(F) \;\;\text{and}\;\;
        (1,2)\dd\,\sInd\big(2\dd\Fr_XV(F)\cap Y,X\big)\leq k-1.  $$
Hence, by inductive hypothesis $(1,2)\dd\,\sInd(2\dd\Fr_XV(F)\cap
Y)\leq k-1$. Since $\tau_1<_{N(1)}\tau_2$ and $Y\in\co\tau_1$, we
have
$$  2\dd\Fr_XV(F)=1\dd\Fr_XV(F) \;\;\text{and}\;\;
        2\dd\Fr_YV'(F)=1\dd\Fr_YV'(F),  $$
where $V'(F)=V(F)\cap Y$. But $1\dd\Fr_YV'(F)$ is $1$-closed in
$1\dd\Fr_XV(F)\cap Y$, and since $\tau_1<_{N(1)}\tau_2$ implies
$\tau_1<_{C(1)}\tau_2$, by the second part of
Pro\-po\-si\-ti\-on~10.2,
\begin{eqnarray*}
    &\ds (1,2)\dd\,\sInd\big(1\dd\Fr_YV'(F)\big)\leq
        (1,2)\dd\,\sInd\big(1\dd\Fr_XV(F)\cap Y\big)= \\
    &\ds =(1,2)\dd\,\sInd\big(2\dd\Fr_XV(F)\cap Y\big)\leq k-1.
\end{eqnarray*}
Hence $(1,2)\dd\,\sInd Y\leq k$ and so $(1,2)\dd\,\sInd Y\leq
(1,2)\dd\,\sInd(Y,X)$. Finally, for $\tau_1<_{N(1)}\tau_2$ it is
evident that $(1,2)\dd\,\sInd(X,X)=   \linebreak
  (1,2)\dd\,\sInd X$.
\end{pf}

\begin{definition}{10.18}
Let $(Y,\tau_1',\tau_2')$ be a $\BsS$ of a $\BS$
$(X,\tau_1,\tau_2)$. Then a pair $(A,B)$, where $A\in\co\tau_i'$,
$B\in\co\tau_i$ and $A\cap B=\vnth$, is said to be a separate
relatively $i$-normal pair. A separate relative partition,
corresponding to $(A,B)$, is a $j$-closed set $T\sbs X$ such that
$X\setminus T$ is not $j$-connected and
$$  (X\setminus T)\cap Y=Y\setminus T\sbs H_1\cup H_2,  $$
where $A\sbs H_i'=H_i\cap Y\in\tau_j'$, $B\sbs H_j\in\tau_j$ and
$H_1\cap H_2=\vnth$.
\end{definition}

\begin{remark}{10.19}
Let $(Y,\tau_1',\tau_2')$ be a $\BsS$ of a $\BS$
$(X,\tau_1,\tau_2)$ and $(A,B)$ be a separate relatively
$i$-normal pair. Then
\begin{enumerate}
\item[(1)] If there is a neighborhood $U(A)\in\tau_j$
$(U(B)\in\tau_j)$ such that $\tau_j\cl U(A)\sbs X\setminus B$
$(\tau_j\cl U(B)\sbs X\setminus A)$, then the set $j\dd\Fr_XU(A)$
$(j\dd\Fr_XU(B))$ is a separate relative partition, corresponding
to $(A,B)$.

\item[(2)] If $T$ is a separate relative partition, corresponding
to $(A,B)$, then $j\dd\Fr_XH_i\sbs T$.
\end{enumerate}
\end{remark}

Indeed, (1) Let us prove the case in brackets. We have:
$$  X\setminus j\dd\Fr_XU(B)=
        \big(X\setminus\tau_j\cl U(B)\big)\cup U(B),  $$
where $A\sbs H_i'=(X\setminus\tau_j\cl U(B))\cap Y\in\tau_j'$,
$B\sbs U(B)=H_j\in\tau_j$ and $H_1\cap H_2=\vnth$.

(2) Since $(X\setminus T,\tau_1'',\tau_2'')$ is not $j$-connected,
i.e., $X\setminus T=H_1\cup H_2$,
$H_i\in\tau_j\setminus\{\vnth\}$, $H_1\cap H_2=\vnth$, $A\sbs
H_i'=H_i\cap Y$ and $B\sbs H_j$, we have
$$  j\dd\Fr_XH_i\cap (X\setminus T)=
        (\tau_j\cl H_i\setminus H_i)\cap (X\setminus T)=\vnth.  $$
Thus $j\dd\Fr_XH_i\sbs T$.

\begin{proposition}{10.20}
Let $(Y,\tau_1',\tau_2')$ be a $\BsS$ of a $\BS$
$(X,\tau_1,\tau_2)$ and $n$ denote a nonnegative integer. If
$(i,j)\dd\,\sInd(Y,X)\leq n$, then to every separate relatively
$i$-normal pair $(A,B)$ there corresponds a separate relative
partition $T$ such that $(i,j)\dd\,\sInd(T\cap Y,X)\leq n-1$.
\end{proposition}

\begin{pf}
Let $(i,j)\dd\,\sInd(Y,X)\leq n$ and $(A,B)$ be a separate
relatively         \linebreak          $i$-normal pair. Then
$U(A)=X\setminus B\in\tau_i$ and since $(i,j)\dd\,\sInd(Y,X)\leq
n$, there is $V(A)\in\tau_j$ such that
$$  \tau_j\cl V(A)\sbs U(A) \;\;\text{and}\;\;
        (i,j)\dd\,\sInd\big(j\dd\Fr_XV(A)\cap Y,X\big)\leq n-1.   $$
Thus, it remains to use (1) of Remark~10.19.
\end{pf}

\begin{corollary}{10.21}
Let $(Y,\tau_1'<\tau_2')$ be a $\BsS$ of a $\BS$
$(X,\tau_1<\tau_2)$ and $n$ denote a nonnegative integer.  Then
$(2,1)\dd\,\sInd(Y,X)\leq n$ if and only if to every separate
relatively $2$-normal pair $(A,B)$ there corresponds a separate
relative partition $T$ such that $(2,1)\dd\,\sInd(T\cap Y,X)\leq
n-1$.
\end{corollary}

\begin{pf}
The proof of the first part is given by Proposition~10.20.

Conversely, let us suppose that the condition is satisfied,
$A\in\co\tau_2'$ and $U(A)\in\tau_2$. Then $(A,B=X\setminus U(A))$
is a separate relatively      \linebreak        $2$-normal pair
and, by condition, there is a separate relative partition $T$,
corresponding to $(A,B)$, such that $(2,1)\dd\,\sInd(T\cap
Y,X)\leq n-1$. But $Y\setminus T\sbs H_1\cup H_2$, where $A\sbs
H_i'=H_i\cap Y\in\tau_1'$, $B\sbs H_j\in\tau_1$ and $H_1\cap
H_2=\vnth$. Clearly,
$$  \tau_1\cl H_i\sbs X\setminus H_j\sbs X\setminus B=
                U(A)\in\tau_2.    $$
Let $H_i=V(A)$. Then $\tau_1\cl V(A)\sbs U(A)$ and by (2) of
Remark~10.19,
$$  1\dd\Fr_XV(A)\cap Y=1\dd\Fr_XH_i\cap Y\sbs T\cap Y. $$
Since $(2,1)\dd\,\sInd(T\cap Y,X)\leq n-1$ and $1\dd\Fr_XV(A)\cap
Y$ is $1$-closed and so, $2$-closed in $T\cap Y$, by
Proposition~10.13,
$$  (2,1)\dd\,\sInd\big(1\dd\Fr_XV(A)\cap Y,X)\leq
        (2,1)\dd\,\sInd(T\cap Y,X)\leq n-1        $$
so that $(2,1)\dd\,\sInd(Y,X)\leq n$.
\end{pf}

\begin{proposition}{10.22}
If a $\BsS$ $(Y,\tau_1',\tau_2')$ of a $\BS$ $(X,\tau_1,\tau_2)$
is          \linebreak        $(i,j)\dd\WS$-separately supernormal
in $X$ and $j$-extremally disconnected in $X$, then
$(i,j)\dd\,\sInd(Y,X)=0$.
\end{proposition}

\begin{pf}
Let $F\in\co\tau_i'$ and $U(F)\in\tau_i$. Then, by (2) of
Proposition~2.31, there is $V'(F)\in\tau_j'$ such that $\tau_j\cl
V'(F)\sbs U(x)$, where by topological version of Definition~2.39,
$\tau_j\cl V'(F)=V(F)\in\tau_j\cap\co\tau_j$, so that
$j\dd\Fr_XV(F)=\vnth$. Therefore, by (2) of Definition~10.11,
$(i,j)\dd\,\sInd(Y,X)=0$.
\end{pf}

\vskip+0.5cm
\section*{\textbf{11. $(i,j)$-Relative Covering Dimension Functions}}
\vskip+0.2cm

\begin{definition}{11.1}
Let $(Y,\tau_1',\tau_2')$ be a $\BsS$ of a $\BS$
$(X,\tau_1,\tau_2)$ and $n$ denote a nonnegative integer. Then
\begin{enumerate}
\item[(1)] $(i,j)\dd\dim(Y,X)=-1$ if and only if $Y=\vnth$.

\item[(2)] $(i,j)\dd\dim(Y,X)\leq n$ if every family
$\boldsymbol{\cU}=\{U_s\}_{s\in S}\sbs\tau_i$ for which there is a
family of sets $\boldsymbol{\boldsymbol{\cF}}=\{F_t\}_{t\in
T}\sbs\co\tau_j'$ which refines $\boldsymbol{\cU}$, can be refined
by a family $\boldsymbol{\cV}=\{V_r\}_{r\in R}\sbs\tau_i$ such
that the family $\boldsymbol{\boldsymbol{\cF}}$ also refines the
family $\boldsymbol{\cV}$ and $\ord\{(j,i)\dd\Fr_XV_r\cap
Y\}_{r\in R}\leq n$.

\item[(3)] $(i,j)\dd\dim(Y,X)=n$ if $(i,j)\dd\dim(Y,X)\leq n$ and \\
$(i,j)\dd\dim(Y,X)\leq n-1$ is false.

\item[(4)] $(i,j)\dd\dim(Y,X)=\infty$ if $(i,j)\dd\dim(Y,X)\leq n$
is false for all $n$.
\end{enumerate}
\end{definition}

As usual,
$$  p\,\dd\,\dim(Y,X)\!\leq\!n\!\llra\!
        \big((1,2)\dd\dim(Y,X)\!\leq\!n\wedge (2,1)\dd\dim(Y,X)\!\leq\!n\big). $$

Clearly, for the topological case we have

\begin{definition}{11.2}
Let $(Y,\tau')$ be a $\TsS$ of a $\TS$ $(X,\tau)$ and $n$ denote a
nonnegative integer. Then
\begin{enumerate}
\item[(1)] $\dim(Y,X)=-1$ if and only if $Y=\vnth$.

\item[(2)] $\dim(Y,X)\leq n$ if every family
$\boldsymbol{\cU}=\{U_s\}_{s\in S}\sbs\tau$ for which there is a
family of sets $\boldsymbol{\boldsymbol{\cF}}=\{F_t\}_{t\in
T}\sbs\co\tau'$ which refines $\boldsymbol{\cU}$, can be refined
by a family $\boldsymbol{\cV}=\{V_r\}_{r\in R}\sbs\tau$ such that
the family $\boldsymbol{\boldsymbol{\cF}}$ also refines
$\boldsymbol{\cV}$ and $\ord\{\Fr_XV_r\cap Y\}_{r\in R}\leq n$.

\item[(3)] $\dim(Y,X)=n$ if $\dim(Y,X)\leq n$ and $\dim(Y,X)\leq
n-1$ is false.

\item[(4)] $\dim(Y,X)=\infty$ if $\dim(Y,X)\leq n$ is false for
all $n$.
\end{enumerate}
\end{definition}

\begin{proposition}{11.3}
If $(Y,\tau_1''',\tau_2''')\sbs (Y_1,\tau_1'',\tau_2'')\sbs
(X_1,\tau_1',\tau_2')\sbs (X,\tau_1,\tau_2)$, $Y\in\co\tau_j''$
and $X_1\in\tau_i$, then
$$  (i,j)\dd\dim(Y,X)\leq (i,j)\dd\dim(Y_1,X_1).    $$
\end{proposition}

\begin{pf}
We shall use the induction under a nonnegative integer
$n=(i,j)\dd\dim(Y_1,X_1)$. Let us suppose that the inequality is
proved for $n\leq k-1$ and prove it for $n=k$. If
$\boldsymbol{\cU}=\{U_s\}_{s\in S}\sbs\tau_i$ is a family for
which there is a family of sets $\boldsymbol{\cF}=\{F_t\}_{t\in
T}\sbs\co\tau_j'''$ which refines $\boldsymbol{\cU}$, then it is
evident that $\boldsymbol{\cU}\cap X_1=\{U_s\cap X_1\}_{s\in
S}\sbs\tau_i'$ is a family of subsets of $X_1$ such that the
family $\boldsymbol{\cF}$ also refines $\boldsymbol{\cU}\cap X_1$.
It follows from $Y\in\co\tau_j''$ that
$\boldsymbol{\cF}\sbs\co\tau_j''$ and since
$(i,j)\dd\dim(Y_1,X_1)=k$, there is a family
$\boldsymbol{\cV}=\{V_r\}_{r\in R}\sbs\tau_i'$ such that
$\boldsymbol{\cF}$ refines $\boldsymbol{\cV}$, $\boldsymbol{\cV}$
refines $\boldsymbol{\cU}\cap X_1$ and
$\ord\{(j,i)\dd\Fr_{X_1}V_r\cap Y_1\}_{r\in R}\leq k$. Since
$Y\sbs Y_1\sbs X_1\sbs X$, it is clear, that
$$  (j,i)\dd\Fr_{X_1}V_r\cap Y_1=(j,i)\dd\Fr_XV_r\cap Y_1   $$
and
$$  (j,i)\dd\Fr_XV_r\cap Y\sbs (j,i)\dd\Fr_XV_r\cap Y_1 $$
for each $r\in R$. Besides, $\boldsymbol{\cV}$ also refines
$\boldsymbol{\cU}$, $X_1\in\tau_i$ implies that
$\boldsymbol{\cV}\sbs\tau_i$ and
\begin{eqnarray*}
    &\ds \ord\big\{(j,i)\dd\Fr_XV_r\cap Y\big\}_{r\in R}\leq \\
    &\ds \leq \ord\big\{(j,i)\dd\Fr_XV_r\cap Y_1\big\}_{r\in R}=
        \ord\big\{(j,i)\dd\Fr_{X_1}\cap Y_1\big\}_{r\in R}\leq k.
\end{eqnarray*}
Thus $(i,\!j)\dd\dim(Y,X)\!\leq\!k$ so that
$(i,\!j)\dd\dim(Y,X)\!\leq\!(i,\!j)\dd\dim(Y_1,X_1)$.~\end{pf}

\begin{corollary}{11.4}
If $(Y,\tau''')\sbs (Y_1,\tau'')\sbs (X_1,\tau')\sbs (X,\tau)$,
$Y\in\co\tau''$ and $X_1\in\tau$, then
$\dim(Y,X)\leq\dim(Y_1,X_1)$.
\end{corollary}

\begin{corollary}{11.5}
If $(Z,\tau_1'',\tau_2'')\sbs (Y,\tau_1',\tau_2')\sbs
(X,\tau_1,\tau_2)$ and $Y\in\tau_i$, then
$(i,j)\dd\,\dim(Z,X)\!\leq\!(i,j)\dd\,\dim(Z,Y)$ $($\textbf{the
anti-monotonicity of the relative bitopological covering dimension
functions}$)$.

If, in addition, $Z\!\in\!\co\tau_j'$,then
$(i,j)\dd\,\dim(Z,X)\!\leq\!(i,j)\dd\,\dim Y$, where
$(i,j)\dd\,\dim Y$ is considered in the sense of
Definition~$3.3.1$ in $[8]$.
\end{corollary}

\begin{pf}
For the first inequality let us suppose in Proposition~11.3 that
$Y=Y_1=Z$ and $X_1=Y$, and for the second inequality -- that $Y=Z$
and $Y_1=X_1=Y$.
\end{pf}

\begin{corollary}{11.6}
If $(Z,\tau'')\sbs (Y,\tau')\sbs (X,\tau)$ and $Y\in\tau$, then
$\dim(Z,X)\leq \dim(Z,Y)$ $($\textbf{the anti-monotonicity of the
relative topological covering dimension}$)$.

If, in addition, $Z\in\co\tau'$, then $\dim(Z,X)\leq \dim Y$.
\end{corollary}

\begin{theorem}{11.7}
The equalities $(i,j)\dd\,\dim(Y,X)=0$ and
   \linebreak         $(i,j)\dd\Ind(Y,X)=0$ are equivalent for every $\BsS$
$(Y,\tau_1',\tau_2')$ of a $\BS$ $(X,\tau_1,\tau_2)$ and either of
them yields the $(i,j)$-supernormality of $Y$ in $X$ together with
other notions of relative bitopological normality.
\end{theorem}

\begin{pf}
We begin by assuming that $(i,j)\dd\,\dim(Y,X)=0$ for a $\BsS$ $Y$
of a $\BS$ $X$. Let $A\!\in\!\co\tau_j'$ and $U(A)\!\in\!\tau_i$.
Since $(i,j)\dd\,\dim(Y,X)\!=\!0$, there is $V(A)\!\in\!\tau_i$
such that $F\!\sbs\!V(A)\!\sbs\!U(A)$ and
$\ord\big((j,i)\dd\Fr_XV(A)\cap Y\big)\!\leq~\!\!\!0$. Hence
$(j,i)\dd\Fr_XV(A)\cap Y\!=\!\vnth$ and by Remark~8.2,
$(i,j)\dd\Ind(Y,X)\!=~\!\!\!0$.

Conversely, let $(i,j)\dd\Ind(Y,X)=0$,
$\boldsymbol{\cU}=\{U_s\}_{s\in S}\sbs\tau_i$ and
$\boldsymbol{\cF}=\{F_t\}_{t\in T}\sbs\co\tau_j'$ be families such
that $F_{t(s)}\sbs U_s$ for each $s\in S$. Since
$(i,j)\dd\Ind(Y,X)=0$, for each $s\in S$ there $V_s\in\tau_i$ such
that $F_{t(s)}\sbs V_s\sbs\tau_j\cl V_s\sbs U_s$ (i.e., the family
$\boldsymbol{\cV}=\{V_s\}_{s\in S}$ refines $\boldsymbol{\cU}$ and
$\boldsymbol{\cF}$ refines $\boldsymbol{\cV})$ and
$(j,i)\dd\Fr_XV_s\cap Y=\vnth$. Hence,
$$  \ord\big\{(j,i)\dd\Fr_XV_s\cap Y\big\}_{s\in S}=-1  $$
so that $(i,j)\dd\,\dim(Y,X)\!\leq\!0$. But $Y\!\neq\!\vnth$ and
so $(i,j)\dd\,\dim(Y,X)\!=\!0$.~\end{pf} \vskip+0.2cm

For the rest see the notice between Definition~8.1 and Remark~8.2.

\begin{corollary}{11.8}
The equalities $\dim(Y,X)=0$ and $\Ind(Y,X)=0$ are equivalent for
every $\TsS$ $(Y,\tau')$ of a $\TS$ $(X,\tau)$ and either of them
yields the supernormality of $Y$ in $X$ together with other
notions of relative topological normality.
\end{corollary}

\begin{corollary}{11.9}
The following conditions are satisfied for any $\BsS$
$(Y,\tau_1',\tau_2')$ of a $\BS$ $(X,\tau_1,\tau_2)$:
\begin{enumerate}
\item[(1)] $p\,\dd\,\dim(Y,X)=0\llra p\,\dd\Ind(Y,X)$ and either
of them yields the $p$-supernormality of $Y$ in $X$ together with
other notions of relative bitopological ``absolute'' normality.

\item[(2)] If $(Y,\tau_1',\tau_2')$ is a $j\dd\TT_1$ $\BsS$ of a
$\BS$ $(X,\tau_1,\tau_2)$, then
$$  (i,j)\dd\,\dim(Y,X)\!=\!0\,\big(\!\llra\!(i,j)\dd\Ind(Y,X)\!=\!0\big)\!\lra\!
        (i,j)\dd\ind(Y,X)\!=\!0        $$
so that if $(Y,\tau_1',\tau_2')$ is a $\RR\dd\,p\,\dd\TT_1$ $\BsS$
of a $\BS$ $(X,\tau_1,\tau_2)$, then
$$  p\,\dd\,\dim(Y,X)=0\,\big(\llra p\,\dd\Ind(Y,X)=0\big)\lra
        p\,\dd\ind(Y,X)=0.        $$
\end{enumerate}
\end{corollary}

\begin{pf}
Follows directly from Theorem~11.7 and Proposition~8.4.
\end{pf}

\begin{corollary}{11.10}
If $(Y,\tau')$ is a $\TT_1$ $\TsS$ of a $\TS$ $(X,\tau)$, then
$$  \dim(Y,X)=0\,\big(\llra \Ind(Y,X)=0\big)\lra \ind(Y,X)=0.  $$
\end{corollary}

\begin{corollary}{11.11}
If $(Y,\tau_1',\tau_2')$ is a $(j,i)\dd\WS$-supernormal $\BsS$ of
a $\BS$ $(X,\tau_1,\tau_2)$ and $Y$ is $(i,j)$-extremally
disconnected in $X$, then $(i,j)\dd\,\dim(Y,X)=0$.
\end{corollary}

\begin{pf}
The condition is an immediate consequence of Theorem~11.7 and
Proposition~8.11.
\end{pf}

\begin{corollary}{11.12}
If $(Y,\tau')$ is $\WS$-supernormal $\TsS$ of a $\TS$ $(X,\tau)$
and $Y$ is extremally disconnected in $X$, then $\dim(Y,X)=0$.
\end{corollary}

In the last Section~12 we introduce and study relative versions of
Baire spaces for both the topological and the bitopological case.

\vskip+0.5cm
\section*{\textbf{12. Relative Baire Spaces and Relative $(i,j)$-Baire Spaces}}
\vskip+0.2cm

A $\TS$ $(X,\tau)$ is said to have the Baire property if any first
category subset of $X$ has an empty interior in $(X,\tau)$ and a
Baire space is a $\TS$ such that every nonempty open set is of
second category. A brief but cognitive story of the origin and
development of the theory of Baire spaces as well as of its
applications can be found in [15]. Most of the results on Baire
spaces are stated in [12]. The bitopological versions of Baire
spaces are introduced and studied, in particular, in [7], [8] and
[1].

\begin{definition}{12.1}
Let $(Y,\tau_1',\tau_2')$ be a $\BsS$ of a $\BS$
$(X,\tau_1,\tau_2)$. Then a subset $A\sbs X$ is said to be
$(i,j)$-nowhere dense $((i,j)$-somewhere dense) relative to $Y$ if
$\tau_i'\nt(\tau_j\cl A\cap Y)=\vnth$
$(\tau_i'\nt(\tau_j\cl A\cap Y)\neq\vnth)$.
\end{definition}

The families of all subsets of $X$ which are $(i,j)$-nowhere dense
$((i,j)$-somewhere dense) relative to $Y$ are denoted by
$(i,j)\dd\,\cN\cD(Y,X)$        \linebreak
$((i,j)\dd\,\cS\cD(Y,X))$.

The bitopological counterparts of these families are introduced
and studied in [8].

\begin{definition}{12.2}
Let $(Y,\tau')$ be a $\TsS$ of a $\TS$ $(X,\tau)$. Then a subset
$A\sbs X$ is said to be nowhere dense (somewhere dense) relative
to $Y$ if $\tau'\nt(\tau\cl A\cap Y)=\vnth$ $(\tau'\nt(\tau\cl
A\cap Y)\neq\vnth)$.
\end{definition}

The family of all subsets of $X$ which are nowhere dense
(somewhere dense) relative to $Y$ is denoted by $\cN\cD(Y,X)$
$(\cS\cD(Y,X))$.

Of considerable importance later on is

\begin{proposition}{12.3}
For a $\BsS$ $(Y,\tau_1',\tau_2')$ of a $\BS$ $(X,\tau_1,\tau_2)$
and a subset $A\sbs X$ the following conditions are satisfied:
\begin{enumerate}
\item[(1)] $A\in(i,j)\dd\,\cN\cD(Y,X)$ if and only if
$$  \tau_j\cl A\cap Y\sbs\tau_i'\cl\big(Y\setminus(\tau_j\cl A\cap Y)\big). $$

\item[(2)] If $A\in(i,j)\dd\,\cN\cD(Y,X)$, then for any set
$U'\in(\tau_1'\cap\tau_2')\setminus\{\vnth\}$ there is a set
$V\in\tau_j\setminus\{\vnth\}$ such that $V\cap Y\sbs U'$ and
$V\cap A=\vnth$.

\item[(3)] If for any set $U'\in\tau_i'\setminus\{\vnth\}$ there
is a set $V\in\tau_j\setminus\{\vnth\}$ such that $V\cap Y\sbs U'$
and $V\cap A=\vnth$, then $A\in(i,j)\dd\,\cN\cD(Y,X)$.
\end{enumerate}
\end{proposition}

\begin{pf}
(1) If $A\in(i,j)\dd\,\cN\cD(Y,X)$, then
$$  Y=Y\setminus\tau_i'\nt(\tau_j\cl A\cap Y)=
        \tau_i'\cl\big(Y\setminus(\tau_j\cl A\cap Y)\big)    $$
and, conversely, if
$$  \tau_j\cl A\cap Y\sbs\tau_i'\cl\big(Y\setminus(\tau_j\cl A\cap Y)\big)=
        Y\setminus \tau_i'\nt(\tau_j\cl A\cap Y),   $$
then
$$  \vnth=(\tau_j\cl A\cap Y)\cap \tau_i'\nt(\tau_j\cl A\cap Y)=
        \tau_i'\nt(\tau_j\cl A\cap Y).      $$

(2) If $A\in(i,j)\dd\,\cN\cD(Y,X)$ and
$U'\in(\tau_1'\cap\tau_2')\setminus\{\vnth\}$ is an arbitrary set,
then $U'\setminus(\tau_j\cl A\cap Y)\neq\vnth$, since the contrary
means that
$$  \vnth=U'\cap\big(X\setminus(\tau_j\cl A\cap Y)\big)=
        U'\cap\big(Y\setminus(\tau_j\cl A\cap Y)\big),  $$
i.e.,
$$  \vnth=U'\cap\tau_i'\cl\big(Y\setminus(\tau_j\cl A\cap Y)\big)=
        U'\setminus\vnth=U'.        $$
Let $V=U\setminus\tau_j\cl A$, where $U\cap Y=U'$, $U\in\tau_j$.
Then
$$  V\cap Y=(U\setminus\tau_j\cl A)\cap Y=
        U'\setminus(\tau_j\cl A\cap Y)\sbs U' $$
and $V\cap A=\vnth$.

(3) If $A\,\ol{\in}\,(i,j)\dd\,\cN\cD(Y,X)$, i.e., if
$\tau_i'\nt(\tau_j\cl A\cap Y)\neq\vnth$, then for any set
$V\in\tau_j\setminus\{\vnth\}$ such that $V\cap
Y\sbs\tau_i'\nt(\tau_j\cl A\cap Y)$, we have $V\cap
Y\sbs(\tau_j\cl A\cap Y)$ so that $V\cap\tau_j\cl A\neq\vnth$.
Hence $V\cap A\neq\vnth$.
\end{pf}

\begin{corollary}{12.4}
If $(Y,\tau_1'<\tau_2')$ is a $\BsS$ of a $\BS$
$(X,\tau_1<\tau_2)$, then the following conditions are satisfied:
\begin{enumerate}
\item[(1)] \ \ \ \\[-0.85cm]
$$  \begin{matrix}
        (2,1)\dd\,\cN\cD(Y,X) & \sbs & 1\dd\,\cN\cD(Y,X) \\
        \bigcap & & \bigcap \\
        2\dd\,\cN\cD(Y,X) & \sbs & (1,2)\dd\,\cN\cD(Y,X)
    \end{matrix}        $$
and hence
$$  \begin{matrix}
        (1,2)\dd\,\cS\cD(Y,X) & \sbs & 1\dd\,\cS\cD(Y,X) \\
        \bigcap & & \bigcap \\
        2\dd\,\cS\cD(Y,X) & \sbs & (2,1)\dd\,\cS\cD(Y,X).
    \end{matrix}        $$

\item[(2)] $A\in (1,2)\dd\,\cN\cD(Y,X)$ if and only if for each
set $U'\in\tau_1'\setminus\{\vnth\}$ there is a set
$V\in\tau_2\setminus\{\vnth\}$ such that $V\cap Y\sbs U'$ and
$V\cap A=\vnth$.

\item[(3)] If $A\in (2,1)\dd\,\cN\cD(Y,X)$, then for any set
$U'\in\tau_2'\setminus\{\vnth\}$ such that $\tau_1'\nt
U'\neq\vnth$ there is a set $V\in\tau_1\setminus\{\vnth\}$ such
that $V\cap Y\sbs U'$ and $V\cap A=\vnth$.

\item[(4)] If for any set $U'\in\tau_2'\setminus\{\vnth\}$ there
is a set $V\in\tau_1\setminus\{\vnth\}$ such that $V\cap Y\sbs U'$
and $V\cap A=\vnth$, then $A\in (2,1)\dd\,\cN\cD(Y,X)$.
\end{enumerate}
\end{corollary}

\begin{pf}
(1) The inclusions are evident.

(2) The condition is an immediate consequence of (2) and (3) of
Proposition~12.3.

(3) Let $U'\!\in\!\tau_2'\setminus\{\vnth\}$ and $\tau_1'\nt
U'\!\neq\!\vnth$. Then $\tau_1'\nt U'\setminus(\tau_1\cl A\cap
Y)\!\neq~\!\!\vnth$, since the contrary means that
\begin{eqnarray*}
    \vnth &&\hskip-0.6cm =\tau_1'\nt U'\cap\big(Y\setminus\tau_1'\nt
                (\tau_1\cl A\cap Y)\big)= \\
    &&\hskip-0.6cm =\tau_1'\nt U'\setminus\tau_1'\nt(\tau_1\cl A\cap Y)=
        \tau_1'\nt U'\setminus\vnth=\tau_1'\nt U',
\end{eqnarray*}
where the equality $\tau_1'\nt(\tau_1\cl A\cap Y)=\vnth$ is
conditioned by inclusion $(2,1)\dd\,\cN\cD(Y,X)\sbs
1\dd\,\cN\cD(Y,X)$. Now, it is clear that the set
$V=W\setminus\tau_1\cl (A\cap Y)$, where $W\in\tau_1$, $W\cap
Y=\tau_1'\nt U'$, is the required set.

(4) Follows directly from (3) of Proposition~12.3.
\end{pf}

\begin{corollary}{12.5}
If $(Y,\tau')$ is a $\TsS$ of a $\TS$ $(X,\tau)$, then
$A\in\cN\cD(Y,X)$ if and only if for any set
$U'\in\tau'\setminus\{\vnth\}$ there is a set
$V\in\tau\setminus\{\vnth\}$ such that $V\cap Y\sbs U'$ and $V\cap
A=\vnth$.
\end{corollary}

\begin{definition}{12.6}
Let $(Y,\tau_1',\tau_2')$ be a $\BsS$ of a $\BS$
$(X,\tau_1,\tau_2)$. Then a subset $A\sbs X$ is said to be
$(i,j)$-boundary $((i,j)$-dense) relative to $Y$ if
$\tau_i'\nt(\tau_j\nt A\cap Y)=\vnth$ $(\tau_i'\cl(\tau_j\cl A\cap
Y)=Y)$.
\end{definition}

The families of all subsets of $X$ which are $(i,j)$-boundary
$((i,j)$-dense) relative to $Y$ are denoted by
$(i,j)\dd\,\Bd(Y,X)$ $((i,j)\dd\,\cD(X,Y))$. Evidently
\begin{eqnarray*}
    &\displaystyle \tau_i'\nt(\tau_j\nt A\cap Y)=\vnth\llra
        Y=Y\setminus \tau_i'\nt(\tau_j\nt A\cap Y)= \\
    &\displaystyle =\tau_i'\cl\big(Y\setminus(\tau_j\nt A\cap Y)\big)=
        \tau_i'\cl\big((X\setminus\tau_j\nt A)\cap Y\big)= \\
    &\displaystyle =\tau_i'\cl\big(\tau_j\cl(X\setminus A)\cap Y\big)
\end{eqnarray*}
so that $A\in (i,j)\dd\,\Bd(Y,X)\llra X\setminus
A\in(i,j)\dd\,\cD(Y,X)$.

\begin{definition}{12.7}
Let $(Y,\tau')$ be a $\TsS$ of a $\TS$ $(X,\tau)$. Then a subset
$A\sbs X$ is said to be boundary (dense) relative to $Y$ if
$\tau\nt A\cap Y=\vnth$ $(\tau\cl A\cap Y=Y)$.
\end{definition}

The family of all subsets of $X$ which are boundary (dense)
relative to $Y$ is denoted by $\Bd(Y,X)$ $(\cD(Y,X))$. Clearly,
$A\in \Bd(Y,X)\llra X\setminus A\in\cD(Y,X)$.

\begin{remark}{12.8}
In contrast to the usual case, when $A\in\cN\cD(X)\llra\tau\cl
A\in\Bd(X)$, for the relative case $A\in\cN\cD(Y,X)\lra\tau\cl
A\in\Bd(Y,X)$, but not conversely as shows the following
elementary example: $X=\{a,b,c,d,e,f\}$,
$\tau=\{\vnth,\{a\},\{a,b\},\{c,d,e\},\{a,c,d,e\},
       \linebreak    \{a,b,c,d,e\},\{a,c,d,e,f\},X\}$,
$Y=\{f\}$ and $A=\{a,b\}$. Then $\tau\cl A\in\Bd(Y,X)$, but
$A\,\ol{\in}\,\cN\cD(Y,X)$.
\end{remark}

The following simple but useful notions are employed frequently
below.

\begin{definition}{12.9}
Let $(Y,\tau_1',\tau_2')$ be a $\BsS$ of a $\BS$
$(X,\tau_1,\tau_2)$. Then a subset $A\sbs X$ is $i$-strongly
boundary $(i$-strongly dense) relative to $Y$ if $\tau_i'\nt(A\cap
Y)=\vnth$ $(\tau_i'\cl(A\cap Y)=Y)$.
\end{definition}

The families of all subsets of $X$ which are $i$-strongly boundary
     \linebreak         $(i$-strongly dense) relative to $Y$ are denoted by
$i\dd\,\sBd(Y,\!X)$ $(i\dd\,\sD(Y,\!X))$.

Therefore,
\begin{eqnarray*}
    A\in i\dd\,\sBd(Y,X) &&\hskip-0.6cm \llra A\cap Y\in i\dd\,\Bd(Y), \\
    A\in i\dd\,\sD(Y,X) &&\hskip-0.6cm \llra A\cap Y\in i\dd\,\cD(Y), \\
    A\in (i,j)\dd\,\cN\cD(Y,X) &&\hskip-0.6cm \llra \tau_j\cl A\in i\dd\,\sBd(Y,X)
\end{eqnarray*}
and
$$  (i,j)\dd\,\cN\cD(Y,X)\sbs i\dd\,\sBd(Y,X)\sbs (i,j)\dd\,\Bd(Y,X).    $$

Hence, if $(Y,\tau')$ is a $\TsS$ of a $\TS$ $(X,\tau)$, then
\begin{eqnarray*}
    A\in\sBd(Y,X) &&\hskip-0.6cm \llra A\cap Y\in\Bd(Y), \\
    A\in\sD(Y,X) &&\hskip-0.6cm \llra A\cap Y\in\cD(Y), \\
    A\in\cN\cD(Y,X) &&\hskip-0.6cm \llra \tau\cl A\in\sBd(Y,X)
\end{eqnarray*}
and $\cN\cD(Y,X)\sbs \sBd(Y,X)\sbs\Bd(Y,X)$.

\begin{proposition}{12.10}
If $(Y,\tau_1'<\tau_2')$ is a $\BsS$ of a $\BS$
$(X,\tau_1<\tau_2)$, $A\in j\dd\,\sBd(Y,X)$ and
$B\in(i,j)\dd\,\cN\cD(Y,X)$, then $A\cup B\in 1\dd\,\sBd(Y,X)$.
\end{proposition}

\begin{pf}
Since $\tau_j'\cl(B\cap Y)\sbs\tau_j\cl B\cap Y$, we have
\begin{eqnarray*}
    &\ds Y\setminus(\tau_j\cl B\cap Y)=
        \tau_j'\cl(Y\setminus(A\cap Y))\setminus(\tau_j\cl B\cap Y)\sbs \\
    &\ds \sbs \tau_j'\cl(Y\setminus(A\cap Y))\setminus\tau_j'\cl(B\cap Y)\sbs
        \tau_j'\cl(Y\setminus(A\cup B)).
\end{eqnarray*}
Since $B\in(i,j)\dd\,\cN\cD(Y,X)$, by (2) of Lemma~0.2.1 in [8],
\begin{eqnarray*}
    Y &&\hskip-0.6cm =Y\setminus\tau_i'\nt(\tau_j\cl B\cap Y)= \\
    &&\hskip-0.6cm =\tau_i'\cl\big(Y\setminus(\tau_j\cl B\cap Y)\big)\sbs
        \tau_i'\cl\tau_j'\cl(Y\setminus(A\cup B))= \\
    &&\hskip-0.6cm =\tau_1'\cl(Y\setminus(A\cup B))=
        \tau_1'\cl\big((Y\setminus(A\cup B)\cap Y)\big).
\end{eqnarray*}
Thus $\tau_1'\nt((A\cup B)\cap Y)=\vnth$, i.e., $A\cup B\in
1\dd\,\sBd(Y,X)$.
\end{pf}

\begin{corollary}{12.11}
If $(Y,\tau')$ is a $\TsS$ of a $\TS$ $(X,\tau)$, $A\in\sBd(Y,X)$
and $B\in\cN\cD(Y,X)$, then $A\cup B\in\sBd(Y,X)$.
\end{corollary}

\begin{corollary}{12.12}
If $(Y,\tau_1'<\tau_2')$ is a $\BsS$ of a $\BS$
$(X,\tau_1<\tau_2)$ and $A\in2\dd\,\sBd(Y,X)$, then $A\cup B\in
1\dd\,\sBd(Y,X)$ if anyone of the following conditions is
satisfied:
\begin{enumerate}
\item[(1)] $B\in (2,1)\dd\,\cN\cD(Y,X)$.

\item[(2)] $B\in 1\dd\,\cN\cD(Y,X)$.

\item[(3)] $B\in 2\dd\,\cN\cD(Y,X)$.

\item[(4)] $B\in (1,2)\dd\,\cN\cD(Y,X)$.
\end{enumerate}
\end{corollary}

\begin{pf}
Follows directly from Proposition~12.10 taking into account (1) of
Corollary~12.4.
\end{pf}

\begin{corollary}{12.13}
For a $\BsS$ $(Y,\tau_1'<\tau_2')$ of a $\BS$ $(X,\tau_1<\tau_2)$
the following conditions are satisfied:
\begin{enumerate}
\item[(1)] $A,\,B\!\in\!i\dd\,\cN\cD(Y,X)\!\lra\!A\cup B\!\in\!
i\dd\,\cN\cD(Y,X)\!\sbs\!(1,2)\dd\,\cN\cD(Y,X)$.

\item[(2)] $A\in 1\dd\,\cN\cD(Y,X)$, $B\in
(2,1)\dd\,\cN\cD(Y,X)\lra A\cup B\in       \linebreak
1\dd\,\cN\cD(Y,X)\sbs (1,2)\dd\,\cN\cD(Y,X)$.

\item[(3)] $A\in 2\dd\,\cN\cD(Y,X)$, $B\in
(1,2)\dd\,\cN\cD(Y,X)\lra A\cup B\in \linebreak
      (1,2)\dd\,\cN\cD(Y,X)$.

\item[(4)] $A\in 2\dd\,\cN\cD(Y,X)$, $B\in
(2,1)\dd\,\cN\cD(Y,X)\lra A\cup B\in        \linebreak
2\dd\,\cN\cD(Y,X)\sbs (1,2)\dd\,\cN\cD(Y,X)$.

\item[(5)] $A,\,B\in (2,1)\dd\,\cN\cD(Y,X)\lra A\cup B\in
(2,1)\dd\,\cN\cD(Y,X)\sbs 1\dd\,\cN\cD(Y,X)\cap
2\dd\,\cN\cD(Y,X)\sbs (1,2)\dd\,\cN\cD(Y,X)$.

\item[(6)] $A\in (2,1)\dd\,\cN\cD(Y,X)$, $B\in
(1,2)\dd\,\cN\cD(Y,X)\lra A\cup B\in (1,2)\dd\,\cN\cD(Y,X)$.
\end{enumerate}
\end{corollary}

\begin{pf}
(1) $A\in i\dd\,\cN\cD(Y,X)\lra\tau_i\cl A\in i\dd\,\sBd(Y,X)$ and
$B\in i\dd\,\cN\cD(Y,X)\lra\tau_i\cl B\in i\dd\,\cN\cD(Y,X)$.
Hence, by Corollary~12.11, $\tau_i\cl A\cup\tau_i\cl
B=\tau_i\cl(A\cup B)\in i\dd\,\sBd(Y,X)$ so that $A\cup B\in
i\dd\,\cN\cD(Y,X)$. The inclusion follows from (1) of
Corollary~12.4.

(2) By (1) of Corollary~12.4, $(2,1)\dd\,\cN\cD(Y,X)\sbs
1\dd\,\cN\cD(Y,X)$ and hence, it remains to use (1) above for
$i=1$.

(3) $A\in 2\dd\,\cN\cD(Y,X)\lra \tau_2\cl A\in 2\dd\,\sBd(Y,X)$
and $B\in        \linebreak      (1,2)\dd\,\cN\cD(Y,X)\lra
\tau_2\cl B\in (1,2)\dd\,\cN\cD(Y,X)$. Therefore, by (4) it
Corollary~12.12, $\tau_2\cl A\cup\tau_2\cl B=\tau_2\cl(A\cup B)\in
1\dd\,\sBd(Y,X)$ so that $A\cup B\in (1,2)\dd\,\cN\cD(Y,X)$.

(4) By (1) of Corollary~12.4, $(2,1)\dd\,\cN\cD(Y,X)\sbs
2\dd\,\cN\cD(Y,X)$ and thus, it remains to use (1) above for
$i=2$.

(5) We have $\tau_2'\nt(\tau_1\cl A\cap
Y)=\vnth=\tau_2'\nt(\tau_1\cl B\cap Y)$. Then $\tau_1\cl A\in
2\dd\,\sBd(Y,X)$ and (2) of Lemma~0.2.1 in [8],
$\tau_2'\nt(\tau_2\cl\tau_1\cl B\cap Y)=\tau_2'\nt(\tau_1\cl B\cap
Y)=\vnth$ so that $\tau_1\cl B\in 2\dd\,\cN\cD(Y,X)$. Hence, by
Corollary~12.11, $\tau_1\cl A\cup\tau_1\cl B\in 2\dd\,\sBd(Y,X)$.
Since $\tau_1\cl A\cup\tau_1\cl B=\tau_1\cl(A\cup B)\in\co\tau_1$,
we have $A\cup B\in (2,1)\dd\,\cN\cD(Y,X)$.

The rest follows from (1) of Corollary~12.4.

(6) Follows directly from (3) above and (1) of Corollary~12.4.
\end{pf}

\begin{definition}{12.14}
Let $(Y,\tau_1',\tau_2')$ be a $\BsS$ of a $\BS$
$(X,\tau_1,\tau_2)$. Then a subset $A\sbs X$ is said to be an
$i\dd\,\cF_\sg$-set $(i\dd\,\cG_\dl$-set) relative to $Y$, if
$A=\bigcup\limits_{n=1}^\infty A_n$
$(A=\bigcap\limits_{n=1}^\infty A_n)$ and $A_n\cap Y\in\co\tau_i'$
$(A_n\cap Y\in\tau_i')$ for each $n=\ol{1,\infty}$.
\end{definition}

The families of all subsets of $X$ which are $i\dd\,\cF_\sg$
$(i\dd\,\cG_\dl)$ relative to $Y$, are denoted by
$i\dd\,\cF_\sg(Y,X)$ $(i\dd\,\cG_\dl(Y,X))$.

Therefore, $A\in i\dd\,\cF_\sg(Y,X)\llra A\cap Y\in
i\dd\,\cF_\sg(Y)$ and
$$  A\in i\dd\,\cG_\dl(Y,X)\llra A\cap Y\in i\dd\,\cG_\dl(Y).   $$

It is evident that $i\dd\,\cF_\sg(X)\sbs i\dd\,\cF_\sg(Y,X)$
$(i\dd\,\cG_\dl(X)\sbs i\dd\,\cG_\dl(Y,X))$ for any $Y\sbs X$ and
$$  A\in i\dd\,\cF_\sg(Y,X)\llra X\setminus A\in i\dd\,\cG_\dl(Y,X).  $$

The topological case is obvious and the notions above are of
independent interest which is not our objective at the present
time. Note only here, that by analogy with Lemma~2.21 one can
prove that if $(Y,\tau_1',\tau_2')\sbs(X,\tau_1,\tau_2)$, $Y$ is
$p$-supernormal in $X$, $P\in 1\dd\,\cF_\sg(Y,X)$, $Q\in
2\dd\,\cF_\sg(Y,X)$ and $(\tau_1\cl P\cap Q)\cup (P\cap \tau_2\cl
Q)=\vnth$, then there are disjoint sets $U\in\tau_2$, $V\in\tau_1$
such that $P\cap Y\sbs U$, $Q\cap Y\sbs V$.

\begin{definition}{12.15}
Let $(Y,\tau_1',\tau_2')$ be a $\BsS$ of a $\BS$
$(X,\tau_1,\tau_2)$. Then a subset $A\sbs X$ is said to be of
$(i,j)$-first category (also called $(i,j)$-meager,
$(i,j)$-exhaustible) relative to $Y$ if
$A=\bigcup\limits_{n=1}^\infty A_n$, where $A_n\in
(i,j)\dd\,\cN\cD(Y,X)$ for each $n=\ol{1,\infty}$, and $A$ is said
to be of $(i,j)$-second category (also called $(i,j)$-nonmeager,
$(i,j)$-inexhaustible) relative to $Y$ if it is not of
$(i,j)$-first category relative to $Y$.
\end{definition}

A subset $A\sbs X$ is said to be of $(i,j)$-first category in
itself relative to $Y$ if $A=\bigcup\limits_{n=1}^\infty A_n$,
where $A_n\in(i,j)\dd\,\cN\cD(Y\cap A,A)$ for each
$n=\ol{1,\infty}$ and hence, $A$ is of $(i,j)$-second category in
itself relative to $Y$ if it is not of $(i,j)$-first category in
itself relative to $Y$.

The families of all subsets of $X$ which are of $(i,j)$-first
$((i,j)$-second) categories relative to $Y$, are denoted by
$(i,j)\dd\,\Catg_{{}_{\I}}(Y,X)$         \linebreak
$((i,j)\dd\,\Catg_{{}_{\II}}(Y,X))$, while
$$  A\in (i,j)\dd\,\Catg_{{}_{\I}}(Y\cap A,A) \;\;\;
        \big(\,A\in(i,j)\dd\,\Catg_{{}_{\II}}(Y\cap A,A)\,\big)     $$
is abbreviated to $A$ is of $(i,j)\dd\,\Catg_{{}_{Y}}\I$
$((i,j)\dd\,\Catg_{{}_{Y}}\II)$.

Clearly
$$  (i,j)\dd\,\cN\cD(Y,X)\sbs
        (i,j)\dd\,\Catg_{{}_{\I}}(Y,X)=2^X\setminus(i,j)\dd\,\Catg_{{}_{\II}}(Y,X).  $$
Now, if $(Y,\tau')$ is a $\TsS$ of a $\TS$ $(X,\tau)$ and $A\sbs
X$, then $A\in \Catg_{{}_{\I}}(Y,X)\llra
A=\bigcup\limits_{n=1}^\infty A_n$, where $A_n\in \cN\cD(Y,X)$ for
each $n=\ol{1,\infty}$ and $\Catg_{{}_{\II}}(Y,X)=2^X\setminus
\Catg_{{}_{\I}}(Y,X)$. Therefore, $A\in\Catg_{{}_{\I}}(Y\cap A,A)$
$(A\in\Catg_{{}_{\II}}(Y\cap A,A))$ is abbreviated to $A$ is of
$\Catg_{{}_{Y}}\I$ $(\Catg_{{}_{Y}}\II)$.

\begin{theorem}{12.16}
For a $\BsS$ $(Y,\tau_1',\tau_2')$ of a $\BS$ $(X,\tau_1,\tau_2)$
the following conditions are satisfied:
\begin{enumerate}
\item[(1)] The families $(i,j)\dd\,\Catg_{{}_{\I}}(Y,X)$ are
$\sg$-ideals so that if $A_n\in
  \linebreak       (i,j)\dd\,\Catg_{{}_{\I}}(Y,X)$ for each $n=\ol{1,\infty}$, then
$$  \bigcup\limits_{n=1}^\infty A_n\in (i,j)\dd\,\Catg_{{}_{\I}}(Y,X),  $$
and if $B\!\in\!(i,j)\dd\,\Catg_{{}_{\I}}(Y,X)$, $A\!\sbs\!B$,
then $A\!\in\!(i,j)\dd\,\Catg_{{}_{\I}}(Y,X)$.

\item[(2)] If $A=\bigcup\limits_{n=1}^\infty A_n\in
j\dd\,\cF_\sg(X)$ and $A_n\in i\dd\,\sBd(Y,X)$ for each
$n=\ol{1,\infty}$, then $A\in (i,j)\dd\,\Catg_{{}_{\I}}(Y,X)$.

\item[(3)] For every set $A\in (i,j)\dd\,\Catg_{{}_{\I}}(Y,X)$
there is a set $B\in j\dd\,\cF_\sg(X)\cap
(i,j)\dd\,\Catg_{{}_{\I}}(Y,X)$ such that $A\sbs B$.

\item[(4)] The families $(i,j)\dd\,\Catg_{{}_{\II}}(Y,X)$ are
closed under arbitrary unions and $A\in
(i,j)\dd\,\Catg_{{}_{\II}}(Y,X)$, $A\sbs B$ imply that $B\in
\linebreak        (i,j)\dd\,\Catg_{{}_{\II}}(Y,X)$.

\item[(5)] If $X$ is of $(i,j)\dd\,\Catg_Y\II$ and for any subset
$A\sbs X$ there is a set $B\sbs A$ such that
$B=\bigcap\limits_{n=1}^\infty B_n\in j\dd\,\cG_\dl(X)$ and
$B_n\in i\dd\,\sD(Y,X)$ for each $n=\ol{1,\infty}$, then $A\in
(i,j)\dd\,\Catg_{{}_{\II}}(Y,X)$.

\item[(6)] $X$ is of $(i,j)\dd\,\Catg_{{}_Y}\II$ if and only if
the intersection of any sequence $\{A_n\}_{n=1}^\infty$, where
$A_n\in\tau_j\cap i\dd\,\sD(Y,X)$ for each $n=\ol{1,\infty}$, is
nonempty.
\end{enumerate}

Moreover, for a $\BsS$ $(Y,\tau_1'<\tau_2')$ of a $\BS$
$(X,\tau_1<\tau_2)$ we have
\begin{enumerate}
\item[(7)] \ \ \ \\[-0.85cm]
$$  \begin{matrix}
        (2,1)\dd\,\Catg_{{}_{\I}}(Y,X) & \sbs & 1\dd\,\Catg_{{}_{\I}}(Y,X) \\
        \bigcap & & \bigcap \\
        2\dd\,\Catg_{{}_{\I}}(Y,X) & \sbs & (1,2)\dd\,\Catg_{{}_{\I}}(Y,X)
    \end{matrix}        $$
and hence,
$$  \begin{matrix}
        (1,2)\dd\,\Catg_{{}_{\II}}(Y,X) & \sbs & 1\dd\,\Catg_{{}_{\II}}(Y,X) \\
        \bigcap & & \bigcap \\
        2\dd\,\Catg_{{}_{\II}}(Y,X) & \sbs & (2,1)\dd\,\Catg_{{}_{\II}}(Y,X)
    \end{matrix}        $$
\end{enumerate}
\end{theorem}

\begin{pf}
The conditions (1) and (4) are obvious, (7) follows directly from
(1) of Corollary~12.4.

(2) Let $A=\bigcup\limits_{n=1}^\infty A_n\in j\dd\,\cF_\sg(X)$,
i.e., $A_n\in\co\tau_j$ for each $n=\ol{1,\infty}$. By condition
$\tau_i'\nt(A_n\cap Y)=\vnth$ so that $\tau_i'\nt(\tau_j\cl
A_n\cap Y)=\vnth$ and thus, $A_n\in(i,j)\dd\,\cN\cD(Y,X)$ for each
$n=\ol{1,\infty}$. Hence $A\in (i,j)\dd\,\Catg_{{}_{\I}}(Y,X)$.

(3) Let $A\in (i,j)\dd\,\Catg_{{}_{\I}}(Y,X)$. Then
$A=\bigcup\limits_{n=1}^\infty A_n$, where $A_n\in
(i,j)\dd\,\cN\cD(Y,X)$ for each $n=\ol{1,\infty}$. But $A_n\in
(i,j)\dd\,\cN\cD(Y,X)$ implies that $\tau_j\cl A_n\in
(i,j)\dd\,\cN\cD(Y,X)$ for each $n=\ol{1,\infty}$. Let
$B=\bigcup\limits_{n=1}^\infty \tau_j\cl A_n$. Then $A\sbs B$ and
$B\in j\dd\,\cF_\sg(X)\cap (i,j)\dd\,\Catg_{{}_{\I}}(Y,X)$.

(5) Let $A\sbs X$ be any subset and $B\sbs A$, where
$B=\bigcap\limits_{n=1}^\infty B_n\in j\dd\,\cG_\dl(X)$, i.e.,
$B_n\in\tau_j$ for each $n=\ol{1,\infty}$. By condition $B_n\cap
Y\in i\dd\,\cD(Y)$ so that $F_n=Y\setminus B_n\in i\dd\,\Bd(Y)$
for each $n=\ol{1,\infty}$. Since $X\setminus
B=\bigcup\limits_{n=1}^\infty (X\setminus B_n)\in
j\dd\,\cF_\sg(X)$ and $(X\setminus B_n)\cap Y=Y\setminus B_n$, we
have
$$  \vnth=\tau_i'\nt(Y\setminus B_n)=\tau_i'\nt((X\setminus B_n)\cap Y)=
        \tau_i'\nt(\tau_j\cl (X\setminus B_n)\cap Y)     $$
and so $(X\setminus B_n)\in (i,j)\dd\,\cN\cD(Y,X)$ for each
$n=\ol{1,\infty}$. Hence $X\setminus B\in
(i,j)\dd\,\Catg_{{}_{\I}}(Y,X)$. But $B\sbs A$ implies $X\setminus
A\sbs X\setminus B$ and by (1), $X\setminus A\in
(i,j)\dd\,\Catg_{{}_{\I}}(Y,X)$. Then $A\in
(i,j)\dd\,\Catg_{{}_{\II}}(Y,X)$ since, by virtue of (1), the
contrary means that $X$ is of $\Catg_{{}_Y}\I$.

(6) First, let $\{A_n\}_{n=1}^\infty$ be a sequence of subsets of
$X$ such that $A_n\in\tau_j\cap i\dd\,\sD(Y,X)$ for each
$n=\ol{1,\infty}$ and $\bigcap\limits_{n=1}^\infty A_n=\vnth$.
Then
$$  X=X\setminus \bigcap\limits_{n=1}^\infty A_n=
        \bigcup\limits_{n=1}^\infty (X\setminus A_n),       $$
where $X\setminus A_n\in\co\tau_j\cap i\dd\,\sBd(Y,X)$ for each
$n=\ol{1,\infty}$. Hence, by (2), $X$ is of
$(i,j)\dd\,\Catg_{{}_Y}\I$.

Conversely, if $X$ is of $(i,j)\dd\,\Catg_{{}_Y}\I$, then by (3),
$X\in j\dd\,\cF_\sg(X)$, i.e., $X=\bigcup\limits_{n=1}^\infty
F_n$, where $F_n\in\co\tau_j\cap i\dd\,\sBd(Y,X)$ for each
$n=\ol{1,\infty}$. Let us consider the sequence $\{X\setminus
F_n\}_{n=1}^\infty$, where $X\setminus F_n\in\tau_j\cap
i\dd\,\sD(Y,X)$ for each $n=\ol{1,\infty}$. It is clear that
$\bigcap\limits_{n=1}^\infty (X\setminus F_n)=\vnth$.
\end{pf}

\begin{corollary}{12.17}
For a $\TS$ $(Y,\tau')$ of a $\TS$ $(X,\tau)$ the following
conditions are satisfied:
\begin{enumerate}
\item[(1)] The family $\Catg_{{}_{\I}}(Y,X)$ is a $\sg$-ideal.

\item[(2)] If $A=\bigcup\limits_{n=1}^\infty A_n\in\FF_\sg(X)$ and
$A_n\in\sBd(Y,X)$ for each $n=\ol{1,\infty}$, then $A\in
\Catg_{{}_{\I}}(Y,X)$.

\item[(3)] For every set $A\in \Catg_{{}_{\I}}(Y,X)$ there is a
set $B\in\FF_\sg(X)\cap \Catg_{{}_{\I}}(Y,X)$ such that $A\sbs B$.

\item[(4)] The family $\Catg_{{}_{\II}}(Y,X)$ is closed under
arbitrary unions and $A\in \Catg_{{}_{\II}}(Y,X)$, $A\sbs B$ imply
that  $B\in \Catg_{{}_{\II}}(Y,X)$.

\item[(5)] If $X$ is of $\Catg_{{}_Y}\II$ and for any subset
$A\sbs X$ there is a set $B\sbs A$ such that
$B=\bigcap\limits_{n=1}^\infty B_n\in\cG_\dl(X)$ and
$B_n\in\sD(Y,X)$ for each $n=\ol{1,\infty}$, then $A\in
\Catg_{{}_{\II}}(Y,X)$.

\item[(6)] $X$ is of $\Catg_{{}_Y}\II$ if and only if the
intersection of any sequence $\{A_n\}_{n=1}^\infty$, where
$A_n\in\tau\cap \sD(Y,X)$ for each $n=\ol{1,\infty}$, is nonempty.
\end{enumerate}
\end{corollary}

\begin{corollary}{12.18}
If $(Y,\tau_1',\tau_2')$ is a $\BsS$ of a $\BS$
$(X,\tau_1,\tau_2)$, $X$ is of $(i,j)\dd\,\Catg_{{}_Y}\II$,
$A=\bigcap\limits_{n=1}^\infty A_n\in j\dd\,\cG_\dl(X)$ and
$A_n\in i\dd\,\sD(Y,X)$ for each $n=\ol{1,\infty}$, then $A\in
(i,j)\dd\,\Catg_{{}_{\II}}(Y,X)$.
\end{corollary}

\begin{pf}
Follows directly from (5) of Theorem~12.16 for $B=A$.
\end{pf}

\begin{corollary}{12.19}
If $(Y,\tau')$ is a $\TsS$ of a $\TS$ $(X,\tau)$, $X$ is of
$\Catg_{{}_Y}\II$, $A=\bigcap\limits_{n=1}^\infty A_n\in
\cG_\dl(X)$ and $A_n\in \sD(Y,X)$ for each $n=\ol{1,\infty}$, then
$A\in \Catg_{{}_{\II}}(Y,X)$.
\end{corollary}

\begin{corollary}{12.20}
If $(Y,\tau_1<\tau_2')$ is a $\BsS$ of a $\BS$ $(X,\tau_1<\tau_2)$
and $X$ is of $(1,2)\dd\,\Catg_{{}_Y}\II$, then for
$A=\bigcap\limits_{n=1}^\infty A_n$ we have: {\small
$$  \xymatrix{
        A_n\!\in\!\tau_1\cap 2\dd\,\sD(Y,X) \ar@{=>}[r] \ar@{=>}[d] \!&\!
            A_n\!\in\!\tau_2\cap 2\dd\,\sD(Y,X) \ar@{=>}[d] & \\
        A_n\!\!\in\!\tau_1\!\cap\!1\dd\,\sD(Y\!,\!X)\!\ar@{=>}[r] \!&\!
            A_n\!\!\in\!\tau_2\!\cap\!1\dd\,\sD(Y\!,\!X)\!\ar@{=>}[r] \!&\! A\!\!\in\!
                (1,2)\dd\,\Catg_{{}_{\II}}(Y\!,\!X). }        $$ }
\end{corollary}

\begin{pf}
Follows directly from Corollary~12.18 and the inclusion
   \linebreak $\tau_1\sbs\tau_2$.
\end{pf}

\begin{definition}{12.21}
A $\BsS$ $(Y,\tau_1',\tau_2')$ of a $\BS$ $(X,\tau_1,\tau_2)$ is
$(i,j)$-quasi regular in $X$ if for any set
$U'\in\tau_i'\setminus\{\vnth\}$ there is a set
$V\in\tau_i\setminus\{\vnth\}$ such that $\tau_j'\cl(V\cap Y)\sbs
U'$.
\end{definition}

\begin{theorem}{12.22}
If a $\BsS$ $(Y,\tau_1<\tau_2')$ of a $\BS$ $(X,\tau_1<\tau_2)$ is
$(2,1)$-quasi regular in $X$ and $1$-compact, then
$(2,1)\dd\,\Catg_{{}_{\I}}(Y,X)\sbs 2\dd\,\sBd(Y,X)$.
\end{theorem}

\begin{pf}
Let $A\in (2,1)\dd\,\Catg_{{}_{\I}}(Y,X)$ be any set. Then
$A=\bigcup\limits_{n=1}^\infty A_n$, where $A_n\in
(2,1)\dd\,\cN\cD(Y,X)(\llra \tau_2'\nt(\tau_1\cl A_n\cap
Y)=\vnth)$ for each $n=\ol{1,\infty}$. We shall prove that
$\tau_2'\nt(A\cap Y)=\vnth$. Clearly, it suffices to prove that
for any set $V\in\tau_2'\setminus\{\vnth\}$ we have $V\cap
(Y\setminus A)\neq\vnth$. Suppose that
$V\in\tau_2'\setminus\{\vnth\}$ is an arbitrary set. Then
$V\setminus(\tau_1\cl A_1\cap Y)\neq\vnth$ since the contrary
means that $V\sbs \tau_1\cl A_1\cap Y$ which is impossible since
$\tau_2'\nt(\tau_1\cl A_1\cap Y)=\vnth$. By condition, there is
$U_1\in\tau_2\setminus\{\vnth\}$ such that
$$  \tau_1'\cl(U_1\cap Y)\sbs V\setminus(\tau_1\cl A_1\cap Y)   $$
so that $\tau_1'\cl(U_1\cap Y)\sbs V$ and $\tau_1'\cl(U_1\cap
Y)\cap A_1=\vnth$. By analogy, there is
$U_2\in\tau_2\setminus\{\vnth\}$ such that
$$  \tau_1'\cl(U_2\cap Y)\sbs U_1\cap Y\sbs \tau_1'\cl(U_1\cap Y)   $$
and $\tau_1'\cl(U_2\cap Y)\cap A_2=\vnth$. Therefore, one can
construct a sequence
$\{U_n\}_{n=1}^\infty\sbs\tau_2\setminus\{\vnth\}$ such that
$\tau_1'\cl(U_n\cap Y)\sbs U_{n-1}$ and $\tau_1'\cl(U_n\cap Y)\cap
A_n=\vnth$. Now, $Y$ is $1$-compact implies that
$\bigcap\limits_{n=1}^\infty \tau_1'\cl(U_n\cap Y)\neq\vnth$,
i.e., there is a point $x\in\tau_1'\cl(U_n\cap Y)$ for each
$n=\ol{1,\infty}$. Therefore, $x\in
V\setminus\bigcup\limits_{n=1}^\infty A_n=V\setminus A$ so that
$V\cap (X\setminus A)\neq\vnth$ and since $V\sbs Y$, $V\cap
(Y\setminus A)\neq\vnth$.
\end{pf}

\begin{corollary}{12.23}
If a $\TsS$ $(Y,\tau')$ of a $\TS$ $(X,\tau)$ is quasi regular in
$X$ and compact, then $\Catg_{{}_{\I}}(Y,X)\sbs \sBd(Y,X)$.
\end{corollary}

\begin{theorem}{12.24}
For a $\BsS$ $(Y,\tau_1',\tau_2')$ of a $\BS$ $(X,\tau_1,\tau_2)$
the following conditions are equivalent:
\begin{enumerate}
\item[(1)] $U\in\tau_i\setminus\{\vnth\}\lra
U\in(i,j)\dd\,\Catg_{{}_{\II}}(Y,X)$.

\item[(2)] If $\{U_n\}_{n=1}^\infty$ is a countable family of
subsets of $X$ such that $U_n\in\tau_j\cap i\dd\,\sD(Y,X)$ for
each $n=\ol{1,\infty}$, then $\bigcap\limits_{n=1}^\infty U_n\in
i\dd\,\cD(X)$.

\item[(3)] $A\in(i,j)\dd\,\Catg_{{}_{\I}}(Y,X)\lra X\setminus A\in
i\dd\,\cD(X)$.

\item[(4)] If $\{F_n\}_{n=1}^\infty$ is a countable family of
subsets of $X$ such that $F_n\in\co\tau_j\cap i\dd\,\sBd(Y,X)$ for
each $n=\ol{1,\infty}$, then $\bigcup\limits_{n=1}^\infty F_n\in
i\dd\,\Bd(X)$.
\end{enumerate}
\end{theorem}

\begin{pf}
$(1)\lra (2)$. Let $\{U_n\}_{n=1}^\infty\sbs 2^X$, where $U_n\in
\tau_j\cap i\dd\,\sD(Y,X)$ for each $n=\ol{1,\infty}$. Then
$X\setminus U_n\in\co\tau_j\cap i\dd\sBd(Y,X)$ for each
$n=\ol{1,\infty}$ and so
$$  V=X\setminus\bigcap\limits_{n=1}^\infty U_n=
            \bigcup\limits_{n=1}^\infty (X\setminus U_n)\in
        (i,j)\dd\,\Catg_{{}_{\I}}(Y,X)      $$
as $X\setminus U_n\in (i,j)\dd\,\cN\cD(Y,X)$ for each
$n=\ol{1,\infty}$. If $\bigcap\limits_{n=1}^\infty
U_n\,\ol{\in}\,i\dd\,\cD(X)$, then there is
$U\in\tau_i\setminus\{\vnth\}$ such that $U\cap
\bigcap\limits_{n=1}^\infty U_n=\vnth$.

Hence $U\sbs X\setminus\bigcap\limits_{n=1}^\infty U_n=V$ and
since $V\in(i,j)\dd\,\Catg_{{}_{\I}}(Y,X)$, by (1) of Theorem
12.16, $U\in(i,j)\dd\,\Catg_{{}_{\I}}(Y,X)$.

$(2)\lra (3)$. Let $A\in(i,j)\dd\,\Catg_{{}_{\I}}(Y,X)$. Then, by
(3) of Theorem~12.16, there is a set $B\in j\dd\,\cF_\sg(X)\cap
(i,j)\dd\,\Catg_{{}_{\I}}(Y,X)$ such that $A\sbs B$. Therefore
$B=\bigcup\limits_{n=1}^\infty B_n$, where $B_n\in\co\tau_j\cap
i\dd\,\sBd(Y,X)$ for each $n=\ol{1,\infty}$, and $X\setminus
B=\bigcap\limits_{n=1}^\infty (X\setminus B_n)$. Clearly
$X\setminus B_n\in\tau_j\cap i\dd\,\sD(Y,X)$ for each
$n=\ol{1,\infty}$ and by (2), $X\setminus B\in i\dd\,\cD(X)$. Now,
it is evident that $X\setminus A\in i\dd\,\cD(X)$.

$(3)\lra (4)$. If $F=\bigcup\limits_{n=1}^\infty F_n$, where
$F_n\in\co\tau_j\cap i\dd\,\sBd(Y,X)$ for each $n=\ol{1,\infty}$,
then $F_n\in (i,j)\dd\,\cN\cD(Y,X)$ for each $n=\ol{1,\infty}$ and
so $F\in (i,j)\dd\,\Catg_{{}_{\I}}(Y,X)$. Hence, by (3),
$X\setminus F\in i\dd\,\cD(X)$ and thus $F\in i\dd\,\Bd(X)$.

$(4)\lra (1)$. Let $U\in(\tau_i\setminus\{\vnth\})\cap
(i,j)\dd\,\Catg_{{}_{\I}}(Y,X)$. Then
$U=\bigcup\limits_{n=1}^\infty U_n$, where
$U_n\in(i,j)\dd\,\cN\cD(Y,X)$ for each $n=\ol{1,\infty}$.
Therefore, $\tau_j\cl U_n\in (i,j)\dd\,\cN\cD(Y,X)$ for each
$n=\ol{1,\infty}$ and so $\bigcup\limits_{n=1}^\infty \tau_j\cl
U_n\in (i,j)\dd\,\Catg_{{}_{\I}}(Y,X)$. But $\{\tau_j\cl
U_n\}_{n=1}^\infty$ is a countable family, where       \linebreak
$\tau_j\cl U_n\in\co\tau_j\cap i\dd\,\sBd(Y,X)$ for each
$n=\ol{1,\infty}$. Therefore, by (4),
$\tau_i\nt\bigcup\limits_{n=1}^\infty \tau_j\cl U_n=\vnth$ and
hence,
$$  U=\bigcup\limits_{n=1}^\infty U_n\sbs
        \tau_i\nt\bigcup\limits_{n=1}^\infty \tau_j\cl U_n      $$
implies that $U=\vnth$. A contradiction.
\end{pf}

\begin{corollary}{12.25}
For a $\TsS$ $(Y,\tau')$ of a $\TS$ $(X,\tau)$ the following
conditions are equivalent:
\begin{enumerate}
\item[(1)] $U\in\tau\setminus\{\vnth\}\lra
U\in\Catg_{{}_{\II}}(Y,X)$.

\item[(2)] If $\{U_n\}_{n=1}^\infty$ is a countable family of
subsets of $X$ such that $U_n\in\tau\cap s\dd\,\cD(Y,X)$ for each
$n=\ol{1,\infty}$, then $\bigcap\limits_{n=1}^\infty
U_n\in\cD(X)$.

\item[(3)] $A\in \Catg_{{}_{\I}}(Y,X)\lra X\setminus A\in\cD(X)$.

\item[(4)] If $\{F_n\}_{n=1}^\infty$ is a countable family of
subsets of $X$ such that $F_n\in\co\tau\cap s\dd\,\Bd(Y,X)$ for
each $n=\ol{1,\infty}$, then $\bigcup\limits_{n=1}^\infty
F_n\in\Bd(X)$.
\end{enumerate}
\end{corollary}

\begin{definition}{12.26}
A $\BsS$ $(Y,\tau_1',\tau_2')$ $(\TsS$ $(Y,\tau'))$ of a $\BS$
$(X,\tau_1,\tau_2)$ $(\TS$ $(X,\tau))$ is a almost $(i,j)$-Baire
space (almost Baire space) in $X$ (briefly, $Y$ is an
$A\dd\,(i,j)\dd\BrS$ in $X$ $(Y$ is an $A\dd\BrS$ in $X))$ if
anyone of the equivalent conditions (1)--(4) of Theorem~12.24
(Corollary~12.25) is satisfied.
\end{definition}

\begin{definition}{12.27}
A $\BsS$ $(Y,\tau_1',\tau_2')$ $(\TsS$ $(Y,\tau'))$ of a $\BS$
$(X,\tau_1,\tau_2)$ $(\TS$ $(X,\tau))$ is an $(i,j)$-Baire space
(Baire space) in $X$ (briefly, $Y$ is an $(i,j)\dd\BrS$ in $X$
$(Y$ is a $\BrS$ in $X))$ if $U\in\tau_i\setminus\{\vnth\}$
$(U\in\tau\setminus\{\vnth\})$ implies that $U$ is of
$(i,j)\dd\,\Catg_{{}_Y}\II$ $(\Catg_{{}_Y}\II)$.
\end{definition}

\begin{lemma}{12.28}
If $(Y,\tau_1',\tau_2')$ is a $\BsS$ of a $\BS$
$(X,\tau_1,\tau_2)$ and $U\in\tau_i\setminus\{\vnth\}$, then
\begin{enumerate}
\item[(1)] $U$ is of $(i,j)\dd\,\Catg_{{}_Y}\II$ implies that
$U\in (i,j)\dd\,\Catg_{{}_{\II}}(Y,X)$.
\end{enumerate}

If, in addition, $\tau_1\sbs\tau_2$, $Y\in 2\dd\,\cD(X)$ and
$U\in\tau_1\setminus\{\vnth\}$, then
\begin{enumerate}
\item[(2)] $U$ is of $(1,2)\dd\,\Catg_{{}_Y}\II\llra U\in
(1,2)\dd\,\Catg_{{}_{\II}}(Y,X)$.
\end{enumerate}
\end{lemma}

\begin{pf}
(1) Let $U$ be of $(i,j)\dd\,\Catg_{{}_Y}\II$ and $U\in
(i,j)\dd\,\Catg_{{}_{\I}}(Y,X)$. Then
$U=\bigcup\limits_{n=1}^\infty U_n$, where $\tau_i'\nt(\tau_j\cl
U_n\cap Y)=\vnth$, i.e., $U_n\in (i,j)\dd\,\cN\cD(Y,X)$ for each
$n=\ol{1,\infty}$. Since $(Y\cap U,\tau_1''',\tau_2''')\sbs
(U,\tau_1'',\tau_2'')\sbs (X,\tau_1,\tau_2)$, we have
$\tau_i'''\nt(\tau_j''\cl U_n\cap(U\cap Y))=\vnth$, i.e.,
$U_n\in(i,j)\dd\,\cN\cD(Y\cap U,U)$ for each $n=\ol{1,\infty}$ so
that $U$ is of $(i,j)\dd\,\Catg_{{}_Y}\I$.

(2) By (1), it suffices to prove the implication from the right to
the left, i.e., $U$ is of $(1,2)\dd\,\Catg_{{}_Y}\I$ implies that
$U\in (1,2)\dd\,\Catg_{{}_{\I}}(Y,X)$, where $U$ is of
$(1,2)\dd\,\Catg_{{}_Y}\I\llra U\in
(1,2)\dd\,\Catg_{{}_{\I}}(Y\cap U,U)$.

Since $U\in (1,2)\dd\,\Catg_{{}_{\I}}(Y\cap U,U)$, we have
$U=\bigcup\limits_{n=1}^\infty A_n$, where $A_n\in
(1,2)\dd\,\cN\cD(Y\cap U,U)$ for each $n=\ol{1,\infty}$.

Let $V'\in\tau_1'\setminus\{\vnth\}$ be any set. If $V'\cap
U=\vnth$, i.e., $V'\cap(Y\cap U)=\vnth$, then
$V'\cap\tau_2'\cl(Y\cap U)=\vnth$ as $V'\in\tau_1'\sbs\tau_2'$, so
that $V'\cap\tau_2\cl(Y\cap U)=\vnth$. Let $V=W\setminus\tau_2\cl
U$, where $W\in\tau_1\sbs\tau_2$ and $W\cap Y=V'$. Since $Y\in
2\dd\,\cD(X)$, we have $\tau_2\cl U=\tau_2\cl(U\cap Y)$ and so
$$  V\cap Y=(W\setminus \tau_2\cl U)\cap Y=V'\setminus \tau_2\cl U=
        V'\setminus \tau_2\cl(U\cap Y)=V'.       $$
Since $V\cap \tau_2\cl U=\vnth$, we have $V\cap A_n=\vnth$ and by
(3) of Proposition~12.3, $A_n\in (1,2)\dd\,\cN\cD(Y,X)$ for each
$n=\ol{1,\infty}$. Thus $U\in (1,2)\dd\,\Catg_{{}_{\I}}(Y,X)$.

Now, let $V\in\tau_1'\setminus\{\vnth\}$ and $V'\cap
U=V'\cap(Y\cap U)\neq\vnth$. Then $V'\cap
U\in\tau_1'''\setminus\{\vnth\}$ and since $A_n\in
(1,2)\dd\,\cN\cD(Y\cap U,U)$, there is
$U'\in\tau_2''\setminus\{\vnth\}$ such that $U'\cap (Y\cap
U)=U'\cap Y\sbs V'$ and $U'\cap A_n=\vnth$ for each
$n=\ol{1,\infty}$. But $U'\in\tau_2''$ and $U\in\tau_2$ imply that
$U'\in\tau_2\setminus\{\vnth\}$. Therefore, for
$V'\in\tau_1'\setminus\{\vnth\}$ there is
$U'\in\tau_2\setminus\{\vnth\}$ such that $U'\cap Y\sbs V'$ and
$U'\cap A_n=\vnth$ for each $n=\ol{1,\infty}$. Once more applying
(3) of Proposition~12.3 gives that $A_n\in (1,2)\dd\,\cN\cD(Y,X)$
for each $n=\ol{1,\infty}$ and thus $U\in
(1,2)\dd\,\Catg_{{}_{\I}}(Y,X)$.
\end{pf}

\begin{corollary}{12.29}
Let $(Y,\tau_1',\tau_2')$ be a $\BsS$ of a $\BS$
$(X,\tau_1,\tau_2)$. Then
\begin{enumerate}
\item[(1)] If $Y$ is an $(i,j)\dd\BrS$ in $X$, then $Y$ is an
$A\dd\,(i,j)\dd\BrS$ in $X$.
\end{enumerate}

Moreover, if $\tau_1\sbs \tau_2$ and $Y\in2\dd\,\cD(X)$, then
\begin{enumerate}
\item[(2)] $Y$ is a $(1,2)\dd\BrS$ in $X$ if and only if $Y$ is an
$A\dd\,(1,2)\dd\BrS$ in $X$.
\end{enumerate}
\end{corollary}

\begin{corollary}{12.30}
If $(Y,\tau')$ be a $\TsS$ of a $\TS$ $(X,\tau)$ and $Y\in\cD(X)$,
then
\begin{enumerate}
\item[(1)] If $U\in\tau\setminus\{\vnth\}$, then $U$ is of
$\Catg_{{}_Y}\II$ if and only if $U\in\Catg_{{}_{\II}}(Y,X)$.

\item[(2)] $Y$ is a $\BrS$ in $X$ if and only if $Y$ is an
$A\dd\BrS$ in $X$.
\end{enumerate}
\end{corollary}

\begin{proposition}{12.31}
For a $\BsS$ $(Y,\tau_1'<\tau_2')$ of a $\BS$ $(X,\tau_1<\tau_2)$
the following conditions are satisfied:
\begin{enumerate}
\item[(1)] If $Y\in 2\dd\,\cD(X)$ and $Y$ is a $(1,2)\dd\BrS$ in
$X$, then $Y$ is a $1\dd\BrS$ in $X$.

\item[(2)] If $Y$ is a $2\dd\BrS$ in $X$, then $Y$ is a
$(2,1)\dd\BrS$ in $X$.
\end{enumerate}
\end{proposition}

\begin{pf}
(1) Let $U\in\tau_1\setminus\{\vnth\}$ be any set. Then $U$ is of
$(1,2)\dd\,\Catg_{{}_Y}\II$ and by (2) of Lemma~12.28,
$U\in(1,2)\dd\,\Catg_{{}_{\II}}(Y,X)$. Hence, following (7) of
Theorem~12.16, $U\in 1\dd\,\Catg_{{}_{\II}}(Y,X)$. Hence, by the
topological version of (2) of Lemma~12.28, $U\in
1\dd\,\Catg_{{}_Y}\II$.

(2) If $U\in\tau_2\setminus\{\vnth\}$ is any set, then by the
topological part of Definition~12.27, $U$ is of
$2\dd\,\Catg_{{}_Y}\II$, i.e., $U\in 2\dd\,\Catg_{{}_{\II}}(Y\cap
U,U)$. But, by (7) of Theorem~12.16, $2\dd\,\Catg_{{}_{\II}}(Y\cap
U,U)\sbs (2,1)\dd\,\Catg_{{}_{\II}}(Y\cap U,U)$ and so $U$ is of
$(2,1)\dd\,\Catg_{{}_Y}\II$. Thus, according to Definition~12.27,
$Y$ is a $(2,1)\dd\BrS$ in $X$.
\end{pf}

\begin{proposition}{12.32}
If $(Z,\tau_1''<\tau_2'')\sbs (Y,\tau_1'<\tau_2')\sbs
(X,\tau_1<\tau_2)$, then the following conditions are satisfied:
\begin{enumerate}
\item[(1)] $Y\in\tau_1$ and $Z$ is a $(1,2)\dd\BrS$ in $X$ imply
that $Z$ is also a $(1,2)\dd\BrS$ in $Y$.

\item[(2)] $Y\in 2\dd\,\cD(X)$ and $Z$ is a $(1,2)\dd\BrS$ in $Y$
imply that $Z$ is also a $(1,2)\dd\BrS$ in $X$.
\end{enumerate}
\end{proposition}

\begin{pf}
(1) If $U\in\tau_1'\setminus\{\vnth\}$, then
$U\in\tau_1\setminus\{\vnth\}$ and so $U$ is of
$(1,2)\dd\,\Catg_{{}_Z}\II$.

(2) Contrary: there is $U\in\tau_1\setminus\{\vnth\}$ such that
$U$ is of $(1,2)\dd\,\Catg_{{}_Z}\I$. Hence, by (2) of
Lemma~12.28, $U\in (1,2)\dd\,\Catg_{{}_{\I}}(Z,X)$ and so
$U=\bigcup\limits_{n=1}^\infty A_n$, where
$A_n\in(1,2)\dd\,\cN\cD(Z,X)$, i.e., $\tau_1''\nt(\tau_2\cl
A_n\cap Z)=\vnth$ for each $n=\ol{1,\infty}$. Therefore,
$\tau_1''\nt(\tau_2'\cl(A_n\cap Y)\cap Z)=\vnth$ for each
$n=\ol{1,\infty}$ so that $U\cap Y\in
(\tau_1'\setminus\{\vnth\})\cap (1,2)\dd\,\Catg_{{}_{\I}}(Z,Y)$.
But this is impossible, since by (2) of Corollary~12.29, $Z$ is a
$(1,2)\dd\BrS$ in $Y$ if and only if $Z$ is an
$A\dd\,(1,2)\dd\BrS$ in $Y$.
\end{pf}

\begin{corollary}{12.33}
If $(Z,\tau)\sbs (Y,\tau)\sbs (X,\tau)$, then the following
conditions are satisfied:
\begin{enumerate}
\item[(1)] $Y\in\tau$ and $Z$ is a $\BrS$ in $X$ imply that $Z$ is
also a $\BrS$ in $Y$.

\item[(2)] $Y\in\cD(X)$ and $Z$ is a $\BrS$ in $Y$ imply that $Z$
is also a $\BrS$ in~$X$.
\end{enumerate}
\end{corollary}

It is of interest and very important for further investigations to
study various types of relative properties of relative notions as
for the topological case so for the bitopological case.

\begin{proposition}{12.34}
If $(Z,\tau_1'',\tau_2'')\sbs (Y,\tau_1',\tau_2')\sbs
(X,\tau_1,\tau_2)$, $Z\in\tau_i'$ and $A\sbs X$, then the
following conditions are satisfied:
\begin{enumerate}
\item[(1)] $A\in (i,j)\dd\,\cN\cD(Y,X)$ $(A\in
(i,j)\dd\,\Catg_{{}_{\I}}(Y,X))$ implies that $A\cap Y\in
(i,j)\dd\,\cN\cD(Z,Y)$ $(A\cap Y\in
(i,j)\dd\,\Catg_{{}_{\I}}(Z,Y))$ and so $A\cap Y\in
(i,j)\dd\,\cS\cD(Z,Y)$ $(A\cap Y\in
(i,j)\dd\,\Catg_{{}_{\II}}(Z,Y))$ implies that $A\in
(i,j)\dd\,\cS\cD(Y,X)$ $(A\in (i,j)\dd\,\Catg_{{}_{\II}}(Y,X))$.

\item[(2)] $A\in (i,j)\dd\,\cN\cD(Y,X)$ $(A\in
(i,j)\dd\,\Catg_{{}_{\I}}(Y,X))$ implies that $A\in
(i,j)\dd\,\cN\cD(Z,X)$ $(A\in (i,j)\dd\,\Catg_{{}_{\I}}(Z,X))$ and
so $A\in       \linebreak        (i,j)\dd\,\cS\cD(Z,X)$ $(A\in
(i,j)\dd\,\Catg_{{}_{\II}}(Z,X))$ implies that $A\in \linebreak
(i,j)\dd\,\cS\cD(Y,X)$ $(A\in (i,j)\dd\,\Catg_{{}_{\II}}(Y,X))$.
\end{enumerate}
\end{proposition}

\begin{pf}
It is evident that in both cases it suffices to prove only the
first implication.

(1) Using (1) of Proposition~12.3, let us prove that if
$$  \tau_j\cl A\cap Y\sbs\tau_i'\cl\big(Y\setminus(\tau_j\cl A\cap Y)\big),  $$
then
$$  \tau_j'\cl(A\cap Y)\cap Z\sbs
        \tau_i''\cl\big(Z\setminus(\tau_j'\cl(A\cap Y)\cap Z)\big).  $$
Since $Z\in\tau_i'$, we have
\begin{eqnarray*}
    &\ds \tau_j'\cl(A\cap Y)\cap Z\sbs \tau_j\cl A\cap Y\cap Z\sbs
        \tau_i'\cl\big(Y\setminus \tau_j'\cl(A\cap Y)\big)\cap Z\sbs \\
    &\ds \sbs \tau_i'\cl\big(\big(Y\setminus \tau_j'\cl(A\cap Y)\big)\cap Z\big)=
        \tau_i'\cl\big(Z\setminus\tau_j'\cl(A\cap Y)\big).
\end{eqnarray*}
Thus
\begin{eqnarray*}
    &\ds \tau_j'\cl(A\cap Y)\cap Z\sbs
        \tau_i'\cl\big(Z\setminus\tau_j'\cl(A\cap Y)\big)\cap Z= \\
    &\ds =\tau_i''\cl\big(Z\setminus(\tau_j'\cl(A\cap Y)\cap Z)\big).
\end{eqnarray*}

(2) By analogy with (1), we use (1) of Proposition~12.3. Since
$Z\in\tau_i'$ and $A\in (i,j)\dd\,\cN\cD(Y,X)$ we have
\begin{eqnarray*}
    &\ds (\tau_j\cl A\cap Y)\cap Z=\tau_j\cl A\cap Z\sbs
        \tau_i'\cl\big(Y\setminus(\tau_j\cl A\cap Y)\big)\cap Z\sbs \\
    &\ds \sbs \tau_i'\cl\big(Z\setminus(\tau_j\cl A\cap Z)\big),
\end{eqnarray*}
so that
$$  \tau_j\cl A\cap Z\sbs
        \tau_i'\cl\big(Z\setminus(\tau_j\cl A\cap Z)\big)\cap Z=
            \tau_i''\cl\big(Z\setminus(\tau_j\cl A\cap Z)\big). $$
\end{pf}

\begin{corollary}{12.35}
If $(Z,\tau_1''<\tau_2'')\sbs (Y,\tau_1'<\tau_2')\sbs
(X,\tau_1<\tau_2)$, then the conditions (1) and (2) of
Proposition~12.34 are satisfied for $Z\in\tau_1'$.
\end{corollary}

\begin{corollary}{12.36}
If $(Z,\tau_1'',\tau_2'')\sbs (Y,\tau_1',\tau_2')\sbs
(X,\tau_1,\tau_2)$, $Z\in\tau_i'$ and $Y\in i\dd\,\cD(X)$, then
$Z$ is an $A\dd\,(i,j)\dd\BrS$ in $Y$ implies that $Y$ is an
$A\dd\,(i,j)\dd\BrS$ in $X$.
\end{corollary}

\begin{pf}
Let $U\in\tau_i\setminus\{\vnth\}$ be any set. Since $Y\in
i\dd\,\cD(X)$, we have $U'=U\cap Y\in\tau_i'\setminus\{\vnth\}$
and so $U'\in(i,j)\dd\,\Catg_{{}_{\II}}(Z,Y)$ as $Z$ is an
$A\dd\,(i,j)\dd\BrS$ in $Y$. Since $Z\in\tau_i'$, by (1) of
Proposition~12.34, $U\in(i,j)\dd\,\Catg_{{}_{\II}}(Y,X)$ and it
remains ro recall Definition~12.26.
\end{pf}

\begin{corollary}{12.37}
If $(Z,\tau_1''<\tau_2'')\sbs (Y,\tau_1'<\tau_2')\sbs
(X,\tau_1<\tau_2)$ and $Y\in 2\dd\,\cD(X)$, then
\begin{enumerate}
\item[(1)] $Z\in\tau_1'$ and $Z$ is a $(1,2)\dd\BrS$ in $Y$ imply
that $Y$ is a $(1,2)\dd\BrS$ in $X$.

\item[(2)] $Z\in\tau_2'$ and $Z$ is an $A\dd\,(2,1)\dd\BrS$ in $Y$
imply that $Y$ is an $A\dd\,(2,1)\dd\BrS$ in $X$.
\end{enumerate}
\end{corollary}

\begin{pf}
(1) If $Z$ is a $(1,2)\dd\BrS$ in $Y$, then by (1) of
Corollary~12.29, $Z$ is an $A\dd\,(1,2)\dd\BrS$ in $Y$. Hence, by
Corollary~12.36, $Y$ is an $A\dd\,(1,2)\dd\BrS$ in $X$ and so, it
remains to use (2) of Corollary~12.29.

(2) Follows directly from Corollary~12.36.
\end{pf}

\begin{corollary}{12.38}
If $(Z,\tau'')\sbs (Y,\tau')\sbs (X,\tau)$, $Z\in\tau'$ and $A\sbs
X$, then the following conditions are satisfied:
\begin{enumerate}
\item[(1)] $A\!\in\!\cN\cD(Y,X)$ $(A\!\in\!\Catg_{{}_{\I}}(Y,X))$
implies that $A\cap Y\!\in\!\cN\cD(Z,Y)$ $(A\cap
Y\in\Catg_{{}_{\I}}(Z,Y))$ and so $A\cap Y\in\cS\cD(Z,Y)$ $(A\cap
Y\in\Catg_{{}_{\II}}(Z,Y))$ implies that $A\in\cS\cD(Y,X)$
$(A\in\Catg_{{}_{\II}}(Y,X))$.

\item[(2)] $A\in\cN\cD(Y,X)$ $(A\in\Catg_{{}_{\I}}(Y,X))$ implies
that $A\in\cN\cD(Z,X)$ $(A\in\Catg_{{}_{\I}}(Z,X))$ and so
$A\in\cS\cD(Z,X)$ $(A\in\Catg_{{}_{\II}}(Z,X))$ implies that
$A\in\cS\cD(Y,X)$ $(A\in\Catg_{{}_{\II}}(Y,X))$.
\end{enumerate}
\end{corollary}

\begin{corollary}{12.39}
If $(Z,\tau'')\sbs (Y,\tau')\sbs (X,\tau)$, $Z\in\tau'$ and
$Y\in\cD(X)$, then $Z$ is a $\BrS$ in $Y$ implies that $Y$ is a
$\BrS$ in $X$.
\end{corollary}

\begin{proposition}{12.40}
If $(Z,\tau_1''<\tau_2'')\sbs (Y,\tau_1'<\tau_2')\sbs
(X,\tau_1<\tau_2)$, $Z\in 2\dd\,\cD(Y)$ and $A\sbs X$, then $A\in
(2,1)\dd\,\cN\cD(Z,X)\llra A\in (2,1)\dd\,\cN\cD(Y,X)$ $(A\in
(2,1)\dd\,\Catg_{{}_{\I}}(Z,X)\llra A\in
(2,1)\dd\,\Catg_{{}_{\I}}(Y,X))$ and so
\begin{eqnarray*}
    A\in (2,1)\dd\,\cS\cD(Z,X) &&\hskip-0.6cm \llra
            A\in (2,1)\dd\,\cS\cD(Y,X) \\
    \big(\,A\in (2,1)\dd\,\Catg_{{}_{\II}}(Z,X) &&\hskip-0.6cm \llra
        A\in (2,1)\dd\,\Catg_{{}_{\II}}(Y,X)\,\big).
\end{eqnarray*}
\end{proposition}

\begin{pf}
Evidently, it suffices to prove only the first equivalence.

Let $A\in (2,1)\dd\,\cN\cD(Z,X)$ and
$A\,\ol{\in}\,(2,1)\dd\,\cN\cD(Y,X)$. Then          \linebreak
$\tau_2'\nt(\tau_1\cl A\cap Y)\neq\vnth$ and $Z\in 2\dd\,\cD(Y)$
imply that
$$  \vnth\neq\tau_2'\nt(\tau_1\cl A\cap Y)\cap Z\sbs
        \tau_2''\nt(\tau_1\cl A\cap Y)      $$
which is impossible.

Conversely, let $A\in (2,1)\dd\,\cN\cD(Y,X)$ and
$A\,\ol{\in}\,(2,1)\dd\,\cN\cD(Z,X)$. Then $\tau_2''\nt(\tau_1\cl
A\cap Y)\neq\vnth$ and so there is $U\in\tau_2'\setminus\{\vnth\}$
such that $\vnth\neq U\cap Z=\tau_2''\nt(\tau_1\cl A\cap Y)$.
Since $Z\in 2\dd\,\cD(Y)$, we have
$$  \tau_2'\cl U=\tau_2'\cl(U\cap Z)=
        \tau_2'\cl\tau_2''\nt(\tau_1\cl A\cap Y).       $$
Hence
$$  \tau_1'\cl\tau_2'\cl U=
        \tau_1'\cl \tau_2'\cl\tau_2''\nt(\tau_1\cl A\cap Y)    $$
and by (2) of Lemma~0.2.1 in [8], $\tau_1'\cl
U=\tau_1'\cl\tau_2''\nt(\tau_1\cl A\cap Y)$. Therefore,
$U\in\tau_2'\setminus\{\vnth\}$ gives that
\begin{eqnarray*}
    &\ds \vnth\neq\tau_2'\nt\tau_1'\cl U=
        \tau_2'\nt\tau_1'\cl\tau_2''\nt(\tau_1\cl A\cap Y)\sbs \\
    &\ds \sbs \tau_2'\nt\tau_1'\cl(\tau_1\cl A\cap Y)=
        \tau_2'\nt(\tau_1\cl A\cap Y),
\end{eqnarray*}
i.e., $A\,\ol{\in}\,(2,1)\dd\,\cN\cD(Y,X)$ which is a
contradiction.
\end{pf}

\begin{corollary}{12.41}
If $(Z,\tau_1''<\tau_2'')\sbs (Y,\tau_1'<\tau_2')\sbs
(X,\tau_1<\tau_2)$ and $Z\in 2\dd\,\cD(Y)$, then $Z$ is an
$A\dd\,(2,1)\dd\BrS$ in $X$ if and only if $Y$ is an
$A\dd\,(2,1)\dd\BrS$ in $X$.
\end{corollary}

\begin{pf}
Indeed, if $U\in\tau_2\setminus\{\vnth\}$ is any set, then by
Proposition~12.40, $U\in (2,1)\dd\,\Catg_{{}_{\II}}(Z,X)$ if and
only if $U\in (2,1)\dd\,\Catg_{{}_{\II}}(Y,X)$.
\end{pf}

\begin{corollary}{12.42}
If $(Z,\tau'')\sbs (Y,\tau')\sbs (X,\tau)$, $Z\in\cD(Y)$ and
$A\sbs X$, then $A\in\cN\cD(Z,X)\llra A\in\cN\cD(Y,X)$
$(A\in\Catg_{{}_{\I}}(Z,X)\llra A\in\Catg_{{}_{\I}}(Y,X))$ and so
$A\in\cS\cD(Z,X)\llra A\in\cS\cD(Y,X)$
$(A\in\Catg_{{}_{\II}}(Z,X)\llra A\in\Catg_{{}_{\II}}(Y,X))$.
\end{corollary}

\begin{corollary}{12.43}
If $(Z,\tau'')\sbs (Y,\tau')\sbs (X,\tau)$ and $Z\in\cD(Y)$, then
$Z$ is a $\BrS$ in $X$ if and only if $Y$ is a $\BrS$ in $X$.
\end{corollary}

\begin{proposition}{12.44}
If $(Y,\tau_1',\tau_2')\sbs (X,\tau_1,\tau_2)$ and $Y\in\tau_i$,
then $(i,j)\dd\,\cN\cD(X)\sbs (i,j)\dd\,\cN\cD(Y,X)$
$((i,j)\dd\,\Catg_{{}_{\I}}(X)\sbs
(i,j)\dd\,\Catg_{{}_{\I}}(Y,X))$ and so
\begin{eqnarray*}
    (i,j)\dd\,\cS\cD(Y,X) &&\hskip-0.6cm \sbs (i,j)\dd\,\cS\cD(X) \\
    \big(\,(i,j)\dd\,\Catg_{{}_{\II}}(Y,X) &&\hskip-0.6cm \sbs
        (i,j)\dd\,\Catg_{{}_{\II}}(X)\,\big).
\end{eqnarray*}
\end{proposition}

\begin{pf}
If $A\in (i,j)\dd\,\cN\cD(X)$, i.e., if $\tau_i\nt\tau_j\cl
A=\vnth$, then $\tau_i\nt(\tau_j\cl A\cap Y)=\vnth$ so that
$\tau_i'\nt(\tau_j\cl A\cap Y)=\vnth$ as $Y\in\tau_i$.
\end{pf}

\begin{corollary}{12.45}
If $(Y,\tau_1',\tau_2')\sbs (X,\tau_1,\tau_2)$ and $Y\in\tau_i$,
then $Y$ is an $A\dd\,(i,j)\dd\BrS$ in $X$ implies that $X$ is an
$A\dd\,(i,j)\dd\BrS$.
\end{corollary}

\begin{pf}
If $U\in\tau_i\setminus\{\vnth\}$ is an arbitrary set, then
$U\in(i,j)\dd\,\Catg_{{}_{\II}}(Y,X)$ and so
$U\in(i,j)\dd\,\Catg_{{}_{\II}}(X)$. Hence, it remains  to use
Definition~4.1.5 in [8].
\end{pf}

\begin{corollary}{12.46}
If $(Y,\tau_1'<\tau_2')\sbs (X,\tau_1<\tau_2)$ and $Y\in\tau_1$,
then the conditions of Proposition~$12.44$ are satisfied.
\end{corollary}

\begin{corollary}{12.47}
If $(Y,\tau_1'<\tau_2')\sbs (X,\tau_1<\tau_2)$ and $Y\in\tau_1$,
then
\begin{enumerate}
\item[(1)] $Y\in 2\dd\,\cD(X)$ and $Y$ is a $(1,2)\dd\BrS$ in $X$
imply that $X$ is a $(1,2)\dd\BrS$.

\item[(2)] $Y$ is an $A\dd\,(2,1)\dd\BrS$ in $X$ implies that $X$
is an $A\dd\,(2,1)\dd\BrS$.
\end{enumerate}
\end{corollary}

\begin{pf}
(1) It remains to use (2) of Lemma~12.28 and (4) of Theorem~4.1.6
in [8].

(2) Follows directly from Corollary~12.45.
\end{pf}

\begin{corollary}{12.48}
If $(Y,\tau')\sbs (X,\tau)$ and $Y\in\tau$, then $\cN\cD(X)\sbs
\cN\cD(Y,X)$ $(\Catg_{{}_{\I}}(X)\sbs\Catg_{{}_{\I}}(Y,X))$ and so
$\cS\cD(Y,X)\sbs \cS\cD(X)$
$(\Catg_{{}_{\II}}(Y,X)\sbs\Catg_{{}_{\II}}(X))$.
\end{corollary}

\begin{corollary}{12.49}
If $(Y,\tau')\sbs (X,\tau)$ and $Y\in\tau\cap\cD(X)$, then $Y$ is
a $\BrS$ in $X$ implies that $X$ is a $\BrS$.
\end{corollary}

\begin{pf}
The condition is an immediate consequence of (2) of
Corollary~12.30, Corollary~12.48 and the well-known topological
fact, following which

$U$ is of second category in itself if and only if
$U\in\Catg_{{}_{\II}}(X)$.
\end{pf}

\begin{proposition}{12.50}
If $(Z,\tau_1'',\tau_2'')\sbs (Y,\tau_1',\tau_2')\sbs
(X,\tau_1,\tau_2)$ and $A\sbs Y$ is any set, then $A\in
(i,j)\dd\,\cN\cD(Z,Y)\llra A\in (i,j)\dd\,\cN\cD(Z,X)$ $(A\in
(i,j)\dd\,\Catg_{{}_{\I}}(Z,Y)\llra A\in
(i,j)\dd\,\Catg_{{}_{\I}}(Z,X))$ and so
\begin{eqnarray*}
    A\in (i,j)\dd\,\cS\cD(Z,Y) &&\hskip-0.6cm \llra
            A\in (i,j)\dd\,\cS\cD(Z,X) \\
    \big(\,A\in (i,j)\dd\,\Catg_{{}_{\II}}(Z,Y) &&\hskip-0.6cm \llra
        A\in (i,j)\dd\,\Catg_{{}_{\II}}(Z,X)\,\big).
\end{eqnarray*}
\end{proposition}

\begin{pf}
Since $Z\sbs Y$ and $A\sbs Y$ we have
$$  \tau_i''\nt(\tau_j'\cl A\cap Z)=
        \tau_i''\nt(\tau_j\cl A\cap Y\cap Z)=
            \tau_i''\nt(\tau_j\cl A\cap Z).       $$
Hence $A\in(i,j)\dd\,\cN\cD(Z,Y)\llra A\in(i,j)\dd\,\cN\cD(Z,X)$.
The rest is obvious.
\end{pf}

\begin{corollary}{12.51}
If $(Z,\tau_1''<\tau_2'')\sbs (Y,\tau_1'<\tau_2')\sbs
(X,\tau_1<\tau_2)$ and $Y\in\tau_1$, then $Z$ is a $(1,2)\dd\BrS$
in $X$ implies that $Z$ is a $(1,2)\dd\BrS$ in~$Y$.
\end{corollary}

\begin{pf}
Let $U\in\tau_1'\setminus\{\vnth\}$ be any set. Then
$U\in\tau_1\setminus\{\vnth\}$ and by condition, $U$ is of
$(1,2)\dd\,\Catg_{{}_Z}\II$.
\end{pf}

\begin{corollary}{12.52}
If $(Z,\tau'')\sbs (Y,\tau')\sbs (X,\tau)$ and $A\sbs Y$ is any
set, then $A\in\cN\cD(Z,Y)\llra A\in\cN\cD(Z,X)$ $(A\in
\Catg_{{}_{\I}}(Z,Y)\llra A\in\Catg_{{}_{\I}}(Z,X))$ and so
$A\in\cS\cD(Z,Y)\llra A\in\cS\cD(Z,X)$
$(A\in\Catg_{{}_{\II}}(Z,Y)\llra A\in\Catg_{{}_{\II}}(Z,X))$.
\end{corollary}

\begin{corollary}{12.53}
If $(Z,\tau'')\sbs (Y,\tau')\sbs (X,\tau)$ and $Y\in\tau$, then
$Z$ is a $\BrS$ in $X$ implies that $Z$ is a $\BrS$ in $Y$.
\end{corollary}

\begin{proposition}{12.54}
If $(Z,\tau_1'',\tau_2'')\sbs (Y,\tau_1',\tau_2')\sbs
(X,\tau_1,\tau_2)$ and $A\sbs X$ is any set, then $A\in
(i,j)\dd\,\cN\cD(Z,X)$ $(A\in (i,j)\dd\,\Catg_{{}_{\I}}(Z,X))$
implies that $A\cap Y\in (i,j)\dd\,\cN\cD(Z,Y)$ $(A\cap Y\in
(i,j)\dd\,\Catg_{{}_{\I}}(Z,Y))$ and so $A\cap Y\in
(i,j)\dd\,\cS\cD(Z,Y)$ $(A\cap Y\in
(i,j)\dd\,\Catg_{{}_{\II}}(Z,Y))$ implies that $A\in
(i,j)\dd\,\cS\cD(Z,X)$ $(A\in (i,j)\dd\,\Catg_{{}_{\II}}(Z,X))$.
\end{proposition}

\begin{pf}
It is evident that $\tau_i''\nt(\tau_j\cl A\cap Z)=\vnth$ implies
that      \linebreak        $\tau_i''\nt(\tau_j'\cl(A\cap Y)\cap
Z)=\vnth$. The rest is obvious.
\end{pf}

\begin{corollary}{12.55}
If $(Z,\tau_1'',\tau_2'')\sbs (Y,\tau_1',\tau_2')\sbs
(X,\tau_1,\tau_2)$, $Y\in i\dd\,\cD(X)$ and $Z$ is an
$A\dd\,(i,j)\dd\BrS$ in $Y$, then $Z$ is an $A\dd\,(i,j)\dd\BrS$
in $X$.
\end{corollary}

\begin{pf}
Let $U\in\tau_i\setminus\{\vnth\}$ be any set. Then $U'=U\cap
Y\in\tau_i'\setminus\{\vnth\}$ and so
$U'\in(i,j)\dd\,\Catg_{{}_{\II}}(Z,Y)$. Hence
$U\in(i,j)\dd\,\Catg_{{}_{\II}}(Z,X)$.
\end{pf}

\begin{corollary}{12.56}
Let $(Z,\tau_1''<\tau_2'')\sbs (Y,\tau_1'<\tau_2')\sbs
(X,\tau_1<\tau_2)$ and $Y\in 2\dd\,\cD(X)$. Then
\begin{enumerate}
\item[(1)] $Z$ is a $(1,2)\dd\BrS$ in $Y$ implies that $Z$ is a
$(1,2)\dd\BrS$ in $X$.

\item[(2)] $Z$ is an $A\dd\,(2,1)\dd\BrS$ in $Y$ implies that $Z$
is an $A\dd\,(2,1)\dd\BrS$ in~$X$.
\end{enumerate}
\end{corollary}

\begin{corollary}{12.57}
If $(Z,\tau'')\sbs (Y,\tau')\sbs (X,\tau)$ and $A\sbs X$ is any
set, then $A\in\cN\cD(Z,X)$ $(A\in\Catg_{{}_{\I}}(Z,X))$ implies
that $A\cap Y\in\cN\cD(Z,Y)$ $(A\cap Y\in\Catg_{{}_{\I}}(Z,Y))$
and so $A\cap Y\in\cS\cD(Z,Y)$ $(A\cap Y\in\Catg_{{}_{\II}}(Z,Y)$
implies that $A\in\cS\cD(Z,X)$ $(A\in\Catg_{{}_{\II}}(Z,X))$.
\end{corollary}

\begin{corollary}{12.58}
If $(Z,\tau'')\sbs (Y,\tau')\sbs (X,\tau)$ and $Y\in\cD(X)$, then
$Z$ is a $\BrS$ in $Y$ implies that $Z$ is a $\BrS$ in $X$.
\end{corollary}

\begin{theorem}{12.59}
Let $(Y,\tau_1'<\tau_2')$ be a $\BsS$ of a $\BS$
$(X,\tau_1<\tau_2)$. Then the following conditions are satisfied:
\begin{enumerate}
\item[(1)] If $Y$ is a $(1,2)\dd\BrS$ in $X$ and
$\{A_n\}_{n=1}^\infty$ is a sequence of subsets of $X$, where
$A_n\in 2\dd\,\cG_\dl(X)\cap 1\dd\,\sD(Y,X)$ for each
$n=\ol{1,\infty}$, then $\bigcap\limits_{n=1}^\infty A_n\in
2\dd\,\cG_\dl(X)\cap 1\dd\,\cD(X)$.

\item[(2)] If $Y$ is an $A\dd\,(2,1)\dd\BrS$ in $X$ and
$\{A_n\}_{n=1}^\infty$ is a sequence of subsets of $X$, where
$A_n\in 1\dd\,\cG_\dl(X)\cap 2\dd\,\sD(Y,X)$ for each
$n=\ol{1,\infty}$, then $\bigcap\limits_{n=1}^\infty A_n\in
1\dd\,\cG_\dl(X)\cap 2\dd\,\cD(X)$.
\end{enumerate}
\end{theorem}

\begin{pf}
(1) Clearly, $X\setminus A_n=\bigcup\limits_{k=1}^\infty
(X\setminus A_n^k)$, where $X\setminus A_n^k\in\co\tau_2$ for each
$k=\ol{1,\infty}$, $n=\ol{1,\infty}$ and $X\setminus A_n\in
2\dd\,\cF_\sg(X)\cap 1\dd\,\sBd(Y,X)$ for each $n=\ol{1,\infty}$.
Then $X\setminus A_n^k\in 1\dd\,\sBd(Y,X)$ for each
$k=\ol{1,\infty}$, $n=\ol{1,\infty}$ and by (2) of Theorem~12.16,
$X\setminus A_n\in (1,2)\dd\,\Catg_{{}_{\I}}(Y,X)$. Hence, by (1)
of Theorem~12.16, $\bigcup\limits_{n=1}^\infty (X\setminus A_n)\in
(1,2)\dd\,\Catg_{{}_{\I}}(Y,X)$. But by (1) of Corollary~12.29,
$Y$ is an $A\dd\,(1,2)\dd\BrS$ in $X$ and hence, by (3) of
Theorem~12.24 with Definition~12.26 taken into account,
$$  \tau_1\cl\Big(X\setminus
            \bigcup\limits_{n=1}^\infty (X\setminus A_n)\Big)=
        \tau_1\cl \bigcap\limits_{n=1}^\infty A_n=X.    $$
Thus $\bigcap\limits_{n=1}^\infty A_n\in 2\dd\,\cG_\dl(X)\cap
1\dd\,\cD(X)$.

(2) Evidently, $X\setminus A_n=\bigcup\limits_{k=1}^\infty
(X\setminus A_n^k)$, where $X\setminus A_n^k\in\co\tau_1$ for each
$k=\ol{1,\infty}$, $n=\ol{1,\infty}$ and $X\setminus A_n\in
1\dd\,\cF_\sg(X)\cap 2\dd\,\sBd(Y,X)$ for each $n=\ol{1,\infty}$.
Then $X\setminus A_n^k\in 2\dd\,\sBd(Y,X)$ for each
$k=\ol{1,\infty}$, $n=\ol{1,\infty}$ and by (2) of Theorem~12.16,
$X\setminus A_n\in (2,1)\dd\,\Catg_{{}_{\I}}(Y,X)$ for each
$n=\ol{1,\infty}$. Hence, by (1) of Theorem~12.16,
$\bigcup\limits_{n=1}^\infty (X\setminus A_n)\in
(2,1)\dd\,\Catg_{{}_{\I}}(Y,X)$. Thus, by (3) of Theorem~12.24 in
conjunction with Definition~12.26,
$$  \tau_2\cl\Big(X\setminus
            \bigcup\limits_{n=1}^\infty (X\setminus A_n)\Big)=
        \tau_2\cl \bigcap\limits_{n=1}^\infty A_n=X,    $$
i.e., $\bigcap\limits_{n=1}^\infty A_n\in 1\dd\,\cG_\dl(X)\cap
2\dd\,\cD(X)$.
\end{pf}

\begin{corollary}{12.60}
If $(Y,\tau')$ is a $\TsS$ of a $\TS$ $(X,\tau)$, $Y$ is a $\BrS$
in $X$ and $\{A_n\}_{n=1}^\infty$ is a sequence of subsets of $X$,
where $A_n\in\cG_\dl(X)\cap \sD(Y,X)$ for each $n=\ol{1,\infty}$,
then $\bigcap\limits_{n=1}^\infty A_n\in\cG_\dl(X)\cap\cD(X)$.
\end{corollary}

Note here that the theorem below characterizes all subsets of $X$,
belonging to the family$(1,2)\dd\,\Catg_{{}_{\I}}(Y,X)$
$((2,1)\dd\,\Catg_{{}_{\I}}(Y,X))$ for a      \linebreak
$(1,2)\dd\BrS$ (an $A\dd\,(2,1)\dd\BrS)$ $Y$ in $X$.

\begin{theorem}{12.61}
Let $(Y,\tau_1'<\tau_2')$ be a $\BsS$ of a $\BS$
$(X,\tau_1<\tau_2)$. Then for any subset $A\sbs X$ the following
conditions are satisfied:
\begin{enumerate}
\item[(1)] If $Y\!\in\!\tau_1$ and $Y$ is a $(1,2)\dd\BrS$ in $X$,
then $A\!\in\!(1,2)\dd\,\Catg_{{}_{\I}}(Y,X)$ if and only if
$X\setminus A$ contains a set $B\in 2\dd\,\cG_\dl(X)\cap
1\dd\,\cD(X)$.

\item[(2)] If $Y\in\tau_2$ and $Y$ is an $A\dd\,(2,1)\dd\BrS$ in
$X$, then $A\in         \linebreak
    (2,1)\dd\,\Catg_{{}_{\I}}(Y,X)$ if and only if $X\setminus A$
contains a set $B\in 1\dd\,\cG_\dl(X)\cap 2\dd\,\cD(X)$.
\end{enumerate}
\end{theorem}

\begin{pf}
(1) First, let $A\in (1,2)\dd\,\Catg_{{}_{\I}}(Y,X)$, i.e.,
$A=\bigcup\limits_{n=1}^\infty A_n$, where $A_n\in
(1,2)\dd\,\cN\cD(Y,X)$ for each $n=\ol{1,\infty}$. Then
$\{X\setminus\tau_2\cl A_n\}_{n=1}^\infty$ is a countable family
of $2$-open subsets of $X$ and
\begin{eqnarray*}
    &\ds \tau_1'\cl\big((X\setminus\tau_2\cl A_n)\cap Y\big)=
            \tau_1'\cl (Y\setminus\tau_2\cl A_n)= \\
    &\ds =\tau_1'\cl\big(Y\setminus(\tau_2\cl A_n\cap Y)\big)=
        Y\setminus\tau_1'\nt(\tau_2\cl A_n\cap Y)=Y
\end{eqnarray*}
so that $X\setminus\tau_2\cl A_n\in\tau_2\cap1\dd\,\sD(Y,X)$ for
each $n=\ol{1,\infty}$. Since $Y$ is a $(1,2)\dd\BrS$ in $X$, by
(1) of Corollary~12.29 in conjunction with (2) of Theorem~12.24
$$  B=\bigcap\limits_{n=1}^\infty
        (X\setminus\tau_2\cl A_n)\in 2\dd\,\cG_\dl(X)\cap 1\dd\,\cD(X)  $$
and $B\sbs X\setminus A$.

Conversely, let $B=\bigcap\limits_{n=1}^\infty A_n\in
2\dd\,\cG_\dl(X)\cap 1\dd\,\cD(X)$ and $B\sbs X\setminus A$. Since
$\tau_1\cl B=X$, we have $\tau_1\cl A_n=X$ for each
$n=\ol{1,\infty}$. Hence,
$$  \tau_1\nt(X\setminus A_n)=\tau_1\nt\tau_2\cl(X\setminus A_n)=
        X\setminus\tau_1\cl\tau_2\nt A_n=\vnth      $$
so that $X\setminus A_n\in (1,2)\dd\,\cN\cD(X)$ for each
$n=\ol{1,\infty}$. Since $Y\in\tau_1$, taking into account
Proposition~12.44,
$$  \bigcup\limits_{n=1}^\infty (X\setminus A_n)\in
        (1,2)\dd\,\Catg_{{}_{\I}}(X)\sbs
                (1,2)\dd\,\Catg_{{}_{\I}}(Y,X).     $$
Moreover, $B=\bigcap\limits_{n=1}^\infty A_n\sbs X\setminus A$
implies that $A\sbs X\setminus B$ and thus, by (1) of
Theorem~12.16, $A\in (1,2)\dd\,\Catg_{{}_{\I}}(Y,X)$.

(2) Let $A\in (2,1)\dd\,\Catg_{{}_{\I}}(Y,X)$, i.e.,
$A=\bigcup\limits_{n=1}^\infty A_n$, where $A_n\in      \linebreak
(2,1)\dd\,\cN\cD(Y,X)$ for each $n=\ol{1,\infty}$. Then
$\{X\setminus\tau_1\cl A_n\}_{n=1}^\infty$ is a countable family
of $1$-open subsets of $X$ and
\begin{eqnarray*}
    &\ds \tau_2'\cl\big((X\setminus\tau_1\cl A_n)\cap Y\big)=
            \tau_2'\cl\big(Y\setminus(\tau_1\cl A_n\cap Y)\big)= \\
    &\ds =Y\setminus\tau_2'\nt(\tau_1\cl A_n\cap Y)=Y
\end{eqnarray*}
so that $X\setminus\tau_1\cl A_n\in \tau_1\cap 2\dd\,\sD(Y,X)$ for
each $n=\ol{1,\infty}$. Since $Y$ is an $A\dd\,(2,1)\dd\BrS$ in
$X$, by (2) of Theorem~12.24, $B=\bigcap\limits_{n=1}^\infty
(X\setminus\tau_1\cl A_n)\in 1\dd\,\cG_\dl(X)\cap 2\dd\,\cD(X)$
and $B\sbs X\setminus A$.

Conversely, let $B=\bigcap\limits_{n=1}^\infty A_n\in
1\dd\,\cG_\dl(X)\cap 2\dd\,\cD(X)$ and $B\sbs X\setminus A$. Since
$\tau_2\cl B=X$, we have $\tau_2\cl A_n=X$ for each
$n=\ol{1,\infty}$. Hence,
$$  \tau_2\nt(X\setminus A_n)=\tau_2\nt\tau_1\cl(X\setminus A_n)=
        X\setminus\tau_2\cl\tau_1\nt A_n=\vnth      $$
so that $X\setminus A_n\in (2,1)\dd\,\cN\cD(X)$ for each
$n=\ol{1,\infty}$. Therefore, since $Y\in\tau_2$, by
Proposition~12.44,
$$  \bigcup\limits_{n=1}^\infty (X\setminus A_n)\in
        (2,1)\dd\,\Catg_{{}_{\I}}(X)\sbs (2,1)\dd\,\Catg_{{}_{\I}}(Y,X)  $$
and since $A\!\sbs\!X\setminus B$, by (1) of Theorem~12.16,
$A\!\in\!(2,1)\dd\,\Catg_{{}_{\I}}(Y,X)$.~\end{pf}

\begin{corollary}{12.62}
If $(Y,\tau')$ is a $\BrS$ in $(X,\tau)$, $Y\in\tau$ and $A\sbs X$
is any set, then $A\in\Catg_{{}_{\I}}(Y,X)$ if and only if
$X\setminus A$ contains a set $B\in\cG_\dl(X)\cap\cD(X)$
\end{corollary}

Thus, Corollary~12.62 characterizes all subsets of $X$, belonging
to the family $\Catg_{{}_{\I}}(Y,X)$ for a $\BrS$ $Y$ in $X$.

\begin{theorem}{12.63}
Let $\{(Y_\al,\tau_1^\al,\tau_2^\al)\}_{\al\in D}$ be a family of
bitopological subspaces of a $\BS$ $(X,\tau_1,\tau_2)$ and
$Y=\bigcup\limits_{a\in D} Y_\al\in i\dd\,\cD(X)$. If there is
$\al_0\in D$ such that $Y_{\al_0}\in \tau_i$ and $Y_{\al_0}$ is an
$A\dd\,(i,j)\dd\BrS$ in $Y$, then $Y$ is an $A\dd\,(i,j)\dd\BrS$
in $X$.
\end{theorem}

\begin{pf}
Let $\al_0\in D$ and $Y_{\al_0}$ be an $A\dd\,(i,j)\dd\BrS$ in
$(Y,\tau_1,\tau_2')$. Then $Y_{\al_0}\in\tau_i$ and $Y_{\al_0}\sbs
Y$ imply that $Y_{\al_0}\in\tau_i'$. Hence, it remains to use
Corollary~12.36.
\end{pf}

\begin{corollary}{12.64}
Let $\{(Y_\al,\tau_1^\al<\tau_2^\al)\}_{\al\in D}$   be a family
of bitopological subspaces of a $\BS$ $(X,\tau_1<\tau_2)$ such
that $Y=\bigcup\limits_{\al\in D} Y_\al\in 2\dd\,\cD(X)$. Then the
following conditions are satisfied:
\begin{enumerate}
\item[(1)] If $Y_{\al_0}$ is a $(1,2)\dd\BrS$ in $Y$ and
$Y_{\al_0}\in\tau_1$ for some $\al_0\in D$, then $Y$ is a
$(1,2)\dd\BrS$ in $X$.

\item[(2)] If $Y_{\al_0}$ is an $A\dd\,(2,1)\dd\BrS$ in $Y$ and
$Y_{\al_0}\in\tau_2$ for some $\al_0\in D$, then $Y$ is an
$A\dd\,(2,1)\dd\BrS$ in $X$.
\end{enumerate}
\end{corollary}

\begin{pf}
The conditions follows directly from Corollary~12.37.
\end{pf}

\begin{corollary}{12.65}
Let $\{(Y_\al,\tau_\al)\}_{\al\in D}$ be a family of topological
subspaces of a $\TS$ $(X,\tau)$, $Y_{\al_0}$ be a $\BrS$ in $Y$
and $Y_{\al_0}\in\tau$ for some $\al_0\in D$, then
$Y=\bigcup\limits_{\al\in D} Y_\al\in\cD(X)$ implies that $Y$ is a
$\BrS$ in $X$.
\end{corollary}

\begin{theorem}{12.66}
If $(Y,\tau_1'<\tau_2')\sbs (X,\tau_1<\tau_2)$, $Y\in\tau_2$ and
$Y$ is a $(1,2)\dd\BrS$ in $X$, then
$$  2\dd\,\cG_\dl(X)\cap (1,2)\dd\,\Catg_{{}_{\I}}(Y,X)\sbs
            (1,2)\dd\,\cN\cD(Y,X).        $$
\end{theorem}

\begin{pf}
Contrary: there is $A\in 2\dd\,\cG_\dl(X)\cap
(1,2)\dd\,\Catg_{{}_{\I}}(Y,X)$ such that
$A\,\ol{\in}\,(1,2)\dd\,\cN\cD(Y,X)$. Hence, there is a set $V'\in
\tau_1'\setminus\{\vnth\}$ such that $V'\sbs\tau_2\cl A\cap Y$.
Let $V\in\tau_1$, $V\cap Y=V'$ and let us prove that $V$ is of
$(1,2)\dd\,\Catg_{{}_Y}\I$, i.e., $V\in
(1,2)\dd\,\Catg_{{}_{\I}}(V',V)$. Since $V=(V\cap A)\cup
(V\setminus A)$, by (1) of Theorem~12.16 it suffices to prove that
$V\cap A,\,V\setminus A\in (1,2)\dd\,\Catg_{{}_{\I}}(V',V)$. Since
$A\in (1,2)\dd\,\Catg_{{}_{\I}}(Y,X)$, $V'\sbs Y\sbs X$ and
$V'\in\tau_1'$, by (2) of Proposition~12.34, $A\in
(1,2)\dd\,\Catg_{{}_{\I}}(V',X)$. Hence, by (1) of Theorem~12.16,
A$\cap V\in (1,2)\dd\,\Catg_{{}_{\I}}(V',X)$. Since $A\cap V\sbs
V$ and $V'\sbs V\sbs X$, by Proposition~12.50,
$$  A\cap V\in (1,2)\dd\,\Catg_{{}_{\I}}(V',X)\llra A\cap V\in
        (1,2)\dd\,\Catg_{{}_{\I}}(V',V).        $$
Now, let us prove that $V\setminus A\in
(1,2)\dd\,\Catg_{{}_{\I}}(V',V)$. For this purpose, taking into
account (2) of Theorem~12.16, it suffices to prove $V\setminus
A=\bigcup\limits_{n=1}^\infty F_n\in 2\dd\,\cF_\sg(V)$ and $F_n\in
1\dd\,\sBd(V',V)$ for each $n=\ol{1,\infty}$. Evidently, $A\in
2\dd\,\cG_\dl(X)$ implies that $V\setminus A=(X\setminus A)\cap
V\in 2\dd\,\cF_\sg(V)$, so that $V\setminus
A=\bigcup\limits_{n=1}^\infty F_n$, where $F_n\in\co\tau_2''$ in
$(V,\tau_1''<\tau_2'')$. Let us prove that $\tau_1'''\nt(F_n\cap
V')=\vnth$ in $(V',\tau_1'''<\tau_2''')$, or equivalently,
$\tau_1'''\cl(V'\setminus F_n)=V'$ for each $n=\ol{1,\infty}$. We
have:
$$  V'\cap F_n\sbs V'\cap (V\setminus A)=V'\setminus A      $$
and hence, $V'\setminus(V'\setminus A)\sbs V'\setminus(V'\cap
F_n)$ so that $V'\cap A\sbs V'\setminus F_n$. Therefore, it
suffices to prove that $V'=\tau_1'''\cl(V'\cap A)$. Since
$V\in\tau_1\sbs\tau_2$, $V'\sbs\tau_2\cl A\cap Y$ and
$Y\in\tau_2$, we have
\begin{eqnarray*}
    V' &&\hskip-0.6cm =(V'\cap V)\cap Y\sbs
        (\tau_2\cl A\cap Y\cap V)\cap Y\sbs \tau_2\cl (A\cap Y\cap V)\cap Y= \\
    &&\hskip-0.6cm =\tau_2'\cl(V'\cap A)
\end{eqnarray*}
and so $V'=\tau_1'''\cl(V'\cap A)$. Thus $V'\cap A\in
1\dd\,\sD(V',V)$ so that $F_n\in 1\dd\,\sBd(V',V)$ for each
$n=\ol{1,\infty}$. Thus, $V\setminus A\in
(1,2)\dd\,\Catg_{{}_{\I}}(V',V)$ and so $V\in
(1,2)\dd\,\Catg_{{}_{\I}}(V',V)$, i.e., $V$ is of
$\Catg_{{}_Y}\I$, which is impossible, since $Y$ is a
$(1,2)\dd\BrS$ in $X$.
\end{pf}

\begin{corollary}{12.67}
If $(Y,\tau')\sbs (X,\tau)$, $Y\in\tau$ and $Y$ is $\BrS$ in $X$,
then $\cG_\dl(X)\cap\Catg_{{}_{\I}}(Y,X)\sbs \cN\cD(Y,X)$.
\end{corollary}

Furthermore, let us prove that thanks to (5) of Corollary~12.13, a
countable family of $1$-open $2$-strongly dense relative to $Y$
sets in (2) of Theorem~12.24 can be essentially narrowed for an
$A\dd\,(2,1)\dd\BrS$       \linebreak $(Y,\tau_1'<\tau_2')$ in a
$\BS$ $(X,\tau_1<\tau_2)$. Namely, take place

\begin{theorem}{12.68}
A $\BsS$ $(Y,\tau_1'<\tau_2')$ of a $\BS$ $(X,\tau_1<\tau_2)$ is
an $A\dd\,(2,1)\dd\BrS$ in $X$ if and only if the intersection of
any monotone decreasing sequence of $2$-strongly dense relative to
$Y$ and $1$-open sets is $2$-dense in $X$.
\end{theorem}

\begin{pf}
Evidently, it suffices to prove that if the intersection of any
monotone decreasing sequence of $2$-strongly dense relative to $Y$
and $1$-open sets is $2$-dense in $X$, then (3) of Theorem~12.24
is satisfied.

Let $A\in (2,1)\dd\,\Catg_{{}_{\I}}(Y,X)$ be any set. Then, by
analogy with the proof of the implication $(2)\lra(3)$ in
Theorem~12.24, one can assume that $A\in 1\dd\,\cF_\sg(X)\cap
(2,1)\dd\,\Catg_{{}_{\I}}(Y,X)$ so that
$A=\bigcup\limits_{n=1}^\infty A_n$, where $A_n\in\co\tau_1\cap
(2,1)\dd\,\cN\cD(Y,X)$ for each $n=\ol{1,\infty}$. Suppose that
$X\setminus A\,\ol{\in}\,2\dd\,\cD(X)$ that is,
$\bigcap\limits_{n=1}^\infty (X\setminus
A_n)\,\ol{\in}\,2\dd\,\cD(X)$. Let $U_n=\bigcap\limits_{k=1}^n
(X\setminus A_k)$ for each $n=\ol{1,\infty}$. It is obvious that
$\{U_n\}_{n=1}^\infty$ is a monotone decreasing sequence.
Moreover, $X\setminus U_n=\bigcup\limits_{k=1}^n A_k$, where
$A_k\in\co\tau_1\cap (2,1)\dd\,\cN\cD(Y,X)$ for each
$k=\ol{1,\infty}$ and $n=\ol{1,\infty}$. According to (5) of
Corollary~12.13, $X\setminus
U_n\in\co\tau_1\cap(2,1)\dd\,\cN\cD(Y,X)$ and so $U_n\in \tau_1$
and $\tau_2'\nt(\tau_1\cl(X\setminus U_n)\cap Y)=\vnth$ for each
$n=\ol{1,\infty}$. Hence
\begin{eqnarray*}
    &\ds Y\setminus\tau_2'\nt\big(\tau_1\cl(X\setminus U_n)\cap Y\big)=
        Y\setminus\tau_2'\nt(Y\setminus U_n)= \\
    &\ds =\tau_2'\cl(U_n\cap Y)=Y
\end{eqnarray*}
so that $U_n\in\tau_1\cap 2\dd\,\sD(Y,X)$ for each
$n=\ol{1,\infty}$. Hence, by condition
$\bigcap\limits_{n=1}^\infty U_n\in 2\dd\,\cD(X)$, where
$\bigcap\limits_{n=1}^\infty U_n=\bigcap\limits_{n=1}^\infty
(X\setminus A_n)$ and we obtain a contradiction to
$\bigcap\limits_{n=1}^\infty (X\setminus
A_n)\ol{\in}\,2\dd\,\cD(X)$.
\end{pf}

\begin{corollary}{12.69}
A $\TsS$ $(Y,\tau')$ of a $\TS$ $(X,\tau)$ is an $A\dd\BrS$ in $X$
if and only if the intersection of any monotone decreasing
sequence of strongly dense relative to $Y$ open sets is dense in
$X$.
\end{corollary}

Clearly, if $Y\in\cD(X)$, then Corollary~12.69 gives the
characterization of a $\BrS$ $Y$ in $X$.

\begin{theorem}{12.70}
Let $(Y,\tau_1'<\tau_2')$ be a $\BsS$ of a $\BS$
$(X,\tau_1<\tau_2)$. Then the following conditions are equivalent:
\begin{enumerate}
\item[(1)] $Y$ is an $A\dd\,(i,j)\dd\BrS$ in $X$.

\item[(2)] If $A\in
(i,j)\dd\,\Catg_{{}_{\I}}(Y,X)\setminus\{\vnth\}$, $X\setminus
A\in i\dd\,\cD(X)$ and $Y\sbs X\setminus A$, then $Y$ is an
$A\dd\,(i,j)\dd\BrS$ in $X\setminus A$.
\end{enumerate}
\end{theorem}

\begin{pf}
$(1)\lra (2)$ Let $Y$ be an $A\dd\,(i,j)\dd\BrS$ in $X$,
$A\!\in\!(i,j)\dd\,\Catg_{{}_{\I}}(Y,X)\setminus \{\vnth\}$ and
$Y\sbs X\setminus A$. If $B\in
(i,j)\dd\,\Catg_{{}_{\I}}(Y,X\setminus A)$, then by
Proposition~12.50, $B\in (i,j)\dd\,\Catg_{{}_{\I}}(Y,X)$ and so,
by (1) of Theorem~12.16 and (3) of Theorem~12.24,
$X\setminus(A\cup B)\in i\dd\,\cD(X)$. But $X\setminus(A\cup
B)=(X\setminus A)\setminus B$ and so $(X\setminus A)\setminus B\in
i\dd\,\cD(X\setminus A)$. Now, once more applying (3) of
Theorem~12.24, we obtain that $Y$ is an $A\dd\,(i,j)\dd\BrS$ in
$X\setminus A$.

$(2)\lra (1)$ First, let $X$ contains nonempty subset $A\in
\linebreak     (i,j)\dd\,\Catg_{{}_{\I}}(Y,X)$. Then
$U\in\tau_i\setminus\{\vnth\}$ implies that $U\in
(i,j)\dd\,\Catg_{{}_{\II}}(Y,X)$ and by (1) of Theorem~12.24, $Y$
is an $A\dd\,(i,j)\dd\BrS$ in $X$. Now, if there is $A\in
(i,j)\dd\,\Catg_{{}_{\I}}(Y,X)\setminus\{\vnth\}$ such that
$X\setminus A\in i\dd\,\cD(X)$, $Y\sbs X\setminus A$ and $Y$ is an
$A\dd\,(i,j)\dd\BrS$ in $X\setminus A$, then by Corollary~12.55,
$Y$ is an $A\dd\,(i,j)\dd\BrS$ in $X$.
\end{pf}

\begin{corollary}{12.71}
For a $2$-dense $\BsS$ $(Y,\tau_1'\!<\!\tau_2')$ of a $\BS$
$(X,\tau_1\!<~\!\!\!\tau_2)$ the following conditions are
equivalent:
\begin{enumerate}
\item[(1)] $Y$ is a $(1,2)\dd\BrS$ in $X$.

\item[(2)] If $A\in
(1,2)\dd\,\Catg_{{}_{\I}}(Y,X)\setminus\{\vnth\}$, $Y\sbs
X\setminus A$, then $Y$ is a $(1,2)\dd\BrS$ in $X\setminus A$.
\end{enumerate}
\end{corollary}

\begin{pf}
Evidently, $Y\in 2\dd\,\cD(X)$ implies, on the one hand, that
$Y\in 2\dd\,\cD(X\setminus A)$ and, on the other hand, that
$X\setminus A\in 2\dd\,\cD(X)$. Hence, it remains to use (2) of
Corollary~12.29 together with Theorem~12.70.
\end{pf}

\begin{corollary}{12.72}
Let $(Y,\tau')\sbs (X,\tau)$. Then the following conditions are
equivalent:
\begin{enumerate}
\item[(1)] $Y$ is an $A\dd\BrS$ in $X$.

\item[(2)] If $A\in\Catg_{{}_{\I}}(Y,X)\setminus\{\vnth\}$,
$X\setminus A\in \cD(X)$ and $Y\sbs X\setminus A$, then $Y$ is an
$A\dd\BrS$ in $X\setminus A$.
\end{enumerate}
\end{corollary}

Clearly, if $Y\in\cD(X)$, then Corollary~12.72 remains valid for
Baire spaces.

Our last aim of Section~12 is to show that the notion of
concentration of a set on a subspace is a link connecting the
relative property of this subspace to be Baire in the whole space
with the property of the same subspace to be Baire in itself for
both the topological and the bitopological case. As we shall see
below, this idea is essentially based on the density of a subspace
in the whole space.

\begin{definition}{12.73}
Let $(Y,\tau_1',\tau_2')$ $((Y,\tau'))$ be a $\BsS$ $(\TsS$) of a
$\BS$ $(\TS)$ $(X,\tau_1,\tau_2)$ $((X,\tau))$. Then a subset
$A\sbs X$ is $(i,j)\dd\WW$-nowhere dense $(\WW$-nowhere dense)
relative to $Y$ if $\tau_i'\nt\tau_j'\cl(A\cap Y)=\vnth$
\linebreak     $(\tau'\nt\tau'\cl(A\cap Y)=\vnth)$.
\end{definition}

The families of all subsets of $X$ which are $(i,j)\dd\WW$-nowhere
$(\WW$-now\-he\-re) dense relative to $Y$, are denoted by
$(i,j)\dd\,\wND(Y,X)$    \linebreak        $(\wND(Y,X))$. It is
evident that
$$  A\in (i,j)\dd\,\wND(Y,X)\llra A\cap Y\in (i,j)\dd\,\cN\cD(Y)    $$
so that if $(Y,\tau')\sbs (X,\tau)$, then $A\in \wND(Y,X)\llra
A\cap Y\in\cN\cD(Y)$.

Moreover, we say that a subset $A$ of a $\BS$ $(X,\tau_1,\tau_2)$
$(\TS$ $(X,\tau))$ is $j\dd\WW$-concentrated $(\WW$-concentrated)
on a $\BsS$ $(Y,\tau_1',\tau_2')$ $(\TsS$ $(Y,\tau'))$ if
$\tau_j\cl A\cap Y\sbs\tau_j\cl(A\cap Y)$ $(\tau\cl A\cap
Y\sbs\tau\cl(A\cap Y))$. Clearly, if $A$ is $j$-concentrated
(concentrated) on $Y$ or $Y\in\tau_j$ $(Y\in\tau)$, then $A$ is
$j\dd\WW$-concentrated $(\WW$-concentrated) on $Y$.

Now, observe that if $(Y,\tau_1',\tau_2')\sbs (X,\tau_1,\tau_2)$
$((Y,\tau')\sbs (X,\tau))$, $A\sbs X$ and $A$ is
$j\dd\WW$-concentrated $(\WW$-concentrated) on $Y$, then
\begin{eqnarray*}
    A\in (i,j)\dd\,\cN\cD(Y,X) &&\hskip-0.6cm \llra A\in(i,j)\dd\,wND(Y,X) \\
    \big(\,A\in\cN\cD(Y,X) &&\hskip-0.6cm \llra A\in\wND(Y,X)\,\big).
\end{eqnarray*}

Indeed, if $A$ is $j\dd\WW$-concentrated $(\WW$-concentrated) on
$Y$, then        \linebreak   $\tau_j\cl A\cap Y\sbs
\tau_j\cl(A\cap Y)$ $(\tau\cl A\cap Y\sbs \tau\cl(A\cap Y))$ and
therefore,
\allowdisplaybreaks
\begin{eqnarray*}
    \tau_j\cl A\cap Y &&\hskip-0.6cm =\tau_j\cl(A\cap Y)\cap Y=\tau_j'\cl(A\cap Y) \\
    \big(\,\tau\cl A\cap Y &&\hskip-0.6cm =\tau'\cl (A\cap Y)\,\big).
\end{eqnarray*}
Thus
\begin{eqnarray*}
    \tau_i'\nt(\tau_j\cl A\cap Y) &&\hskip-0.6cm =
        \tau_i'\nt\tau_j'\cl(A\cap Y) \\
    \big(\,\tau'\nt(\tau\cl A\cap Y) &&\hskip-0.6cm =\tau'\nt\tau'\cl(A\cap Y)\,\big)
\end{eqnarray*}
so that
\begin{eqnarray*}
    A\in(i,j)\dd\,\cN\cD(Y,X) &&\hskip-0.6cm \llra A\cap Y\in(i,j)\dd\,\cN\cD(Y) \\
    \big(\,A\in\cN\cD(Y,X) &&\hskip-0.6cm \llra A\cap Y\in\cN\cD(Y)\,\big), \\
    A\in(i,j)\dd\,\cS\cD(Y,X) &&\hskip-0.6cm \llra A\cap Y\in(i,j)\dd\,\cS\cD(Y) \\
    \big(\,A\in\cS\cD(Y,X) &&\hskip-0.6cm \llra A\cap Y\in\cS\cD(Y)\,\big).
\end{eqnarray*}

Therefore, if $A$ is $j\dd\WW$-concentrated $(\WW$-concentrated)
on $Y$, then
\begin{eqnarray*}
    A\in (i,j)\dd\,\Catg_{{}_{\I}}(Y,X) &&\hskip-0.6cm \llra
        A\cap Y\in (i,j)\dd\,\Catg_{{}_{\I}}(Y) \\
    \big(\,A\in\Catg_{{}_{\I}}(Y,X) &&\hskip-0.6cm \llra
        A\cap Y\in\Catg_{{}_{\I}}(Y)\,\big)
\end{eqnarray*}
and so
\begin{eqnarray*}
    A\in (i,j)\dd\,\Catg_{{}_{\II}}(Y,X) &&\hskip-0.6cm \llra
        A\cap Y\in (i,j)\dd\,\Catg_{{}_{\II}}(Y) \\
    \big(\,A\in\Catg_{{}_{\II}}(Y,X) &&\hskip-0.6cm \llra
        A\cap Y\in\Catg_{{}_{\II}}(Y)\,\big).
\end{eqnarray*}

\begin{theorem}{12.74}
If $(Y,\tau_1'<\tau_2')\sbs (X,\tau_1<\tau_2)$ and $Y\in
2\dd\,\cD(X)$, then $Y$ is an $A\dd\,(i,j)\dd\BrS$ in $X$ if and
only if $Y$ is an $A\dd\,(i,j)\dd\BrS$.
\end{theorem}

\begin{pf}
First of all, let us prove that if $Y\in 2\dd\,\cD(X)$, then
$\tau_j\cl U=\tau_j\cl(U\cap Y)$ for each set
$U\in\tau_i\setminus\{\vnth\}$. Indeed, if
$U\in\tau_2\setminus\{\vnth\}$, then $\tau_2\cl U=\tau_2\cl(U\cap
Y)$. Hence, by (2) of Lemma~0.2.1 in [8], $\tau_1\cl
U=\tau_1\cl(U\cap Y)$. Now, if $U\in\tau_1\setminus\{\vnth\}$,
then $U\in\tau_2\setminus\{\vnth\}$ and so $\tau_2\cl
U=\tau_2\cl(U\cap Y)$.

Therefore, for any set $U\in\tau_i\setminus\{\vnth\}$ we have
$$  \tau_j\cl U\cap Y\sbs \tau_j\cl U=\tau_j\cl(U\cap Y)    $$
so that, every set $U\in\tau_i\setminus\{\vnth\}$ is
$j\dd\WW$-concentrated on $Y$. By reasonings before Theorem~12.74,
$$  U\in(i,j)\dd\,\Catg_{{}_{\II}}(Y,X)\llra
        U\cap Y\in(i,j)\dd\,\Catg_{{}_{\II}}(Y).      $$

Now, let $Y$ be an $A\dd\,(i,j)\dd\BrS$ in $X$ and
$U'\in\tau_i'\setminus\{\vnth\}$ be any set. Then there is a set
$U\in\tau_i\setminus\{\vnth\}\cap(i,j)\dd\,\Catg_{{}_{\II}}(Y,X)$
such that $U\cap Y=U'$. Thus $U'\in(i,j)\dd\,\Catg_{{}_{\II}}(Y)$
and so $Y$ is an $A\dd\,(i,j)\dd\BrS$. Conversely, if $Y$ is an
$A\dd\,(i,j)\dd\BrS$ and $U\in\tau_i\setminus\{\vnth\}$ is any
set, then the set $U'=U\cap Y\in(\tau_i'\setminus\{\vnth\})\cap
(i,j)\dd\,\Catg_{{}_{\II}}(Y)$ and so $U\in
(i,j)\dd\,\Catg_{{}_{\II}}(Y,X)$. Thus $Y$ is an
$A\dd\,(i,j)\dd\BrS$ in $X$.
\end{pf}

\begin{corollary}{12.75}
Under the hypotheses of Theorem~12.74 the following conditions are
satisfied:
\begin{enumerate}
\item[(1)] $Y$ is a $(1,2)\dd\BrS$ in $X$ if and only if $Y$ is a
$(1,2)\dd\BrS$.

\item[(2)] $Y$ is an $A\dd\,(2,1)\dd\BrS$ in $X$ if and only if
$Y$ is an $A\dd\,(2,1)\dd\BrS$.
\end{enumerate}
\end{corollary}

\begin{pf}
(1) Follows directly from (2) of Corollary~12.29 and (4) of
Theorem~4.1.6 in [8].

(2) The condition is obvious.
\end{pf}

\begin{corollary}{12.76}
Let $(Y,\tau')\sbs (X,\tau)$ and $Y\in\cD(X)$. Then $Y$ is a
$\BrS$ in $X$ if and only if $Y$ is a $\BrS$.
\end{corollary}

\begin{corollary}{12.77}
If $(Y,\tau_1'<\tau_2')$ is a $\BsS$ of a $\BS$
$(X,\tau_1<_S\tau_2)$ and $Y\in\tau_1\cap 2\dd\,\cD(X)$, then the
following conditions are equivalent:
\begin{enumerate}
\item[(1)] $Y$ is a $1\dd\BrS$.

\item[(2)] $Y$ is an $A\dd\,(1,2)\dd\BrS$.

\item[(3)] $Y$ is a $(1,2)\dd\BrS$.

\item[(4)] $Y$ is a $2\dd\BrS$.

\item[(5)] $Y$ is an $A\dd\,(2,1)\dd\BrS$.

\item[(6)] $Y$ is a $(2,1)\dd\BrS$.

\item[(7)] $Y$ is a $1\dd\BrS$ in $X$.

\item[(8)] $Y$ is an $A\dd\,(1,2)\dd\BrS$ in $X$.

\item[(9)] $Y$ is a $(1,2)\dd\BrS$ in $X$.

\item[(10)] $Y$ is a $2\dd\BrS$ in $X$.

\item[(11)] $Y$ is an $A\dd\,(2,1)\dd\BrS$ in $X$.

\item[(12)] $Y$ is a $(2,1)\dd\BrS$ in $X$.

\item[(13)] $X$ is a $1\dd\BrS$.

\item[(14)] $X$ is an $A\dd\,(1,2)\dd\BrS$.

\item[(15)] $X$ is a $(1,2)\dd\BrS$.

\item[(16)] $X$ is a $2\dd\BrS$.

\item[(17)] $X$ is an $A\dd\,(2,1)\dd\BrS$.

\item[(18)] $X$ is a $(2,1)\dd\BrS$.
\end{enumerate}
\end{corollary}

\begin{pf}
First of all, let us note that if $(Y,\tau_1'<\tau_2')\sbs
(X,\tau_1<\tau_2)$, $\tau_1<_S\tau_2$ and $Y\in 2\dd\,\cD(X)$,
then by (3) of Corollary~2.1.6 in [8], $\tau_1'<_S\tau_2'$.
Therefore on the one hand, $\tau_1<_S\tau_2$ and (8) of
Theorem~4.1.6 in [8] give that $(13)\llra (14)\llra (15)\llra
(16)\llra (17)\llra (18)$ and, on the other hand,
$\tau_1'<_S\tau_2'$ and the same (8) of Theorem~4.16 give that
$(1)\llra (2)\llra (3)\llra (4)\llra (5)\llra (6)$.

Therefore, since $Y\in 2\dd\,\cD(X)$, by (1) and (2) of
Corollary~12.75, $(3)\llra (9)$, $(5)\llra (11)$ and since
$2\dd\,\cD(X)\sbs 1\dd\,\cD(X)$, Corollary~12.76 gives that
$(1)\llra (7)$ and $(4)\llra (10)$. Moreover, by (2) of
Corollary~12.29, $(8)\llra (9)$. Therefore, we obtain that
$(7)\llra (8)\llra (9)\llra (10)\llra (11)$. Since $(10)\llra
(11)$ and by (2) of Proposition~12.31, $(10)\lra (12)$, it remains
to prove only that $(12)\llra (11)$. But this implication is given
by (1) Corollary~12.29.

Finally, since $Y\in \tau_1\cap 2\dd\,\cD(X)\lra Y\in\tau_1\cap
1\dd\,\cD(X)$, by Proposition~1.14 and Theorem~1.15 in [12], we
have $(1)\llra (13)$. This completes the proof.
\end{pf}

\begin{corollary}{12.78}
If $(Y,\tau')$ is a $\TsS$ of a $\TS$ $(X,\tau)$ and $Y\in
      \linebreak    \tau\cap \cD(X)$, then the following conditions are equivalent:
\begin{enumerate}
\item[(1)] $Y$ is a $\BrS$.

\item[(2)] $Y$ is a $\BrS$ in $X$.

\item[(3)] $X$ is a $\BrS$.
\end{enumerate}
\end{corollary}

\vskip+0.8cm

\end{document}